%% file: all.tex
% Evaluation on 2/20/04

\input defv7.tex

%%%%%%%%%%%%%%%%%%%%%%%%%%%%
% 0. Introduction
%%%%%%%%%%%%%%%%%%%%%%%%%%%%

\input 0int.tex

%%%%%%%%%%%%%%%%%%%%%%%%%%%%%%%%%%%%%%%%%%%%%%%%%%%
%%%%%%%%%%%%%%%%%%%%%%%%%%%%%%%%%%%%%%%%%%%%%%%%%%%

%%%%%%%%%%%%%%%%%%%%%%%%%%%%
% I.1.Definition of log-Riemann surfaces (*****)
%%%%%%%%%%%%%%%%%%%%%%%%%%%%

\input i1.tex

%%%%%%%%%%%%%%%%%%%%%%%%%%%%
%I.2 Euclidean metric (****) finish I.2.3.1 and I.2.3.2
%%%%%%%%%%%%%%%%%%%%%%%%%%%%

\input i2.tex

%%%%%%%%%%%%%%%%%%%%%%%%%%%%
%I.3. Topology (****) minimal skeleton?
%%%%%%%%%%%%%%%%%%%%%%%%%%%%

\input i3.tex

%%%%%%%%%%%%%%%%%%%%%%%%%%%%
%I.4. Ramified coverings (****) review
%%%%%%%%%%%%%%%%%%%%%%%%%%%%

\input i4.tex

%%%%%%%%%%%%%%%%%%%%%%%%%%%%
%I.5 Ramification surgery (***) complete and enlarge
%%%%%%%%%%%%%%%%%%%%%%%%%%%%

\input i5.tex

%%%%%%%%%%%%%%%%%%%%%%%%%%%%%%%%%%%%%%%%%%%%%%%%%%%%%%%
%%%%%%%%%%%%%%%%%%%%%%%%%%%%%%%%%%%%%%%%%%%%%%%%%%%%%%%

%%%%%%%%%%%%%%%%%%%%%%%%%%%%
%II.1.Type (***) add more survey
%%%%%%%%%%%%%%%%%%%%%%%%%%%%

\input ii1.tex

%%%%%%%%%%%%%%%%%%%%%%%%%%%%
%II.2.Boundary behaviour (***) finish proof
%%%%%%%%%%%%%%%%%%%%%%%%%%%%

\input ii2.tex

%%%%%%%%%%%%%%%%%%%%%%%%%%%%
%II.3. Caratheodory theorem (***) convergence for completed surfaces
%%%%%%%%%%%%%%%%%%%%%%%%%%%%

\input ii3.tex

%%%%%%%%%%%%%%%%%%%%%%%%%%%%
%II.4.Quasi-conformal theory (*) To Do
%%%%%%%%%%%%%%%%%%%%%%%%%%%%

\input ii4.tex

%%%%%%%%%%%%%%%%%%%%%%%%%%%%%
%II.5. Uniformization (***) Proof by Schwarz-Christ. and review
%%%%%%%%%%%%%%%%%%%%%%%%%%%%%

\input ii5.tex

%%%%%%%%%%%%%%%%%%%%%%%%%%%%%
%II.6. Cyclotomic log-Riemann surfaces (*) To Do
%%%%%%%%%%%%%%%%%%%%%%%%%%%%%

\input ii6.tex

%%%%%%%%%%%%%%%%%%%%%%%%%%%%%
%II.7. Uniformization of infinite (*) To Do
%%%%%%%%%%%%%%%%%%%%%%%%%%%%%

\input ii7.tex

%%%%%%%%%%%%%%%%%%%%%%%%%%%%%%%%%%%%%%%%%%%%%%%%%%%%%%%%%%%%%%%%%%%%
%%%%%%%%%%%%%%%%%%%%%%%%%%%%%%%%%%%%%%%%%%%%%%%%%%%%%%%%%%%%%%%%%%%%

%%%%%%%%%%%%%%%%%%%%%%%%%%%%%%%%
%III.1.Ring of special functions (*****)
%%%%%%%%%%%%%%%%%%%%%%%%%%%%%%%%

\input iii1.tex

%%%%%%%%%%%%%%%%%%%%%%%%%%%%%%%%%%
%III.2. Refined analytic estimates (****) correct Stolz
%%%%%%%%%%%%%%%%%%%%%%%%%%%%%%%%%%

\input iii2.tex

%%%%%%%%%%%%%%%%%%%%%%%%%%%%%%%%%%
%III.3. Liouville Thm (****) correct Stolz
%%%%%%%%%%%%%%%%%%%%%%%%%%%%%%%%%%

\input iii3.tex

%%%%%%%%%%%%%%%%%%%%%%%%%%%%%%%%%%
%III.4. Structural ring (**) Enlarge and complete
%%%%%%%%%%%%%%%%%%%%%%%%%%%%%%%%%%

\input iii4.tex

%%%%%%%%%%%%%%%%%%%%%%%%%%%%%%%%%%%%%%%%%%%%%%%%%%%%%%%%%%%%%%%%%%%%%%%%
%%%%%%%%%%%%%%%%%%%%%%%%%%%%%%%%%%%%%%%%%%%%%%%%%%%%%%%%%%%%%%%%%%%%%%%%

%%%%%%%%%%%%%%%%%%%%%%%%%%%%%%%%%%
% Bibliography
%%%%%%%%%%%%%%%%%%%%%%%%%%%%%%%%%%

\input bib.tex

\end

%% file: defv7.tex
\def\D{\Delta} 
\def\a{\alpha}

\def\d{\delta}

\def\G{\Gamma}
\def\g{\gamma}
\def\eps{\varepsilon}
\def\t{\theta}
\def\O{\Omega}
\def\o{\omega}
\def\w{\omega}

\def \dd#1{{\bf#1}}

\def\cl#1{{\cal#1}}

%Simbolos matematicos.

\def\mod{\mathop{\rm mod}\nolimits}
\def\Max{\mathop{\rm Max}}

\def\det{\mathop{\rm det}}

\def\Arg{\mathop{\rm Arg}\nolimits}

\def\Im{\mathop{\rm Im}\nolimits}
\def\Re{\mathop{\rm Re}\nolimits}
\def\id{\mathop{\rm id}\nolimits}

%Macros.

\def\ouv#1{\smash{\mathop{#1}\limits^{\lower 1pt\hbox
{$\scriptscriptstyle\circ$}}}}

\def\apl#1#2#3#4#5{{#1}:\matrix{
		\hfill#2 & \longrightarrow & #3\hfill \cr
		\hfill#4 & \longmapsto     & #5\hfill \cr}}

\def\ap#1#2#3#4{\matrix{
		\hfill#1 & \longrightarrow & #2\hfill \cr
		\hfill#3 & \longmapsto     & #4\hfill \cr}}

\def\hfl#1#2{\smash{\mathop{\hbox to 12mm{\rightarrowfill}}
\limits^{\scriptstyle#1}_{\scriptstyle#2}}}

%Titulos, enunciados.

\long\def\eno#1#2{\par\smallskip{\bf{#1}}{\it\ {#2}}\par\medskip}

\def\tit#1{\vskip 5mm plus 1mm minus 2mm {\tir #1}
		\vskip 3mm plus 1mm minus 2mm}

\def\stit#1{\vskip 3mm plus 1mm minus 2mm {\bf{#1}}
		\smallskip}

\font\tir=cmbx10 at 12pt

\long\def\ctitre#1{\vbox{\leftskip=0pt plus 1fil
\rightskip=0ptplus 1fil \parfillskip=0pt#1}}

\def\ref#1#2#3#4{{\bf #1}{\ #2}{\it ,\ #3}{,\ #4}\medskip}

%Dibujos

\def \picture #1 by #2 (#3){\midinsert \centerline 
{\vbox to #2{\hrule width #1 heigth 0pt 
depth 0pt \null \vfill \special {picture #3}}}\endinsert }

\def\scaledpicture #1 by #2 (#3 scaled #4) {{
\dimen0 =#1 \dimen1 =$2
\divide \dimen0 by 1000 \multiply \dimen0 by #4
\divide \dimen1 by 1000 \multiply \dimen1 by #4
\picture \dimen0 by \dimen1 (#3 scaled $4)}}

\def\figure #1 #2 #3 {\midinsert \vglue 3mm 
{\vbox to #3 {\hrule width 6cm height 0cm depth 0cm \vfill
{\special {picture #1 scaled #2}}}}\vglue 2mm \endinsert}

\magnification=1200

%% file: 0int.tex
\input defv7.tex
\input psfig.sty
\font\gothic=eufm10

\overfullrule=0pt

\input defv7.tex
\long\def\ctitre#1{\vbox{\leftskip=0pt plus 1fil
\rightskip=0ptplus 1fil \parfillskip=0pt#1}}

\vskip 2cm plus 5mm minus 5mm
\bigskip
\ctitre{\tir LOG-RIEMANN SURFACES}
\bigskip
\bigskip
{\centerline {{{\bf K. Biswas, R. P\'erez-Marco} \footnote{*}{CNRS,
LAGA, UMR 7539, Universit\'e Paris 13, France.}}}

%{\centerline {{{\bf K. Biswas, R. P\'erez-Marco} \footnote{*}{UCLA,
%Dept. Mathematics, 405, Hilgard Ave., Los Angeles CA-90024, USA.}}}

\bigskip
\bigskip

%{\centerline {Preliminary version}}

\bigskip
\bigskip

{\centerline {\tir {Table of contents.}}}

\bigskip
\bigskip

{\bf  Introduction: Motivations and contents.}\par

\bigskip

    {\bf I. Geometric theory of log-Riemann surfaces.}\par
    \medskip

    {\parindent=2cm

        \item{}I.1) Definition of log-Riemann surfaces.\par
        \smallskip

        {\parindent=2.5cm

            \item{}I.1.1) Definition.\par
            \smallskip

                    \item{}I.1.2) Examples\par
            \smallskip
                      }

        \item{}I.2) Euclidean metric and ramification points.\par
        \medskip

        {\parindent=2.5cm

            \item{}I.2.1) Definition.\par
            \smallskip

                    \item{}I.2.2) Ramification points.\par
            \smallskip

            \item{}I.2.3) Log-euclidean geometry and charts.\par
            \smallskip
            {\parindent=3cm

                \item{}I.2.3.1) Preliminaries.\par
                \smallskip

                        \item{}I.2.3.2) Constructions of charts.\par
                \smallskip

                \item{}I.2.3.3) Other applications.\par
                \smallskip

                \item{}I.2.3.4) The Kobayashi-Nevanlinna net.\par
                \smallskip
                    }
                 }
        \medskip

        \item{}I.3) Topology of log-Riemann surfaces.\par
        \medskip

        {\parindent=2.5cm

            \item{}I.3.1) The skeleton.\par
            \smallskip

                    \item{}I.3.2) Skeleton and fundamental group. \par
            \smallskip
                      }

        \item{}I.4) Ramified coverings.\par
        \medskip

        {\parindent=2.5cm

            \item{}I.4.1) Ramified coverings and formal Riemann surfaces.\par
            \smallskip

                    \item{}I.4.2) Ramified coverings and log-Riemann surfaces.\par
            \smallskip

                    \item{}I.4.3) Ramified coverings of formal Riemann surfaces.\par
            \smallskip

            \item{}I.4.4) Universal covering of log-Riemann surfaces.\par
            \smallskip

                    \item{}I.4.5) Ramified coverings and degree.\par
            \smallskip

                      }
        \medskip

        \item{}I.5) Ramification surgery.\par
        \medskip

        {\parindent=2.5cm

            \item{}I.5.1) Grafting of ramification points.\par
            \smallskip

                    \item{}I.5.2) Pruning of ramification points.\par
            \smallskip
                      }
        \medskip
    }

\null
\vfill
\eject

    {\bf II. Analytic theory of log-Riemann surfaces.}\par
    \medskip

        {\parindent=2cm

        \item{}II.1) Type of log-Riemann surfaces.\par
        \medskip

        {\parindent=2.5cm

            \item{}II.1.1) General facts.\par
            \smallskip

            \item{}II.1.2) Kobayashi-Nevanlinna criterium.\par
            \smallskip

                      }
        \medskip
        \item{}II.2) Boundary behaviour of the universal cover.\par
        \medskip

        \item{}II.3) Caratheodory theorem for log-Riemann surfaces.\par
        \medskip

        {\parindent=2.5cm

            \item{}II.3.1) Kernel convergence.\par
            \smallskip

                    \item{}II.3.2) Caratheodory theorem.\par
            \smallskip

                    \item{}II.3.3) Conformal radius.\par
            \smallskip

            \item{}II.3.4) Closure of algebraic log-Riemann surfaces.\par
            \smallskip

                    \item{}II.3.5) Stolz continuity.\par
            \smallskip

                    \item{}II.3.6) Functions holomorphic at ramification points.\par
            \smallskip

            \item{}II.3.7) General Weierstrass theorem.\par
            \smallskip

                      }
        \medskip

        \item{} II.4) Quasi-conformal theory of log-Riemann surfaces.\par

        \medskip

        {\parindent=2.5cm

            \item{}II.4.1) Complex automorphisms of log-Riemann surfaces.\par
            \smallskip

                    \item{}II.4.2) Teichm\"uller distance.\par
            \smallskip

                    \item{}II.4.3) Teichm\"uller space of a finite log-Riemann surface.\par
            \smallskip

                      }
        \medskip

        \item{} II.5) Uniformization of finite log-Riemann surfaces.\par

        \medskip

        {\parindent=2.5cm

            \item{}II.5.1) A fundamental example.\par
            \smallskip

                    \item{}II.5.2) On the general case.\par
            \smallskip

                    \item{}II.5.3) Uniformization theorem.\par
            \smallskip

                    \item{}II.5.4) Uniformization and Schwarz-Christoffel formula.\par
            \smallskip

                    \item{}II.5.5) General uniformization theorem.\par
            \smallskip
                      }

        \medskip

        \item{} II.6) Cyclotomic log-Riemann surfaces.\par

        \medskip

        {\parindent=2.5cm

            \item{}II.6.1) Definition.\par
            \smallskip

                    \item{}II.6.2) Ramification values .\par
            \smallskip

                    \item{}II.6.3) Subordination of cyclotomic log-Riemann surfaces.\par
            \smallskip

                    \item{}II.6.4) Caratheodory limits of cyclotomic log-Riemann surfaces.\par
            \smallskip

                    \item{}II.6.5) Continued fraction expansion of the uniformizations.\par
            \smallskip

                    \item{}II.6.6) Relation to the incomplete Gamma function.\par
            \smallskip

                    \item{}II.6.7) Relation to Hermite polynomials.\par
            \smallskip

                      }

        \medskip

        \item{} II.7) Uniformization of infinite log-Riemann surfaces.\par

        \medskip

    }

\null
\vfill
\eject

    {\bf III. Algebraic theory of log-Riemann surfaces.}\par
    \medskip

        {\parindent=2cm

        \item{}III.1) A ring of special functions.\par
        \smallskip

        {\parindent=2.5cm

            \item{}III.1.1) Definition.\par
            \smallskip

                    \item{}III.1.2) Asymptotics. \par
            \smallskip

                    \item{}III.1.3) Linear independence. \par
            \smallskip

                    \item{}III.1.4) Algebraic independence. \par
            \smallskip

                    \item{}III.1.5) Integrals. \par
            \smallskip

                    \item{}III.1.6) Differential ring structure. \par
            \smallskip

                    \item{}III.1.7) Picard-Vessiot extensions.\par
            \medskip

                    \item{}III.1.8) Liouville classification.\par
            \medskip
                      }

        \item{}III.2) Refined analytic estimates.\par
        \smallskip

        {\parindent=2.5cm

            \item{}III.2.1) Decomposition of the end at infinity.\par
            \smallskip

                    \item{}III.2.2) Analytic estimates.\par
            \medskip
                      }

                \item{}III.3) Liouville theorem on log-Riemann surfaces.\par
                \smallskip

        \item{}III.4) The structural ring.\par
        \smallskip

        {\parindent=2.5cm

            \item{}III.4.1) Definition.\par
            \smallskip

            \item{}III.4.2) Points of $\cl S^*$ as maximal ideals.\par
            \smallskip

            \item{}III.4.3) The ramificant determinant.\par
            \smallskip

            \item{}III.4.4) Infinite ramification points.\par
            \medskip
                      }

}

\bigskip

{\bf Bibliography.}

\bigskip

\bigskip

\bigskip


\vfill

date of compilation : \the\day -\the\month -\the\year .

\null
\eject

\bigskip
\ctitre{\tir INTRODUCTION}
\bigskip

\tit {Historical motivation.}

Since the revolutionary ideas of Bernhard Riemann in the XIXth
century the notion of Riemann surface has experienced dramatic
changes. The motivation of Riemann's fundamental memoir [Ri1] is the
study of Abelian Integrals which are in the XIXth century the "new"
transcendental functions attracting the attention since the work of
Abel, Galois and Jacobi. The arithmetico-geometric history of
transcendental functions is very instructive. The first
transcendental functions are associated to the geometry of the
circle. These have been studied since the origins of Mathematics
(the famous Babylonian clay tablet Plimpton 322, of the Plimpton
collection, is a tabulation of arctangents for Pythagorean
triangles, see [Fr] and [Va]). These elementary trigonometric
functions, are also obtained by integration of elementary algebraic
differentials. For instance,
$$
\arcsin z=\int^z {dx \over \sqrt { 1-x^2}} \ .
$$
These functions satisfy addition formulae as
$$
\arcsin z+\arcsin w =\arcsin \left ( z\sqrt {1-w^2}+w\sqrt{1-z^2}\right ) \ .
$$
The close relation of trigonometric functions to the complex
exponential was unveiled by L. Euler. J. Wallis (1655) attempted
the computation of the arc-length of ellipses leading to Elliptic
Integrals of the form
$$
\int ^z {dx \over \sqrt {P(x)}} \ ,
$$
where $P$ is a polynomial of degree $3$ or $4$.
Elliptic
Integrals form a new family of transcendental functions that are associated to
the geometry of elliptic curves or genus $1$ algebraic curves.

Giulio Fagnano (1716) and L. Euler (1752, 1757) discovered
addition theorems for them, as
$$
\int_0^{z} {dx \over \sqrt {(1-x^2)(1-k^2x^2)}}+\int_0^{w} {dx
\over \sqrt {(1-x^2)(1-k^2x^2)}}=\int_0^{\xi} {dx \over \sqrt
{(1-x^2)(1-k^2x^2)}} \ ,
$$
where $k$ is the parameter of the elliptic integral and $\xi$ is
determined by
$$
\xi (1-k^2 z^2 w^2)=w\sqrt {(1-z^2)(1-k^2z^2)} +z\sqrt
{(1-w^2)(1-k^2w^2)} \ .
$$

%$$
%\int^z {dx \over \sqrt {P(x)}} + \int^w {dx \over \sqrt {P(x)}}
%=\int^{F(z,w)} {dx \over \sqrt {P(x)}} \ ,
%$$
%where
%$$
%F(z,w)={1\over 4} \left ( \sqrt {P(z)} -\sqrt {P(w)} \over z-w \right )^2 -z-w \ .
%$$

Elliptic Integrals were later studied by A.M. Legendre and C.F. Gauss, and later
by N.H. Abel, C. Jacobi, Ch. Hermite,... who also studied more general Hyperelliptic
Integrals of the form
$$
\int^z {R(x) \over \sqrt {P(x)} } \ dx \ ,
$$
where $R$ is a rational function and $P$ is a polynomial of arbitrary degree.

\medskip

The next important progress was achieved by N.H. Abel who
considered general Abelian Integrals of the general form
$$
\int^z R(x,y) \ dx \ ,
$$
where $y$ is an algebraic function of $x$, i.e. satisfies an algebraic equation
with polynomial coefficients in $x$ or
$$
P(x,y)=0
$$
for $P(x,y)\in \dd C[x,y]$. Abel proved his famous general
addition theorem for these Abelian integrals. The sum
$$
\int^{x_1} R(x,y) \ dx + \int^{x_2} R(x,y) \ dx +\ldots
+\int^{x_n} R(x,y) \ dx
$$
taken with extremes at $(x_1, \ldots , x_n)$, which are the
intersection points of the curve
$$
P(x,y)=0
$$
and a family of algebraic curves
$$
Q(x,y,a_1,\ldots , a_n) =0\ ,
$$
is a rational function plus logarithmic terms of $(a_1, \ldots ,
a_m)$, the parameters that parameterize the intersecting family of
algebraic curves (the usual modern formulation of Abel theorem is
only a weaker particular case of the original result).

\medskip

It is understandable the sensation that this result caused (even
if it was ignored for some time): Abel's result shows that the
algebraic theory of these very general new transcendental
functions is very rich. Abel's result seems to have been also
discovered independently by \'E. Galois, as can be found in the
"brouillons" left by Galois (see [Gal] p.187 and p.518). The
corresponding geometry of general Abelian Integrals are general
algebraic curves. These new transcendental functions motivated
several fundamental theories.

\medskip

The close inspection of the manuscripts of Galois shows that his
motivation to build what is now called Galois theory was well beyond
the problem of resolution of algebraic equations. His ultimate goal was
a full classification of transcendental functions. He did made important
progress for Abelian Integrals dividing them into three kinds and studying
their periods (see [Pi] volume III p.472).

\medskip

Riemann discovered that Abelian Integrals live naturally on
Riemann surfaces spread over the Riemann sphere. Well before
Riemann, Euler was well aware of the natural multivaluedness of
algebraic and other important functions special functions (see
[Eu] chapter I where he defines "Functiones multiformes", and his
famous writings on the logarithm of negative numbers). The
audacious idea of Riemann is to pass at once from the
"multivaluedness" of the function to a geometric dissociation of
the space were the function lives by imagining new sheets were
each branch of the function is univalued. Riemann surfaces, as
Riemann understood them, are abstract manifolds, in the sense that
they are not embedded in any ambient space. Riemann had a clear
understanding of the notion of abstract manifold as is
demonstrated by his Inaugural Dissertation on the foundations of
geometry (see [Ri2]). One cannot explain otherwise the
consideration of non-Riemannian metrics. But Riemann surfaces, in
the view of Riemann, are always spread over the complex plane $\dd
C$ or the Riemann sphere ${\overline {\dd C}}$. These are called
today Riemann domains. They come equipped with canonical
coordinates or charts. This is why some schools, as the German or
the Russian around the middle of the XXth century, gave them the
name of "concrete Riemann surfaces". The equipment with canonical
charts or coordinates enriches the Riemann surface structure. In
particular it enables the link between the geometry and the
transcendental functions. The modern notion of Riemann surface (we
should say "abstract Riemann surface") does not come equipped with
a preferred set of charts. This modern notion was conceived by T.
Rado and H. Weyl (see [Wey]) and marks the origin of intrinsic
differential geometry.

\medskip

The influence of Riemann's ideas was deep in the Mathematics of
the XIXth century. R. Dedekind assisted in 1855-1856, to Riemann's
lectures in G\"ottingen on Abelian Functions. They were close
friends, and one may wonder how much of his achievements in
Algebraic Number Theory were influenced by Riemann's theory of
algebraic curves. Dedekind's theory of ideals, extending E. Kummer
"ideal numbers", marks the birth point of Commutative Algebra. It
allows the unification of the theory of Number Fields and that of
Function Fields on algebraic curves (or compact Riemann surfaces).
The culmination of this unification is his article with H. Weber
published in 1882 ([De-We]). For a modern and faithful account of
this theory the reader can read the excellent exposition of J.
Mu\~noz D\' \i az in [Mu], the book of C. Chevalley [Che], and the
original memoir of Dedekind and Weber. The idea of unification in
Science was in the mood of times. Maybe better known, or better
popularized, in Physics by Maxwell's theory of Electrodynamics.
Dedekind-Weber arithmetico-geometrical unification is very much in
the spirit of subsequent work in Algebraic Geometry and Number
Theory during the XXth century. Dedekind-Weber theory provides a
dictionary between Number Fields and Function Fields. From the
algebraic structure of the function field of meromorphic functions
on the compact Riemann surface, Dedekind-Weber theory recovers
algebraically the points of the Riemann surface, they are
identified with the valuation sub-rings in the field. In the
affine model we recover points as prime ideals, and ramification
points correspond to prime ideals that ramify on the Number Field.
The spectral reconstruction of the space is now a well known
important idea that penetrates Mathematics of the XXth century
(for example, I.M. Gel'fand theory of normed algebras (1940) is
based on it), as well as Physics (Quantum Mechanics).

\medskip

Despite these early success, during the XXth century, the
intrinsic and coordinate independent inclination in differential
geometry, with a total aversion to preferred coordinate systems,
erased completely Riemann's original notion of Riemann surface. It
is easy to check that an important number of contemporary
mathematicians have problems telling precisely what is the
difference between the Riemann surface of the logarithm and $\dd
C$. Riemann's notion was progressively replaced by Weyl's
intrinsic, coordinate independent, notion. In that way the direct
link to transcendental functions was broken, and this original,
historical and fruitful motivation through Abelian Integrals
(which were progressively degrated to Abelian Differentials) was
lost.

\medskip

One of the main goals of this article is to go back to these
origins and show how much is missing through this modern point of
view. In particular to pursue the link between new transcendental
functions and new geometries necessitates Riemann's original point
of view. This article is the first of a series where we enlarge
the class of Abelian Integrals
$$
\int R(x,y) \ dx
$$
to a larger class of integrals leading to new transcendental
functions of the form
$$
\int R_1(x,y) \ e^{R_2(x,y)} \ dx
$$
where $R_1$ and $R_2$ and rational functions and $y$ is an algebraic
function of $x$. In the restricted situation that we consider in
this article, $R_1$ and $R_2$ are polynomials and $y=x$. The
corresponding geometry is a class of Riemann domains with a finite
number of ramification points, some of which can be logarithmic
ramification points (also called infinite ramification points). Thus
these complex curves are no longer algebraic. Our aim is to extend
Dedekind-Weber theory to this setting so that we can include Riemann
surfaces (in Riemann sense) spread over $\dd C$ (or ${\overline {\dd
C}}$) with some infinite ramification points. In Dedekind-Weber's
dictionary, points on the algebraic curve correspond to prime ideals
(or maximal ideals since these coincide in dimension $1$), and
finite ramification points do correspond to ramified primes in
Number Fields. The extension of the geometric picture in order to
include Riemann surfaces with infinite ramification points should
correspond to a certain type of transcendental extensions of $\dd Q$
of finite transcendental degree, probably not unrelated to a
non-abelian Iwasawa theory ([Iwa], [Was]). Today this Transalgebraic
Number Theory remains largely unexplored, but it remains one of our
motivations. We refer to [PMBBJM] for a historical introduction, for
the exposition of the philosophy governing this research, and a few
steps into this unknown territory. This point of view was deeply
rooted in Galois mind, as the second author has noticed after
reading in repeated occasions Galois' memoirs. As Galois writes, his
meditations in that subject did occupy him during his last year of
life while in prison. As he announced with clairvoyance in his
posthumous letter to his friend Auguste Chevalier (see [Ga2] p.185,
and [Ga1] p.32),

\medskip

{\it "...Mais je n'ai pas le temps, et mes id\'ees ne sont pas
encore bien d\'evelopp\'ees sur ce terrain qui est immense..."}
\footnote {$\dagger$}{"...But I don't have time, and my ideas are
not well developed in this immense domain..."}

\medskip

This article is a step into that direction. Our aim is to develop
the geometric side of the Dedekind-Weber dictionary that we believe
does extend to the Transalgebraic world. This should shed new light
on Transalgebraic Number Theory, the counter-part of the dictionary.
In particular, into the main problem of determining which
transcendental extensions are transalgebraic. This is a central
problem that demands to be elucidated. This transalgebraic
extensions should be transcendental and maybe of finite
transcendental degree, but not all transcendental extensions of
finite degree are transalgebraic. It is natural to expect that
generators of these finite transcendental extensions should be
provided by values obtained from the special functions appearing in
the geometric Dedekind-Weber theory. The location of finite
ramification points for algebraic log-Riemann surfaces with
algebraic uniformizations defined over a Number Field, generate
algebraic extensions. Examples of transalgebraic numbers should
include (in order of increasing complexity):
$$
\pi, e, \log (p/q), \Gamma (p/q) , \gamma , \zeta (3),\ldots
$$
It is natural to expect that the location of infinite ramification
points for transalgebraic log-Riemann surfaces with uniformizations
defined over a Number Field (i.e. with rational coefficients for
example) should define transalgebraic extensions of this Number
Field. For example, for cyclotomic log-Riemann surfaces studied in
section II.6, we obtain the values of the $\Gamma$-function at
rational arguments. This philosophy can be linked to Kronecker's
"Judgendtraum" and Hilbert's twelfth problem, which seems to have
remained largely misunderstood. This is all about the generation of
(trans)algebraic extensions by analytic means. A close parallel
philosophy comes also from deep intuitions of J. Mu\~noz D\'\i az
about the possibility of generating points of algebraic curves
defined over a number field by means of divisors of certain types of
transcendental functions. A remarkable result from the Salamanca
school is the Thesis of P. Cutillas ([Cu]), where this author proves
Mu\~noz conjecture on the existence and uniqueness of a canonical
field of transcendental functions with finite order fixed essential
singularities, that determines completely the Riemann surface in
Dedekind-Weber style.

\medskip

 The goals of this first article are modest. Only the affine model
and the genus $0$ case are considered here. As said before, this
corresponds to log-abelian integrals of the form
$$
\int P_1 (x)e^{P_2(x)} \ dx \ .
$$
We develop in this preliminary setting the different angles through
which we can view the theory: Geometric, Analytic and Algebraic.

%\vfill \null \eject

\medskip

\stit {Meccano motivation.}

\medskip

What lies behind Dedekind-Weber motivation is the correspondence
between a geometric meccano and an algebraic meccano. Under this
diccionary simple geometric operations correspond to intricate
arithmetic operations and conversely. We conceive the geometric
meccano as a lego box containing some sort of pieces or building
blocks, and a set of construction rules. Using these we can build
a class of geometric objects. These geometric objects have an
algebraic counterpart.

To fix the ideas we can consider a very simple geometric meccano:
We are allowed to cut and paste by the identity a finite set of
complex planes without creating any topology (i.e. the geometric
manifold thus constructed is supposed to be simply connected).
This geometric meccano corresponds to the algebraic meccano build
up with finite operations with a free variable $z$: These are all
polynomial expression generated by $z$, that is the ring of
polynomials $\dd C[z]$. Via the uniformization, we establish the
identity of these two meccanos: A uniformization from $\dd C$ to
the Riemann surface constructed has polynomial expression in the
canonical charts. Observe that the simple operation of addition of
polynomials is a mysterious binary operation on the corresponding
Riemann surfaces. An open question is: Construct geometrically the
Riemann surface corresponding to the sum. Conversely, the grafting
of one such Riemann surface onto another, by cutting and pasting
through slits, is an even more mysterious algebraic binary
operation at the level of the ring of polynomials.

\medskip

 The richness of this point of view consists in the possibility of
enlarging successively the geometric meccano by adding new
building blocks or new construction rules. Thus in the precedent
meccano we could allow to paste an infinite number of planes and
allow infinite ramification points, but keep the total number of
ramification points finite. The enlarged class of geometric
objects becomes the class of transalgebraic log-Riemann surfaces
that are studied in this article. The corresponding enlargement of
the algebraic meccano consist in allowing not only primitives of
polynomials but also primitives of products of polynomials and
exponentials of polynomials.

\medskip

Another strong point of the meccano intuition is that we are led
naturally to consider sub-meccanos. For example, in the precedent
meccano we may allow only to paste the planes through slits ending
at algebraic points, i.e. branch points lie only above algebraic
points. This restricted construction rule defines a sub-class of
the precedent class. This sub-class corresponds to polynomials in
$\bar \dd Q [z]$, i.e. polynomials with algebraic coefficients.
Belyi theorem states that, up to birational equivalence, the same
geometric sub-meccano is obtained by allowing only ramification
points over $0$ and $1$ (and also $\infty$), i.e. by using only
cuts ending at $0$ and $1$ in the sheets.

\medskip

The possibilities for enlargement of the meccano are endless. For
instance, we may want to have as uniformization the integral of a
rational function without simple poles. Then we obtain the
projective model of transalgebraic log-Riemann surfaces. But if we
allow to integrate arbitrary rational functions, then we need to
enlarge the geometric building blocks by allowing not only complex
planes $\dd C$ but also "tubes" $\dd C/\dd Z$ and euclidean
polygons (or more precisely "log-polygons"). This generates the
class of tube-log Riemann surfaces that will be studied in future
articles.

\medskip

\stit {Dynamical System motivation.}

\medskip

As just mentioned tube-log Riemann surfaces are similar to
log-Riemann surfaces but more general: We allow to cut and paste
not only planes $\dd C$ but also tubes $\dd C/\dd Z$. There is the
subclass of those with only a finite number of tubes and a finite
number of ramification points. Originally the second author used
the tube-log Riemann surface similar to the log-Riemann surface of
the logarithm except that one plane was replaced by one cylinder,
in order to solve several open problems in holomorphic dynamics
(see figure below and [PM1],[PM2]). The tube-log Riemann surface
of the figure has an important special function as uniformization:
The logarithmic integral
$$
\int^z_{-\infty} {e^z \over z} \ dz \ .
$$
This is just a particular illustration of the general algebraic
meccano correspondence.

\medskip

This geometry, and not another, proves the optimality of the
diophantine condition ($(p_n/q_n)$ is the sequence of convergents
of the rotation number appearing in the problem)
$$
\sum_{n=1}^{+\infty } {\log \log q_{n+1} \over q_n} < +\infty \ ,
$$
in Siegel problem of linearization of holomorphic dynamics with no
strict periodic orbits (see [PM1]).

\vfill\null\eject

\tit {Description of the article.}

\medskip

 The article is divided into three main sections. In section
I we define log-Riemann surfaces. That are the proper formalization
of Riemann surfaces from Riemann's original point of view. We
discuss there topological and geometrical aspects. In section II we
study various analytic aspects of log-Riemann surfaces. Finally in
section III we develop the algebraic theory and study the natural
function fields defined on log-Riemann surfaces. As in
Dedekind-Weber theory we can recover algebraically from this field
of functions the points of the log-Riemann surface, including the
ramification points finite or infinite. We are able to distinguish
algebraically finite from infinite ramification points.

\medskip

\stit {Contents of part I: Geometric Theory of log Riemann
surfaces.}

\medskip

In section I.1 we define log-Riemann surfaces which is the proper
formalization of Riemann's classical notion of Riemann surfaces. The
definition restricts the class of modern Riemann surfaces by
imposing the existence of an atlas with trivial (identity) change of
charts. The defines a preferred coordinate in the charts. We give
numerous examples of "classical" log-Riemann surfaces. These can be
constructed by isometrically gluing together by the identity complex
planes through half-line slits. In particular, classical algebraic
curves over $\dd C$, once represented over the complex plane as
Riemann domains, are examples of log-Riemann surfaces. We prove in
section II.2 that all such algebraic Riemann surfaces can be
obtained in such a way, i.e. we can build them using half-line cuts
and not just segment cuts as it is classically done. In section I.2
we develop the metrical theory. A log-Riemann surface inherits a
natural flat conformal metric coming from its preferred coordinate:
The log-euclidean metric. Log-euclidean geometry is at the same time
an elementary locally euclidean geometry, but very rich globally.
This metric is non-smooth but numerous results from Riemannian
geometry subsist. We develop some of these, and some results extend
euclidean geometry to this setting. We have a rich convex geometry.
Also log-euclidean geometry is useful in order to construct "minimal
atlases" which play an important role in the applications. The
completion of a log-Riemann surface $\cl S$ is a completed space
$$
\cl S^*=\cl S \cup \cl R \ ,
$$
where $\cl R$ is a closed set: The ramification set. Ramification
points are defined as isolated points in $\cl R$. Ramification
points are of two kinds: Finite ramification points, of finite
order $n<+\infty$, and infinite ramification points, of order
$n=+\infty$. We restrict the study to those log-Riemann surfaces
having a discrete ramification set $\cl R$. It is natural to
consider more general ramification sets as those appearing in the
study of entire functions, but these will be discussed elsewhere
(see the forthcoming [Bi-PM1]). Even when $\cl R$ is discrete, the
completed log-Riemann surface $\cl S^*$ does not inherit in
general of a Riemann surface structure, not even the structure of
a topological surface. They are not even locally compact when
there are infinite ramification points. This completion $\cl S^*$
is by definition a formal Riemann surface. These formal Riemann
surfaces are natural objects that deserve a study by themselves.
In section I.4 we develop a general theory of these objects and
their ramified coverings. The natural notion of ramified covering
extends the classical one. It is a more sophisticated notion that
unexpectedly has better behavior in the category of formal Riemann
surfaces. In section I.3 we study the topology of log-Riemann
surfaces. A combinatorial object, a skeleton, is associated to
log-Riemann surface. Such skeleton contains all the information
about the topology of $\cl S^*$, in particular the fundamental
group of $\cl S^*$ can be read in the fundamental group of a
skeleton associated to it.

In the last subsection I.5 we introduce some useful surgeries
involving ramification points. In particular the grafting of
ramification points plays an important role in the Algebraic
Theory. This surgery for simply connected Riemann surfaces has
been studied recently by M. Taniguchi [Ta3].

\medskip

\stit {Contents of part II: Analytic Theory of log-Riemann
surfaces.}

\medskip

 In section II.1 we discuss the type problem, a topic that
has occupy most of the work on the related field of entire
functions. When the log-Riemann surface is simply connected, it is
important to recognize if it is of parabolic or hyperbolic type. A
geometric criterium of Z. Kobayashi [Ko] and R. Nevanlinna [Ne2] is
adapted for log-Riemann surfaces. In section II.2 we study radial
limits of the uniformization in the spirit of the theory of boundary
behavior of conformal representations. We prove that for hyperbolic
log-Riemann surfaces, infinite ramification points correspond to a
countable set in the boundary ofthe universal cover. In section II.3
we generalize Caratheodory kernel convergence for planar domains to
log-Riemann surfaces and domains in log-Riemann surfaces. We
determine the closure of algebraic log-Riemann surfaces with a
bounded number of finite ramification points. This yields the class
of transalgebraic log-Riemann surfaces, the simplest class of
log-Riemann surfaces just after the algebraic ones. In section II.4
we start a quasi-conformal theory of log-Riemann surfaces. This
theory is richer than just the quasi-conformal theory of the
underlying Riemann surfaces. Local lipschitz behaviour at the
ramification points is critical. The space of quasi-conformal
deformations of log-Riemann surfaces is larger, i.e. has more
parameters, than the one of the underlying Riemann surface. We
define the Teichm\"uller distance and generalize the classical
convergence theorems. In section II.5 we give formulas for the
uniformization of transalgebraic log-Riemann surfaces. Their
uniformizations are of the form
$$
F_0(z)=\int_0^z P_1(t) e^{P_0(t)} \ dt \ ,
$$
where $P_1$ and $P_0$ are polynomials. Conversely any log-Riemann
surface with such uniformization is a transalgebraic log-Riemann
surface. The number of finite ramification points is $d_1=\deg
P_1$ and the number of infinite ramification points is $d_0=\deg
P_0$. Thus we can identify the space of transalgebraic log-Riemann
surfaces with
$$
\dd C [z]^*\times \dd C [z] \ .
$$
 This result can be attributed to R. Nevanlinna [Ne1] and
has been rediscovered since then (this happen to us and to others
[Ta1]). In recent work, M. Taniguchi studies entire functions
which are uniformizations of log-Riemann surfaces from a geometric
point of view and overlaps with some parts of our study. In
section II.6 we study a particular class of transalgebraic
log-Riemann surfaces, cyclotomic log-Riemann surfaces, with
uniformizations of the form
$$
F_{j,d}(z)=\int_0^z t^j e^{t^d}\ dt \ .
$$
This uniformization has a remarkable continued fraction expansion
and also asymptotic expansions in some sectors that we study in
detail.

\medskip

\stit {Contents of part III: Algebraic theory of log-Riemann
surfaces.}

\medskip

The Algebraic study is very much in the spirit of Mathematics of the
XIXth century. In that time complex analysis and algebra were
working hand to hand. Algebraic theories, as the one of elliptic
functions and theta functions, were developed by analysts, more
precisely, the path of research was set by intuition driven by
complex analysis. The right algebraic objects are located thanks to
analytic results of Liouville type. An example of this is Liouville
theorem. The elements of the basic ring of polynomials $\dd C[z]$
are characterized as those entire functions with at most polynomial
growth at infinite. The algebraic theory of section III follows the
same path. The main difficulty has been to guess the right ring of
functions on which the extended Dedekind-Weber theory builds upon.
For a transalgebraic log-Riemann surface $\cl S$ with only
$d<+\infty$ infinite ramification points and having as
uniformization
$$
F_0(z)=\int_0^z e^{P_0(t)} \ dt \ ,
$$
the basic special functions generating our ring are,
$j=0,1,\ldots, d-1$, $d=\deg P_0$,
$$
\eqalign{F_0(z)&=\int_0^z \ e^{P_0(t)} \ dt \cr F_1(z)&=\int_0^z t \
e^{P_0(t)} \ dt  \cr \ldots \cr F_{d-1}(z)&=\int_0^z t^{d-1} \
e^{P_0(t)} \ dt \ . \cr}
$$
These $d$ special functions are algebraically independent over the
field of rational functions $\dd C(z)$. They define a
Piccard-Vessiot extension of the simplest kind: It is a Liouville
extension. Coming back to the classical motivation of mathematicians
from the XIXth century, we prove that these special functions are
exactly the new transcendentals needed in order to be able to
compute all integrals of the form
$$
\int Q e^{P_0}  \ ,
$$
where $Q$ is an arbitrary polynomial. These integrals form the
vector space,
$$
V_{P_0}=z\dd C[z] e^{P_0} \oplus \dd C \oplus \dd C F_0 \oplus
\ldots \oplus \dd C F_{d-1} \ .
$$
Observe that this vector space is finite dimensional modulo known
functions (polynomials and exponential). The ring $A_0$ generated by
this vector space and its field of fractions $K_0$ are the
fundamental objects for building Dedekind-Weber theory. More
precisely, composing these functions with the inverse of the
uniformization $F_0$ they define functions on $\cl S$. These
functions enjoy the remarkable property of having Stolz limits at
infinite ramification points, thus they are indeed remarkably well
defined in $\cl S^*$. A reason corroborating that this is the right
choice of ring of functions is provided by a Liouville theorem on
$\cl S$ that we discovered. The functions in the vector space
$V_{P_0}$ can be characterized by their growth at infinite in the
log-Riemann surface $\cl S$. The meaning of "growth at infinite" has
to be made precise in $\cl S$. One can escape to infinite in
different ways, one being the classical one in the plane. The other
two being the convergence to or spiraling around infinite
ramification points. Once this is understood, we can establish the
estimates for the functions in $V_{P_0}$ and prove that conversely
any function holomorphic in $\cl S^*$ fulfilling these estimates is
a function in $V_{P_0}$.

\medskip

After establishing these results we proceed to identify
algebraically points of the log-Riemann surface $\cl S^*$ from the
ring $A_0$ in Dedekind-Weber style. We prove that distinct points
determine distinct maximal ideals of $A_0$. This is straightforward
except for the separation of infinite ramification points which is a
non-trivial result. It is based on the non-vanishing of a
determinant: The ramificant determinant. Normalize the polynomial
$P_0$ as
$$
P_0(z)=-{1\over d} z^d+a_{d-1} z^{d-1}+a_{d-2}z^{d-2}+\ldots +a_1
z+a_0 \ .
$$
Consider the $d$ $d$-roots of unity
$$
\omega_1=1 ,  \omega_2=e^{2\pi
i\over d}, \ldots  \omega_d=e^{2\pi i(d-1)\over d} .
$$
The ramificant determinant is
$$
\Delta (a_0,a_1,\ldots ,a_{d-1})=\left | \matrix{F_0(+\infty .\o_1)
& F_1(+\infty .\o_1) & \ldots &F_{d-1}(+\infty .\o_1) \cr
F_0(+\infty .\o_2) &F_1(+\infty .\o_2) & \ldots &F_{d-1}(+\infty
.\o_2^*) \cr \vdots &\vdots &\ddots & \vdots \cr F_0(+\infty .\o_d)
& F_1(+\infty .\o_d) & \ldots &F_{d-1}(+\infty .\o_d) \cr} \right |
$$
or more explicitly
$$
\Delta (a_0,a_1,\ldots ,a_{d-1})=\left | \matrix{\int_0^{+\infty .
\omega_1} e^{P_0(z)}\ dz & \int_0^{+\infty .\omega_1} z e^{P_0(z)}\
dz & \ldots &\int_0^{+\infty .\omega_1} z^{d-1} e^{P_0(z)} \ dz \cr
\int_0^{+\infty .\omega_2} e^{P_0(z)} \ dz &\int_0^{+\infty
.\omega_2} z e^{P_0(z)} \ dz & \ldots &\int_0^{+\infty .\omega_2}
z^{d-1} e^{P_0(z)} \ dz \cr \vdots &\vdots &\ddots & \vdots \cr
\int_0^{+\infty .\omega_d} e^{P_0(z)} \ dz & \int_0^{+\infty
.\omega_d} z e^{P_0(z)} & \ldots &\int_0^{+\infty .\omega_d}
z^{d-1}e^{P_0(z)} \cr} \right | \ \ .
$$
The non-vanishing of the ramificant is equivalent to the separation
of the ramification points by the ring $A_0$. More than the
non-vanishing, we have the remarkable result that the ramificant can
actually be explicitly computed (even if the entries cannot!). We
have the remarkable formula:
$$
\Delta (a_0,a_1,\ldots , a_d)= {1 \over \sqrt {\pi}} \left ( \pi
d\right )^{d\over 2}  \ e^{\Pi_d (a_0 ,a_1 ,\ldots , a_{d-1}) } \ \
,
$$
where $\Pi_d$ is a universal polynomial with positive rational
entries. For example we compute:
$$
\eqalign{ \Pi_1 (X_0)&=X_0 \ ,\cr \Pi_2 (X_0,X_1) &= 2X_0+{1\over 2}
X_1^2 \ ,\cr \Pi_3 (X_0,X_1,X_2) &= 3X_0 + 2 X_1 X_2 +{4\over 3}
X_2^3 \ ,\cr \Pi_4(X_0,X_1,X_2,X_3)&=4X_0 +3X_3 X_1+ 2
X_2^2+9X_3^2X_2+\ldots \cr}
$$

\medskip

The computation of the ramificant relies on analytic tools. The
first observation is that $\Delta$ is an entire function in
$(a_0,a_1,\ldots , a_{d-1})\in \dd C^d$. The second observation is
that it satisfies a system of linear PDE's of the form
$$
\partial_{a_k} \Delta= c_k \Delta \ ,
$$
where $c_k$ is a polynomial on the $a_k$'s. These two facts prove
that the vanishing of $\Delta$ one point implies that $\Delta$ is
identically $0$. Also it proves that
$$
\Delta = C.e^{\Pi_d (a_0,\ldots , a_{d-1})} \ ,
$$
where $C$ is independent of the $a_k$'s. We only need to show that
$C\not=0$. But we can compute explicitly $\Delta (0,0,\ldots , 0)$
(essentially the only place where we know to do that!) by the
explicit computations for cyclotomic Riemann surfaces from section
II.6 which gives a non-zero Vandermonde determinant.

\medskip
 This result gives other significant corollaries. For instance,
if we assume $a_0=P_0(0)=0$, then the polynomial $P_0$ is uniquely
determined by the asymptotic values $(F_j(+\infty . \o_k))$ for
$j=0,1,\ldots d-1$ and $k=1,\ldots ,d$. The coefficients of $P_0$
are universal polynomial functions on these asymptotic values. Also
the locus ramification mapping
$$
\Upsilon: \dd C^d \to \dd C^d
$$
defined by
$$
\Upsilon (a_0,a_1 ,\ldots , a_{d-1})=(F_0(+\infty
.\o_1^*),F_0(+\infty .\o_2^*),\ldots , F_0(+\infty .\o_d^*)) \ ,
$$
is a local diffeomorphism everywhere.

\medskip

Finally we distinguish algebraically finite from infinite
ramification points. For a point $w_0 \in \cl S^*$, let ${\hbox
{\gothic M}} = {\hbox {\gothic M}}_{w_0}$ be the associated maximal
ideal in the ring $A_0$.

Let $A_{0,{\hbox {\gothic M}}}$ be the localization of $A_0$ at the
maximal ideal ${\hbox {\gothic M}}$, and let $\widehat{{\hbox
{\gothic M}}} \subset A_{0,{\hbox {\gothic M}}}$ be the image of
${\hbox {\gothic M}}$ in $A_{0,{\hbox {\gothic M}}}$. Then $w_0$ is
an infinite ramification point if and only if $\widehat{{\hbox
{\gothic M}}} / \widehat{{\hbox {\gothic M}}}^2$ is an infinite
dimensional $\dd C$-vector space.

\bigskip

\stit {Acknowledgements.}

We are grateful to the numerous people with whom we had discussed
the topics developed on this article. In particular with the group
of Mathematicians from Salamanca that have patiently listen to
earlier versions of this manuscript during 2004: R. Alonso Blanco,
A. \'Alvarez V\'azquez, D. Bl\'azquez Sanz, P. Cutillas Ripoll, S.
Jim\'enez Verdugo, J. Lombardero, J. Mu\~noz D\'\i az. We are
grateful to D. Barsky for his comments and his interest, and to F.
Marcell\'an Espa\~nol for helpful discussions.

\medskip

We would like to thank the Institut des Hautes \'Etudes
Scientifiques (I.H.E.S.) and its Director J.-P. Bourguignon for its
support and hospitality. This work was started there during a visit
of both authors in 2002. We thank also the Institute for Mathematical
Sciences at Stony Brook and J. Milnor  for its support and
hospitality that enable both authors to meet there in February 2006
and work towards the final version of this work.

\medskip

{\bf Addendum (December 2015).}

This manuscript was finished in March 2007 and posted in ArXiv in December 2015. 
There are now 3 other related papers posted in ArXiv, two of them 
published:

\medskip

\noindent {\it Caratheodory convergence of log-Riemann surfaces and Euler's formula}, ArXiv:1011.0535, 
Contemporary Mathematics, Volume 639, 2015, pg 197-203.

\noindent DOI: http://dx.doi.org/10.1090/conm/639

\medskip

\noindent {\it Uniformization of simply connected finite type log-Riemann surfaces }, ArXiv:1011.0812, 
Contemporary Mathematics, Volume 639, 2015, pg 205-216.

\noindent DOI: http://dx.doi.org/10.1090/conm/639

\medskip

\noindent {\it Uniformization of higher genus finite type log-Riemann surfaces}, ArXiv:1305.2339.

\bigskip

The first two papers contain results from this manuscript while the third contains new material.

\null
\vfill
\eject

%% file: i1.tex
\input defv7.tex
\bigskip

{\centerline {\tir I. Geometric theory of log-Riemann surfaces.}}

\bigskip

\stit {I.1) Definition of log-Riemann surfaces.}

\stit {I.1.1) Definition.}

\eno {Definition I.1.1.1}{A cut $\g$ with base point $w\in \dd C$ is
a path homeomorphic to $[0,+\infty [$ starting at $w$ and
tending to $\infty$. A straight cut is a cut which is a metric
half line in $\dd C$.}

\eno { Definition I.1.1.2 (log-Riemann surface).}
{The surface $\cl S$ is a log-Riemann surface if we
have:

\item {$(1)$} $\cl S $ is a Riemann surface.

\item {$(2)$} $\cl S$ is equipped with an atlas $\cl A=
\{ (U_i, \varphi_i )\}$ where $\varphi_i: U_i\to \dd C$ are
charts such that
$$
\varphi_i (U_i)=\dd C-\G_i
$$
where $\G_i$ is a discrete, i.e. locally finite, union of
disjoint straight cuts with base points
forming a discrete set $F_i$. We call such charts log-charts.

\item {$(3)$} For each point $z$ in a cut $\g_i$, not an endpoint, the map
$\varphi^{-1}$ extends to a local holomorphic diffeomorphism into the
surface. We have two extensions, one from each side, that we assume do not
coincide.

\item {$(4)$} The changes of charts in the atlas are
the identity
$$
\varphi_{ij}=\varphi_i\circ \varphi_j^{-1} =\id \ .
$$

\medskip

We do identify log-Riemann surface structures for which there
is a homeomorphism from the underlying surfaces that is the
identity on charts.
}

\medskip

\stit {Observations.}

\medskip

{\bf 1.} Condition (3) ensures that the "cuts" do not belong to the
geometry of the surface, i.e. there is no "boundary" at these cuts.

\medskip

{\bf 2.} In condition (3), we do want distinct extensions, otherwise we could just remove the cut, keeping the
endpoint, in order to get a chart into a slit pointed plane.

\medskip

{\bf 3.} We can define log-Riemann surface structures using non-straight
cuts. This introduces technical difficulties when for example the cuts
spiral. We will show later that this more general definition is equivalent
to the one given here, i.e. we can always find charts with straight cuts.
In condition (3) we have to be careful to use Jordan
theorem to define the two sides of the cuts.

\medskip

{\bf 4.} Staying within the same homotopy class for the cuts does not
change the log-Riemann surface structure.

\medskip

{\bf 5.} Riemann surfaces are assumed to be connected.

\medskip

{\bf 6.} The most general definition of log-Riemann surfaces would allow
cuts that are homeomorphic to segments with two finite end-points. As
we see in the example 5 below this is the classical view of algebraic
curves. In fact as we will prove in section I.2.3.3 our definition
covers the general case (assuming as we will that the ramification set
defined in section I.2 is discrete.)

\medskip

{\bf 7.} In condition (2) the discretness of the endpoints does not
ensure the non-accumulation of cuts onto a point of another cut. This
is the reason why we must add that the cuts themselves form a discrete set,
that is, any ball of finite radius intersects a finite number of cuts.

\medskip

\eno {Definition I.1.1.3 (Affine class).}{Two log-Riemann surfaces $\cl S_1$
and $\cl S_2$ are affine equivalent if there exists a holomorphic
diffeomorphism $\varphi : \cl S_2 \to \cl S_1$ and an affine automorphism
$l:\dd C \to \dd C$ such that $\varphi$ on all log-charts is equal to $l$.

The affine class of a log-Riemann surface $\cl S$ is the set of all
log-Riemann surfaces that are affine equivalent.
}

\eno {Definition I.1.1.4 (Projection mapping).}{
The change of charts being the identity,
there is a well defined map $\pi : \cl S \to \dd C$
given by the charts called the projection mapping.

The fiber of (or above) a point $z\in \dd C$ is the discrete set
$\pi^{-1} (z)\subset \cl S$.}

The projection mapping $\pi$ is a local holomorphic
diffeomorphism. It can be used as
a canonical coordinate for the log-Riemann surface structure.

\medskip

Given a log-Riemann surface $\cl S$ and an automorphism $l$ of
$\dd C$, $l(z)=az+b$ with $a\in \dd C-\{0 \}$, $b\in C$, we denote
by $a\cl S +b$ the log-Riemann surface affine equivalent to $\cl
S$ by $l$, i.e. there is a complex diffeomorphism $\varphi: \cl S
\to a\cl S+b$ such that we have the commutative diagram
$$
\pi_{a\cl S+b} \circ \varphi =l\circ \pi_{\cl S} \ .
$$

%% Commutative diagram here

\stit {I.1.2) Examples.}

\medskip

\stit {1. Planes glued together.}

Given a collection of slit planes with their euclidean structure, if the slits can be pasted isometrically by the
identity, we get a Riemann surface with a canonical log-Riemann surface structure inherited from the surgery.
Conversely, any log-Riemann surface structure can be realized in that way. We refer to such a structure of slit
pasted planes as "the log-Riemann surface". Notice that it is a Riemann surface with a set of distinguished
charts.

\stit {2. The complex plane as log-Riemann surface.}

The identity  map ${\hbox {\rm id}} : \dd C \to \dd C$ defines a canonical
log-Riemann surface structure with only one chart, or without cuts.
We denote it by
$\dd C_{\hbox {\rm id}}$.
Any other log-Riemann surface structure on $\dd C$ with only one
chart is given by an affine automorphism $l:\dd C\to \dd C$. We denote
this structure by $\dd C_l$. Observe that if $l_1 \not= l_2$ then the
log-Riemann surface structures $\dd C_{l_1}$ and $\dd C_{l_2}$
are not equivalent but are affine equivalent.

\bigskip

\stit {3. Log-surface of $\root n \of{z}$.}

For $k=0, \ldots n-1$ let $$ \eqalign { U_k &= \{ z\in \dd C^* ; {2\pi k\over n} < \Arg z < {2\pi (k+1) \over n}
\} \ ,\cr \tilde U_k &=  \{ z\in \dd C^* ; {2\pi k\over n}+{\pi \over n} < \Arg z < {2\pi (k+1) \over n}+{\pi
\over n} \} \ .\cr } $$ Let $\g =[0, +\infty[$ and $\tilde \g=]-\infty, 0]$ and $$ \eqalign { &\varphi_k: U_k \to
\dd C-\g \cr &\tilde \varphi_k: \tilde U_k \to \dd C-\tilde \g \cr } $$ defined by $\varphi_k (z)=z^n$ and $\tilde
\varphi_k(z)=z^n$. The atlas $\{(U_k, \varphi_k), (\tilde U_k, \tilde \varphi_k)\}$ defines a log-Riemann surface
structure on $\dd C^*$, denoted by $\cl S_n$. This log-Riemann surface structure is equivalent to the one given by
$n$ complex planes slit and pasted together along $[0, +\infty[$. We visualize $\cl S_n$ in that way. These planes
correspond to the domain of definition of the charts $\varphi_k$. Note that this representation is equivalent to
the one given by $n$ complex planes slit and pasted together along $]-\infty , 0]$. Observe that there is a well
defined $n$-th root $\root n \of{\  }: \cl S_n \to \dd C^*$ that satisfies for $w\in \cl S_n$ $$ \left (\root n
\of{w} \right )^n =\pi (w) \ . $$

\bigskip
\bigskip

{\hfill {\centerline {\psfig {figure=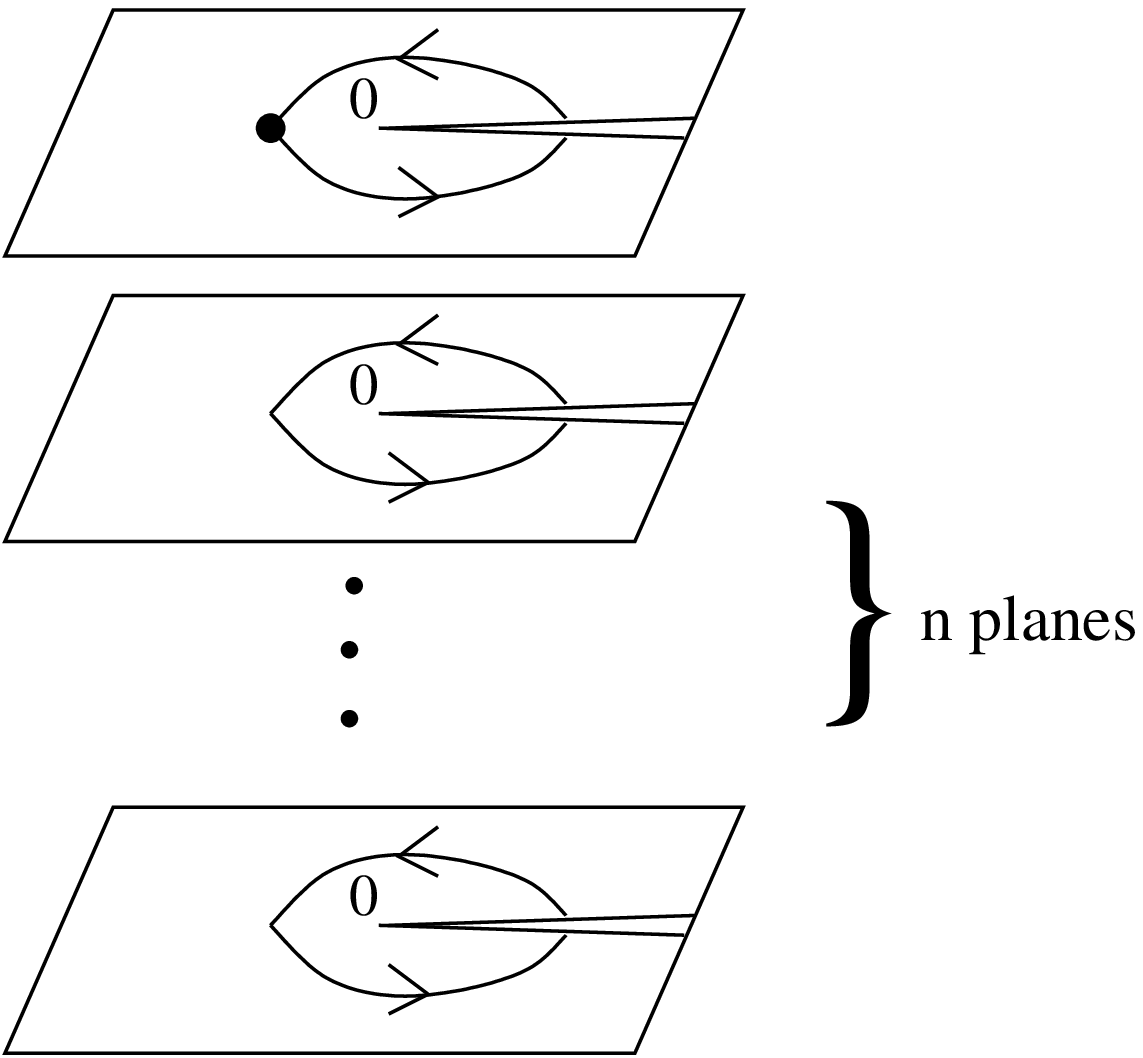,height=5cm}}}}

\bigskip

{\centerline {\bf Figure I.1.1} }

\bigskip
\bigskip

Let $\cl S'_n$ be the log-Riemann surface structure defined on
$\dd C-\{-n\}$ using the atlas $\{ (U_{n,k}, \varphi_{n,k}),
(\tilde U_{n,k}, \tilde \varphi_{n,k})\}$ where
$$
\eqalign {
U_{n,k} &= \{ z\in \dd C-\{-n\} ;
{2\pi k\over n} < \Arg (z+n) < {2\pi (k+1) \over n} \} \ ,\cr
\tilde U_{n,k} &=  \{ z\in \dd C-\{ -n\} ; {2\pi k\over n}+{\pi \over n}
< \Arg (z+n) < {2\pi (k+1) \over n}+{\pi \over n} \} \ .\cr
}
$$
and $\varphi_{n,k}=(1+z/n)^n$ and $\tilde
\varphi_{n,k}=(1+z/n)^n$. Then $\cl S'_n$ is affine equivalent (by
the affine map $l(z)=1+z/n$) to the previous log-Riemann surface
structure:
$$
\cl S'_n = 1+{1\over n} \cl S_n \ .
$$

\medskip

\stit {4. Log-Riemann surfaces associated to polynomials.}

Given a polynomial $Q_0(z)\in \dd C[z]$ we can construct a
pointed log-Riemann surface $(\cl S, z_0)$, $\pi(z_0)=0$,
such that we have a holomorphic diffeomorphism
$$
\apl {F_0}{\dd C-Q_0^{-1}(0)} {\cl S} {z}
{F_0(z)}
$$
such that $F_0(0)=z_0$ and
$$
\pi \circ F_0 =\int_0^z Q_0(t)\ dt
$$
is the polynomial integral of $Q_0$.
This generalizes the previous example where $Q_0(z)=n z^{n-1}$.

The log-Riemann surface $\cl S$ is bi-holomorphic to the complex plane
minus a finite set and can be built with a finite number of log-charts
and a finite number of cuts.

\medskip

Conversely, any such pointed log-Riemann surface has a uniformization
$F_0 : \dd C-\{z_1,\ldots, z_k\} \to \cl S$ which in log-charts
is a polynomial $z\mapsto \int_0^z Q_0(t) \ dt$ and $z_1, \ldots
, z_k$ are the zeros of the polynomial $Q_0$.

\medskip

Observe that in this way we get a correspondence between an algebraic structure,
the ring of polynomials $\dd C[z]$, and a geometric structure, the space of
such log-Riemann surfaces. Elementary algebraic operations on the algebraic
side are not simple operations on the geometric counter-part, and vice versa,
simple geometric surgeries do not correspond to simple algebraic operations.
This philosophy is a guiding idea of the whole theory.

\medskip

\stit {5. Algebraic curves over $\dd C$: Algebraic log-Riemann surfaces.}

This is a classical example that generalizes the previous one. If we glue together a finite number of planes with
a finite number of slits we obtain an algebraic curve spread over $\dd C$. It is bi-holomorphic to a compact
Riemann surface (not necessarily the Riemann sphere as before) minus a finite number of points (those at $\infty$
and those corresponding to finite ramification points). If the projection map $\pi : \cl S \to \dd C$ is taken as
the $z$-variable then the surface can be identified to an algebraic curve over $\dd C$ given by an algebraic
equation $$ P(w,z)=0 \ , $$ where $P$ is a polynomial.

These log-Riemann surfaces are named algebraic log-Riemann surfaces. Classically algebraic curves are defined by
using log-Riemann surfaces constructed with segment cuts with disjoint end-points as well (see remark 6 above.)
Indeed we will prove that the class of algebraic log-Riemann surfaces constructed with only infinite cuts gives
all algebraic curves (this appears missing in the classical literature).

A particular case occurs when we only use two plane sheets.
These are called hyper-elliptic log-Riemann surfaces and are
hyper-elliptic curves whose equation is of the form
$$
w^2=P(z)\ .
$$

This log-Riemann surface description is the classical point of
view that Mathematicians had of algebraic curves in the XIXth
century after B. Riemann's celebrated memoir on Abelian Integrals
([Ri], [Ab]). The abstract modern definition of Riemann surfaces
is due to H. Weyl and T. Rad\'o ([We]). Classical references on
the XIXth century theory of algebraic curves are [Ap-Go], [Jo],
[Pi] or [Va]. For a historical survey see [Ho].

\stit {6. Belyi log-Riemann surfaces.}

We consider a log-Riemann surface build up with a finite number of
sheets with only possible slits  $]-\infty, 0]$ and $[1,
+\infty[$. This defines a Riemann surface with projection mapping
branched only over $0$, $1$ and $\infty$. We define this to be a
Belyi log-Riemann surface. The associated algebraic curve is a
compact Riemann surface defined over the algebraic closure
${\overline {\dd Q}}$ of $\dd Q$. Conversely, Belyi's theorem
states that we get in this way all compact Riemann surfaces
defined over ${\overline {\dd Q}}$ (any ramified cover of a
compact Riemann surface over the Riemann sphere ${\overline {\dd
C}}$ branched only over $0$, $1$ and $\infty$ is the projection
mapping for a log-Riemann surface structure with log-charts having
only  $]-\infty, 0]$ and $[1, +\infty[$ as cuts). Thus we can
characterize compact Riemann surfaces define over ${\overline {\dd
Q}}$ as those possessing a Belyi log-Riemann surface structure.
This is equivalent to be able to tile the surface by flat
congruent equilateral triangles isometrically pasted along the
sides. We refer to [Be] and [Bo] p.99.

Note that in this example we have a geometric sub-meccano of the
general geometric meccano yielding algebraic log-Riemann surfaces.
This sub-meccano corresponds to the arithmetic sub-meccano of
algebraic equations defining the algebraic curves with polynomials
with algebraic coefficients.

\stit {7. Log-Riemann surface of the logarithm.}

For $k\in \dd Z$ we define $$ \eqalign { U_k &=\{z\in \dd C; 2\pi k < \Im z < 2\pi (k+1) \} \ ,\cr \tilde U_k &=\{
z\in \dd C; 2\pi k  +\pi< \Im z < 2\pi (k+1) +\pi \} \ .\cr } $$ Let $\g =[0, +\infty[$ and $\g'=]-\infty, 0]$ and
$$ \eqalign { &\varphi_k: U_k \to \dd C-\g \cr &\tilde \varphi_k: \tilde U_k \to \dd C-\g' \cr } $$ be defined
$\varphi_k (z)=e^z$ and $\tilde \varphi_k (z)=e^z$. This defines a log-Riemann surface structure on $\dd C$. This
log-Riemann surface structure is the same as the one obtained by considering a countable number of copies of $\dd
C$ slit along $\g$ and pasted together to form the log-Riemann surface of the logarithm. Note that a logarithm
function $\log: \cl S \to \dd C$  is well defined on $\cl S$ so that for $w\in \cl S$ $$ \exp (\log w)=\pi (w) \ .
$$

\bigskip
\bigskip

{\hfill {\centerline {\psfig {figure=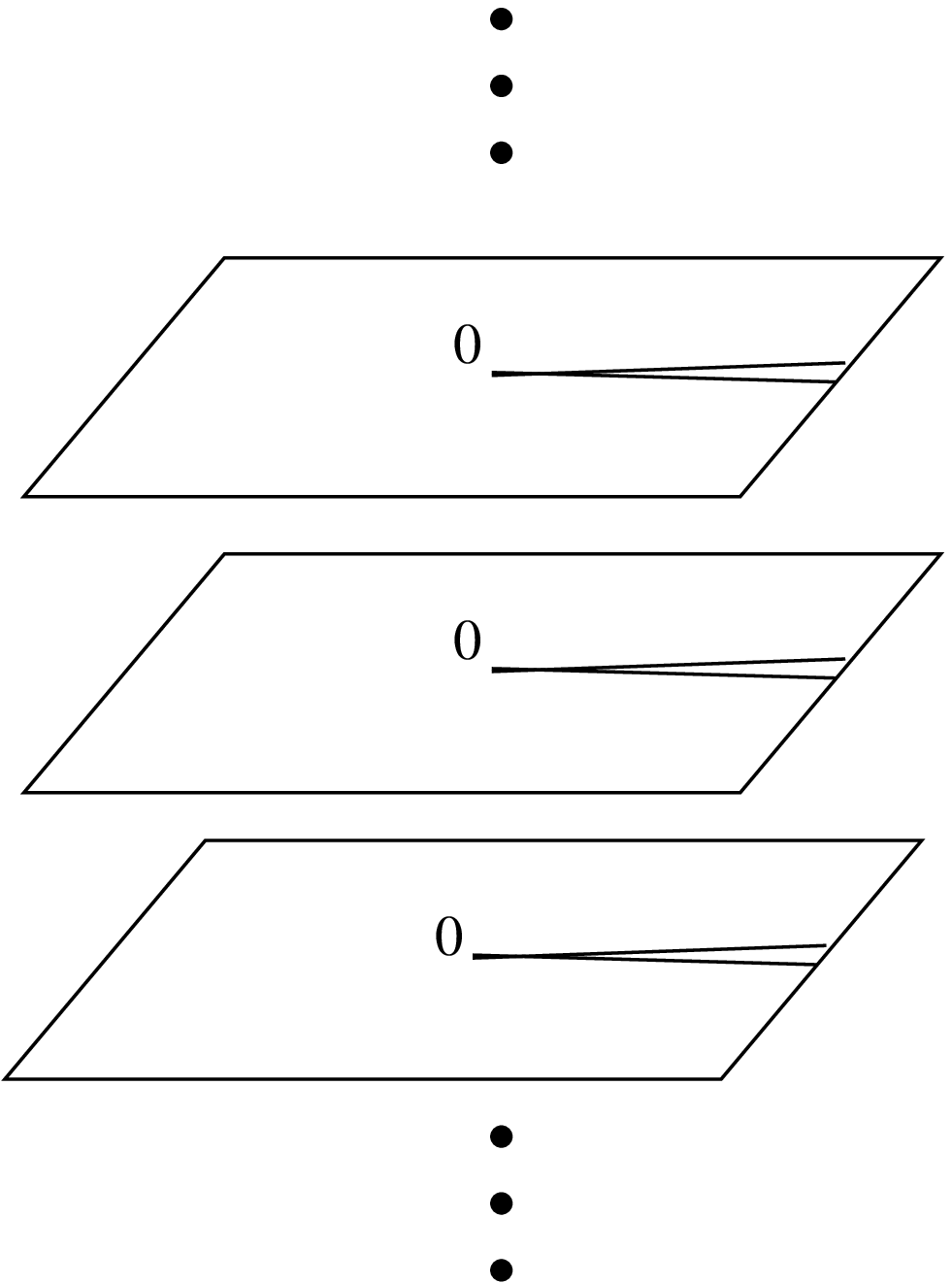,height=5cm}}}}

\bigskip

{\centerline {\bf Figure I.1.2} }

\bigskip
\bigskip

Note that in some sense this log-Riemann surface is the limit of the log-Riemann surfaces $\cl S_n$ of example 3
when $n\to +\infty$. More precisely, we can observe that for a fixed $k\in \dd Z$, when $n\to +\infty$, $U_{n,k}
\to U_k$ in Caratheodory kernel topology (choosing $i(2\pi k+\pi)$ as base point for example) as well as $\tilde
U_{n,k}\to \tilde U_k$. Moreover we also have $$ \eqalign { \varphi_{n,k} &\to \varphi_k \cr \tilde \varphi_{n,k}
&\to \tilde \varphi_k \cr } $$ uniformly on compact sets of $U_k$ and $\tilde U_k$ since $$ \lim_{n\to +\infty}
\left ( 1+ {z\over n} \right )^n =e^z\ . $$ Notice that the charts $\varphi_{n,k}$ (and $\tilde \varphi_{n,k}$)
are uniformly normalized such that $$ \eqalign { \varphi_{n,k} (0) &=1  \ , \cr \varphi'_{n,k} (0) &=1 \ , \cr }
$$ thus the general theory of univalent functions (see [Du] for example) shows that they form a normal family on
the kernel of their domain of definition. In section II.3 we define and discuss the notion of convergence of
log-Riemann surfaces.

\stit {8. Gauss log-Riemann surface.}

In this example we just describe the construction of the log-Riemann surface by pasting slit planes. We consider
the cuts $\g=[1, +\infty[$ and $\g'=]-\infty , -1]$. To the slit plane $\dd C-(\g\cup \g')$ we paste a countable
family of copies of $\dd C-\g$ and another countable family of copies of $\dd C-\g'$. We graft these families on
each cut in the same way that we do for the log-Riemann surface of the logarithm. This defines the Gauss
log-Riemann surface. The reason for this terminology is that, as we prove later, this Riemann surface is
bi-holomorphic to $\dd C$ and the Gauss integral $$ z\mapsto {2\over \sqrt{\pi}}\int_0^z e^{-t^2}\ dt $$ defines a
uniformization from $\dd C$ into this log-Riemann surface.

\bigskip
\bigskip

{\hfill {\centerline {\psfig {figure=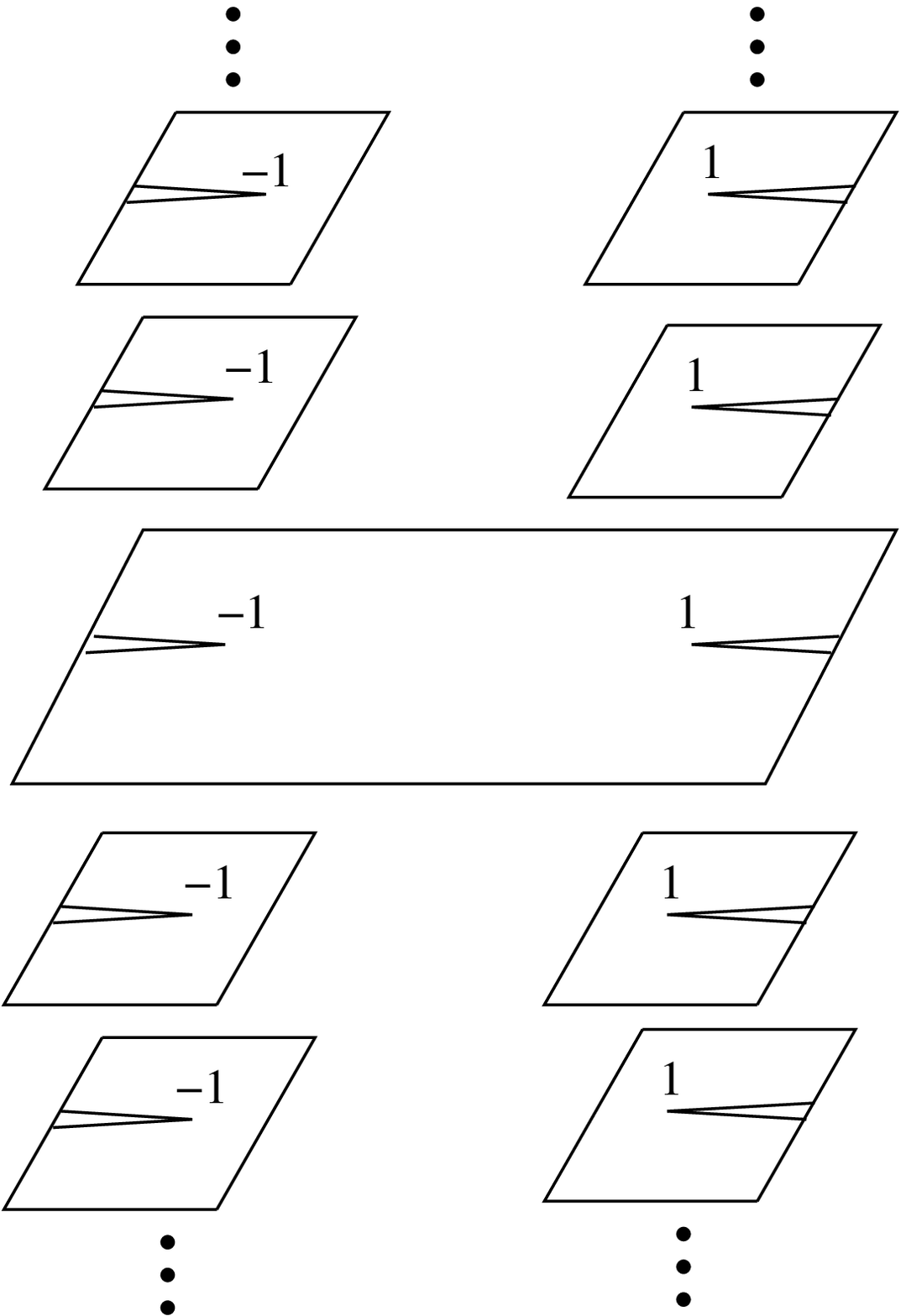,height=7cm}}}}

\bigskip

{\centerline {\bf Figure I.1.3} }

\bigskip
\bigskip

There is a natural generalization of this example. Let $d\geq 2$ and consider a base sheet $\dd C$ with $d$ radial
cuts with end-points at the $d$th-roots of unity. We paste on these cuts distinct families of planes as in the
construction of the surface of the logarithm. We call this log-Riemann surface the Gauss log-Riemann surface of
log-degree $d$. We will show in section II.6 that this log-Riemann surface is bi-holomorphic to the complex plane
and the integral, $$ z\mapsto {d \over \G (1/d)} \ \int_0^{z} e^{-t^d} \ dt \ $$ defines the uniformization from
$\dd C$ to this log-Riemann surface mapping $0$ into $0$ of the base sheet.

\medskip

R. Nevanlinna uses these entire functions as first examples of
entire functions with $d$ exceptional values in the sense of
Nevanlinna theory ([Ne3] p.20 and p.90). The exceptional values
are the $d$-th roots of unity, whose fiber contains the end-points
of the cuts (i.e. these points have an abnormal fiber above them.)

\bigskip

\stit {9. Log-Riemann surfaces associated to entire functions.}

Let $\cl S$ be a log-Riemann surface holomorphically embedded into a
simply connected parabolic Riemann surface $\cl S^\times$ (that is bi-holomorphic
to the complex plane $\dd C$) such that $\cl S^\times -\cl S$
is discrete (once we define ramification points in section I.2, we can say
that $\cl S$ is simply connected and parabolic once we add the finite ramification
points.) We consider the uniformization
$$
F_0 : \dd C \to \cl S^\times \ .
$$
Then we get an entire function $F=\pi \circ F_0$. Note that $\cl
S$ is what is classically called the Riemann surface of the
multivalued inverse function $F^{-1}$ (note that this gives a
log-Riemann surface structure and not just a Riemann surface
structure despite the slightly confusing classical terminology).

\medskip

Conversely, given an arbitrary entire function $F$ we ask whether
we can associate to it a log-Riemann surface $\cl S$. At a
non-critical value image point we can choose an inverse branch of
$F$ in its neighborhood. We can then build the Riemann surface of
this germ of univalent function. This Riemann surface $\cl S$
comes equipped with a canonical chart $\pi : \cl S \to \dd C$ such
that $F$ lifts into a biholomorphic map $F_0: \dd C \to \cl S$
such that $F=\pi\circ F_0$ (see [Ma] volume II, chapter VIII,
section 5 p. 502-540). In general $\cl S$ is not endowed with a
log-Riemann surface structure as defined in section I.1. It is not
always possible to fulfill the requirement to have a locally
finite cuts in each log-chart (according to condition (2) of the
definition.) For example, as when we have a chart with parallel
cuts converging to a cut (indeed a half cut). As we will see
elsewhere [Bi-PM] there are entire functions which give rise to
such structures that are not equivalent to a log-Riemann surface
structure with locally finite cuts. But a large and natural class
of entire functions have log-Riemann surfaces associated to them.

\stit {10. Modular log-Riemann surface.}

Consider a countable family of copies of $\dd C-(]-\infty , 0]
\cup [1, +\infty [)$. We start with one copy and we paste four
distinct copies, two in each slit. Next we paste $12$ distinct
copies in the free slit boundaries. Next $36$ distinct copies in
the free slits, and so on. In such a way we build the modular
log-Riemann surface. It is simply connected and bi-holomorphic
to the unit disk. The classical modular function $\lambda$ (see
for example [Ah1] p. 281) is a uniformization from the upper
half plane into this modular log-Riemann surface.

\stit {11. Polylogarithm log-Riemann surface.}

We consider a complex plane slit along $[1, +\infty [$, $\dd C-[1, +\infty [$. We paste two copies of slit planes
$\dd C-(]-\infty , 0]\cup [1, +\infty [)$. In the $6=3\times 2$ remaining slits we paste copies of the slit plane
$\dd C-(]-\infty , 0]\cup [1, +\infty [)$. We keep pasting this slit plane on the remaining free slits and so on.
All the sheets of these log-Riemann surface are the same except the first one. Polylogarithm functions are defined
by its holomorphic germ at $0$, $$ {\hbox {\rm Li}}_k (z)=\sum_{n=1}^{+\infty } {z^n \over n^k} \ . $$ All
polylogarithm functions ${\hbox {\rm Li}}_k$ extend holomorphically and are well defined (i.e. single valued) on
the above log-Riemann surface (see for example [Oe]). We call this log-Riemann surface the polylogarithm
log-Riemann surface.

\stit {12. Billiard log-Riemann surfaces.}

There is a classical construction that associates to a polygonal billiard dynamical system (see for example [Bi])
a log-Riemann surface. We start with the original polygon and we attach all possible reflections across the
boundary segments. We continue reflecting the new copies. This generates a log-Riemann surface. This construction
is usually done for billiards with rational angles (i.e. commensurable with $\pi$.) At the vertices, after a
finite number of reflections, the last reflected polygon is glued to the first one. This construction generates
algebraic curves. A particular case occurs for a rectangular or equilateral triangle billiard that gives the
log-complex plane. In the general case, for incommensurable angles, we obtain log-Riemann surfaces with infinitely
many sheets.

\medskip

%% file: i2.tex
\input defv7.tex
\stit {I.2) Euclidean metric and ramification points.}

\stit {I.2.1) Definition.}

\eno {Definition I.2.1.1}{
Pulling back the Euclidean metric
on $\dd C$ by the projection mapping $\pi$ we get a flat conformal
metric on $\cl S$.
We call this metric the Euclidean metric.
}

\eno {Definition I.2.1.2}{Associated to the Euclidean metric we have
a metric space by defining a distance as, for $z_1, z_2 \in \cl S$,
$$
d(z_1, z_2)=\inf_{z_1, z_2 \in \g} l(\g ) \ ,
$$
where the infimum runs over all rectifiable paths containing the two
points and $l(\g )$ denotes the Euclidean length of $\g$.
}

The projection mapping $\pi$ is a local isometry, a global
contraction (not strict),
and an open map.
The Euclidean metric space is never complete when we have charts
with cuts (just construct a Cauchy sequence that converges on the
chart to the end-point of a cut.)

\stit {I.2.2) Ramification points.}

\eno {Definition I.2.2.1 (Ramification set).}{Let $\cl S^*=\cl S\cup \cl R$ be
the completion of $\cl S$ in the Euclidean metric. The set $\cl R$
is closed and is called the ramification set.
The projection mapping $\pi$ extends continuously uniquely to $\cl R$.
We keep the same notation $\pi$ for the extension.
}

\eno {Definition I.2.2.2 (Ramification point).}{ A ramification point is an isolated point of $\cl R$.

The log-Riemann surface $\cl S$ is called finite or transalgebraic if $\cl R$ is a finite set. } The terminology
"structurally finite" instead of "finite" is used by M. Taniguchi ([Ta1], [Ta2], [Ta3].)

Note that "ramification point" means isolated point in $\cl R$ and not
point in $\cl R$.

\eno {Lemma I.2.2.3}{For any ramification point $z^*\in \cl R$ there is a ball
$B(z^*,r)$ centered  at $z^*$ with no other points of $\cl R$
such that $\pi (B(z^*,r))=B(\pi(z^*),r)$,
$\pi (B(z^*,r)-\{z^*\})=B(\pi(z^*),r)-\{\pi(z^*)\}$ is a pointed disk,
and $\pi: B(z^*,r)-\{z^*\}\to \pi (B(z^*,r)-\{z^*\})$ is a
covering.}

\stit {Proof.}
Choose $r>0$ small enough so that there are no points of $\cl R$
 in $B(z^*,r)$ apart from $z^*$. Pick a $z \in B(z^*,r) - \{z^*\}$;
then $\pi(z) \in B(\pi(z^*),r) - \{\pi(z^*)\}$, and the local inverse of $\pi$ at $\pi(z)$ satisfying
$\pi^{-1}(\pi(z)) = z$ can be analytically continued to all points of $B(\pi(z^*),r) - \{\pi(z^*)\}$, since the
only possible obstruction to the continuation is encountering points in $\cl R$ above, but by the choice of $r$
this is not possible. This shows that $\pi$ maps $B(z^*,r) - \{z^*\}$ onto $B(\pi(z^*),r) - \{\pi(z^*)\}$.

The proof that $\pi: B(z^*,r)-\{z^*\}\to \pi (B(z^*,r)-\{z^*\})$
is a covering is similar. For each point
$z_0 \in \pi (B(z^*,r)-\{z^*\})=B(\pi(z^*),r)-\{\pi(z^*)\}$,
we choose $\rho >0$ such that $B(z_0,\rho) \subset B(\pi(z^*), r)-\{
\pi(z^*)\}$. Let $U$ be a connected component of the preimage
$\pi^{-1}(B(z_0, \rho))$; we can pick a $z_1 \in U$ and as before
continue without obstruction the local inverse of $\pi$ satisfying
$\pi^{-1}(\pi(z_1)) = z_1$ to all of the disk $B(z_0,\rho)$. Since
this disk is simply connected, the continuation of $\pi^{-1}$ to it
is single-valued, so $\pi_{|U} : U \to B(z_0,\rho)$ has an inverse
and is therefore a diffeomorphism. $\diamond$

\medskip

\eno {Corollary I.2.2.4}{The set of ramification points is at most countable.}

\stit {Examples.}

{\bf 1.} The ramification set $\cl R$ can be uncountable.
Consider the hierarchy of
end-points of segments generating the triadic Cantor set
$$
\eqalign {
F_0 &=\{ 0, 1 \} \ , \cr
F_1 &=\{ 1/3, 2/3 \} \ , \cr
F_2 &=\{ 1/9, 2/9, 7/9, 8/9 \} \ ,\cr
\vdots & \cr
}
$$
We consider a copy of $\dd C$ with two vertical cuts going to $-i\infty$
with set of end-points $F_0$. We paste a single plane sheet at these
cuts, and on each of these new plane sheets we consider vertical cuts
going to $-i\infty$ with set of end-points $F_1$. On each cut we paste a
single plane and we make cuts on each one with set of end-points $F_2$, etc.

In this way we construct a log-Riemann surface with a countable number of
ramification points, and with an uncountable ramification set $\cl R$ that
projects by $\pi$ onto the triadic Cantor set.

\medskip

{\bf 2.}
A simple modification of the previous construction yields a log-Riemann
surface with $\pi (\cl R ) =\dd C$ (just take a sequence of finite sets
$F_n$ such that $\bigcup_n F_n $ is dense on $\dd C$.)

\medskip

{\bf 3.} It is possible for $\cl R$ to be perfect thus no
ramification point exists. We modify the previous example by
pasting not one but a countable number of planes above and below
each slit (as we do in the construction of the log-Riemann surface
of the logarithm, example 7 of section I.2.) The log-Riemann
surface thus constructed has a perfect ramification set $\cl R$
and moreover $\pi (\cl R )=\dd C$.

\bigskip

Using the previous lemma we define the degree or
order of a ramification point.

\eno {Definition I.2.2.5}{For each ramification point $z^*$ for a small disk $B(\pi (z^*), r)$, the connected
component $U$ of $\pi^{-1} (B(\pi(z^*), r))$ containing $z^*$ satisfies that $U-\{z^*\}$ is simply connected or
bi-holomorphic to a pointed disk. The degree $1\leq n=n(z^*)\leq +\infty$ of the covering $$ \pi : U-\{z^*\} \to
B(\pi (z^*),r)-\{\pi (z^*) \} $$ is the degree or order of the ramification point $z^*$.

If the ramification point $z^*$ has finite order we say that
$z^*$ is a finite ramification point. Then the Riemann
surface structure (but not the log-Riemann surface structure)
of $\cl S$ can be extended to $\cl S\cup \{z^* \}$
and $\pi$ is a ramified covering at this point.

The log-degree of a log-Riemann surface is the number of infinite
ramification points.
}

\eno {Definition I.2.2.6}{Let $\cl S$ be a log-Riemann surface. The Riemann
surface obtained by adding the finite ramification points of $\cl S$
and extending the Riemann surface structure to them is called the
finitely completed Riemann surface of $\cl S$ and denoted by
$\cl S^\times$.}

\stit {Examples.}

{\bf 1.} The number and order of the ramification points only depend
on the affine class.

{\bf 2.} We refer to the examples given in section I.2.
Example 2 has no ramification points. It is easy to prove the converse.

\eno {Proposition I.2.2.7}{A log-Riemann surface with no ramification
points is a planar log-Riemann surface $\dd C_l$.}

{\bf 3.} Examples 3 and 7 have a unique ramification point of order $n$ and
$+\infty$ respectively. It is also easy to prove the converse.

\eno {Theorem I.2.2.8}{A log-Riemann surface with only one ramification point is in
the affine class of $\cl S_n$, the log-Riemann surface of $\root n \of z$, or
of the log-Riemann surface of the logarithm, $\cl S_{\log}$.}

The log-Riemann surface $\cl S_n$ has log-degree $0$ and $\cl S_{\log }$
has log-degree $1$.

{\bf 4.} Log-Riemann surfaces associated to polynomials as in example 4
have the property that the Riemann surfaces obtained by adding the finite
number of finite ramification points are simply connected and parabolic.
The converse also holds.

\eno {Theorem I.2.2.9}{Let $\cl S$ be a log-Riemann surface with a finite ramification
set with all ramification points of finite order such that the finitely completed
Riemann surface $\cl S^\times$ is simply connected. Then it is
parabolic (i.e. bi-holomorphic to $\dd C$) and the uniformization map
$$
F_0: \dd C \to \cl S^\times
$$
is such that $\pi \circ F_0: \dd C \to \dd C$ is a polynomial map.
}

{\bf 5.} Algebraic log-Riemann surfaces in example 5 have a finite ramification
set, all ramification points are of finite order. The converse also holds.

\eno {Theorem I.2.2.10}{A log-Riemann surface with a finite ramification set and
with all ramification points of finite order is an algebraic log-Riemann surface.}

{\bf 6.} Belyi log-Riemann surfaces as defined in example 6 have a
finite ramification set, all ramification points are of finite order and project
only onto $0$, $1$. The converse also holds. Note that if we render projective
invariant the definition of ramification points then we can talk of ramification
points over $\infty$.

\medskip

{\bf 7.} The Gauss log-Riemann surface has a ramification set composed of
two ramification points of
infinite order. Any other simply connected log-Riemann surface with
this property is in the affine class of the Gauss log-Riemann
surface.

\medskip

{\bf 8.} The Modular log-Riemann surface has an infinite number of ramification
points all of infinite order projecting only onto $0$ and $1$.

\medskip

{\bf 9.} Billiard log-Riemann surfaces associated to polygonal billiards with at least three sides and mutually
incommensurable angles with $\pi$ give examples of log-Riemann surfaces with a countable ramification set composed
by an infinite number of infinite ramification points with a dense projection on $\dd C$.

\bigskip

{\it From now on and for the rest of the article we only consider log-Riemann
surfaces with a discrete ramification set $\cl R$. Thus all points of $\cl R$
are ramification points.}

\bigskip

The Euclidean metric on the log-Riemann surface $\cl S$ essentially
characterizes the log-Riemann surface structure on $\cl S$ when $\cl S^\times$
is simply connected. More precisely we have:

\eno {Theorem I.2.2.11}{ Let $\cl S$ be a Riemann surface endowed with
a flat conformal metric. We can define as before the ramification
set $\cl R$ as the points added to $\cl S$ in the completion for
this metric. We assume that $\cl R$ is a discrete set. We define
finite ramification points as those for which the Riemann surface
structure of $\cl S$ extends to them. Then $\cl S^\times$ is well
defined and we assume that $\cl S^\times $ is simply connected.
We assume that the metric at the finite
ramification points is of the form $|dw|=|z|^n \ |dz|$ for some $n\geq 1$.
Then there is a unique, up to translation and rotation (thus in the same
affine class),
log-Riemann surface structure compatible with the given metric.}

\stit {Proof.}

Pick a point and a local isometric chart at this point mapping the point to $0
\in \dd C$. This defines
a germ of holomorphic diffeomorphism $\pi$ at this point. We extend
$\pi$ by holomorphic continuation to all of $\cl S^\times$.
This is possible because there is no monodromy in a neighborhood of the finite
ramification points since we assume that the metric has the given normal form
at these points, and also globally because $\cl S^\times$ is simply connected.
We can build log-charts using $\pi^{-1}$. For each $z_0 \in \cl S$ we choose
the branch of $\pi^{-1}$  such that $\pi^{-1}\circ \pi (z_0)=z_0$. We extend
$\pi^{-1}$ radially from $\pi(z_0)$.
Each time we encounter
a ramification point we draw a radial cut starting at its image by $\pi$.
This procedure builds a log-chart on a plane. Such charts define
a log-Riemann surface structure.$\diamond$

\stit {Observations.}

{\bf 1.} If the metric does not have the stated form at the finite ramification points then the continuation of
$\pi$ may have a monodromy locally around a ramification point and the construction of a global mapping $\pi$ is
impossible. An example of an inadmissible conformal metric would be $|dw|=|z|^\a \ |dz|$ where $\a >-1$ and is not
an integer.

\medskip

{\bf 2.} The assumption that $\cl S^\times $ is simply connected is not
superfluous. For example $\cl S =\dd C /\dd Z$ carries a flat metric inherited
from the euclidean metric on $\dd C$ by isometric quotient but has no log-chart
compatible with this metric. Indeed $\cl S=\dd C/\dd Z$ is complete, thus
$\cl R=\emptyset$ and a log-chart can not have a cut.

A more general structure that allows these cylindrical ends are
the tube-log Riemann surface structure that we will study
elsewhere [Bi-PM2].

\medskip

{\bf 3.} Also when $\cl S^\times$ is not simply connected and without
tubular ends, $\pi$ can still have a non-trivial monodromy. For example,
consider two planes $\dd C-(]-\infty , 0]\cup [1, +\infty[)$ and
$\dd C-(]-\infty, 0]\cup [2, +\infty[)$. We paste the left cuts by the identity
and the right cuts by $z\to z+1$. The Riemann surface obtained inherits
a flat metric. Continuing $\pi$ along the loop in the figure gives a monodromy
$+1$ for $\pi$.

%% Figure
\bigskip
\bigskip

{\hfill {\centerline {\psfig {figure=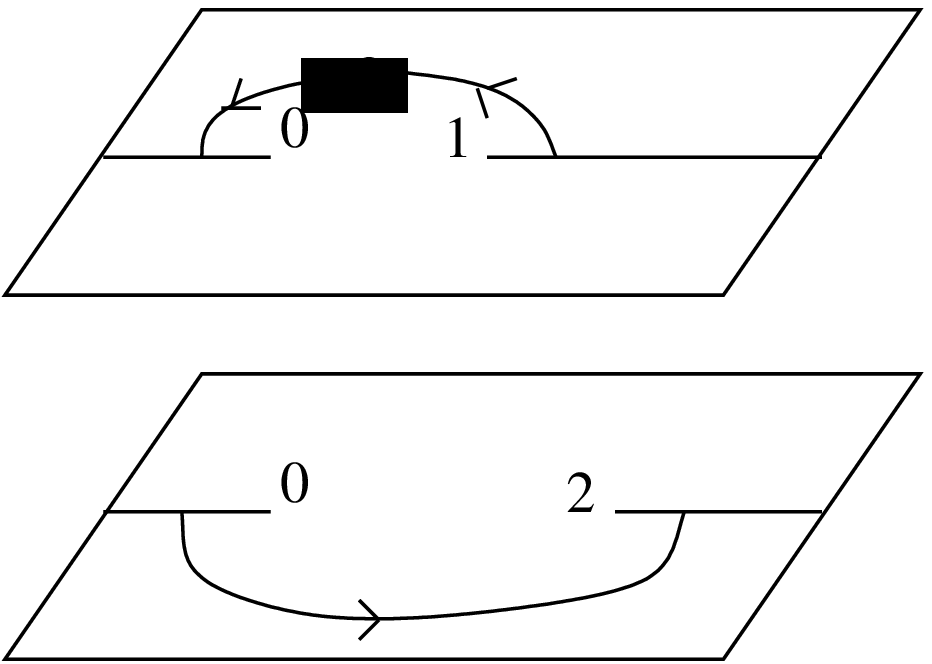,height=3cm}}}}

\bigskip

{\centerline {\bf Figure I.2.1} }

\bigskip

\stit {I.2.3) Log-euclidean geometry and charts.}

\stit {I.2.3.1) Preliminaries.}

In order to do some geometry on $\cl S$ it is convenient to complete
the space into $\cl S^*$. The space $\cl S^*$ endowed with the euclidean
metric is a path metric space in the sense of Gromov (see [Gr].)
Many properties of smooth riemannian geometry do extend to our setting.

A geodesic is a rectifiable path $\g$ such that for any two points close
enough in $\g$, the distance between these two points is the length of the segment
in $\g$ that joins them.

The space is geodesically complete in the sense that every geodesic can be continued
indefinitely, and geodesics starting at any point cover a full
neighborhood of the point. Indeed we have a replacement for the classical
exponential mapping at a point $m\in M$,
$$
\exp_m^t: T^1_m M \to M
$$
such that $(\exp_m^t)_{0\leq t\leq t_0(m)}$ is a geodesic segment
parametrized by length. We have to replace the unit tangent bundle
$T^1 M$ by a covering of degree $n$ (finite or infinite)
$T^1_{z^*} \cl S^* \to \dd R/\dd Z$ at any ramification point
$z^*\in  \cl S^* -\cl S$ of order $n$.

Even if our space is not compact or locally compact, we see below
that some of the conclusions of the Hopf-Rinow theorem hold: Any
two points can be joined by a minimal geodesic. Some of the other
properties of geodesically complete spaces do not hold due to the
non-smoothness of the metric at the ramification points. Thus for
instance it is not true that bounded sets do have compact closure.
We refer to [Mi1] p. 62, [Pe] chapter 7 (or [Ho-Ri]) for these
questions.

Observe that closed and bounded sets are not necessarily compact if they contain a neighborhood of an infinite
ramification point. But even worse: We can have a closed bounded set with no ramification point that is not
compact as for example a spiraling strip centered around the ramification point on the log-Riemann surface of the
logarithm.

\medskip

\eno {Theorem I.2.3.1}{Given a log-Riemann surface $\cl S^*$ and any compact
set $K\subset \cl S^*$ there exist $\eps=\eps (K)>0$ depending on $K$
such that for any two points $z_1, z_2 \in K$ with $d(z_1, z_2) <\eps$,
there exists a unique geodesic segment joining them.}

\eno {Lemma I.2.3.2}{Let $z^*\in \cl S^*$ be a ramification point and consider $r>0$ such that $B(z^*, r)\cap \cl
R =\{ z^* \}$. Then for $z_1, z_2 \in B(z^*, r/2)$ there exists a unique geodesic segment joining the two points.
This segment is either an euclidean segment $[\pi(z_1) , \pi(z_2)]$ in a log-chart or composed of two euclidean
segments $[\pi (z_1), \pi (z^*)]$ and $[\pi(z^*), \pi(z_2)]$ in two log-charts.}

\stit {Proof of the lemma.}

If there exists a log-chart containing $z_1$ and $z_2$ we have two
possibilities. First, the euclidean segment $[\pi(z_1), \pi(z_2)]$
lifts into a segment entirely contained in the log-chart. Then we
are in the first situation. Second, this does not hold and then it
is elementary to prove that the geodesic joining the two points is
composed by two euclidean segments contained on the log-chart
having end-points at $z_1$, $z_2$ and $z^*$.

If no log-chart contains both points, then they cannot be seen from
$z^*$ through an angle less than $2\pi$. In that case the union of the euclidean
segment from $z_1$ to $z^*$, $[z_1, z^*]\subset \cl S^*$, and
then $z^*$ to $z_2$, $[z^*, z_2]\subset \cl S^*$ is the shortest path joining
the two points. To prove this we can choose two disjoint log-charts containing
respectively $z_1$ and $z_2$ such that the part of any other path joining $z_1$ and $z_2$
in each log-chart is strictly longer than the segment $[z_1, z^*]$ or $[z^*, z_2]$
respectively.$\diamond$

\medskip

\stit {Proof of Theorem I.2.3.1.}
 By contradiction, take a sequence of pairs $(z_1^{(n)}, z_2^{(n)}) \in K^2$ with
 $d(z_1^{(n)}, z_2^{(n)}) \to 0$ and having distinct geodesics joining them.
 Extracting a converging subsequence
 we can assume that $z_1^{(n)}\to z_0$ and $z_2^{(n)}\to z_0$, with $z_0\in K$.
 If $z_0 $ is a ramification point, then using the lemma we get a contradiction.
 If $z_0$ is not a ramification point, then it has a neighborhood which is a
 euclidean disk and again we get a contradiction.$\diamond$

\medskip

We can now give a full description of geodesics.

\eno {Theorem I.2.3.3}{Let $\g$ be a geodesic path in $\cl S^*$. Then $\g$ is a polygonal line
made up with euclidean segments belonging to log-charts with vertices at the
ramification points of $\cl S^*$.}

\stit {Proof.}

The result follows from the local description provided by
the previous result.$\diamond$

\medskip

Since we can always construct a polygonal line with vertices at ramification
points joining two arbitrary points (the surface is path connected), we get:

\eno {Theorem I.2.3.4}{Any two points can be joined by a geodesic.}

The existence of minimal geodesics is false. See the counter-example below. A minimizing sequence of paths joining
two points does not need to have a convergent subsequence (the classical argument that uses Ascoli-Arzela theorem
only works in locally compact spaces).

\stit {Counter-example.}

Consider a countable family of complex planes indexed by $\dd Z$ with two horizontal slits $]-\infty , -1]$ and
$[1, +\infty [$, and a third vertical slit $[-i/|n+1| , +i\infty [$ in the $n$-th sheet. We glue together both
horizontal slits simultaneously as for the surface of the logarithm. Over the vertical slits we paste independent
planes. There is one infinite ramification point over $-1$ and another over $+1$. As is easy to see there is no
minimal geodesic joining them, and any minimizing sequence of geodesics escapes to infinite. Note that if we pick
two regular points $-2+i\eps$ and $2+i\eps$, and the central cuts to be $[-i(2+1/|n+1|) , +i\infty [$, then again
there is no minimal geodesic joining these two points.

%% figure
\bigskip
\bigskip

{\hfill {\centerline {\psfig {figure=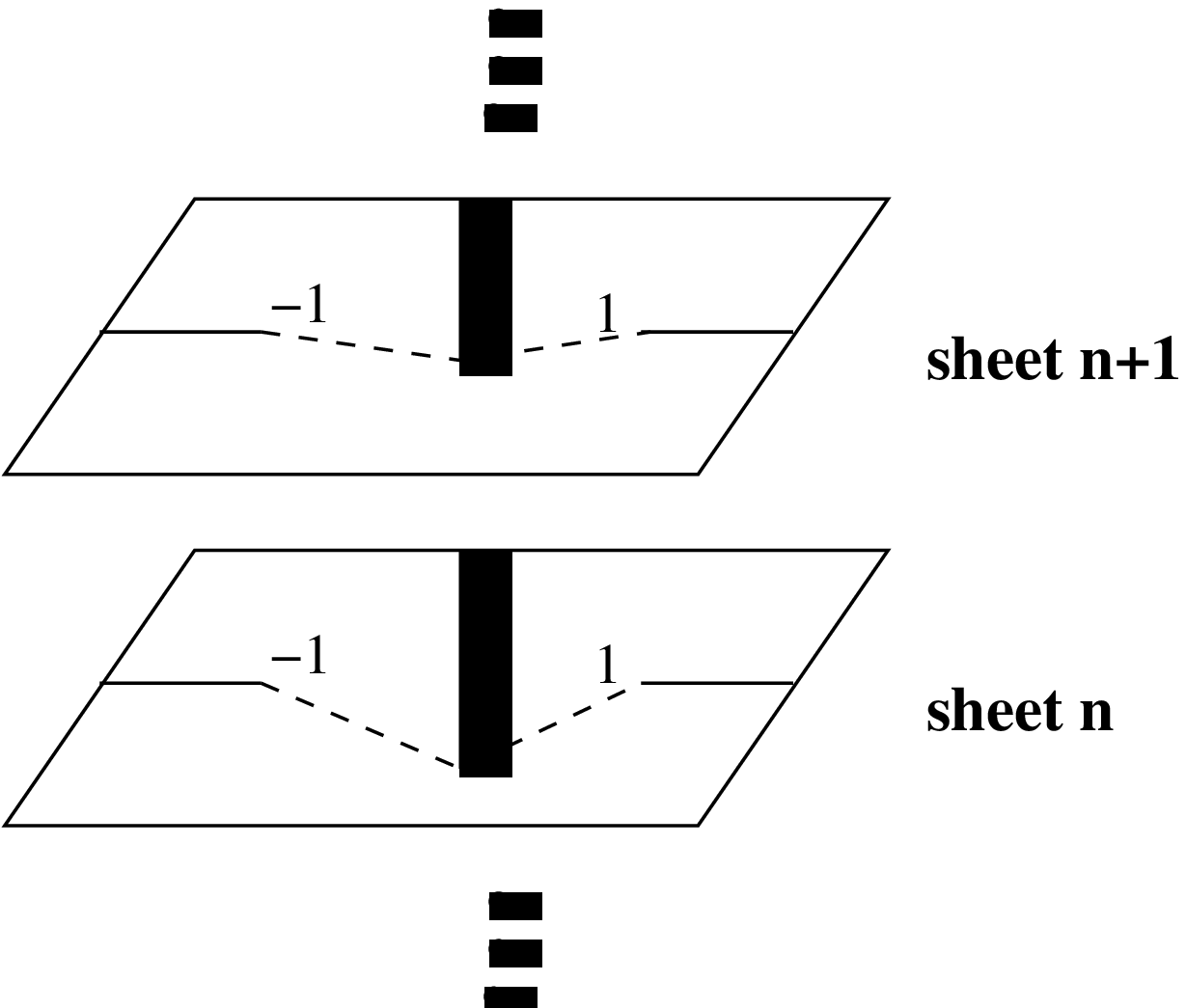,height=5cm}}}}

\bigskip

{\centerline {\bf Figure I.2.2} }

\bigskip
\bigskip

On the other hand, when $\pi(\cl R)$ is finite (in particular when $\cl R$ is
finite) the result holds.

\eno {Theorem I.2.3.5}{Let $\cl S^*$ be a log-Riemann surface with $\pi (\cl R )$
finite. Then for any two points in $\cl S^*$ there exists a minimizing geodesic.}

\stit {Proof.}

Any geodesic $\g$ joining the two points has a projection $\pi (\g )$ which is
a polygonal line of the same length with vertices contained in the finite
set $\pi (\cl R)$. The set of such polygonal lines with a uniform upper bound
on their length is finite. Thus given any minimizing sequence of geodesics their
length would be eventually constant.$\diamond$

\medskip

A minimizing geodesic joining two points is not necessarily unique as the
following counter-example shows.

\stit {Counter-example.}

Consider two complex planes slitted along $]-\infty , -1]$ and $[1, +\infty [$.
We glue them together in order to create two ramification points of order $2$.
There are two minimal geodesics joining them (one in each sheet.)

%% figure
\bigskip
\bigskip

{\hfill {\centerline {\psfig {figure=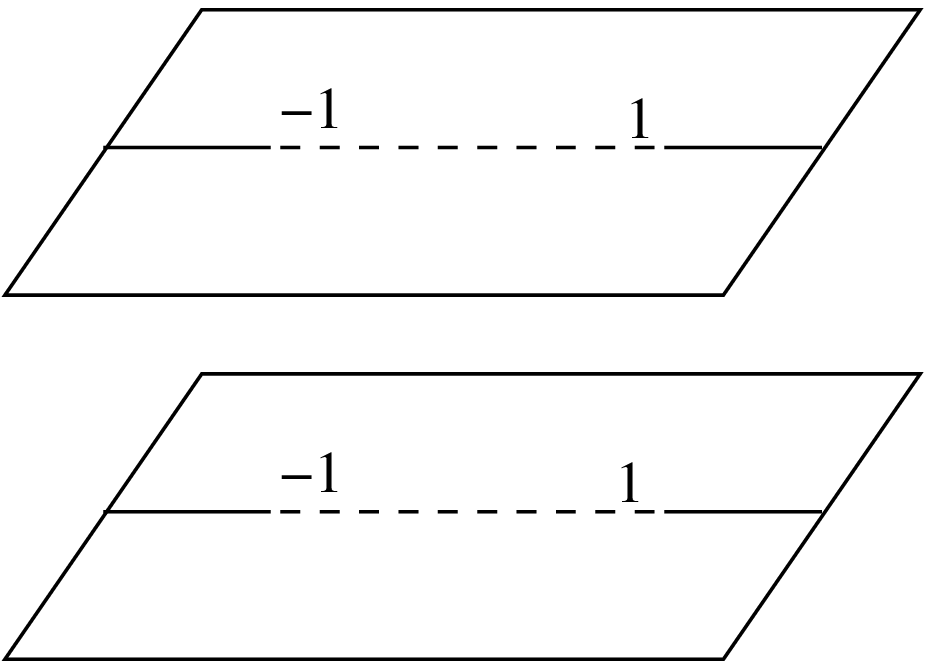,height=3cm}}}}

\medskip

{\centerline {\bf Figure I.2.3} }

\bigskip

\bigskip

A classical result in Riemannian geometry is that on a compact
riemannian manifold each homology class can be realized by
a geodesic. From the previous description of geodesics we can prove
the same result in log-euclidean geometry.

\eno{Proposition I.2.3.6}{Each homology class of $\cl S$ can be realized
by a geodesic. If $\pi (\cl R)$ is finite a minimal geodesic
realizes the homology class.}

{\it We assume that $\pi (\cl R)$ is finite in what follows.}

We can describe circles.

\eno {Theorem I.2.3.7}{A circle centered at a point of $w_0\in \cl S^*$
of positive radius
is composed by arcs of euclidean circles on charts with
centers at $z_0$ and at some ramification points.}

\stit {Proof.} For any point $w$ on this circle, there is a minimal geodesic $\g$ joining this point to the center
$w_0$. Let $w_1^*\in \cl R$ be the first vertex on this geodesic from $w$. Let us consider the two germs of circle
arcs to the left and to the right centered at $w_1^*$, of radius $d(w, w_1^*)$, and starting from $w$.

\eno {Lemma I.2.3.8}{One of these two germs of circle arcs is
part of the circle $C(w_0,r)$ and the circle $C(w_1^*, r-d(w_0, w_1^*))$.}

\stit {Proof of the lemma.}

Consider the angle of $\pi (\g)$ at $\pi(w_1^*)$. If the angle is not flat, let $C'$ be the germ of circle arc
that enters the big sector. Then for $C'$ small and for $w'\in C'$ we can construct a minimal geodesic joining
$w'$ with $w_0$ of length $r$ by just pivoting the segment $[w_1^*, w]$ to $[w_1^*, w']$ and leaving untouched the
rest of the geodesic. When the angle is $\pi$ (resp. $-\pi$), we choose the germ of circle arc by rotating
positively (resp. negatively). $\diamond$

\medskip

\stit {Proof of Theorem (continued).} The same arguments show that the other part is an arc of circle centered at
$w_1^*$, at another ramification point, or at $w_0$.$\diamond$

%picture figure I.2.4

\bigskip
\bigskip

{\hfill {\centerline {\psfig {figure=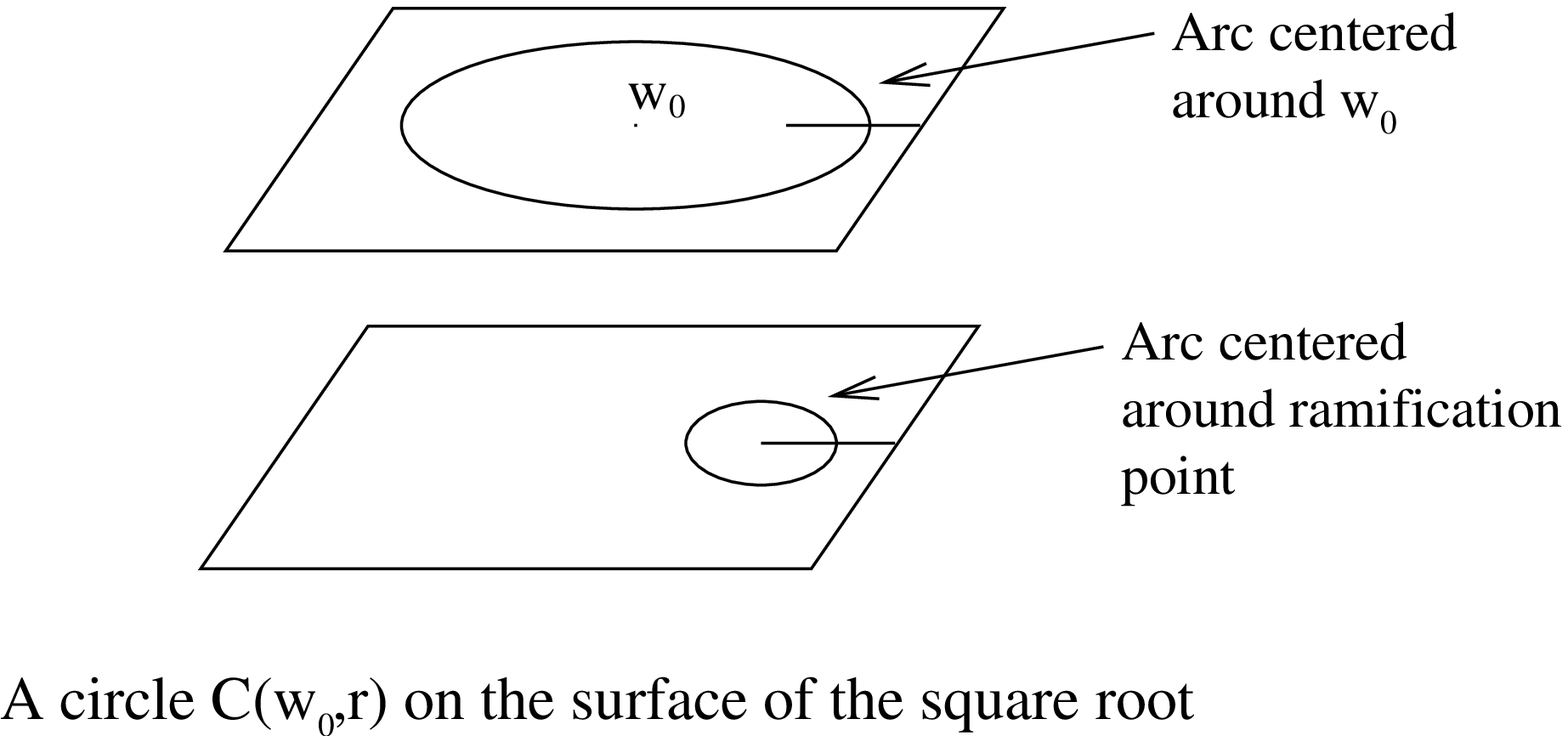,height=5cm}}}}

\bigskip

{\centerline {\bf Figure I.2.4} }

\bigskip

The set of points equidistant to two given points can have an interior.
For example pick two points in the same fiber of the $\cl S_{\log }$ log-Riemann
surface. Using the existence of minimizing geodesics we get the following
description of the boundary of these regions.

\eno {Theorem I.2.3.9}{The locus of points equidistant from two given points
$w_0, w_1 \in \cl S^*$ has a boundary composed by pieces of hyperbolas and
segments.}

We can define angles between geodesics.

\eno {Definition I.2.3.10}{The angle formed by two half geodesics meeting at a regular point $w_0 \in \cl S$ is
the usual euclidean angle (defined modulo $2\pi$.) If the two half geodesics meet at a ramification point $w_0 \in
R$, there is a local argument function that defines the angle: Modulo $2\pi n$ if the ramification point is of
order $1\leq n<+\infty$, or an angle in $\dd R$ if the ramification point is an infinite ramification point.}

We can define convex sets in log-euclidean geometry.

\eno {Definition I.2.3.11}{A subset $U\subset \cl S^*$ is convex if for each pair of
points $w_0, w_1 \in \cl S^*$ all minimizing
geodesics $[w_0, w_1]$ in $\cl S^*$  joining them are
entirely contained in $U$. A point in $w_0\in U$ is extremal if $w_0$ does not
belong to the interior of any geodesic segment contained in $U$.
The convex hull of a set in $\cl S^*$ is the minimal convex set containing it.}

\stit {Examples.}

{\bf 1.} The full space $\cl S^*$ is convex with no extremal points.

{\bf 2.} A single point is convex with one extremal point.

{\bf 3.} The empty set is convex with no extremal points.

\medskip

\eno {Theorem I.2.3.12}{The intersection of convex sets is convex.}

\stit {Proof.}
 If $(C_i)$ is a family of convex sets, and $w_1$ and $w_2$ belong
 to their intersection, then any minimizing geodesic joining them
 is contained in each $C_i$ and hence in their intersection.$\diamond$

\eno {Definition I.2.3.13}{A full geodesic is an unbounded geodesic $\g$ in $\cl S^*$
homeomorphic to $\dd R$ such
that for any two distinct points $w_1, w_2\in \g$ there is a unique
minimizing geodesic joining them and it is the segment
$[w_1,w_2]\subset \g$ defined by these two points.}

\stit {Remark.}

A full geodesic is isometric to the real line, and
its intrinsic metric structure on $\g$ coincides with
the induced metric by the embedding in $\cl S^*$.

\medskip

\eno {Theorem I.2.3.14}{The closure of a component of the complement of a
full geodesic is convex. These are called half-spaces.
}

\stit {Proof.}

Let $w_1$ and $w_2$ be two distinct points in the closure ${\overline  H}$
of a component $H$ of $\cl S^*-\g$. Given a minimizing geodesic $\eta$
joining $w_1$ and $w_2$, if $\eta$ does not intersect $\g$ then it is entirely
contained in ${\overline H}$. If $\eta$ does intersect $\g$, let $w_1^*$ the
first intersection point from $w_1$, and $w_2^*$ the last intersection point.
Let $[w_1^*, w_2^*]_\g\subset \g$ (resp. $[w_1^*, w_2^*]_\eta \subset \eta$)
be the geodesic segment in $\g$ (resp. $\eta$) determined
by $w_1^*$ and $w_2^*$. Since $\g$ is a full geodesic we have that
$$
[w_1^*, w_2^*]_\g=[w_1^*, w_2^*]_\eta
$$
and hence $\eta \subset {\overline H}$.
$\diamond$

\stit {Remark.}

There may be only one half-space as the example of the non-simply connected
elliptic log-Riemann surface below shows (Figure I.2.5).

\bigskip

{\hfill {\centerline {\psfig {figure=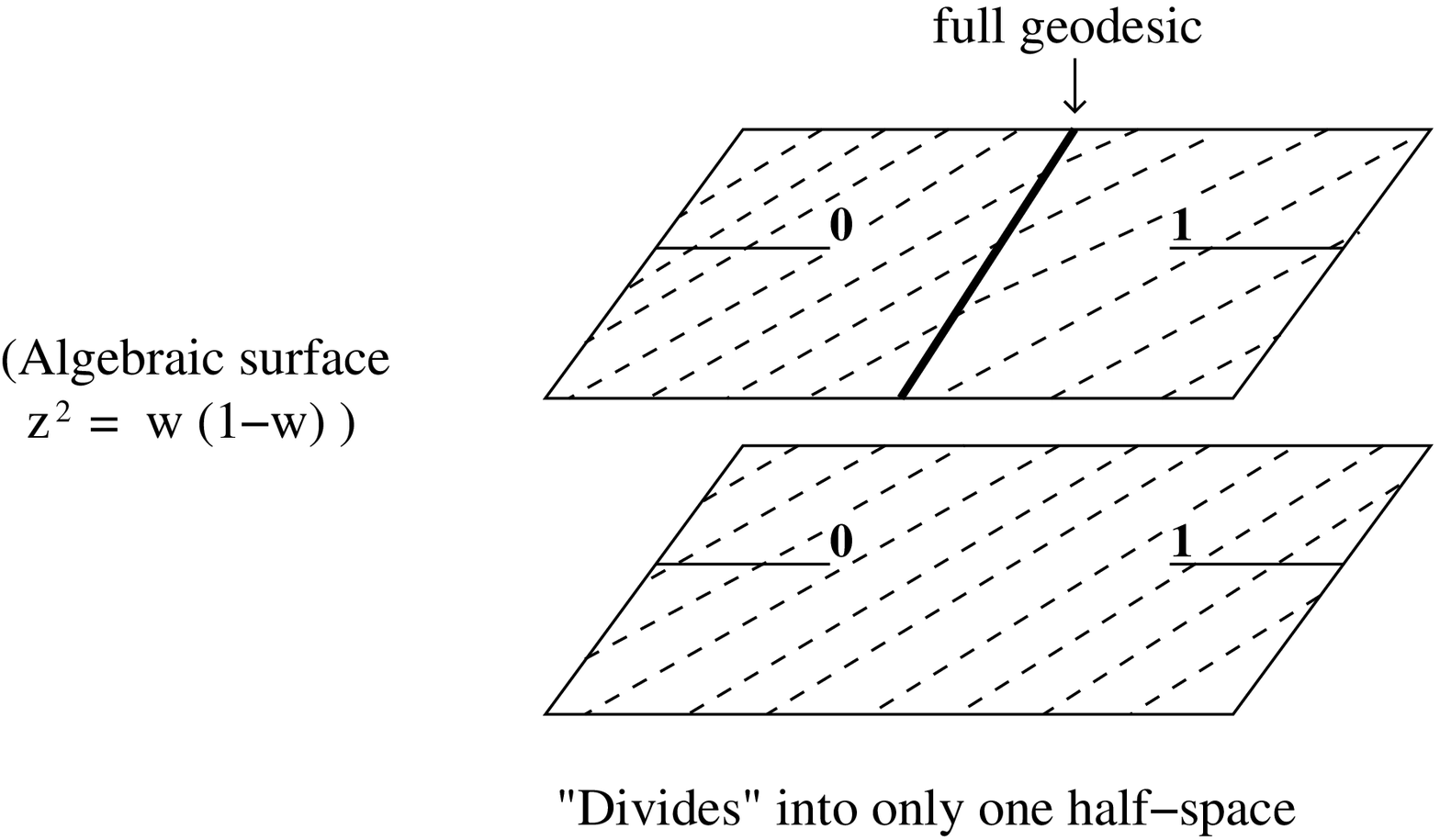,height=5cm}}}}

\bigskip

{\centerline {\bf Figure I.2.5} }

\bigskip

There may be
also more than two half-spaces as a stright line geodesic passing through
the infinite ramification point in $\cl S_{\log }$ shows (three half spaces).

\bigskip
\bigskip

{\hfill {\centerline {\psfig {figure=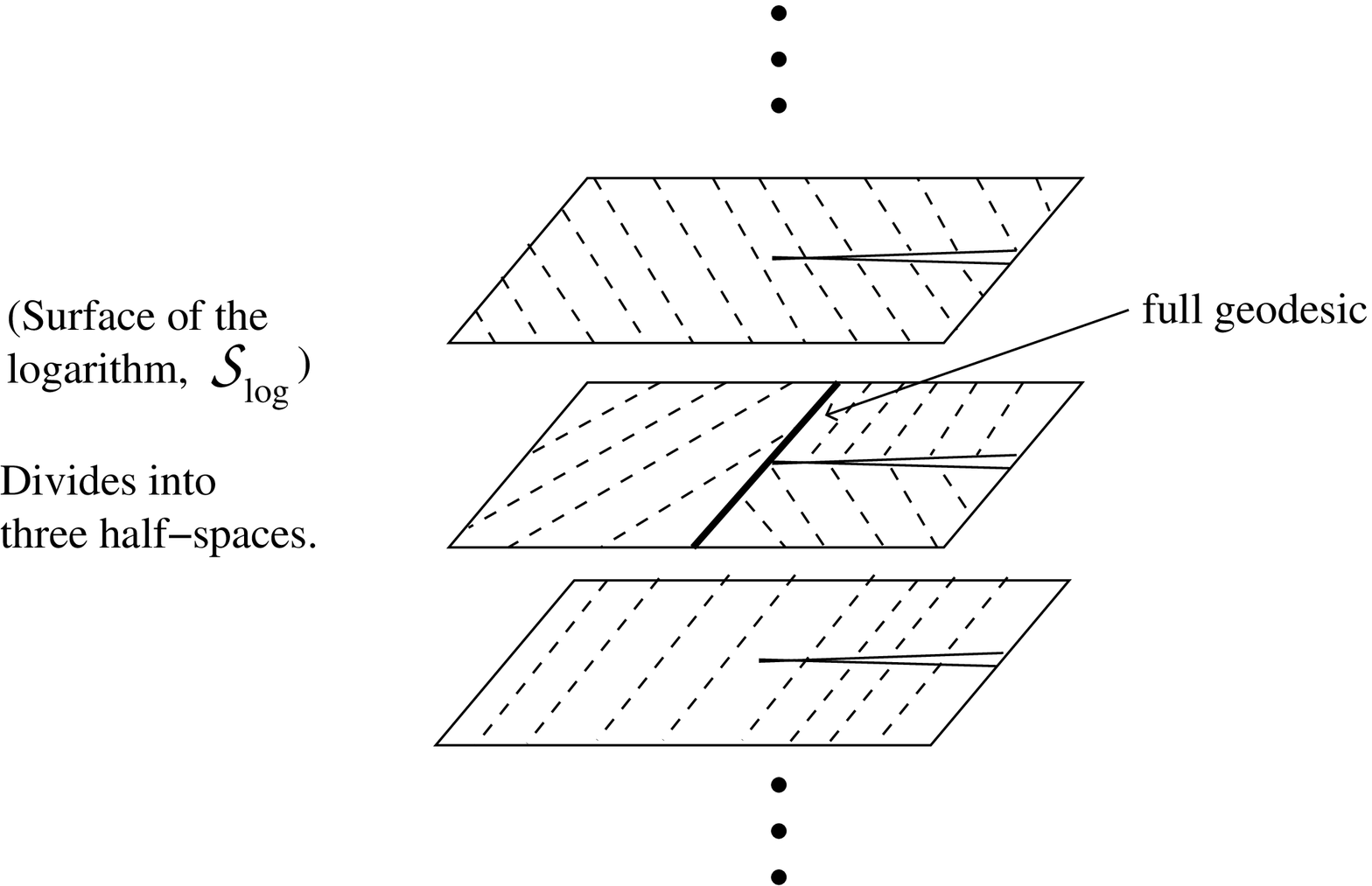,height=5.5cm}}}}

\bigskip

{\centerline {\bf Figure I.2.6} }

\bigskip
\bigskip
\medskip

\eno {Proposition I.2.3.15}{If $\g$ is a full geodesic containing $n\geq 0$
infinite ramification points then there are at most $2+n$ half-spaces
for $\g$.}

\stit {Proof.}
Each point of $\g$ belongs to two or three half spaces. This last possibility
can only happen at infinite ramification points. Hence following a point
along $\g$ we can only meet $2+n$ distinct half-spaces.$\diamond$

The following proposition is clear from the definition of full geodesic.

\eno {Proposition I.2.3.16}{The intersection of two full geodesics is either
empty, a point or a geodesic segment (compact or not).}

\eno {Definition I.2.3.17}{Full geodesics in $\cl S^*$ separate points
from compact sets if for any point in $\cl S^*$ and a compact convex set
$K\subset \cl S^*$ not containing the point there exists a full geodesic
such that the point and the compact set are in distinct half-spaces.}

\eno {Theorem I.2.3.18}{If full geodesics separate points from compact sets then
a compact convex set is the intersection of half-spaces containing it.}

\stit {Proof.}

The intersection of these half-spaces is a closed convex set $C$ containing
the compact set $K$. If $w_0 \notin K$ then there exists a full geodesic
separating $w_0$ from $K$, thus $w_0 \notin C$.$\diamond$

\medskip

The proof of the next theorem is clear.

\eno {Theorem I.2.3.19}{The metric space $\cl S^*$ is locally convex, i.e. any point
has a convex neighborhood. Indeed for any point $w\in \cl S^*$, $B(w, r)$ is
convex for small $r>0$.}

%\item {5.} A compact convex set is the convex hull of its extremal points.

%\eno {Definition }{
%A convex polygon is the convex hull of a finite number of points.
%A triangle is the convex hull of three points (a degenerate triangle
%if the three points
%belong to the same geodesic.)}

\eno {Definition I.2.3.20}{A polygon is an oriented loop (maybe self-intersecting)
formed by a finite
number of geodesic segments.}

Observe that given the orientation, we can talk about the
internal angle
formed at a vertex of the polygon, or at a ramification point on
the polygon.
The sum of angles of a triangle does not add up to $180^o$.
Euclid's axiom of parallels does not hold in log-euclidean geometry.
Nevertheless we can prove.

\eno {Theorem I.2.3.21}{We consider a polygon $\O$ with all internal
angles in $[0,\pi]$. Then the  sum of the
internal angles $(\a_i )$ at
the vertices and at the ramification points
add up to
$$
\sum_i \a_i=\pi (k- 2) \ ,
$$
where $k$ is the number of vertices and ramification points on the polygon.
}

\stit {Proof.}

Since the internal angles are in $[0,\pi]$ then they are
the same as those of the $\pi$ projection of the
polygon. Thus it is enough to prove the result for planar
oriented polygons. This is straightforward by induction: Each time
that we remove one vertex linking the two adjacent ones, the total
sum of angles decreases by $\pi$. Thus at the end we end up with a
triangle and the result holds.$\diamond$

\stit {I.2.3.2) Construction of charts.}

We only assume in this section that $\cl R$ is discrete.
We are going to describe how to construct log-charts just
using the log-euclidean metric and $\pi$.

\eno {Definition I.2.3.22}{Let $w_0 \in \cl S$. The star of $w_0$ is the
union of segment geodesics with one endpoint at $w_0$ and not
meeting ramification points. We denote the star of $w_0$ by $V(w_0)$.}

\eno {Theorem I.2.3.23}{The star $V(w_0)$ of a point $w_0 \in \cl S$ is an
open connected set such that the restriction of $\pi$ to $V(w_0)$
is a holomorphic diffeomorphism and
$$
\pi(V(w_0))=\dd C-\g
$$
where $\g$ is a locally finite union of disjoint straight cuts with base points forming
a locally finite set. Thus $(V(w_0), \pi)=(U_i, \varphi_i)$ is an allowable
chart for the log-Riemann surface structure.}

\medskip

\stit {Proof.}

First observe that it is enough in this case to prove the local finiteness of the endpoints of the cuts. Indeed if
we have a sequence of cuts accumulating at a point $z_1 \in \dd C$, then necessarily the end-points of these cuts
must accumulate some point in the segment $[\pi(w_0), z_1]$ (if a neighborhood is free of them, then a cone
centered around $[\pi(w_0), z_1]$ is free of cuts at finite distance).

Next, an accumulation of end-points of cuts will give after lifting by $\pi$ a
non-isolated point in $\cl R$.$\diamond$

Observe that $\pi (\cl R )\subset \dd C$ is countable.
Choose $z_0\in \dd C-\pi (\cl R )$ and consider its
countable fiber $\pi^{-1}(z_0)=(w_i)$.

\eno {Definition I.2.3.24}{The cell of $w_j$ relative to the fiber $(w_i)$
is the set
$$
U(w_j)=\{ w\in \cl S ; \forall i\ \ d(w,w_j)<d(w, w_i) \} \ .
$$
}

\eno {Theorem I.2.3.25}{Let $(w_i)$ be a fiber of $\pi$.
The cell $U(w_j)$ coincides with the
star of $w_j$,
$$
U(w_j)=V(w_j)
$$
}

\stit {Proof.}

Let $(w_i)=\pi^{-1} (z)$ be this fiber.

We first show that $V(w_j) \subset U(w_j)$.

Consider a point $w$ in the star $V(w_j)$. The point $w_j$
is joined to $w$ by a straight segment $[w_j, w]$ not
containing any ramification points, and $d(w, w_j) = |\pi(w) - z|$.
Thus this segment has an $\epsilon_0$-neighborhood $W(\epsilon_0)$ on
which $\pi$ is univalent. If for some $i$ we have
$d(w, w_i) \leq d(w, w_j)$, then taking $\epsilon << \epsilon_0$ we can find a curve $\gamma \subset \cl S$
joining $w$ to $w_i$ of length less than $|\pi(w) - z| + \epsilon$.
Then $\pi(\gamma)$ has the same length and
hence $\gamma$ must be contained in $W(\epsilon_0)$,
hence $w_i \in W(\epsilon_0)$, a contradiction since
$\pi$ is univalent in $W(\epsilon_0)$.

\medskip

Now we prove the other inclusion.

Let $w$ be a point in the cell $U(w_j)$.
From $\pi(w)$, we lift the segment $[\pi(w), z]$ using $\pi$
to a path $\gamma \subset \cl S^*$. The path $\gamma$ ends
at a point $w_i$ of the fiber, and
$d(w, w_i) = |\pi(w) - z| \leq d(w, w_j)$,
hence since $w$ belongs to the cell $U(w_j)$ we must
have $w_i = w_j$.

Now if $\gamma$ does not meet any ramification points then $w$ belongs to the star $V(w_j)$ and we are done. If it
does meet ramification points, then we can choose at the first ramification point that we meet from $w$ a
different lift. We continue this new lift and at any other ramification points we take care to choose a lift
disjoint from $\gamma$. Then the endpoint will be a point of the fiber distinct from $w_j$, a contradiction.
$\diamond$

\eno {Theorem I.2.3.26}{The cells $(U(w_i))$ form a disjoint collection of open sets and their union is dense in
$\cl S$. Each connected component of the boundary of each cell $U^*(w_i)$ is a geodesic segment but not a full
geodesic, ending in at least one ramification point.}

\medskip

\stit{Proof.} The cells are disjoint by definition.
That their union is dense is true because this holds
for the collection of stars of points in a fiber $\pi^{-1}(z)$.
Consider the at most countable collection $\D_z$ of lines in $\dd C$
passing through $z$ and points of $\pi(\cl R)$. The pre-image
$\pi^{-1}(\dd C - \D_z)$ is dense in $\cl S$ and covered by the
stars of points in the fiber. $\diamond$

The following is a straightforward corollary.

\eno{Corollary I.2.3.27}{Consider the union $\D \subset \dd C$ of
all lines passing through pairs of distinct points of
$\pi(\cl R)$. Since $\cl R$ is countable the set $\D$ has
Lebesgue measure $0$ and the complement of $\D$ is a
$G_\d$-dense. If $z_0\in \dd C-\D$ then the connected components of
the boundaries of
the cells $U(w_i)$ are half lines ending at a ramification
point.}

We can construct an atlas of log-charts using cells.

\eno{Theorem I.2.3.28}{Consider a point $z_1 \in \dd C - \D$. We take a
second point $z_2$ not in $\D_{z_1}$. We can
define $\D_{z_1,z_2}$ as the collection of all lines passing through
pairs of distinct points of $\pi(\cl R)$ and the intersections of
lines in $\D_{z_1}$ and $\D_{z_2}$. Now $\dd C - \D_{z_1,z_2}$ is
a $G_\d$-dense, and picking a point $z_3$ in this set, the stars
of the points of the fibers $\pi^{-1}(z_1), \pi^{-1}(z_2), \pi^{-1}(z_3)$
form an atlas of log-charts for $\cl S$.}

\stit{Proof:} Just observe that $\D_{z_1, z_2}$ is still
an at most countable collection of lines. The conclusion follows
from the choices made. $\diamond$

This Theorem and the previous Corollary
has an important application.

Note that all the objects defined and the
results established up to now depend only on the
euclidean metric. In particular, the definition of ramification
points, geodesics, stars, cells,...
If we consider a log-Riemann surface defined not by straight
cuts but by arbitrary path cuts, we can still define the euclidean
metric and build all the theory {\it mutatis mutandis}.
Now we can use the previous Theorem and the charts given
by the cells of the fibers of the three generic points.

This defines an equivalent log-Riemann surface structure with
only straight cuts. We have proved:

\eno{Theorem I.2.3.29}{Any log-Riemann surface structure defined with path cuts
is equivalent to a log-Riemann surface structure defined with straight
cuts.}

Note that this is not totally trivial. Consider for instance a path in $\dd C$ converging to $0$ and tending to
infinite spiraling. Build the log-Riemann surface using a countable number of copies slit along this path and
pasted as in the construction of the log-Riemann surface of the logarithm. Then this log-Riemann surface is still
the log-Riemann surface of the logarithm.

\stit {I.2.3.3) Other applications.}

Consider a generalized log-Riemann surface structure $\cl S$ defined using the definition stated in section I.1
but allowing finite cuts with two finite end-points, all the other hypotheses being the same. As before we can
define $\pi$ and lift the euclidean metric in order to define the euclidean metric on $\cl S$ and the completion
$\cl S^* =\cl S \cup \cl R$.

\eno {Theorem I.2.3.30}{If $\cl R$ is discrete then the generalized log-Riemann
surface structure is equivalent to a classical log-Riemann structure
defined with only infinite cuts. Namely, there exists a holomorphic
diffeomorphism between the underlying Riemann surfaces which is the
identity on the charts defining the generalized and the classical log-Riemann
surface structure.}

The next Corollary should be a classical result, but
we don't know of a reference for it in the literature.

\eno {Corollary I.2.3.31}{Any algebraic curve over $\dd C$ defined by an algebraic
equation
$$
P(w,z)=0 \ ,
$$
has a classical log-Riemann surface structure defined only using
infinite cuts.
}

\stit {I.2.3.4) The Kobayashi-Nevanlinna net.}

We consider a log-Riemann surface $\cl S$. We define a new cellular decomposition of $\cl S$ that is due to Z.
Kobayashi and R. Nevanlinna (see [Ko] and [Ne2] chapter XII) and is used for purposes of determining the type of
finitely completed log-Riemann surfaces.

\eno {Definition I.2.3.32}{Let $w^* \in \cl S^*-\cl S$. We define the Kobayashi-Nevanlinna cell of $w^*$ as $$
W(w^*)=\{ w\in \cl S ; d(w,w^*) < d(w, \cl R-\{w^* \} ) \} \ . $$ }

 Notice that the star $V(w)$ of a point $w\in \cl S$ is also well defined
when $w\in \cl S^*$ is a ramification point (with the same definition).

\eno {Theorem I.2.3.33}{ We have the following properties of Kobayashi-Nevanlinna cells

\item {$\bullet$} $W(w^*)$ is open and path connected.

\item {$\bullet$} $W(w^*) \subset V(w^*)$, more precisely,
$$
[w,w^*]\subset W(w^*) \ .
$$

\item {$\bullet$} The boundary of $W(w^*)$ is composed of euclidean segments.

\item {$\bullet$} The disjoint union
$$
\bigcup_{w^* \in \cl R} W(w^*)
$$
is a dense open set in $\cl S$.
}

\stit {Proof.}

The Kobayashi-Nevanlinna cell is obviously an open set. We prove that for any $w\in W(w^*)$
we have
$$
[w,w^*]\subset W(w^*) \ .
$$
From $w$ follow the geodesic euclidean segment in the direction of $w^*$. We cannot
hit another ramification point before $w^*$ and the result follows.
This implies path connectedness and that the Kobayashi cell is contained in the
star of $w^*$.

We now study the structure of the cell boundaries. Consider a point $w\in \partial W(w^*) \subset \cl S$. Consider
all ramification points $w_1^*, \ldots , w_n^*$ at minimal distance $r_0>0$ from $w$ (thus $n\geq 2$). The disk
$B(w,r_0)\subset \cl S$ is an euclidean disk. Label $w_1^*, \ldots , w_n^*$ in cyclic order and modulo $n$. Draw
the angular bisectors to the sectors $[w, w_i^*]\cup [w, w_{i+1}^*]$. Since for small $\eps >0$, $B(w, r_0+\eps )$
contains no new ramification point, the local structure at $w$ of $\partial W(w^*)$ is formed by small segments of
these bisectors starting at $w$ (see figure). The generic case corresponds to $n=2$ and the boundary is locally a
segment at $w$.

Finally it is clear that Kobayashi-Nevanlinna cells are disjoint and their union covers all of
$\cl S$ minus the boundaries which have empty interior.$\diamond$

\eno {Definition I.2.3.34}{The Kobayashi-Nevanlinna net is the union of the boundaries of the
Kobayashi-Nevanlinna cells, thus it is a union of euclidean segments.}

%% file: i3.tex
\input defv7.tex
\stit {I.3) Topology of log-Riemann surfaces.}

\stit {I.3.1) The skeleton.}

\eno {Definition I.3.1.1 (Minimal atlas).}{ A minimal atlas $\cl A=\{(U_i, \varphi_i)\}$
is a collection of charts as in definition I.1.1.2 such that the open sets
$(U_i)$ are disjoint, $\bigcup_i U_i$ is dense in $\cl S$, and that can be
completed into a log-Riemann surface atlas as in definition I.1.1.2. The log-Riemann
surface structure is constructed gluing together by the identity on charts the
cuts of the sheets $U_i$.}

Note that a "minimal atlas" is not strictly speaking an atlas since
the open sets $(U_i)$ do not cover completely the surface. Nevertheless
this is irrelevant since it is trivial to add charts covering the cuts
in order to have a complete atlas.

\medskip

It is not difficult to construct minimal atlases.

\eno {Proposition I.3.1.2}{Given a fiber $(w_i)=\pi^{-1}(z_0)$ for a generic point
$z_0 \in \dd C$, the
cells $(U(w_i))$ form a minimal atlas.}

\medskip

Each chart (or sheet) of a minimal atlas of a log-Riemann surface, distinct from
the one sheet log-Riemann surface $\dd C_l$, contains end-points of cuts.
These can be thought of as the trace of ramification
points on charts.

\eno {Definition I.3.1.3 (Clean sheet).}{A sheet with only one trace of
ramification point (that is only with
one cut) is called a clean sheet.}

\eno {Definition I.3.1.4 (Skeleton).}{The skeleton $\G_{\cl S}(U)$ of a
log-Riemann surface $\cl S$ is a connected
graph constructed from a minimal
atlas $U=(U_i)$ as follows:

\item {$\bullet$} Each vertex corresponds to a sheet $U_i$.

\item {$\bullet$} We put an edge between
two vertices for each boundary cut joining the two corresponding $U_i$'s
that are glued together through this cut.

At each vertex the edges occur in pairs corresponding to the same cut.
We call such edges associated.
The graph endowed with the extra information of association
is called the skeleton with articulations and is denoted by $\G'_{\cl S}(U)$.}

\stit {Remark.}

The skeleton does depend on the choice of the  minimal atlas.

%\eno {Proposition.}{The skeleton with articulations (and without)
%$\G'_{\cl S}$ does not depend on the
%choice of the minimal atlas used to construct it.}

%\stit {Proof.}

%Let $(U_i)$ and $(V_j)$ be two minimal atlas. Consider a generic
%base point $z_0\in \cup_i U_i \cap  \cup_j V_j \subset \cl S$ and
%$$
%\pi (z_0)\in \bigcap_i \pi (U_i) \cap \bigcap_j \pi (V_j) \ .
%$$
%We construct an isomorphism between the graphs $\G_{\cl S}^{U}$ and
%$\G_{\cl S}^{V}$. Each point in $\pi^{-1} (\pi (z_0))$ belongs
%to exactly one $U_i$ and one $V_j$. We put them in correspondence.
%We only need to check that this bijection is a graph map.
%Let  $\hat z_0\in \cl S$ with $\pi (\hat z_0)=\pi (z_0)$
%belonging to  $U_{i_0}$ and
%$V_{j_0}$.  For each vertex $U_i$ linked to $U_{i_0}$ we can find
%an oriented  loop in $\dd C$
%with endpoints at $\pi(z_0)$ which lifts to $\cl S$ with end-points
%$\hat z_0$ and $\hat z_1\in U_i\cap \pi^{-1}(\pi (z_0))$.
%Then $\hat z_1$ belongs
%to exactly one $V_j$ that is linked to $V_{j_0}$ through the same
%path and is in correspondence with $U_i$. Conversely any vertex
%$V_j$ linked to $V_{j_0}$ is in correspondence with $U_i$ linked
%to $U_{i_0}$.$\diamond$

\stit {Examples of skeletons.}

With the minimal atlases given in the examples in section I.1.2 we
have the following skeletons.

\bigskip

%{\hfill {\centerline {\psfig {figure=figures/skel1.eps,height=5cm}}}}

\medskip

\bigskip

{\hfill {\centerline {\psfig {figure=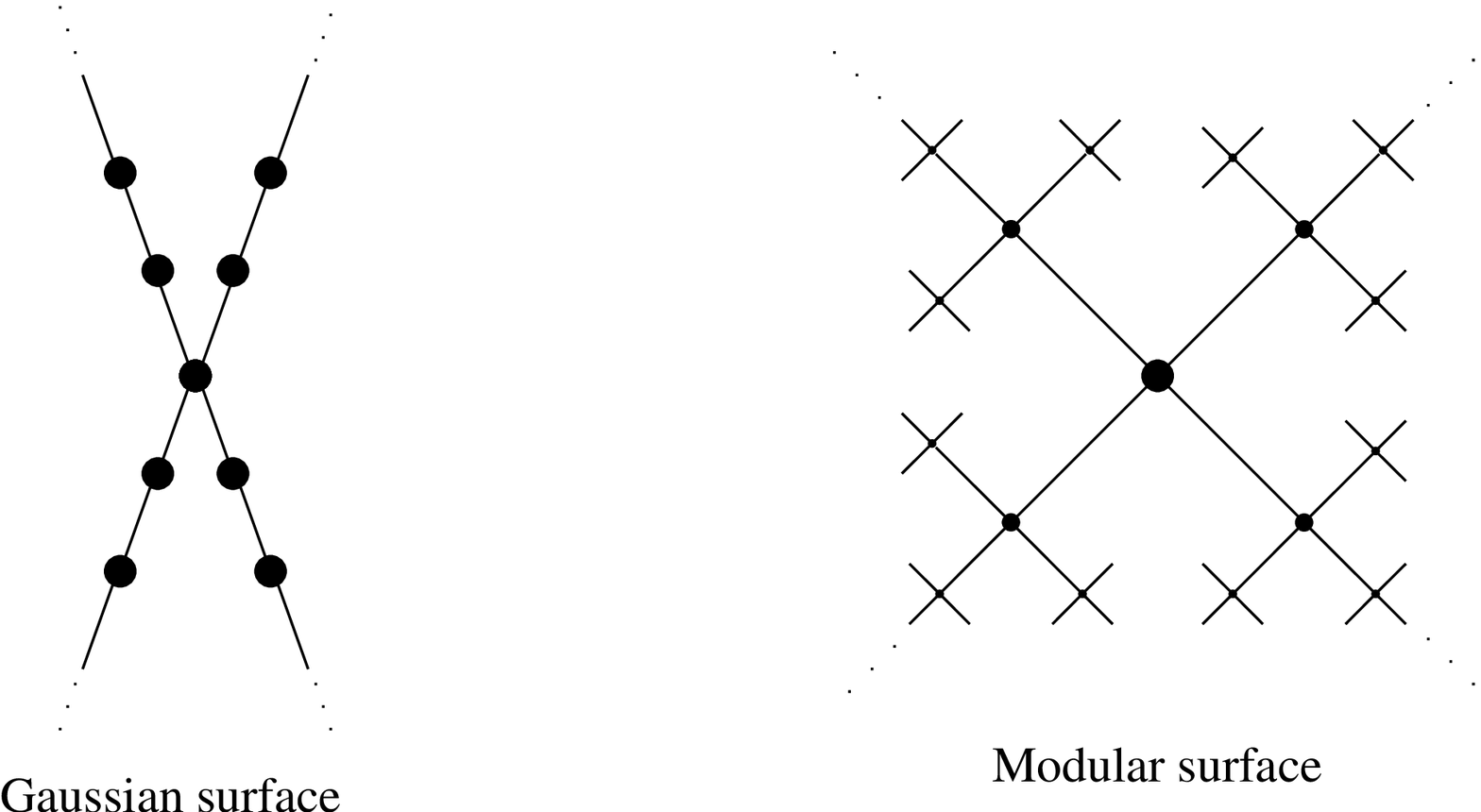,height=5cm}}}}

\medskip

\medskip

A skeleton encodes how the sheets forming the log-Riemann surface
are glued together. Notice that it contains the same information
as the Speiser graph also known as line complex,  which is
classically defined only for those log-Riemann
surfaces having a finite projection of the ramification set (see
[Er2] and [Ne2].) Such log-Riemann surfaces seem to have attracted
most of the classical work related to type
problems (see section II.1 and [BMS].)

We can read on a skeleton many features of the log-Riemann surface.
The following are easy observations.

\eno {Proposition I.3.1.5}{ We have

\item {$\bullet$} Each vertex with only two edges corresponds to a clean sheet.

\item {$\bullet$} Each finite ramification point of order $n$ gives a primitive
cycle of length $n$ in the skeleton.

\item {$\bullet$} Finite ramification points are in one-to-one correspondence
with cycles in the skeleton with articulations formed by edges which are
consecutively associated.

\item {$\bullet$} Infinite ramification points are in one-to-one correspondence
with bi-infinite paths in the skeleton with articulations formed by edges which
are consecutively associated.
}

\medskip

\stit {I.3.2) Skeleton and fundamental group.}

The definition of log-Riemann surface does not imply
that $\cl S$ is simply connected. We can read the fundamental
group of $\cl S$ in the fundamental group of any skeleton
$\G_{\cl S}(U)$.
More precisely,

\eno {Proposition I.3.2.1}{We have
$$
\pi_1(\cl S)\approx \pi_1(\G_{\cl S}(U)) \ .
$$
In particular, the log-Riemann surface $\cl S $ is
simply connected if and only
if the skeleton $\G_S$ is a tree.
}

\stit {Proof.}
Given a maximal atlas $(U_i)$ and chosing a generic point $z_0 \in \cl S$ for
this maximal atlas as before, each loop with base point $z_0 \in U_{i_0}$ and
having a discrete intersection with the cuts, defines
a loop in $\G_{\cl S}(U)$ with base point the vertex $U_{i_0}$ by joining the
vertices $U_i$ through which the loop passes. Conversely, given a loop with
base point $U_{i_0}$ in $\G_{\cl S}$, that is a finite sequence
$$
U_{i_0} \to U_{i_1} \to \ldots \to U_{i_n}\to U_{i_0} \ ,
$$
we can find curves $\g_0, \g_1, \ldots \g_{n+1}$
joining respectively  $z_0$ to $z_1$, $z_1$ to $z_2$, $\ldots$,
$z_n$ to $z_0$ where $z_k\in U_{i_k} \cap \pi^{-1}(\pi(z_0))$
and such that $\g_k \subset \overline {U_{i_k}\cup U_{i_{k+1}}}$.
This gives a loop $\g=\g_0\cup \ldots \g_{n+1} \in \pi_1(\cl S)$.
These two constructions define mutually inverse group homomorphisms.
$\diamond$

\medskip

To each log-Riemann surface $\cl S$ with a minimal atlas $U$,
we associate its skeleton $\G_{\cl S}(U)$
which is a graph with all vertices belonging to an even number of edges.
The converse holds. This is straightforward by direct construction of
the log-Riemann surface. We glue sheets with cuts
according to the connections described by the graph.
Note in particular that the graph is not necesarily planar (for example a
K5 graph is a skeleton but it is not planar.) Note also that distinct
log-Riemann surfaces admit the same skeleton. The skeleton contains
no information about the conformal relative position of the ramification
points.

\eno{Proposition I.3.2.2}{Let $\G$ be a graph with vertices belonging to an even
number of edges. Then there is a non-empty class of log-Riemann surfaces
having $\G$ as skeleton.}

\medskip

%% file: i4.tex
\input defv7.tex
\stit {I.4) Ramified coverings.}

\stit {I.4.1) Ramified coverings and formal Riemann surfaces.}

We present a more general definition of ramified covering between Riemann surfaces (not necessarily log-Riemann
surfaces) than the classical one (see for example [FK] p.15, [Ga] p.441, [BBIF] p.233). The definition below may
appear strange at first, since the new type of ramified covering maps do not necessarily "cover" the base surface.
Nevertheless the definition puts in equal footing finite and infinite ramification points. The notion presented
makes possible the definition of these ramification points, but is quite far from the type of ramified coverings
used to define a Riemann surface orbifold (as defined in [Mi] Appendix E.)

\eno {Definition I.4.1.1 (Ramified covering).}{Let $\cl S_1$ and $\cl S_2$ be
two Riemann surfaces. A mapping $\pi : \cl S_2 \to \cl S_1$ is a ramified
covering if $\pi$ is a local holomorphic diffeomorphism and if the
following condition holds. Given $z_1\in \cl S_1$, for each
neighborhood $U$ of $z_1$
we consider the set $\cl C_U$ of connected components of the pre-image
$\pi^{-1} (U-\{z_1\})$. The set of neighborhoods of $z_1$ form a directed set
by the inclusion, and if $U'\subset U$ then we have a natural map
$$
\cl C_{U'} \to \cl C_U \  .
$$
The inverse limit space
$$
\cl C_{z_1}=\lim_{\leftarrow} \cl C_U
$$
is the space of ends over $z_1$.
For each end over $z_1$
$$
c=(c_U)\in \lim_{\leftarrow }\cl C_U \ ,
$$
we assume that there exists a small Jordan neighborhood
$U = U(z_1,c)$ of $z_1$, such that for its corresponding
connected component $c_U$, the restriction $\pi: c_U \to U - \{z_1\}$ is a classical unramified covering.
Thus one of the following two possibilities must hold:

\item {(1)} Either $c_U$ is bi-holomorphic to a pointed disk
and $\pi_{c_U} : c_U \to U-\{z_1\}$
is a covering of degree $1\leq n<+\infty$.

or

\item {(2)} $c_U$ is biholomorphic to the unit disk and
$\pi_{c_U} : c_U \to U_{c}-\{z_1\}$
is a universal covering of infinite degree $n=+\infty$.

In case (1), $c_U$ corresponds to a ramification of order $n\geq 1$
and we may complete $\cl S_2$ preserving its Riemann surface structure by adding a ramification
point $z^*$. If $n=1$ no ramification exists and we talk of a regular point.
We assume that the sets $c_U$ not corresponding to
regular points are pairwise disjoint.

In case (2), $c_U$ corresponds to a ramification of infinite order $n=+\infty$.
We can also associate to this case a ramification point by adding to the surface $\cl S_2$
a formal point $z^*$. The enlarged set is a topological space by declaring the open
 sets $c_{U'}$, where $U' \subset U = U(z_1,c)$ are Jordan neighbourhoods of $z_1$, to be a base of the
neighborhoods of $z^*$. The enlarged topological space is no longer a surface, it is not even locally compact
in the neighborhood of the infinite ramification points $z^*$.

The formal completion of $\cl S_2$ associated to $\pi$ is denoted by
$\cl S_2^*=\cl S_2 \cup \{ z^*\}$ and is obtained by adding all formal points.
The map $\pi$ extends continuously to
$\cl S_2^*$. We still denote by $\pi : \cl S_2^* \to \cl S_1$ this extension.

}

\stit {Remarks.}

{\bf 1.} Observe that this definition of ramified covering enlarges the classical one. It includes ramified
coverings for which there may exist a point $z_1\in \cl S_1$ for which all neighborhoods $U$ have connected
components of the pre-image $\pi^{-1} (U - \{ z_1 \})$ which are not biholomorphic to a disk or a pointed disk
(see the figure below). However we can still define ramification points, and eventually each ramification point in
$\cl S_2^*$ has a neighborhood which projects nicely by $\pi$.

\bigskip

{\hfill {\centerline {\psfig {figure=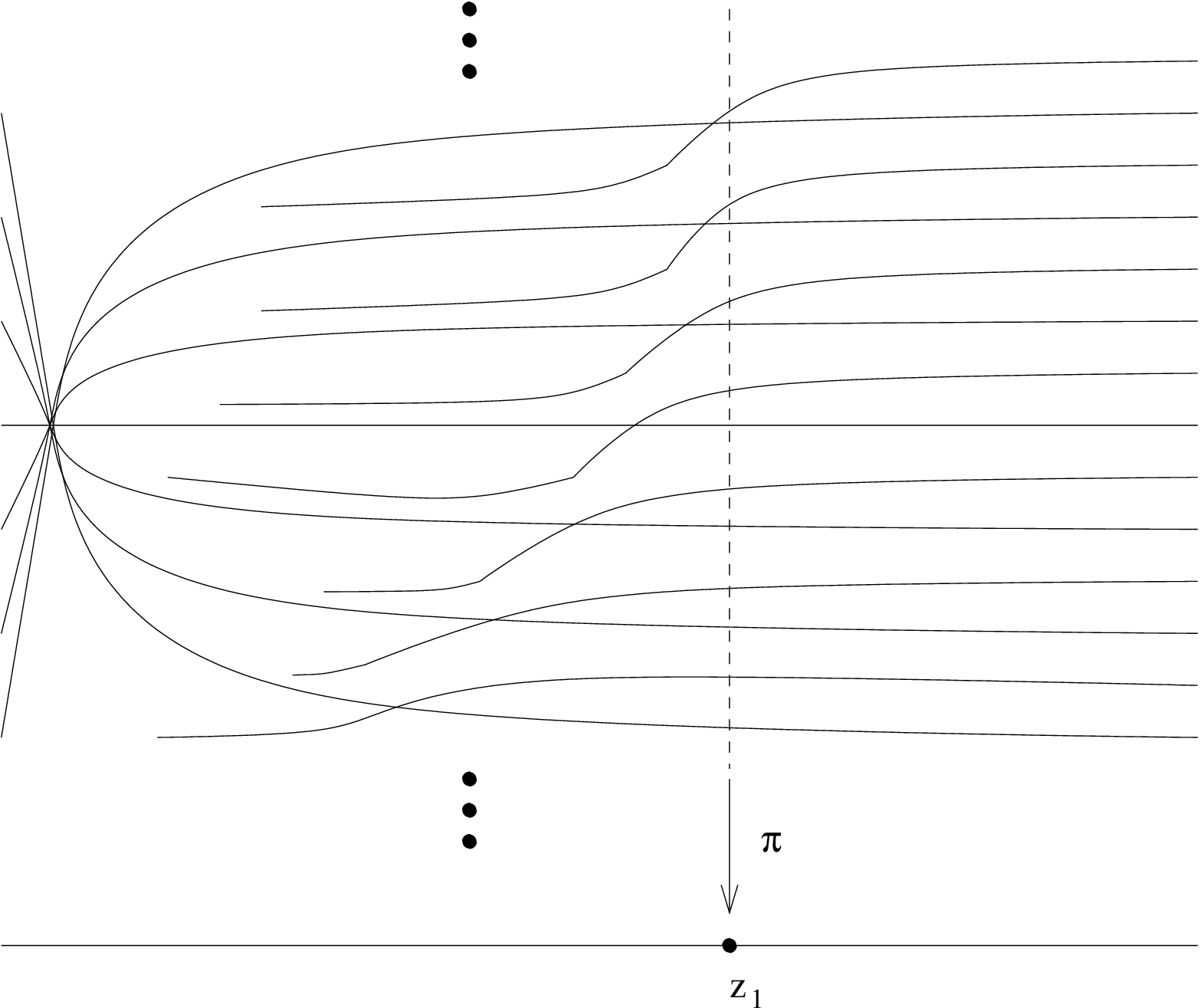,height=6cm}}}}

\bigskip

{\bf 2.} The map $\pi: \cl S_2 \to \cl S_1$ is not necessarily onto, but its extension $\pi : \cl S_2^* \to \cl
S_1$ is.

\medskip

{\bf 3.} The set of ramification points $\cl S_2^*-\cl S_2$ is
discrete.

\medskip

{\bf 4.} For a local holomorphic diffeomorphism $\pi: \cl S_2 \to \cl S_1$
without any extra assumption, the spaces of ends $\cl C_{z_1}$ can
be uncountable as the figure below shows.

\bigskip

{\hfill {\centerline {\psfig {figure=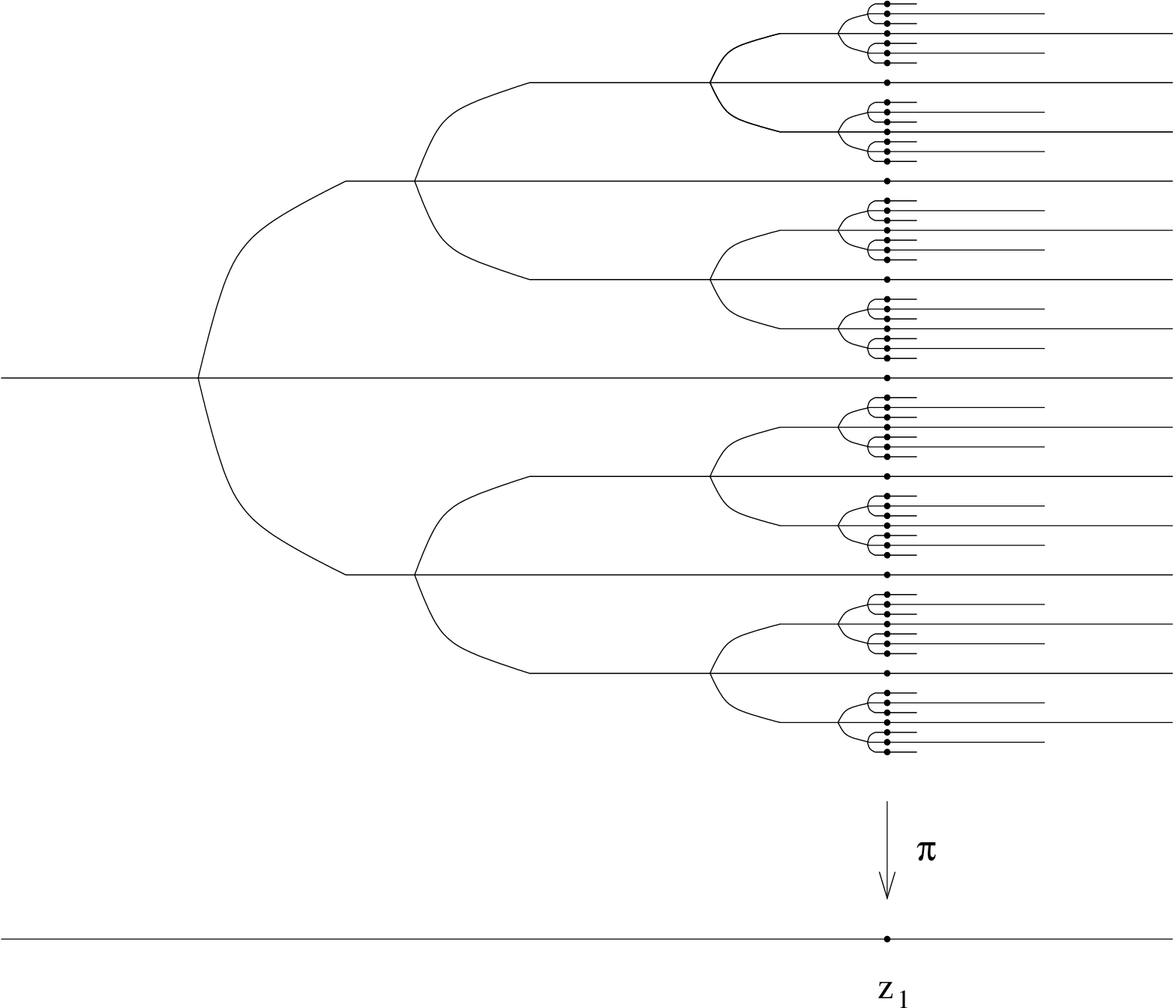,height=6cm}}}}

\bigskip

The formally completed Riemann surface $\cl S^*_2$ is a natural
geometric object and deserves its own terminology.

\eno {Definition I.4.1.2 (Formal Riemann surface).}{A formal Riemann surface
$\cl S^*$ is the formal completion of a Riemann surface $\cl S$ associated
to a ramified covering $\pi : \cl S \to \cl S_1$.}

\eno {Proposition I.4.1.3}{The number of ramification points is at most countable.
Thus a formal Riemann surface is a Riemann surface up to removal of
an at most countable set of points.}

\stit {Proof.}

To each ramification point we can associate a unique open set of $\cl S_2$,
namely $c_U$. These are non-overlapping, and by Rado's theorem a Riemann
surface is $\sigma$-compact. The result follows.$\diamond$

\stit {I.4.2) Ramified coverings and log-Riemann surfaces.}

As our main example of ramified covering we have the projection mapping
$\pi : \cl S \to \dd C$ of a log-Riemann surface $\cl S$ with a
discrete ramification set $\cl R$. In that case the two notions
of ramification points defined so far do coincide.

\eno {Proposition I.4.2.1}{Let $\cl S$ be a log-Riemann surface with discrete
ramification set and
$\pi : \cl S \to \dd C$ its projection mapping.

We have that the projection mapping $\pi$ is a ramified covering of
Riemann surfaces and its ramification points coincide with the ramification
points defined using the Euclidean completion. Both completions $\cl S^*$
are homeomorphic.}

\stit {Proof.}

It has already been shown that for any ramification point $z^*$ in the Euclidean sense, its image $z_1 =
\pi(z^*)$ has a neighborhood $U = \{ |z-z_1| < \varepsilon \}$
enjoying the properties of the definition of
a ramification point for the covering, and these open sets
$c_U$ can be chosen disjoint. Once the union of these,
$\bigcup c_U$, is removed from the surface, all the
end-spaces $\cl C_{z}$ are trivial.
This shows that $\pi$ is a ramification of Riemann surfaces and that
both completions coincide as sets. In a neighborhood of each point and of
ramification
points the topologies coincide, thus both completions are homeomorphic.
$\diamond$

Conversely, we have:

\eno {Theorem I.4.2.2}{If $\cl S$ is a Riemann surface and a ramified covering
$\pi : \cl S \to \dd C$ is given, then this ramified covering
endows $\cl S$ with a log-Riemann
surface structure for which $\pi$ is the projection mapping
and the ramification set $\cl R$ is discrete.}

\stit {Proof.}

We construct the charts using $\pi$. We have to check that
in a given chart the cuts and base points form a discrete set. Note
that each base point of a cut should correspond to a ramification point.
An accumulation of base points at a given point in a given chart
would contradict that $\pi$ is a ramified covering at the $\pi$ image
of that point, because there would be a $c=(c_U)$ as in the definition
with no open set $U(z_1,c)$ associated to $c$. The ramification set $\cl R$
is discrete since ramification points of the covering and ramification
points of the log-Riemann surface as defined in section I.2 do coincide.
$\diamond$

\medskip

Inspired by this last observations it is natural to define a larger class
than that of log-Riemann surfaces. We can take as model of the base
any Riemann surface $\cl S_1$ instead of the complex plane $\dd C$.

\eno {Definition I.4.2.3} {Let $\cl S_1$ be a Riemann surface. A $\cl S_1$-Riemann
surface structure is a Riemann surface $\cl S_2$ endowed with a ramified
covering
$$
\pi : \cl S_2 \to \cl S_1 \ .
$$
}

{\it From now on and for the rest of the article we will only work
with log-Riemann surfaces with a
discrete ramification set $\cl R$.}

\medskip

\stit {I.4.3) Ramified coverings of formal Riemann surfaces.}

\medskip

It is now natural to extend the definition of ramified covering to
formal Riemann surfaces.

\eno{Definition I.4.3.1 (Ramified coverings of formal Riemann surfaces).}
{Let $\pi_{01} :\cl S_1 \to \cl S_{01}$ be a ramified
covering and $\cl S_1^*$ be the associated formal completion
of $\cl S_1$.

We define the notion of ramified covering in the following
two cases:

(1) A mapping $\pi : \cl S_2 \to \cl S_1^*$ from a Riemann surface $\cl S_2$ into $\cl S_1^*$ is defined to be a
ramified covering if it satisfies the conditions of the definition in all the base points, including the points
$z_1^*\in \cl S_1^*-\cl S_1$. For $\pi$-ends $c^{\pi}$ over these last points which are infinite ramification
points of $\cl S_1^*$, the neighborhood $U_1 = U(z_1^*,c^{\pi})$ is not a Jordan neighborhood but a simply
connected neighborhood of the form $U_1 = c_{U_0}^{\pi_{01}}$ where $c^{\pi_{01}}$ is a $\pi_{01}$-end. In this
case $\pi$ adds a formal point $z_2^*$ to $\cl S_2$ lying over $z_1^*$ and the restriction $\pi : c_{U_1}^{\pi}
\to U_1 - \{z_1^*\}$ is univalent, the degree of this ramification point $z_2^*$ for $\pi$ is $1$. As before,
adding all formal points to $\cl S_2$ gives a formal Riemann surface $\cl S_2^* = \cl S_2^*(\pi)$ to which $\pi$
has a continuous extension which we also denote by $\pi$, $\pi : \cl S_2^* \to \cl S_1^*$.

(2) We formulate the same definition
and keep the same terminology when the domain of $\pi$ is also
a formal Riemann surface $\cl S^*_2 = \cl S^*_2(\pi_{02})$,
with formal structure induced by a ramified covering
$\pi_{02} : \cl S_2 \to \cl S_{02}$.
 In this case we impose the additional conditions
that $\pi : \cl S_2^* \to \cl S_1^*$ is continuous
on all of $\cl S_2^*$ (in particular at formal points
introduced by $\pi_{02}$) and also that every $\pi$-end
is equal to an existing $\pi_{02}$-end. This
ensures that $\pi$ does not introduce any new formal points,
and that we must have $\pi^{-1}(\cl
S_1^*-\cl S_1) \subset \cl S_2^*-\cl S_2$.

We observe that in (1) above, the extension
$\pi : \cl S_2^* \to \cl S_1^*$ is trivially a ramified
covering in the sense of (2) above, where $\cl S_2^* = \cl S_2^*(\pi)$
is the formal completion with respect to the original map
$\pi : \cl S_2 \to \cl S_1^*$.
}

\stit {Example.}

{\bf 1.} Let $m|n$ be two positive integers,
and consider the log-Riemann surfaces $\cl S_n$
of the $n$-th root and $\cl S_m$ of the $m$-th root. Then
we have a ramified covering $\cl S_n \to \cl S_m^*$ (and
also $\cl S_n^* \to \cl S_m^*$) preserving
the fibers of the projection.

\medskip

{\bf 2.} If $n$ is a positive integer and $\cl S^*_{\log }$ is the completion
of the log-Riemann surface of the logarithm, we do have a ramified covering
$$
\pi_n: \cl S^*_{\log} \to \cl S^*_{\log } \ ,
$$
with $\pi_n (z)=z^n$.
The image of each plane sheet (corresponding to charts) is
$n$ plane sheets, but the degree of the ramification point above the
infinite ramification point $\cl S_{\log}^* -\cl S_{\log}$ is one.

\medskip

{\bf 3.} We consider the Gauss log-Riemann surface $\cl S_{{\hbox {\rm Gauss}}}$
(example 7 in section I.1.2) and the modular log-Riemann surface
$\cl S_{{\hbox {\rm mod}}}$ (example 8 in section I.1.2) branched at the same
two points as $\cl S_{{\hbox {\rm Gauss}}}$, with projection mappings
$$
\eqalign {
\pi_{{\hbox {\rm Gauss}}} &: \cl S_{{\hbox {\rm Gauss}}} \to \dd C \ ,\cr
\pi_{{\hbox {\rm mod}}} &:\cl S_{{\hbox {\rm mod}}} \to \dd C \ .\cr
}
$$
We have a ramified covering of
formal Riemann surfaces $\pi :\cl S_{{\hbox {\rm mod}}} \to
\cl S_{{\hbox {\rm Gauss}}}$ such that
$$
\pi_{{\hbox {\rm mod}}} = \pi_{{\hbox {\rm Gauss}}} \circ \pi \ .
$$
\medskip

\stit {Remark.}

If $\cl S_1$ and $\cl S_2$ are formal Riemann surfaces and $\pi : \cl S_2 \to \cl S_1$
is a ramified covering of the underlying Riemann surfaces, it is possible that
$\pi$ is not the restriction of a ramified covering of formal Riemann surfaces
$\cl S_2^*\to \cl S_1^*$. For example,
Let $\cl S_2=\dd C$ be endowed with a one-sheet log-Riemann surface structure
(thus $\cl S_2^*=\cl S_2$). Then $\exp: \dd C \to \cl S_{\log}$ is a ramified covering
of Riemann surfaces but is not the restriction of a ramified covering of formal
Riemann surfaces because the ramification point $0^*$ of $\cl S_{\log}^*$
cannot have a pre-image (we denote $\cl S_{\log}^*=\cl S_{\log} \cup \{0^* \}$.)

\medskip

We can observe that, as in the classical case, the composition of ramified coverings is not necessarily a ramified
covering, for example as in the figure below. The Riemann surface $\cl S$ at the top of the figure is built by
pasting planes above and below slits in a base sheet, with the slits converging to a ''half-cut'' at $0$ in this
sheet; this is a slit with only one side, to which slit planes are pasted, so that one can only spiral clockwise
around the point $0$. Note that $\cl S$ is a Riemann surface but {\it not} a log-Riemann surface.

\bigskip

{\hfill {\centerline {\psfig {figure=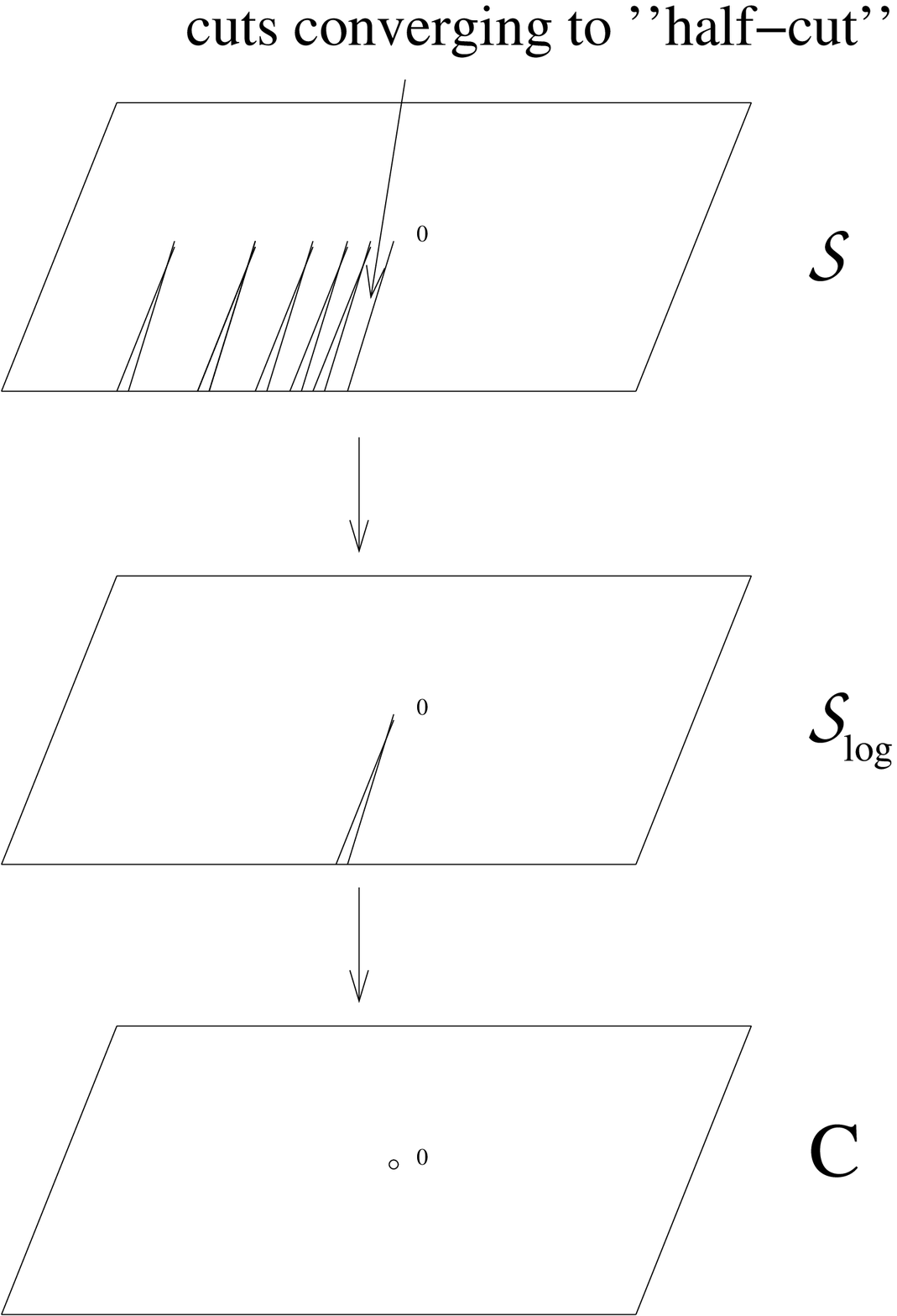,height=9cm}}}}

\bigskip

Nevertheless one of the interests of the new notion of
covering is that it has a better behaviour
under composition once extended to formal Riemann surfaces.

\eno {Theorem I.4.3.2}{Let
$$
\eqalign {
\pi_1 :\cl S_2 \to \cl S_1 & \cr
\pi_2 : \cl S_3 \to \cl S_2 &\cr
}
$$
be ramified coverings of Riemann surfaces, and $\cl S_3^*$ and
$\cl S_2^*$ be the associated formal Riemann surfaces. We assume that
$\pi_2$ is the restriction of a ramified covering of formal Riemann
surfaces
$$
\pi_2 : \cl S_3^* \to \cl S_2^* \ .
$$
Then $\pi_1\circ \pi_2$ is a ramified covering.
}

Another version of this result, staying in the category of
formal log-Riemann surfaces, is the following.

\eno {Theorem I.4.3.3}{Let $\cl S_i^*$, $i=1,2,3$, be formal Riemann surfaces
associated to the ramified coverings
$$
\pi_{0i}:\cl S_i\to \cl S_{0i} \ ,
$$
where $\cl S_{0i}$ are Riemann surfaces.

 Let
$$
\eqalign {
\pi_1 :\cl S^*_2 \to \cl S^*_1 & \cr
\pi_2 : \cl S^*_3 \to \cl S^*_2 &\cr
}
$$
be ramified coverings of these formal Riemann surfaces.

Then
$$
\pi = \pi_1\circ \pi_2: \cl S_3^* \to \cl S_1^*
$$
is a ramified covering of formal Riemann surfaces.
}

\stit {Proof.}

For any $\pi$-end $c^{\pi} = (c_{U}^{\pi})$ over
a point  $z_1 \in \cl S_1^*$, each connected component
$c_{U}^{\pi}$ of ${\pi}^{-1}(U-\{z_1\})$ is a connected
component of ${\pi_2}^{-1}(V)$, where $V \subset \cl S_2^*$
is a connected component of ${\pi_1}^{-1}(U-\{z_1\})$;
thus $c^{\pi}$ determines a $\pi_1$-end $c^{\pi_1}$, which
must then also be a $\pi_{02}$-end, and hence correspond
to a point $z_2 \in \cl S_2^*$. $c^{\pi}$ is then a $\pi_2$-end
over $z_2$, hence also a $\pi_{03}$-end, and so corresponds
to a point $z_3 \in \cl S_3^*$. For $U$ small enough the restrictions
$\pi_1 : c_{U}^{\pi_1} \to U - \{z_1\}$,
$\pi_2 : c_{U}^{\pi} \to c_{U}^{\pi_1}$ are classical coverings, and
considering the different cases when $z_1,z_2$ may be finite or infinite
ramification points of $\pi_1,\pi_2$, it is easily seen that
the composition $\pi : c_{U}^{\pi} \to U - \{z_1\}$ is a classical
covering, as required. $\diamond$

\medskip

Specializing these results to ramified covering between log-Riemann surfaces
we get a notion where infinite ramification points play the
same role as finite ramification points, even if they only exist in the
formal completion.

\eno {Definition I.4.3.4 (Ramified coverings of log-Riemann surfaces).}{Let $\cl S_1$ and
$\cl S_2$ be two log-Riemann surfaces with projection mappings $\pi_1$ and $\pi_2$.

A mapping $\pi: \cl S_2 \to \cl S_1$ is a ramified covering
of log-Riemann surfaces
if $\pi$ is the restriction of a ramified covering of formal Riemann surfaces
$\pi: \cl S_2^* \to \cl S_1^*$ for the formal completions
$\cl S_1^*, \cl S_2^*$ given by the projection mappings $\pi_1,\pi_2$.
}

Observe that in this case $\pi (\cl S_2^*)$ is not necessarily
contained in $\cl S_1$ but we have
$$
\pi (\cl S_2^*)=\cl S_1^* \ .
$$

\medskip

If we denote by $\cl S^*_2 (\pi: \cl S_2 \to \cl S_1)$, resp.
$\cl S_2^*(\pi : \cl S_2 \to \cl S_1^*)$,
$\cl S^*_2 (\pi_1\circ \pi: \cl S_2 \to \dd C)$,
$\cl S_2^*=\cl S_2^*(\pi_2 : \cl S_2 \to \dd C )$, the formal
completion of $\cl S_2$ associated to the ramified coverings
$\pi: \cl S_2 \to \cl S_1$, resp. $\pi: \cl S_2 \to \cl S_1^*$,
$\pi_1\circ \pi: \cl S_2 \to \dd C$, $\pi_2: \cl S_2 \to \dd C$,
we have the continuous embedding
$$
\cl S^*_2 (\pi: \cl S_2 \to \cl S_1) \hookrightarrow
\cl S^*_2 (\pi: \cl S_2 \to \cl S^*_1)=\cl S_2^*(\pi_1\circ \pi: \cl S_2
\to \dd C) = \cl S^*_2(\pi_2: \cl S_2 \to \dd C) \ .
$$

\stit {I.4.4) Universal covering of log-Riemann surfaces.}

We have the following easy, but important, observation.

\eno {Theorem I.4.4.1} {Consider a Riemann surface $\cl S$ endowed with a \
log-Riemann surface structure, and $\tilde {\cl S}$ a universal
covering, $\tilde \pi: \tilde {\cl S} \to \cl S$, of the
underlying Riemann surface.
Then $\tilde {\cl S}$ inherits a log-Riemann surface
structure from $\cl S$ and $\tilde \pi$ is a ramified covering
of log-Riemann surfaces.}

This follows from the next theorem and lemma:

\eno {Theorem I.4.4.2} {Consider a ramified covering between a Riemann surface
$\cl S_2$ and the completed log-Riemann surface $\cl S_1^*$ with
projection mapping $\pi_1 : \cl S_1 \to \dd C$,
$$
\pi : \cl S_2 \to \cl S^*_1 \ .
$$
Then the ramified cover $\pi_2=\pi_1\circ \pi : \cl S_2 \to \dd C$ endows
$\cl S_2 $ with a log-Riemann surface structure, and we have the
homeomorphism
$$
\cl S_2^*(\pi) \approx \cl S^*_2(\pi\circ \pi_1)=\cl S_2^*\ .
$$
}

\eno {Lemma I.4.4.3}{Let $\tilde \pi :\tilde {\cl S } \to \cl S$ be a universal
covering mapping into the Riemann surface $\cl S$ endowed with a log-Riemann
surface structure. Then $\tilde \pi$ defines also a ramified covering
$\tilde \pi :\tilde {\cl S } \to \cl S^*$.
Now $\tilde \cl S^*(\pi)\to \cl S^* \to \dd C$ endows $\tilde \cl S$ with
a log-Riemann surface structure and $\tilde \cl S^*=\tilde \cl S^*(\tilde
\pi)$.}

\stit {Proof of the lemma.}

We can check that the map $\tilde \pi$ defines a
ramified covering $\tilde \pi: {\tilde {\cl S}}
\to \cl S^*$ by choosing a collection of open sets
$\{c_U\}$ associated to the ramification points of
$\cl S^*$ and observing that each connected component of
$\tilde \pi^{-1} (c_U)$ is simply connected and they
give a corresponding set of open neighborhoods associated
to the ramification points of the covering. The rest follows
with similar arguments as before.$\diamond$

\medskip

Therefore we can talk of the universal cover of a log-Riemann surface
structure in the category of log-Riemann surface structures. Since finite
ramification points force a non-trivial fundamental group the following
proposition is obvious.

\eno {Proposition I.4.4.4}{If $\tilde {\cl S}$ is the universal cover of the
log-Riemann surface $\cl S$, the formal completion $\tilde {\cl S}^*$ contains
only infinite ramification points.}

Most universal coverings contain an infinite number of ramification points.

\eno {Proposition I.4.4.5}{If $\tilde {\cl S}$ is the universal cover of the
log-Riemann surface $\cl S$, then $\tilde {\cl S}^*-\cl S$ is infinite
except when $\tilde {\cl S}=\cl S$ and $\cl S$ has a finite number of
ramification points or when $\cl S=\dd C^*$ (then
$\tilde {\cl S}=\cl S_{\log}$ is the Riemann surface of the logarithm.)}

Given a log-Riemann surface $\cl S$ with finite ramification points,
we can complete $\cl S$ adding the finite ramification points in order
to get an enlarged completed Riemann surface denoted by $\cl S^\times
\subset \cl S^*$. We call $\cl S^\times$ the {\it finitely completed
log-Riemann surface}.
The universal covering $\tilde {\cl S^\times}$
of $\cl S^\times$ has no longer
a natural log-Riemann surface structure, but the natural map
$\cl S^\times \to \cl S$ (which is not a covering map according
to our definition since it is not a local diffeomorphism at the points
$\cl S^\times -\cl S$) can be used away from $\cl S^\times -\cl S$ in order
to define infinite ramification points. We denote by $\cl S^{\times *}$ this
Riemann surface with formal points added. Note that we have a ramified covering
$\tilde {\cl S}^* \to \cl S^{\times *}$ (i.e. all points in the base satisfy
the definition.) The universal covering $\tilde {\cl S^\times}$ has more
chances of having a finite number of ramification points and also of being
parabolic. We will discuss in section II.1 the questions related to the type
of log-Riemann surfaces.

\stit {I.4.5) Coverings and degree.}

We restrict our attention to log-Riemann surfaces in this section.

\eno {Proposition-Definition I.4.5.1 (Local and total degree).}{Consider a ramified covering
$\pi: \cl S_2 \to \cl S_1$ of log-Riemann surfaces, and its extension
$\pi: \cl S^*_2\to \cl S^*_1$. Let $z_2 \in \cl S_2^*$ be a point
of order
$1\leq n_2(z_2) \leq +\infty$ and $z_1=\pi (z_2) \in \cl S_1^*$ a
point of order $1\leq n_1(z_1)\leq +\infty$. Then
we have:

\item {(1)} If $n_1(z_1)=+\infty$ then $n_2(z_2)=+\infty$ and $\pi$
has local degree $1$ near $z_2$, we write $d^o_{z_2} \pi =1$.

\item {(2)} If $n_1(z_1) <+\infty$ then $n_2(z_2)=+\infty$ or $n_2(z_2)<+\infty$.
In the first case $\pi$ has local degree $+\infty$ near $z_2$,
$d^o_{z_2} \pi =+\infty$. In the second case, $n_1(z_1)$ divides $n_2(z_2)$
and the local degree of $\pi$ near $z_2$ is
$$
d^o_{z_2} \pi ={n_2(z_2) \over n_1(z_1)} \ .
$$

\item {(3)} The total local degree of $\pi$ does not depend on the point
on the base and is given by
$$
d^o \pi=d^o_{z_1} \pi =\sum_{z_2 \in \pi^{-1} (z_1)} d^o_{z_2} \pi \ .
$$
When $d^o \pi$ is finite we say that $\pi$ is an algebraic
covering.
}

\stit {Proof.}

We construct a loop $\g$ winding $n_2(z_2)$ times around $z_2$. Its image
by $\pi$ winds a multiple of $n_1(z_1)$ times around $z_1$.$\diamond$

\medskip

The local degree is useful in order to determine when a ramified covering
of log-Riemann surfaces is a holomorphic diffeomorphism.

\eno{Theorem I.4.5.2}{Let $\pi: \cl S^*_2 \to \cl S^*_1$ be a ramified covering between
log-Riemann surfaces. If the local degree of $\pi$ is $1$ at all ramification
points of $\cl S_2$ then $\pi$ is a holomorphic
diffeomorphism.}

\stit {Proof.}

If at a ramification point, $\pi$ has degree one then it is a local diffeomorphism into its image. This local
diffeomorphism and its inverse can be continued through the charts. The only possible obstruction to the analytic
univalent continuation are the ramification points in $\cl S_2^*$. But there, $\pi$ is of degree $1$ thus there is
no obstruction. The local diffeomorphisms do match and extend globally. $\diamond$

\medskip

We can prove more.

\eno{Proposition I.4.5.3}{ Let $\pi : \cl S^*_2 \to \cl S^*_1$ as above. Then there exists
an affine automorphism of $\dd C$, $l:\dd C\to \dd C$, such that
$$
\pi_1\circ \pi =l\circ \pi_2 \ .
$$
(see the commutative diagram.)
That is, the log-Riemann surface structure defined on the Riemann surface
$\cl S_2$ by $\pi_1\circ \pi$ is in the affine class of the log Riemann
surface $\cl S_2$.
}

\stit {Proof.}

Note that locally at the $\pi_2$ image of a ramification point of $\cl S_2$ we can define the composition $\pi_1
\circ \pi \circ \pi_2^{-1}$ and this defines a local holomorphic diffeomorphism. By analytic continuation we
extend the domain of definition of this mapping. The only place where we can have an obstruction to a well defined
and univalent continuation is at the $\pi_2$ image of a ramification point. But the assumption ensures that we
still have a univalent continuation across these points. Thus the extension is a local univalent function. The
range is the whole complex plane since $\pi_1\circ \pi (\cl S_2^*)=\dd C$. Thus the inverse is globally well
defined and we get an automorphism of $\dd C$, i.e. an affine diffeomorphism $l$. $\diamond$

\medskip

For a ramified covering of log-Riemann surfaces $\pi :\cl S^*_2 \to \cl S^*_1$,
the formal completion of $\cl S_1^*$ has less
ramification points and of lower order than $\cl S_2^*$.
Also the formal completion of the Riemann surface
$\cl S_2$ associated to $\pi_1\circ \pi$ has less
ramification points and of lower order than $\cl S_2^*$, the formal
completion of $\cl S_2$ associated to $\pi_2$. As pointed out before,
this comes from
the existence of a continuous
embedding  $\cl S_2^*(\pi_1\circ\pi)\hookrightarrow  \cl S_2^*(\pi_2)$.
In some sense it is natural to think of $\cl S_1$ as subordinated
to $\cl S_2$.

\medskip

\eno {Definition I.4.5.4 (Subordination).}{ Let $\cl S_1$ and $\cl S_2$ be log-Riemann
Surfaces with ramified covering maps $\pi_1: \cl S_1 \to \dd C$ and
$\pi_2 :\cl S_2\to \dd C$.

 The log-Riemann surface $\cl S_1$ is subordinate to
$\cl S_2$ if $\cl S_2$ is a log-Riemann Surface
over $\cl S_1$, that is, if there exists a ramified covering of log-Riemann surfaces
$\pi : \cl S^*_2 \to \cl S^*_1$
such that
$$
\pi_2 = \pi_1 \circ \pi \ .
$$
The log-Riemann surface structure defined by $\pi_1\circ \pi$ in the
Riemann surface $\cl S_2$ is weaker than the one defined by $\pi_2$,
$$
\cl S_2^*(\pi_1\circ\pi)\hookrightarrow  \cl S_2^*(\pi_2)=\cl S_2^* \ .
$$

We write $\cl S_2 \geq \cl S_1$.

If $(\cl S^*_1 , z_1)$ and $(\cl S^*_2 , z_2)$ have a distinguished point then
we request that $\pi (z_2)=z_1$, thus the notion of subordination remains well
defined for pointed log-Riemann surfaces.

Note also that we can define the notion of subordination among affine classes.
An affine class is subordinated to another if there are log-Riemann surfaces
in these classes that are subordinated. Then this holds for all the other log-Riemann
surfaces in the affine classes.
}

The next theorem shows that we have defined an order relation.

\eno {Theorem I.4.5.5} {We consider the set $\cl L_0$ of affine classes of
log-Riemann Surfaces with a finite number of infinite ramification points
and finite ramification points of bounded order.

The set $(\cl L_0 , \leq )$ is an ordered set:
{\parindent=2cm

\item {(O1)} $\cl S_1 \leq \cl S_1$.

\item {(O2)}  If $\cl S_1 \leq \cl S_2$ and $\cl S_2 \leq \cl S_1$ then $\cl S_1 =\cl S_2$.

\item {(O3)} If $\cl S_1 \leq \cl S_2$ and $\cl S_2 \leq \cl S_3$ then $\cl S_1 \leq \cl S_3$.
}
}

\stit {Proof.}

Properties (O1) are (O3) are straightforward. To prove (O1)
consider $\pi =\id_{\cl S_1}$.
Property (O3) follows from the observation that a
composition of ramified coverings of log-Riemann surfaces is a
ramified covering of log-Riemann surfaces.
In order to prove (O2) we need the following lemma.

\eno {Lemma I.4.5.6}{Let $\cl S_1, \cl S_2 \in \cl L_0$
with $\cl S_1 \leq \cl S_2$, $\pi_{12}: \cl S^*_2 \to \cl S^*_1$,
and $\cl S_2 \leq \cl S_1$, $\pi_{21}: \cl S^*_1 \to \cl S^*_2$.
Then $\pi= \pi_{12}\circ \pi_{21}$ has degree one at each ramification
point.
}

\stit {Proof.}

Since for both $\pi_{12}$ and $\pi_{21}$ the pre-images of ramification are ramification points of strictly larger
order, the same is true for $\pi: \cl S_1^* \to \cl S_1^*$. But no ramification point can be mapped into a regular
point by $\pi$ because otherwise the total number of infinite ramification points, or finite ones will decrease by
the finiteness assumption in the definition of $\cl L_0$. Thus $\pi$ induces a bijection of the ramification
points. Moreover, $\pi$ must preserve the order of each ramification point. Otherwise the number of ramification
points of infinite order, or finite of a given order will decrease. Thus $\pi$ has order one at each ramification
point.$\diamond$

\medskip

Using this lemma and Proposition I.4.5.3 we get that $\cl S_1$ and $\cl S_2$ are in
the same affine class, thus (O2) holds.

%% file: i5.tex
\input defv7.tex
\stit {I.5) Ramification surgery.}

M. Taniguchi in [Ta3] describes a similar surgery for simply
connected parabolic log-Riemann surfaces (he works with the
entire function uniformizations) and calls it Maskit surgery
by analogy to similar procedures in the theory of Kleinian
groups.

\stit {I.5.1) Grafting of ramification points.}

Given a log-Riemann surface $\cl S_0$,
we can  graft at any point $z_0\in \cl S_0$ a ramification
point of preassigned order (finite or infinite).
This is done by
constructing a new cut  with base point $z_0$
in each chart of $\cl S_0$ containing $z_0$,
$(U_i, \varphi_i)$, and adding a system of
plane sheets associated to this new cut. A finite number of plane
sheets is necessary if the ramification
order is finite. We need an infinite number in order
to create an infinite ramification
point. It is convenient to consider a minimal
atlas containing $z_0$ for carrying out this construction.
The new plane sheets are "clean" in the
sense that they do not contain any other ramification point.
We get in that way a new log-Riemann surface denoted by
$$
\cl S_1=\cl S_0 \sqcup (z_0,n)
$$
where $1 \leq n \leq +\infty$ is the order of the added ramification
point.

The next theorem relies on analytic tools and it will be justified
in part II.

\eno{Theorem I.5.1.1}{We have the subordination
$$
\cl S_0 \leq \cl S_1 \ .
$$
Moreover, $\cl S_1$ has one more ramification point than $\cl S_0$ that
projects into $\pi_0 (z_0)$. We denote this new ramification point by
$z_0^*$.
}

\eno {Proposition I.5.1.2}{Any simply connected log-Riemann surface with a finite number
of ramification points can be obtained from the
complex plane endowed with the one sheet log-Riemann surface structure by
grafting successively a finite number of ramification points.}

% Generalize that to any log-Riemann surface. Prove that the set of
% log-RS is inductive with respect to subordination and do a transfinite
% induction

\eno {Proposition I.5.1.3}{The skeleton $\G_{\cl S_0 \sqcup (z_0, n)}$ is
obtained from $\G_{\cl S_0}$ by adding a loop of length $n$ at the vertex
corresponding to the plane sheet containing $z_0$ if $n$ is finite, or
adding two infinite branches at that vertex if $n=+\infty$.}

We can also graft ramification points at a ramification point
$z_0 \in \cl S^*_0-\cl S_0$. We postpone to the next section this
definition.

\stit {I.5.2) Prunning of ramification points.}

Prunning of ramification point consists in the reverse operation of
grafting.

Consider a log-Riemann surface $\cl S_1$ containing at least one ramification point $z_1^*\in \cl S_1^*$ of order
$\geq 2$. In the plane sheets we can forget about the cuts with base point at $z^*_1$, and add regular points at
the base points of these cuts. In this way we get several connected components in general, each one being a
log-Riemann surface $\cl S_0$ containing  a regular point $z_0$ at the location of $z_1^*$. If all plane sheets
attached to $z_1^*$ are clean planes but one, we obtain only one non-trivial log-Riemann surface, the others being
one-sheeted planes (copies of $\dd C$.) In this last case we denote by $\cl S_0$ the only non-trivial log-Riemann
surface. Note that $$ \cl S_1=\cl S_0\sqcup (z_0, n) \ , $$ where $n$ is the order of $z_1^*$.

\null
\vfill
\eject

%% file: ii1.tex
\input defv7.tex
\bigskip

{\centerline {\tir II. Analytic theory of log-Riemann surfaces.}}

\bigskip

\stit {II.1) Type of log-Riemann surfaces.}

\stit {II.1.1) General facts.}

A large part of the literature on entire functions is about the
problem of the type of the finitely completed Riemann surface $\cl S^\times$.
Since the surface is not compact its universal cover is the unit disk
$\dd D$ (hyperbolic  type) or the complex plane $\dd C$ (parabolic type.)
Most of the literature is devoted to Riemann surfaces branched over a finite set
in the sphere. Classical results are due to R. Nevanlinna, O. Teichm\"uller,
M. Kobayashi, L. V. Ahlfors,...(see [Ne2]) For more recent results the reader
can consult the survey of A. Eremenko [Er2] and one of the latest articles on the
subject [BBIF]. These Riemann surfaces branched over a finite set in the
sphere are combinatorically described by their Speiser graph. The governing
idea is that many ramification points  (or a very arborescent Speiser graph)
favors hyperbolicity. The log-Riemann surfaces we consider are more general, but
the same philosophy holds.

First a trivial remark on the relation between type and subordianation:

\eno {Theorem II.1.1.1}{ Let $\cl S_1$ and $\cl S_2$ be two log-Riemann surfaces with $\cl S_1 \leq \cl S_2$. If
$\cl S_2$ is parabolic, then $\cl S_1$ is parabolic.}

\stit {Proof.}

If $\cl S_1$ was hyperbolic, and $\cl S_2$ parabolic then we could lift
the map $\dd C \to \cl S_2 \to \cl S_1$ into a non constant holomorphic
map from $\dd C$ into $\dd D$.$\diamond$

\medskip

We have seen examples of hyperbolic log-Riemann surfaces such as the modular log-Riemann surface. They seem to
require some important amount of ramification points. On the other hand we have the following theorem by R.
Nevanlinna (that also follows from the uniformization theorem in section II.5). We prove this theorem at the end
of section II.1.2.

\eno {Theorem II.1.1.2}{ If $\cl S^\times$ is a finitely  completed simply connected log-Riemann surface with a
finite number of ramification points then $\cl S$ is parabolic.}

The ramification points of the modular surface all lie in just two
fibers. Thus we don't need a "big" projection set of the ramification
points in order to have a hyperbolic Riemann surface. On the other
hand one may wonder if a "big" projection set implies hyperbolicity.
Note that the projection set is always countable.
The next examples shows that things are not straightforward. We prove
the next theorem at the end of this section.

\eno {Theorem II.1.1.3}{ There exists a parabolic log-Riemann surface $\cl S$ such that the projection set of
ramification points $\pi(\cl S^*-\cl S)$ is dense in $\dd C$.}

We have also.

\eno {Theorem II.1.1.4}{ We consider the log-Riemann surface associated to a generic polygonal billiard. Then the
projection of the ramification set is dense and the Riemann surface is hyperbolic.}

\stit {Proof.}

The composition of two reflections through boundary components define an automorphism of the underlying Riemann
surface. We have at least three distinct pairs of boundary components that define three distinct non-commuting (in
the generic case) automorphisms. Hence the automorphism group is different from the linear group and the surface
cannot be parabolic.$\diamond$

\stit {II.1.2) Kobayashi-Nevanlinna criterium.}

For the rest of this section we consider only log-Riemann surfaces $\cl S$ whose finite completion $\cl
S^{\times}$ is simply connected and having a discrete ramification set $\cl R$.

\medskip

The following Theorem, based on a length-area argument, is to be
found in Nevanlinna ([Ne2] p.317):

\medskip

\eno{Theorem II.1.2.1}{ Let $\cl S$ be a simply connected log-Riemann surface, and $U : \cl S \to \dd R$ a
real-valued function. Suppose that $U$ satisfies the following conditions:
\medskip

{\bf a)} $U$ is continuous on $\cl S$ except for at most isolated
points.

\medskip

{\bf b)} At the points of discontinuity, $U = +\infty$.

\medskip

{\bf c)} The derivatives ${\partial U \over \partial u}$ and
${\partial U \over \partial v} \ (\w = u+iv)$ are continuous
except at most on a family $(\gamma)$ of locally finite smooth
curves.

\medskip

{\bf d)} $\left({\partial U \over \partial u}\right)^2 +
\left({\partial U \over \partial v}\right)^2 > 0$, except for at
most isolated points on the surface.

\medskip

{\bf e)} If $(w_n)$ is an infinite sequence of points with no
accumulation point in $\cl S$, then
$$
U(w_n) \to \infty \hbox{ as } n \to \infty.
$$

Let $\Gamma_{\rho}$ be the union of the curves on $\cl S$ where $U
= \rho$.

\medskip

If the integral
$$
\int^{\infty} {d\rho \over L(\rho)}
$$
is divergent, where
$$
L(\rho) = \int_{\Gamma_{\rho}} |\hbox{grad}_w \ U| \ |dw| \ , \
|\hbox{grad}_w \ U| = \sqrt{\left({\partial U \over
\partial u}\right)^2 + \left({\partial U \over \partial v}\right)^2},
$$
then the surface $\cl S$ is of parabolic type. }

\medskip

Kobayashi,Nevanlinna use the decomposition of the surface into Kobayashi-Nevanlinna cells to derive a more
geometric type criterion which reflects the connection between the type of a surface and its strength of
branching. They work however with surfaces spread over the Riemann sphere, and the cellular decomposition given by
the spherical metric on the sheets, which are slit spheres. Since in this article we only consider the affine
model of log-Riemann surfaces where the sheets are slit planes, and the Kobayashi-Nevanlinna cells as defined in
section I.2.3.4 are defined using the euclidean metric on the sheets, we indicate below how to adapt their
argument to this setting.

\bigskip
\bigskip

{\hfill {\centerline {\psfig
{figure=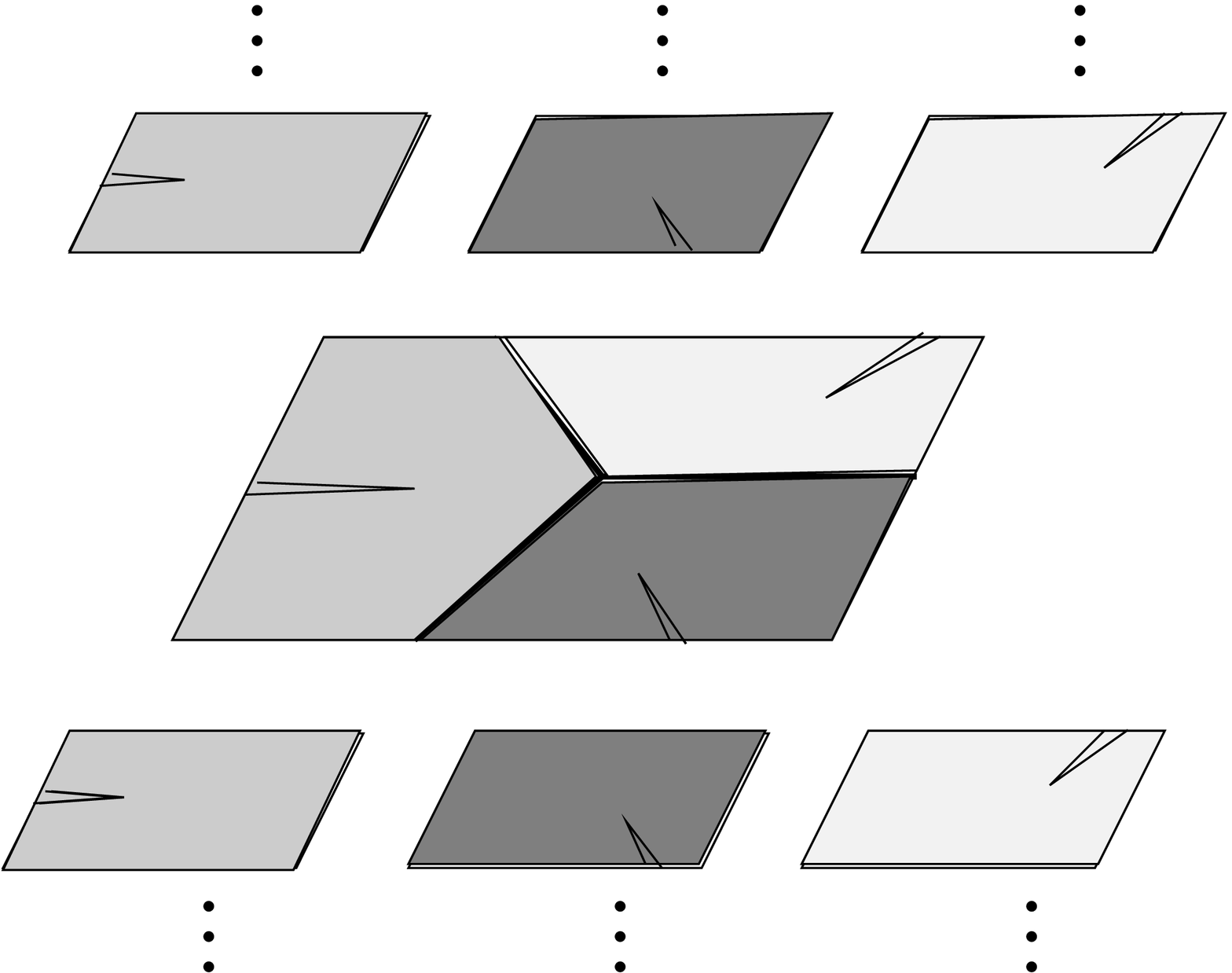,height=5cm}}}}

\bigskip

{\centerline {\bf Figure II.1.1} }

\bigskip
\bigskip

\bigskip

We define on the surface $\cl S$ the continuous differential
$$
d\tau = |d \ \arg(w - a)|
$$
where, for each $w \in \cl S$, $a$ is a ramification point such
that $w$ belongs to the closure of the cell $W(a)$. Fixing a base
point $w_0$ in the net $B$, we define a continuous nonnegative
function $\tau : \cl S \to \dd R$ via
$$
\tau(w) := \inf \int_{w_0}^{w} d\tau
$$
where the infimum is taken over all paths in $\cl S$ joining $w_0$
to $w$.

\medskip

Another continuous nonnegative function $\sigma : \cl S \to \dd R$
is defined by
$$
\sigma(w) := \left |\log |w - a| \right |
$$
where as before $a$ is a ramification point such that $w \in
\overline{W(a)}$.

\medskip

The sum
$$
U(w) := \tau(w) + \sigma(w)
$$
is a function on $\cl S$ that satisfies all the conditions {\bf
(a) - (e) } of the Theorem above.

\medskip

We note that the differential $|\hbox{grad}_w \ U| \ |dw|$ is
conformal and anti-conformal invariant. So making the conformal
change of variables $t = \sigma + i\tau$ gives
$$
|\hbox{grad}_w \ U| \ |dw| = |\hbox{grad}_t \ U| \ |dt| = \sqrt{2}
\, |dt|
$$
So the integral
$$
L(\rho) = \int_{\Gamma_{\rho}} |\hbox{grad}_w \ U| \ |dw|
$$
is equal to the length of the image under $t = \sigma + i\tau$ of
the
 segments $\sigma + \tau = \rho$ of $\Gamma_{\rho}$ multiplied by $\sqrt{2}$.
In order to further estimate this length, we consider for a given
$\tau >0$, the level set $\{ \tau(w) = \tau \}$. This is a union
of line segments which are half-lines or bounded segments. Let
$n(\tau)$ denote the number of such line segments. On each there
lies either {\it one} or {\it no point} of the level set
$\Gamma_{\rho} = \{ \sigma + \tau = \rho \}$. Noting that on
$\Gamma_{\rho}$ we have $t = \rho - \tau + i \tau$ so $|dt| =
\sqrt{2} \, |d\tau|$, we obtain
$$
L(\rho) = \sqrt{2} \ \int_{\Gamma_{\rho}} |dt| = 2 \
\int_{\Gamma_{\rho}} |d\tau| \leq 2 \ \int_{\tau = 0}^{\rho}
n(\tau) \ d\tau
$$
The desired type criterion now follows from the Theorem above:

\medskip

\eno{Theorem II.1.2.2}{ Let $\cl S$ be a log-Riemann surface such that $\cl S^{\times}$ is simply connected and
$\cl R$ is discrete. Let $n(\tau)$ stand for the number of line segments on the surface which are at an angular
distance $\tau$ from a fixed base point $w_0$ in the Kobayashi-Nevanlinna net $B$. If the integral $$
\int^{\infty} {d\rho \over \int_{0}^{\rho} n(\tau) \ d\tau} $$ is divergent, then the surface is of parabolic
type. }

\medskip

We obtain as straightforward corollaries the previously stated Theorems II.1.1.2 and II.1.1.3 :

\medskip

\stit{Proof of Theorem II.1.1.2:} Since $\cl S$ has finitely many ramification points, in this case the counting
function $n(\tau)$ obviously remains bounded, so the integral in the above type criterion is bounded below by an
integral of the form $$ \int^{\infty} {d\rho \over C\rho} $$ and hence divergent.$\diamond$

\medskip

\stit{Proof of Theorem II.1.1.3:} We fix a countable dense set $\{ a_k \}_{k \geq 1} \subset \dd C - \{0\}$, a
strictly increasing sequence of integers $\{ n_k \}_{k \geq 1}$, and an argument function $\arg : {\cl S}_{\log}
\to \dd R$ defined on the log-Riemann surface ${\cl S}_{\log}$ of the logarithm. Consider a log-Riemann surface
$\cl S$ obtained by grafting a ramification point at each point $z_k \in {\cl S}_{\log}, \pi(z_k) = a_k, \arg(z_k)
\in [2\pi n_k, 2\pi (n_k +1) )$. By choosing the sequence $n_k$ growing fast enough, one can make the counting
function $n(\tau)$ for $\cl S$ grow as slowly as one wants, and hence make the integral in the type criterion
divergent.$\diamond$

%% file: ii2.tex
\input defv7.tex
\stit {II.2) Boundary behaviour of the universal cover.}

We study in this section the boundary behaviour of the uniformization when the log-Riemann surface is hyperbolic.
We start by considering a simply connected hyperbolic log-Riemann surface $\cl S$ and a uniformization $$ k:\dd D
\to \cl S \ , $$ and its inverse $$ h: \cl S \to \dd D \ . $$

\medskip

We refer to section II.3.5 for the definitions of Stolz angles and Stolz continuity.

\eno {Theorem II.2.1}{ The uniformization $h$ has a Stolz continuous extension to $\cl S^*$, i.e. $h(w)$ converges
to a limit value on $\dd S^1=\partial \dd D$, denoted by $h(w^*)$, through any Stolz angle at $w^*\in \cl S^*-\cl
S$. In particular, the image by $h$ of any path on $\cl S$ landing at $w^*$ lands at a point on the unit circle
$h(w^*)$.

 The set of points $h(w^*)$ on $\partial \dd D$ corresponding to infinite ramification
 points form a countable set.
 }

 \stit {Proof.}

 The proof is a standard length-area argument. We can
 assume that on a given
 chart containing $w^*$, the point $w^*$ is placed at $0$. Consider
 the circle $C(\rho)$ containing $0$, tangent to the imaginary axes, and
 of diameter $0<\rho \leq r$. We parametrize $C(\rho)$ by the angular
 coordinate and if $w\in C(\rho )$ is the running point
 on $C(\rho )$, $Arg (w)=\t$, we have
 $$
 w=\rho \cos \t e^{i\t }\ ,
 $$
thus
$$
dw=\rho i e^{2i\t} d\t \ ,
$$
and
$$
|dw|=\rho d\t \ .
$$

Now using Cauchy-Schwarz inequality, denoting by $l(\rho)$ the length of $h(C(\rho))$, $$ \eqalign { l(\rho)^2
=\left ( \int_{C(\rho)} |h'(w)| |dw|\right )^2 &\leq \left ( \int_{C(\rho)} |dw| \right ) \left ( \int_{C(\rho)}
|h'(w)|^2 |dw|\right ) \cr &\leq \pi r \int_{-\pi}^\pi  |h'(w)|^2 \rho \ d\t \cr } $$ Now, integrating over
$0<\rho \leq r$ we get, denoting by $A(r)$ the area of $h(D(r))$ where $D(r)$ is the disk bounded by $C(r)$, $$
\int_0^r {l(\rho )^2 \over \rho} \ d\rho \leq \pi \int_0^r \int_{-\pi}^\pi |h'(w)|^2 \ \rho d\t \ d\rho =\pi A(r)
\ . $$ Since $h(D(r))\subset \dd D$ we have that $A(r)<+\infty$. Thus we conclude that $$ \int_0^r {l(\rho )^2
\over \rho} \ d\rho <+\infty \ . $$ From this it follows that there exists a sequence $\rho_n \to 0$ such that
$l(\rho_n)\to 0$. Then ${\hbox {\rm diam}} (D(\rho_n)) \to 0$ and it follows that $h$ has a limit on the Stolz
angle $-\pi +\eps < \t < \pi -\eps$. Rotating this sector we obtain the same result and we conclude that we have a
Stolz limit in any Stolz angle. $\diamond$

\stit {Observation.}

We can extend this theorem when $\cl S$ is not simply connected but still hyperbolic.

\medskip

%% file: ii3.tex
\input defv7.tex
\stit {II.3) Caratheodory theorem for log-Riemann surfaces.}

\stit {II.3.1) Kernel convergence.}

We extend the notion of Kernel convergence of domains in the plane
to log-Riemann surfaces, in view of defining a topology in the
space of log-Riemann surfaces and extending Caratheodory's theorem.

\medskip
Recall that a log-Riemann surface is naturally endowed with its {\it log-Euclidean metric}.

\eno {Definition II.3.1.1}{ A pointed sequence of log-Riemann surfaces $(\cl S_n ,z_n )$ converges to a pointed
log Riemann surface $(\cl S, z)$ if for any compact set in the surface topology  $z\in K\subset \cl S$ there
exists $N=N(K)\geq 1$ such that for $n\geq N$ there is an isometric embedding of $K$ into $\cl S_n$ for the
corresponding log-Euclidean metrics and mapping $z$ into $z_n$. The embeddings of two overlapping such compact
sets are supposed to be compatible. We assume that they are a translation on charts (the translation that maps
$\pi(z)$ into $\pi(z_n)$).}

Note that such a limit must be unique because two such limits
would be isometrically embeddable one into the other, thus will
correspond to the same log-Riemann surface. Such a limit is unique
up to isometry; more precisely we have the following

\medskip

\eno{Proposition II.3.1.2}{ Let $(\cl S, z_0)$ and $(\cl S', z_0')$ be two pointed log-Riemann surfaces both of
which are Caratheodory limits of a sequence of pointed log-Riemann surfaces $(\cl S_n, z_n)$. Suppose $\cl S$ and
$\cl S'$ have discrete ramification sets. Then $\cl S$ and $\cl S'$ are isometric via an isometry $T : \cl S \to
\cl S'$ that takes $z_0$ to $z_0'$, $T(z_0) = T(z_0')$, and whose expression in log-charts is the translation that
maps $\pi(z_0)$ to $\pi'(z_0')$.}

\medskip

\stit{Proof.} Consider the germ of holomorphic diffeomorphism $T$
from $\cl S$ to $\cl S'$ that maps $z_0$ to $z_0'$ and whose
expression in log-charts is the translation that maps $\pi(z_0)$
to $\pi'(z_0')$.

\medskip

We observe the following : Let $\gamma : [a,b] \to \cl S$ be a
curve in $\cl S$ starting from $z_0$, $\gamma(a) = z_0$ along
which $T$ can be continued. Then for $n$ large enough, if $\iota$
and $\iota'$ denote the isometric embeddings of $\gamma$ and
$T(\gamma)$ respectively into $\cl S_n$, then on $\gamma$ we must
have
$$
\iota = \iota' \circ T,
$$
since both maps have derivatives (computed in log-charts) equal to
1, and coincide at $z_0$, $\iota(z_0) = \iota'(T(z_0)) = z_n \in
\cl S_n$.

\medskip

We prove the following lemmas:

\medskip

\eno{Lemma II.3.1.3}{ The germ $T$ can be continued analytically along all paths in $\cl S$.}

\medskip

\stit{Proof.} Suppose there is a path $\gamma : [a,b] \to \cl S,
\, \gamma(a) = z_0, \gamma(b) = z \in \cl S$ along which $T$
cannot be continued, more precisely $T$ can be continued
analytically along $\gamma([a,b))$ but not upto $\gamma(b) = z$.
Since $T$ is a local isometry, as $x \in [a,b)$ tends to $b$ the
following limit must exist in the completion $\cl S'^*$ of $\cl
S$,
$$
z' := \lim_{x \to b} T(\gamma(x)) \ \in \cl S'^*
$$
Since $T$ cannot be continued to $\gamma(b)$, we must have $z'
\notin \cl S'$, so $z'$ is a ramification point of $\cl S'$.

\medskip

Take $\delta > 0$ small enough so that $B(z',\delta)$ contains no
other ramification points and so that $\overline{B(z,\delta)}
\subset \cl S$. Let $b_1 \in [a,b)$ be such that $T(\gamma(b_1))
\in B(z',\delta)$. Let $\alpha : [c,d] \to \cl S$ be a circular
loop in $\cl S$ that winds once around $z$, starting from
$\gamma(b_1)$, $\alpha(c) = \gamma(b_1)$. We note that $T$ can be
continued along $\alpha$, but $T(\alpha(c))$ is not equal to
$T(\alpha(d))$.

\medskip

Now consider $n$ large enough so that the compacts $\gamma([a,b]) \cup \alpha([c,d]) \cup \overline{B(z,\delta)}
\subset \cl S$ and $T(\gamma([a,b_1])) \cup T(\alpha([c,d])) \subset \cl S'$ both embed isometrically into $\cl
S_n$, and denote by $\iota, \iota'$ the respective embeddings. Now, the ball $\iota(B(z,\delta))$ is completely
contained in $\cl S_n$, so for the curve $\alpha$ we have $$ \iota(\alpha(c)) = \iota(\alpha(d)) =
\iota(\gamma(b_1)) \ \hbox{ where } \alpha(c) = \alpha(d) = \gamma(b_1), $$ which implies $$ \iota'(T(\alpha(c)))
= \iota'(T(\alpha(d))) $$ and hence, since $\iota'$ is an isometry, $$ T(\alpha(c)) = T(\alpha(d)), $$ a
contradiction. $\diamond$

\medskip

\eno{Lemma II.3.1.4}{ The continuation of $T$ to all of $\cl S$ is single-valued, ie there is no monodromy from
continuing along closed paths.}

\medskip

\stit{Proof.} Let $\gamma : [a,b] \to \cl S, \, \gamma(a) =
\gamma(b) = z_0$,
 be a closed path in $\cl S$. Consider the curve $T(\gamma) \in \cl S'$ given
by continuing $T$ along $\gamma$. Take $n$ large enough so that
the compacts $\gamma([a,b]) \subset \cl S$ and $T(\gamma([a,b]))
\subset \cl S'$ both embed isometrically into $\cl S_n$ via
isometries $\iota$ and $\iota'$ respectively. As before, along
$\gamma$ we have
$$
\iota = \iota' \circ T
$$
Since
$$
\iota(\gamma(a)) = \iota(\gamma(b)) = \iota(z_0) = z_n,
$$
it follows that
$$
\iota'(T(\gamma(a))) = \iota'(T(\gamma(b)))
$$
and hence
$$
T(\gamma(a))) = T(\gamma(b))),
$$
ie $T$ has no monodromy when continued along $\gamma$.
$\diamondsuit$

\medskip

It follows from the above lemmas that we obtain a globally defined map $T : \cl S \to \cl S'$. Applying the same
arguments to the germ $S = T^{-1}$ given by the inverse of the initial germ $T$ gives a map $S : \cl S' \to \cl
S$, and it is straightforward to check that $T$ and $S$ define global mutual inverses. The conclusions of the
Proposition follow. $\diamond$

\stit {Observation.}

Note that for this convergence notion finite ramification points of increasing
order do converge to infinite ramification points.
An instructive elementary example is the sequence $(\cl S_n ,1)$ of log-Riemann
surfaces of the $\root n \of z$ (branched at $0$), that do converge to the
log-Riemann surface of the logarithm.

\bigskip

We know by Rado's theorem that any Riemann surface is
$\sigma$-compact (i.e. a countable union of compact
sets)\footnote{$\dagger$}{Some texts do include this in the
definition of Riemann surface.}. Thus choosing an exhausting
sequence of compact sets, we can define a base of neighborhoods of
the log-Riemann surface for the above kernel Caratheodory
convergence. This defines a Hausdorff topology in the space of
log-Riemann surfaces.

\bigskip

More generally we can define Caratheodory convergence of domains in log-Riemann surfaces. In the following
definitions and theorems we consider only isometries which are translations in log-charts.

\eno{Definition II.3.1.5}{ A log-domain is a domain $U$ in a log-Riemann surface $\cl S$. We do not distinguish
between log-domains that are isometric for the log-Euclidean metric (even when they belong to different
log-Riemann surfaces).}

\eno{Definition II.3.1.6}{ Let $(U_n, z_n)$ be a sequence of pointed log-domains inside log-Riemann surfaces $\cl
S_n$. A pointed log-domain $(U,z)$ belongs to the kernel of the sequence if for any compact set $K\subset U$ with
$z\in K$ there exists $N=N(K)\geq 1$ such that for $n\geq N$ there exists an isometric embedding $$ K
\hookrightarrow U_n $$ mapping $z$ to $z_n$. }

\eno {Definition II.3.1.7 (Subordination of log-domains).}{The pointed log-domain $(U_1, z_1)$ is subordinated (or
smaller) than the pointed log-domain $(U_2, z_2)$ if we can embed isometrically $U_1$ into $U_2$ mapping $z_1$
onto $z_2$. We write $$ (U_1, z_1) \leq (U_2, z_2) \ . $$ }

\eno {Definition II.3.1.8 (Kernel of a sequence of log-domains).}{Given a sequence of pointed log-domains $(U_n,
z_n)$, we consider all log-domains $(U,z)$ belonging to the kernel of this sequence. If there is one such
log-domain that is maximal in this family, it is the kernel of the sequence.}

\stit {Remark.}

Allowing log-Riemann surfaces constructed with charts with
non-locally finite cuts, and correspondingly enlarging the
definition of log-domains, we would be able to prove the existence
in general of kernels (we refer to [BiPM1]). We cannot avoid a
sequence of log-Riemann surfaces converging to one with
non-discrete ramification set or non-locally finite cuts in the
charts. For now we will prove the following.

\eno {Theorem II.3.1.9}{ Let $(U_n, z_n)$ be a sequence of log-domains belonging to log-Riemann surfaces $\cl S_n$
with finite ramification sets of uniformly bounded cardinality. Then this sequence has a kernel which belongs to
such a log-Riemann surface.}

\stit{Proof.} Consider all pointed log-domains $\{ (U,z) \}$ which
belong to the kernel of the sequence $(U_n, z_n)$. If there are none
such then the kernel is empty and there is nothing to prove. If not,
we can paste the domains together isometrically as follows:
We consider their disjoint union
$$
V = \bigcup_{(U,z)} U.
$$
and quotient $V$ by the following equivalence relation:
$$
\xi \in (U,z) \sim \xi' \in (U',z')
$$
if for all compacts $K, K'$ such that $z,\xi \in K \subset U
, \ z',\xi' \in K' \subset U'$, we have, for $n$ large enough
such that $K, K'$ embed isometrically into $U_n$, that
$$
\iota_n(\xi) = \iota_n'(\xi') \in U_n,
$$
where $\iota_n, \iota_n'$ are the embeddings of $K$ and $K'$ respectively
into $U_n$. This gives a metric space
$$
\overline{V} = V / \sim
$$
which is a Riemann surface with a flat metric. The points $z$ of the
pointed domains $(U,z)$  get identified to a single point
$\overline{z} \in \overline{V}$.

\medskip

Since the log-Riemann surfaces $\cl S_n$ have a uniformly bounded number of
ramification points, it is not hard to see that $\overline{V}$
can be embedded isometrically into a log-Riemann surface $\cl S$ with
a finite number of ramification points. Moreover
$(\overline{V}, \overline{z})$ belongs to the kernel of the
sequence $(U_n, z_n)$, and all the pointed log-domains
$(U,z)$ embed isometrically into $(\overline{V},\overline{z})$,
hence $(\overline{V}, \overline{z})$ is the kernel of the
sequence $(U_n, z_n)$.
$\diamondsuit$

\stit {II.3.2) Caratheodory theorem.}

Now we present a generalization to log-Riemann surfaces of
Caratheodory's kernel convergence theorem.

\eno {Theorem II.3.2.1}{ Let $(\cl S_n , z_n) \to (\cl S, z)$ be a Caratheodory's converging sequence of
log-Riemann surfaces such that the finite completions $\cl S_n^\times, \cl S^\times$ are simply connected. Let
$F_n :\cl S_n^\times \to \dd D_{R_n}$ be the uniformizations of $\cl S_n^\times$ into the complex plane $\dd C$
($R_n=+\infty$) or a finite disk ($R_n<+\infty$), normalized such that $(F_n\circ \pi_n^{-1})'(z_n)=1$. Let $F:\cl
S^\times \to \dd D_R$ be the uniformization of $\cl S^\times$ with $R\in ]0, +\infty ]$ and $(F\circ
\pi^{-1})'(z)=1$.

If
$$
\limsup_{n\to +\infty} R_n \leq R \ ,
$$
then $\lim_{n\to +\infty} R_n =R$ and the sequence $(F_n)$
converges uniformly on compact sets of $\cl S$ to the
uniformization
$$
F: \cl S \to \dd D_R \ .
$$

The uniform convergence on compact sets holds in the following
sense: For each compact set $K\subset \cl S$ in the surface
topology and for $n\geq N(K)$ we have a log-Euclidean isometric
embedding $\iota_n : K \to \cl S_n$ such that the maps $F_n\circ
\iota_n$ are well defined on $K$ for $n\geq N(K)$ and converge
uniformly to $F$. }

\stit {Proof.}

The proof follows the same lines as Caratheodory planar kernel
convergence. The sequence of uniformizations defines on compact
sets of $\cl S$ a family of equicontinuous univalent functions.
Consider a limit point $g:\cl S^\times \to \dd C$ of the sequence of
uniformizations. Then $g$ is univalent and hence $g\circ F^{-1}:
\dd D_R \to \dd D_R$ is an automorphism and by the normalization
it is the identity. Hence the limit is unique and equal to
$F$.$\diamond$

\stit {Remark.}

Note that contrary to the classical Caratheodory theorem, here the
parabolic case is meaningful.

 \stit {II.3.3) Conformal radius.}

The classical definition associates to a simply connected domain
in the plane a conformal radius that depends monotonically on the
domain. We can extend this definition to domains in log-Riemann
surfaces and to log-Riemann surfaces.

\eno {Definition II.3.3.1}{ Let $(\cl S, z_0)$ a pointed log-Riemann surface. We consider a simply connected
domain $U\subset \cl S$ with $z_0 \in U$. The domain $U$ is hyperbolic or parabolic. In the parabolic case we
define its conformal radius to be $+\infty$. In the hyperbolic case there exists a unique $0<R_0 <+\infty$ such
that we have a uniformization into the disk of radius $R_0$, $$ h: U\to \dd D_{R_0} $$ such that $h(z_0)=0$ and
$h'(z_0)=1$, where the derivative of $h$ at $z_0$ is computed using the canonical charts of the log-Riemann
surface. Then we define the conformal radius of $U$ to be $R(U)=R_0$. }

In particular the definition applies to $U=\cl S$ when $\cl S$ is simply connected. The following propositions are
straightforward.

\eno {Proposition II.3.3.2}{ The conformal radius $R(U)$ of a simply connected domain $U$ does not depend on the
log-Riemann surface where it belongs. More precisely, let $(\cl S', z_0')$ be another pointed log-Riemann surface.
Suppose that $U$ can be isometrically embedded into $\cl S'$ by an isometric immersion for the Euclidean metric
mapping $z_0$ to $z_0'$. Then $$ R(U)=R(U') \ . $$ }

\eno{Proposition II.3.3.3}{ Let $U_1\subset U_2\subset \cl S$ be two simply connected domains containing the base
point $z_0$. Then $$ R(U_1)\leq R(U_2) \ , $$ with equality only when $U_1=U_2$. }

We have the following extension of Caratheodory theorem to
log-domains. The proof is the same as before.

\eno {Theorem II.3.3.3}{ Let $(U_n)$ be a sequence of pointed simply connected log-domains converging in the
Caratheodory sense to a log-domain $U$. If $$ \limsup_{n\to +\infty} R(U_n) \leq R(U) \ , $$ then $$ \lim_{n\to
+\infty} R(U_n) =R(U) \ . $$ }

\stit{II.3.4) Closure of algebraic log-Riemann surfaces.}

In this section we define algebraic log-Riemann surfaces and we
determine their closure for the Caratheodory topology.

\eno {Definition II.3.4.1}{ Let $\cl S_1$ and $\cl S_2$ be two log-Riemann surfaces. We say that $\cl S_2$ is
algebraic over $\cl S_1$ if $\cl S_2$ can be obtained from $\cl S_2$ by grafting a finite number of finite
ramification points.

The log-Riemann surface $\cl S_1$ is algebraic if $\cl S_1$ is
algebraic over $\dd C$. }

Obviously a log-Riemann surface is algebraic if and only if it has
a finite number of ramification points of finite order and none of
infinite order.

\eno {Definition II.3.4.2}{ For $n\geq 0$ we define the space $\cl A_n$ of algebraic log-Riemann surfaces  with at
most $n$ ramification points endowed with Caratheodory topology. We define also the space $\cl {TRA}_n$ of
log-Riemann surfaces with at most $n$ ramification points endowed with Caratheodory topology, $$ \cl A_n \subset
\cl {TRA}_n \ . $$

The space
$$
\cl A=\bigcup_{n\geq 0} \cl A_n \ ,
$$
is the space of algebraic log-Riemann surfaces.

The space
$$
\cl {TRA}=\bigcup_{n\geq 0} \cl {TRA}_n \ ,
$$
is the space of transalgebraic (or finite) log-Riemann surfaces.

We also define $\cl L$ as the space of all log-Riemann surfaces,
$$
\cl {A}\subset \cl {TRA}\subset \cl L
$$
}

In [PM] we characterize the functions with a finite number of
exponential rational singularities as the closure of meromorphic
functions with bounded divisor. The following result is closely
related and is the geometric analog.

\eno {Theorem II.3.4.3}{ The space $\cl {TRA}_n$ is closed in $\cl L$. The closure of the space of algebraic
log-Riemann surfaces in $\cl L$ is the space of transalgebraic log-Riemann surfaces, more precisely $$ \overline
{\cl A_n} =\cl {TRA}_n \ . $$ }

Before proving the Theorem we make several observations and prove
a Lemma. We have seen that the order of ramification points can
increase under Caratheodory convergence as shows the example
$$
\cl S_n \to \cl S_{\log} \ .
$$

Also several ramification points can collapse into a ramification point. For example, taking the appropriate base
point, $$ \lim_{\eps \to 0} \eps \  \cl S_{{\hbox {\rm {Gauss}}}} =\cl S_{\log} \ . $$ We also have examples where
the ramification point escapes to $\infty$, such as
 $$ \lim_{a\to \infty} a+\cl S_{\log} = \dd C \ . $$

Nevertheless the number of ramification points cannot increase as
the next lemma shows.

\eno {Lemma II.3.4.4}{ The number of ramification points is lower-semi-continuous for the Caratheodory
convergence.}

\stit {Proof of the Lemma.}

 The proof uses the following lemma that results from
 the definition of Caratheodory convergence.

\eno {Lemma II.3.4.5}{ Let $(\cl S_n , z_n)$ be a sequence of log-Riemann surfaces with projection maps $\pi_n :
\cl S_n \to \dd C$, having $(\cl S, z)$ as Caratheodory limit. If $\pi : \cl S \to \dd C$ is the projection
mapping of $\cl S$ and $$ \lim_{n\to *\infty} \pi_n (z_n)=\pi (z) \ , $$ then the sequence $(\pi_n)$ converges
uniformly on compact sets to $\cl \pi$.}

Recall that the uniform convergence on compact sets mean that for
any compact set $K\subset \cl S$, for $n$ large enough, and
denoting by $\iota_n : K \cl \to S_n$ the isometric embedding, we
have that the sequence $(\pi_n\circ \iota_n)$ converges uniformly
on $K$.

 Now the proof of the first lemma is immediate. By Hurwitz theorem,
if the projection mappings $\pi_n$ are locally one-to-one, then their limit $\pi$ is locally one-to-one. Thus a
ramification point can only be the limit of ramification points.$\diamond$

\stit {Proof of the theorem.}

The proof of the theorem is immediate from the
lower-semi-continuity of ramification points.$\diamond$

We can investigate what happens when we bound the order of the
ramification points instead of their number.

\eno {Definition II.3.4.6}{ We consider the space $\cl {TRA}_{n,m}\subset \cl {TRA}_n$ of transalgebraic
log-Riemann surfaces such that the sum of orders of all finite ramification points is at most $m\geq 0$.}

 \eno {Theorem II.3.4.7} { The order of ramification
points is lower-semi-continuous for the Caratheodory convergence.
More precisely, if several finite critical points collapse into a
ramification point, and $m$ is an upper bound for the sum of their
orders, then the limit ramification point has order $m$ at most.

In particular, we have
$$
\overline {\cl {TRA}_{n,m}} \subset \bigcup_{l=0}^n \cl
{TRA}_{n-l,m}
$$
 }

The theorem follows from the next lemma.

\eno {Lemma II.3.4.8}{ Let $(\cl S_n)$ be a sequence of log-Riemann surfaces with Caratheodory limit $\cl S$. Let
$w^*\in \cl S^*-\cl S$ be a finite ramification of order $m$. For any $\eps
>0$, and for $n$ large enough we embed the compact set
$K_\eps=\bar B(w^*, \eps )-B(w^*, \eps /2 )$ into $\cl S_n$,
$\iota_n: K_\eps \hookrightarrow \iota_n (K_\eps) \subset \cl
S_n$. Then the bounded connected component of $\cl S_n
-\pi_n^{-1}(\pi_n(\iota_n (K_\eps)))$ whose closure intersects
$\iota_n(K_\eps )$, contains ramification points whose orders add
up to $m$ at least.}

\stit {Proof of the lemma.}

Locally there are as many sheets attached locally to a
ramification point as its order. Some of these coincide and they
all together contribute to the $m$ sheets where $\iota(K_\eps)$
lies. $\diamond$

\stit {Example.}

The sum of the orders can actually go down as the following example shows: Take the algebraic elliptic log-Riemann
surface with two ramification points of order $2$ and collapse these two into a single ramification point of order
two.

\stit {II.3.5) Stolz continuity.}

\medskip

Functions that are natural on $\cl S$ do extend in some natural
way to $\cl S^*$.
They do extend continuously to $\cl S^*$ in Stolz angles at the ramification
points.

\eno {Definition II.3.5.1 (Stolz angle).} { Let $w^*\in \cl S^*-\cl S$ be an infinite ramification point. Consider
a log function branched at $w^*$, $\log_{w^*}(w)=\log (w-w^*)$. The argument function $\Arg_{w^*}=\Im \log_{w^*}$
is well defined in a neighborhood of $w^*$.

A Stolz angle at $w^*$  of radius $r >0$, amplitude $M>0$ and
centered at $\t \in \dd R$, is a region of the form
$$
U(M,\t , r)=\{w\in B(w^*,r)-\{w^*\}; |\Arg_{w^*} (w)-\t| < M \}
\cup \{w^* \}\ .
$$
The radius $r>0$ should be small enough
so that $B(w^*, r)$ contains no other
infinite ramification point.

We can define Stolz angles at a finite ramification point $w^*$
in a similar way. When the amplitude is larger than $2\pi n$ where
$n<+\infty$ is the order of the ramification point, the Stolz angle is a
pointed metric ball.
}

\eno {Definition II.3.5.2 (Stolz continuity).}{ A map $f$ defined in $\cl S$ extends Stolz continuously to $\cl
S^*$ if is has limits at any $w^* \in \cl R$ along any Stolz angle with vertex at $w^*$. Then the limit at $w^*$
is unique and is the value of the Stolz continuous extension of $f$. }

\stit {II.3.6) Functions holomorphic at ramification points.}

\eno {Definition II.3.6.1 (Holomorphic function at a ramification point)}{ Let $w^* \in \cl S^* -\cl S$ be a
ramification point. A holomorphic function $f$ defined on a metric neighborhood $U\subset \cl S$ of $w^*$ is
holomorphic at $w^*$ if it has a Stolz continuous extension to $w^*$, i.e. when $w\to w^*$ in a Stolz angle, then
$f(w)$ converges to a well defined value $f(w^*)$. }

\stit {Remark.}

With this definition and using Riemann removability theorem, we see
that the notion of being holomorphic at a finite ramification point
is the classical one.

\stit {II.3.7) General Weierstrass theorem.}

The next theorem generalizes Weierstrass' classical theorem and follows from the fact that a uniform limit of
continuous functions is continuous.

\eno {Theorem II.3.7.1}{ Let $(f_n)$ be a sequence of holomorphic functions defined in an open set for the metric
topology $U \subset \cl S^*$ and converging uniformly on compact sets and on Stolz angles to a function $f: U \to
\dd C$. Then $f$ is holomorphic on $U$, in particular at the ramification point in $U$, $U\cap (\cl S^*-\cl S)$. }

We have even the following stronger result:

\eno{Theorem II.3.7.2}{ We consider a sequence of pointed log-Riemann surfaces $(\cl S^*_n, z_n)$ where $z_n \in
\cl S_n$, and a sequence of holomorphic functions $f_n$ defined on metric open sets $U_n \subset \cl S^*_n$. We
assume that $(\cl S^*_n, z_n) \to (\cl S^*, z)$ and $U_n \to U$ in Caratheodory sense, where $U$ is a metric open
set of the log-Riemann surface $\cl S^*$. We assume that the sequence $(f_n)$ converges uniformly on compact sets
of $U$ (in an obvious sense) to $f$. Then $f$ is holomorphic on $U$.}

%% file: ii4.tex
\input defv7.tex
\stit {II.4) Quasi-conformal theory of log-Riemann surfaces.}

\stit {II.4.1) Complex diffeomorphisms of log-Riemann surfaces.}

%We assume in this section II.4 that the log-Riemann surfaces $\cl S$ considered
%not only do have a discrete ramification set $\cl R$ but that this set is uniformly
%separated in the following sense: There exists $\eta >0$ such that
%for any two distinct ramification points $w^*_1, w^*_2 \in \cl R$ we have
%$$
%d(w_1^*, w_2^*) \geq \eta >0 \ .
%$$
%The minimal number $\eta >0$ is called the separation distance between
%points of $\cl R$.

We assume in this section II.4 that the log-Riemann surfaces that we
consider are finite, that is, the ramification set $\cl R$ is
a finite set. We assume also that the finitely completed Riemann
surface $\cl S^\times$ is simply connected. Then $\cl S^\times$
is parabolic as seen in section II.1.

We define first what it is to be univalent at ramification points.

\eno {Definition II.4.1.1}{ Let $\cl S_1$ and $\cl S_2$ be two log-Riemann surfaces. Let $w_1^* \in \cl S_1^*-\cl
S_1$. Let $\varphi : U_1\subset \cl S_1 \to \cl S_2$ be a germ of complex diffeomorphism defined in a punctured
neighborhood $U_1 \subset \cl S_1$ of $w_1^*$. The map $\varphi$ is said to be univalent at the ramification point
$w_1^*$ if $\varphi$ extends continuously to $w_1^*$ taking values in $\cl S_2^*$, and the extension is
bi-Lipschitz with respect to $w_1^*$, meaning that there exists $L \geq 1$ such that for all $w\in \cl S_1$ in a
neighborhood of $w_1^*$, $$ L^{-1} d(w,w_1^*) \leq d(\varphi(w),\varphi (w_1^*)) \leq L d(w,w_1^*) $$ (the metrics
here being the log-euclidean metrics on $\cl S_1^*$ and $\cl S_2^*$). }

\medskip

 We will prove below that in this case the image $\varphi (w_1^*) =w_2^*$ must be a
ramification point of the same order.

Note that for a complex diffeomorphism $\varphi : \cl S_1 \to \cl S_2$
the derivative $\varphi'=(\pi_2\circ \varphi\circ \pi_1^{-1})'$
computed in charts is well defined and independent
of the choice of the log-charts.

\medskip

\eno{Proposition II.4.1.2}{ If $\varphi$ is univalent at a ramification point then the derivative $\varphi'$ is
bounded in a neighborhood of the ramification point.}

\medskip

\stit{Proof.} We assume for convenience that $\pi_1(w_1^*) = \pi_2(w_2^*) = 0$. Take a neighborhood
$B(w_1^*,\epsilon)$ of $w_1^*$ small enough not to contain any other ramification points. For $w$ in $B(w_1^*, {1
\over 2}\epsilon)$ we can, taking the circle $C$ centered at $w$ and of radius equal to ${1 \over 2} |w|$,
estimate $\varphi'$ using Cauchy's integral formula, $$ \varphi'(w) = \int_C {\varphi(u) \over (u - w)^2} \ {du
\over 2\pi i} \ . $$ The Lipschitz condition on $\varphi$ gives $$\eqalign{ |\varphi(u)| & \leq L|u| \ , \cr
             & \leq L \left(|w| + {1 \over 2}|w| \right) \hbox{ for } u \in C \ , \cr
}
$$
which gives on substituting in the integral the estimate
$$
|\varphi'(w)| \leq 3L \ .
$$
$\diamond$

\medskip

Next we define the same univalence condition at $\infty$ (or at the
ramification points at $\infty$).

\eno {Definition II.4.1.3}{ With the same assumptions as before, assuming $\varphi$ defined in a region $\cl
S_1-\pi_1^{-1} (B(0,R_0))$ for some $R_0 >0$, the map $\varphi$ is univalent at $\infty$ if there is $R>0$ and a
constant $L\geq 1$ such that for $w\in \cl S_1^*$, $|\pi_1(w)| \geq R$, $$ L^{-1}|\pi_1(w)| \leq
|\pi_2(\varphi(w))| \leq L|\pi_1(w)| $$ }

\medskip

We have the weaker notion of continuity.

\eno {Definition II.4.1.4}{ We consider two log-Riemann surfaces $\cl S_1$ and $\cl S_2$ with two base points $z_1
\in \cl S_1$ and $z_2 \in \cl S_2$. A homeomorphism $\varphi: \cl S_1 \to \cl S_2$ is continuous at $\infty$ if
for any $R_2>0$ there exists $R_1 >0$ such that $$ \dd C-\pi_1^{-1} (B(z_1,R_1)) \subset \varphi^{-1} (\dd
C-\pi_2^{-1} (B(z_2, R_2))) \ . $$ }

\stit {Remark.}

The continuity at $\infty$ is independent of the base points chosen $z_1
\in \cl S_1$ and $z_2 \in \cl S_2$.

\medskip

\eno {Proposition II.4.1.5}{ A univalent map at a ramification point $w_1^*$ preserves the order of the
ramification point, that is, $w_1^*$ has the same order as $\varphi (w_1^*)$.}

\stit {Proof.} If $w_1^*$ is an infinite ramification point, then its
image must be an infinite ramification point since pointed neighborhoods
of a finite ramification point and of infinite ramification points are
not homeomorphic (the latter are simply connected whereas the former are not).
Thus assume that $w_1^*$ is order $n<+\infty$ and $\varphi (w_1^*)$
is of order $m <+\infty$. Taking local uniformizations of pointed
neighborhoods the expression of $\varphi$ in these new coordinates
becomes
$$
\psi(z) =\left (\varphi (z^n)\right )^{1/m} \ .
$$
The map $\psi$ is a local diffeomorphism at $0$, $\psi (0)=0$, and
its derivative is
$$
\psi'(z)={n\over m} z^{n-1} \varphi'(z^n) \left ( \varphi (z^n) \right )^{1/m-1} \ .
$$
By the previous Proposition, $\varphi'$ is bounded near $w_1^*$
and by Riemann's removability Theorem it extends analytically to
$w_1^*$. The derivative $\varphi'(0)$ takes a finite and non-zero
value applying the same argument to $(\varphi^{-1})'$. Thus when
$z\to 0$ we have
$$
\psi'(z) \sim {n\over m} (\varphi'(0))^{1/m} z^{{n\over m} -1} \ .
$$
Therefore $n$ must be equal to $m$ in order that $\psi'(0)$ is non-zero and
finite.$\diamond$

\medskip

\eno {Definition II.4.1.6}{ Let $\cl S_1$ and $\cl S_2$ be two finite log-Riemann surfaces. A complex
diffeomorphism between the log-Riemann surfaces $\cl S_1$ and $\cl S_2$ is a complex diffeomorphism $$ \varphi:
\cl S_1 \to \cl S_2 \ , $$ between the underlying Riemann surfaces that is univalent at the ramification points
and at $\infty$. We also call $\varphi$ a univalent map between log-Riemann surfaces. }

\medskip

\eno {Proposition II.4.1.7}{ A complex diffeomorphism between log-Riemann surfaces is globally bi-lipschitz with
respect to the ramification set $\cl R_1$.}

\stit {Proof.}

Just observe that removing a neighbourhood of $\infty$ and
neighbourhoods of each ramification point leaves a bounded set in
the log-Riemann surfaces for the log-euclidean metric that is
bounded away from the ramification sets. Thus $\varphi$ is
bi-Lipschitz with respect to $\cl R_1$ (maybe with a larger
constant than the local constants at the ramification points and
at $\infty$).$\diamond$

From the previous results we get:

\eno {Proposition II.4.1.8}{ A complex diffeomorphism extends continuously to a homeomorphism $\varphi : \cl S_1^*
\to \cl S_2^*$ preserving the order of ramification points.

The inverse of a complex diffeomorphism from $\cl S_1$ to $\cl S_2$
is a complex diffeomorphism from $\cl S_2$ to $\cl S_1$.
}

The main theorem of this section is the following rigidity result.

\eno {Theorem II.4.1.9}{ Let $\varphi: \cl S_1 \to \cl S_2$ be a complex diffeomorphism between finite log-Riemann
surfaces such that the finitely completed Riemann surfaces $\cl S_1^\times$ and $\cl S_2^\times$ are simply
connected. Then $\varphi$ preserves the fibers of $\pi_1$ and $\pi_2$ and $\cl S_1$ and $\cl S_2$ are in the same
affine class. Indeed the expression of $\varphi$ in charts is the same affine map.

Therefore if $\varphi$ is normalized such that $\varphi(z_1)=\varphi(z_2)$
with $z_1\in \cl S_1$ and $z_2\in \cl S_2$ such that
$$
\eqalign {
\pi_1(z_1) &=\pi_2 (z_2) =0 \ , \cr
\varphi' (z_1) &=\varphi'(z_2) =1 \ , \cr
}
$$
then $\varphi=\id$.
}

Thus log-Riemann surfaces have this remarkable rigidity property.

The proof of this Theorem follows from the following Lemma.

\eno {Lemma II.4.1.10}{ Let $\varphi : \cl S_1 \to \cl S_2$ be a complex diffeomorphism of finite log-Riemann
surfaces. Then $\varphi'$, the derivative of $\varphi$, is uniformly bounded $$ ||\varphi'||_{C^0(\cl S_1)}
<+\infty \ . $$ }

\stit {Proof.}

Since $\varphi$ is univalent at the ramification points, $\varphi'$
is bounded near the ramification points, so we choose $\d > 0$
such that $\varphi'$ is bounded in a $\d$-neighbourhood of the
ramification set $\cl R_1$ of $\cl S_1$.

%Let $\eps >0$. Choose $\d >0$ such that for any $w_1^* \in \cl R_1$
%we have,

%$$
%\varphi (B(w_1^* , \d ) ) \subset B(\varphi (w_1^*) , \eps ) \ .
%$$
%and also $0<\eps_1 < \eps$ such that
%$$
%B(\varphi (w_1^*), \eps_1) \subset \varphi (B(w_1^* , \d ) \ .
%$$
%Thus by the subordination principle (comparing with the two linear
%maps near $w_1^*$) for any $w\in \cl S_1^*$ with
%$d(w, w_1^*) < \d$ we have
%$$
%{\eps_1 \over \d } |\pi_1 (w) | \leq |\pi_2 \circ \varphi (w) -\pi_2 \circ
%\varphi (w_1^*) | \leq {\eps \over \d } |\pi_1 (w)| \ .
%$$
%Thus we get that for any $w \in B(w_1^* , \d )$,
%$$
%|\varphi'(w) | \leq {\eps \over \d} \ .
%$$
%We can choose $\d >0$ uniform for the ramification set $\cl R_1$ and
%this
%estimate will hold in the $\d$-neighborhood of $\cl R_1$.

\medskip

We choose $R > 0$ such that $\cl S_1^* -\pi^{-1}(B(0,R))$ does
not contain any points of $\cl R_1$. Then for $|\pi_1(w)| > 2R$
we can take the circle $C$ in $\cl S_1$ with center $w$ and
radius ${1 \over 2}|\pi_1(w)|$,
and estimate $\varphi'(w)$ using Cauchy's formula
$$
\varphi'(w) = {1 \over 2\pi i} \int_C {\varphi(u) \over (u - w)^2} \ du
$$
Since $\varphi$ is univalent at infinity, we have
$$
|\varphi(u)| \leq L|u| \leq L \left(|w| + {1 \over 2}|w| \right)
$$
which gives the estimate
$$
|\varphi'(w)| \leq 3L
$$

\medskip

Finally, on the complement of the $\d$-neighborhood of $\cl R_1$ and $\cl S_1 - \pi_1^{-1}(B(0,2R))$, we have $$
|\pi_2(\varphi(w))| \leq L|\pi_1(w)| \leq L\cdot 2R, $$ so $\pi_2(\varphi(w))$ is bounded. Hence again by Cauchy's
formula, considering circles of uniform radius ${1 \over 2}\d$ around each point avoiding $\cl R_1$, we get the
uniform boundedness of $\varphi'$.$\diamond$

\stit {Proof of the Theorem.}

Since $\varphi'$ is bounded, it extends to a holomorphic bounded function defined on the
finitely completed Riemann surface $\cl S_1^\times$. Since $\cl S_1^\times$ is parabolic
it follows from Liouville's Theorem that
$\varphi'$ is constant, and hence that the expression for $\varphi$ in local charts is an
affine map,
$$
\pi_2 \circ \varphi(w) = a \, \pi_1(w) + b.
$$
The result follows. $\diamond$

\stit {II.4.2) Teicm\"uller distance.}

We continue to make the same assumptions on log-Riemann surfaces.

We define the class of quasi-conformal homeomorphisms  between
completed log-Riemann surfaces.

\eno {Definition II.4.2.1}{ Let $L\geq 1$ and  $\cl S_1$ and $\cl S_2$ be two log-Riemann surfaces. A
homeomorphism $\varphi: \cl S_1^*\to \cl S_2^*$ is $L$-bi-lipschitz if for any $w, w' \in \cl S_1$ $$ L^{-1}
d(w,w') \leq d(\varphi (w), \varphi (w')) \leq L \  d(w,w') \ . $$ We say that $\varphi$ is $L$-bi-lipschitz at
the ends at infinite if for any $w\in \cl S_1$, $$ L^{-1} |\pi_1 (w)| \leq |\pi_2 (\varphi (w))| \leq L \ |\pi_1
(w)| \ . $$ If $\varphi $ extends continuously to the completions, $\varphi : \cl S_1^* \to \cl S_2^*$, we say
that $\varphi$ is $L$-bi-lipschitz at the ramification set if for any $w\in \cl S_1$ and $w^*\in \cl R_1$, $$
L^{-1} d(w,w^*) \leq d(\varphi (w), \varphi (w^*)) \leq L \ d(w,w^*) \ . $$ }

%\eno {Definition.}{ Let $\cl S_1$ and $\cl S_2$ be two log-Riemann surfaces.
%A homeomorphism $\varphi: \cl S_1^*\to \cl S_2^*$ has a derivative at
%$\infty$  if the following limit exists and is non-zero nor infinite,
%when $w\to \infty_i^*$
%$$
%\varphi'(\infty_i^*) =\lim_{w\to w^*} {\pi_2\circ \varphi (w)
%\over \pi_1(w) } \ .
%$$
%}

\eno {Definition II.4.2.2}{ Let $\cl S_1$ and $\cl S_2$ be two log-Riemann surfaces. Given $1\leq K <+\infty$ and
$L\geq 1$ we define a $(K,L)$-quasi-conformal homeomorphism $\varphi : \cl S_1^*\to \cl S_2^*$ as a classical
$K$-quasi-conformal homeomorphism of the underlying Riemann surfaces, such that $\varphi$ extends to a
homeomorphism between the completed log-Riemann surfaces $\varphi : \cl S_1^*\to \cl S_2^*$ which is
$L$-bi-lipschitz at the ramification points and at the ends at $\infty$.

We denote by ${\hbox {\rm QCH}} (\cl S _1, \cl S_2)$ the space of such
quasi-conformal homeomorphisms.
}

\stit {Remark.}

The Lipschitz condition is new compared to classical quasi-conformal theory on Riemann surfaces (where we have no
ramification points and no need of a Lipschitz condition). As we will see it is a fundamental assumption that is
needed in the proof of the main results. Observe, for example, how the Lipschitz condition is used in a
fundamental way in the proof of the quasi-invariance of the degree of finite ramification points that follows
(even if the Lipschitz constant does not appear in the estimate).

The Lipschitz condition introduces some differences with the classical theory. For example, given $L\geq 1$ there
are log-Riemann surfaces that are homeomorphic but not $(K,L)$-quasi-conformal equivalent for any $K\geq 1$ (just
consider two Gauss log-Riemann surfaces of log-degree $3$, one of them having the third ramification point far
away).

\medskip

\eno {Theorem II.4.2.3}{ Consider a germ of $(K,L)$-quasiconformal homeomorphism mapping a finite ramification
point of order $n\geq 1$ into another finite ramification point of order $m\geq 1$. We have $$ K\geq \max \left
({n\over m} , {m\over n} \right ) \geq 1 \ . $$ Moreover this estimate is sharp. }

\eno {Corollary II.4.2.4}{ Finite ramification points and their order are invariant under complex diffeomorphisms
of log-Riemann surfaces.}

\stit {Proof.}

Let $w_1^* \in \cl S_1^*$ and $w_2^* =\varphi (w_1^*)\in \cl S_2$. Let $r>0$ small and consider an annulus
$A_r=B(w_1^*, R)-{\overline {B(w_1^*, r)}}$ with $R\approx 1$ universal constant. Let $\rho =\min_{d(w, w_1^*)=r}
d(\varphi (w), w_1^*)$. Then $$ {1\over m} {1\over 2\pi} \log (1/\rho ) \leq K \mod \varphi (A_r)=K {1\over n}
{1\over 2\pi } \log (1/r) +C \ . $$ Therefore $$ r^{K/n} \leq C_0 \rho^{1/m} \ . $$ But the Lipschitz property
gives $$ \rho \leq L r \ , $$ thus $$ r^{K/n-1/m} \leq C_1=C_0 L^{1/m} \ . $$ This should hold for all $r>0$ small
enough. We conclude that $K\geq n/m$. Similarly considering the inverse mapping we get $K\geq m/n$.

\medskip

We prove now that the estimate is sharp. We consider the local map of the form
(taking as local coordinate $\pi(w)$ and writing $w$ instead)
$$
\varphi (w)=\sqrt {w \bar w } \left ( {w \over \bar w }\right )^{m/2n}
=w^{1/2+m/2n} \bar w^{1/2-m/2n} \ ,
$$
that is
$$
\varphi (re^{i\t })=re^{i{m\t \over n}} \ .
$$
Then we compute
$$
\eqalign {
\bar \partial \varphi &={1\over 2} \left (1-{m\over n}\right )
w^{1/2+m/2n} \bar w^{-1/2-m/2n} \ , \cr
\partial \varphi &={1\over 2} \left (1+{m\over n}\right )
w^{-1/2+m/2n} \bar w^{1/2-m/2n} \ . \cr
}
$$
Finally
$$
\mu (w)={\bar \partial \varphi \over \partial \varphi }
={1-m/n \over 1+m/n} {w\over \bar w} =-{1-n/m \over 1+n/m} {w\over \bar w}\ ,
$$
and $K=\max (m/n, n/m)$.
$\diamond$

As in the classical case given $K,L\geq 1$ we can extract converging subsequences
from a sequence $(\varphi_n)$ of normalized $(K,L)$-quasi-conformal homeomorphisms.

\eno {Theorem II.4.2.5}{ Let $K,L\geq 1$ and $\cl S_1$ and $\cl S_2$ two log-Riemann surfaces. We consider the
space ${\hbox {\rm QCH}}_{K,L} (\cl S _1, \cl S_2)$ of $(K,L)$-quasi-conformal homeomorphisms between $\cl S_1$
and $\cl S_2$. Any sequence $(\varphi_n) \subset {\hbox {\rm QCH}}_{K,L} (\cl S _1, \cl S_2)$ of normalized
homeomorphisms such that $$ \eqalign { \varphi_n (w_1) &=w_2 \cr \varphi (w_1') &=w'_2 \cr } $$ where $w_1,
w'_1\in \cl S_1$ and $w_2, w'_2\in \cl S_2$ are given points, contains subsequences converging uniformly on
compact sets to homeomorphisms in the same class. }

\stit {Proof.}

By the classical result (see [Le-Vi] or [Ahl2]) the $K$-quasi-conformal
homeomorphisms $(\varphi_n)$ have converging subsequences with limits that are
$K$-quasi-conformal. The $L$-bi-lipschitz condition is closed under
pointwise convergence. Moreover it implies that the homeomorphisms do extend
continuously to $\cl S_1^*$. This proves the result.$\diamond$.

\medskip

Now the main problem that we face in order to generalize the classical results is that a sequence minimizing $K$
may have Lipschitz constants $L$ diverging to $+\infty$. For each fixed $L$ however, we may consider the class of
quasi-conformal $L$-bi-Lipschitz homeomorphisms
$$
{\hbox {\rm QCH}}_{L} (\cl S _1, \cl S_2) := \bigcup_{K \geq 1} {\hbox {\rm QCH}}_{K,L} (\cl S _1, \cl S_2)
$$
and define a distance between affine classes as follows:

%\eno {Theorem II.4.2.6}{ Let $\cl S_1$ and $\cl S_2$ be two log-Riemann surfaces. Let $L_0\geq 1$ be such that
%there exists a $L_0$-bi-lipschitz homeomorphism between $\cl S_1^*$ and $\cl S_2^*$. Let $\varphi$ be a $(K,
%L)$-quasi-conformal homeomorphism between $\cl S_1$ and $\cl S_2$. Then there exists a $(K^2,
%2L_0)$-quasi-conformal homeomorphism between $\cl S_1$ and $\cl S_2$. }

%\stit {Proof.}

%The proof consists in modifying $\varphi$ near the ramification points.
%Note that in order to check that a homeomorphism is $2L$ bi-lipschitz
%near a ramification point it is enough to have that it is $L$ lipschitz in a sequence of
%circles centered at this ramification points with radius $2^{-n}$, $n=0,1,2,\ldots $
%Let $(C_n)$ be that sequence of circles and $(\varphi (C_n))$ their images.
%We can construct a local diffeomorphism $\psi$ that is the identity on $\varphi (C_0)$
%and that maps each $\varphi (C_n)$ into an annulus bounded by circles of radius
%$2^{-n}C$ and $2^{-(n+1)}C$ (with $C$ universal constant).

% proof to be completed

\medskip

\eno {Theorem II.4.2.7}{ Let $L \geq 1$. Given two log-Riemann surfaces $\cl S_1$ and $\cl S_2$, we define a
distance between their affine classes by $$ d_L \left ( \cl S_1, \cl S_2 \right ) := \inf_{\varphi \in {\hbox {\rm
QCH}}_{L} (\cl S_1, \cl S_2)} \log K(\varphi ) \ . $$

For each $L$ this defines a distance between affine classes.
}

\stit {Proof.}

The symmetry is obvious. The triangular inequality follows from the fact that a composition of bi-lipschitz
homeomorphisms at the ramification sets and at infinite is a homeomorphism of the same type, and the classical
fact that the composition of a $K_1$-quasi-conformal homeomorphism with a $K_2$-quasi-conformal homeomorphism is a
$K_1K_2$-quasi-conformal homeomorphism. It remains to prove that if $d_L(\cl S_1, \cl S_2)=0$ then $\cl S_1$ and
$\cl S_2$ are in the same affine class. Let $(\varphi_n)$ be a sequence in ${\hbox {\rm QCH}}_{L} (\cl S _1, \cl S_2)$
with $K(\varphi_n) \to 1$. Each $\varphi_n$ is a $(K_n,L)$-quasi-conformal homeomorphism. We can extract a converging
subsequence of $(\varphi_n)$ converging to $\varphi$ that will be $L$-bi-lipschitz at the ramification set and at
$\infty$ and will be a complex
diffeomorphism from $\cl S_1$ to $\cl S_2$. From the main Theorem II.4.1.9 of the previous section we get that
$\varphi$ preserves fibers and the affine class of $\cl S_1$ coincides with the affine class of $\cl
S_2$.$\diamond$

The same arguments show the existence of an $K$-extremal map in the
$L$-bi-lipschitz class.

\eno {Theorem II.4.2.8}{ Let $L \geq 1$ such that there exists a $(K,L)$-bi-lipschitz homeomorphism between
$\cl S_1$ and $\cl S_2$. Then there exists a $L$-bi-lipschitz homeomorphism $\varphi$ with minimal dilatation
$K(\varphi )$. }

% \stit {II.4.3) Teichm\"uller space of a finite log-Riemann surface.}

%% file: ii5.tex
\input defv7.tex
\input psfig.sty

\stit {II.5) Uniformization of parabolic log-Riemann surfaces of finite type.}

\stit {II.5.1) A fundamental example.}

\medskip

Fix $d \geq 2$ and let
$$
F(z) = \int_0^z e^{\xi^d} \, d\xi, z\in \dd{C}.
$$
In this section we show how to associate to $F$ a log-Riemann
surface $\cl S$ so that the uniformization of $\cl S$ is
realised by $F$, ie the expression for the uniformization in
log-charts is given by $F$.

\medskip

We recall the Gauss log-Riemann surface of log-degree $d$,
introduced in example 8 of section I.1.2, with $d$ ramification
points of infinite order placed at the $d$th roots of unity
in a common base sheet.

\medskip

Let the log-Riemann surface $\cl S$ be the log-Riemann surface
affine equivalent to the Gauss log-Riemann surface of
log-degree $d$ via a dilatation, so that the ramification points
of $\cl S$ are placed in the base sheet at the $d$ points
$a_1, \, a_2, \, \dots , \, a_d$ given by
$$\eqalign{
a_1 & = e^{i \pi/d} \int_0^\infty e^{-t^d} dt \cr
a_j & = (e^{2 \pi i /d})^{j-1} a_1 \ , \ j=2, \dots ,d \cr
}$$
Let $\pi : \cl S \to \dd C$ be the projection mapping.
Then we have

\medskip

\eno {Theorem II.5.1.1} { The mapping $F : \dd{C} \to \dd{C}$ 'lifts' to a biholomorphic map $\tilde{F} : \dd{C}
\to \cl S$ such that $F = \pi \circ \tilde{F} \,$. The lift $\tilde{F}$ maps $0 \in \dd C$ to the point $0$ in the
base sheet of $\cl S$. }

\bigskip

{\hfill {\centerline {\psfig {figure=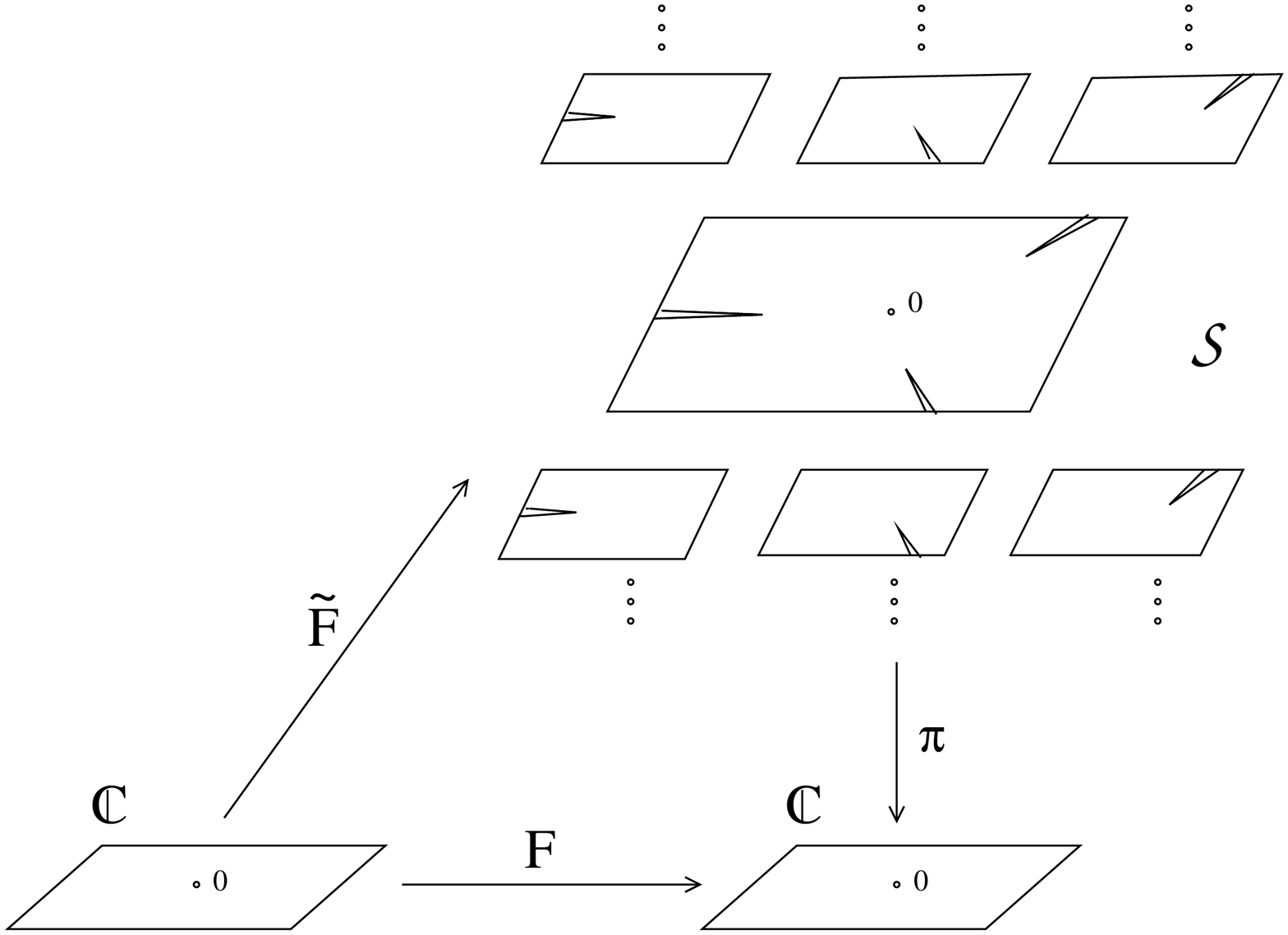,height=6cm}}}}

\bigskip

The proof consists of partitioning the plane into disjoint simply
connected domains such that $F$ maps each univalently to the trace
of a sheet of the minimal atlas, and their boundaries to the cuts
along which the sheets are joined.

\medskip

We consider the level curves $\{ $ Im $F =$ constant $\}$. These
are integral curves for the vector field
$$ X(z) := e^{-i \, \hbox{\sevenrm Im} \, z^d} \, \, \, ,z \in \dd{C}$$
since for any integral curve $Z = Z(t)$,
$$\eqalign{
{d \over dt} F(Z(t)) & = F'(Z(t))\,Z'(t) \cr
                       & = e^{Z(t)^d} \, e^{-i\, \hbox{\sevenrm Im} \, (Z(t)^d)} \cr
                       & = e^{\hbox{\sevenrm Re} \, (Z(t)^d)} \cr
                       & \in \dd{R}_+ .\cr
}$$
We make the following

\medskip

{\bf Observations :}

\medskip

{\bf 1.} $|X| = 1$, so $X$ is nonsingular and
integral curves of $X$ through any initial point
in the plane are defined for all time (they cannot
explode in finite time).

\medskip

{\bf 2.} Im $F$ is constant along $(Z(t))_{t \in \dd{R}}$,
while Re $F$ is strictly increasing in $t$.

\medskip

{\bf 3.} Since $X$ is nonsingular on the whole plane,
$X$ cannot have any limit cycles, so every integral
curve $(Z(t))_{t \in \dd{R}}$ is simple,
and $|Z(t)| \to \infty$ as $|t| \to \infty$.

\medskip

{\bf 4.} $F(\overline{z}) = \overline{F(z)}$

\medskip

{\bf 5.} $F$ commutes with the rotation around 0 by an angle $2 \pi /d$,
since, denoting $\omega = e^{2\pi /d}$, we have
$$\eqalign{
F(\omega z) & = \int_0^{\omega z} e^{\xi^d} \, d\xi \cr
              & = \int_0^z e^{(\omega \tau)^d} \omega \, d\tau \hbox{ (putting } \omega \tau = \xi \hbox{ )}
\cr
              & = \omega  \int_0^z e^{\tau^d} \, d\tau \cr
              & = \omega F(z) \cr
}$$

\medskip

{\bf 6.} Observations {\bf 4.} and {\bf 5.} above imply that
$F$ also commutes with the reflections through
each of the lines $\{ $ arg $z = j \pi /d \} \, , j = 1, \dots , 2d$.
So it suffices to understand how $F$
maps the sector $\Pi$ = $\{ 0 \leq$ arg $z \leq \pi/d \}$.
We consider the foliation given by integral curves
to $X$ in this sector.

\medskip

We define :

\medskip

{\bf 1.} For $z_0 \in \dd{C}$, $(Z(t;z_0))_{t\in \dd{R}}$
to be the integral curve of $X$ starting at $z_0$,
ie $Z(0;z_0) = z_0$.

\medskip

{\bf 2.} For $k \geq 1$, the curves

$$ \Gamma_k = \{ z \in \Pi : \hbox{Im } z^d = k\pi \} $$

\medskip

{\bf 3.} The domains :

$$ D_0 = \{ Z(t;z_0) : z_0 \in \Pi_1, \hbox{Im } {z_0}^d = 0, t > 0 \} $$
$$ D_k = \{ Z(t;z_k(\theta)) : 0 < \theta < \pi/d, -\infty < t < \infty \} \, (\hbox{ integral curves
starting from points on } \Gamma_k) , \, k \geq 1 $$
$$ E_k = \{ z \in \Pi : k\pi < \hbox{Im } z^d < (k+1)\pi \} \, (\hbox{ the region in between the curves
}\Gamma_k \hbox{ and } \Gamma_{k+1}) \, , \, k \geq 1 $$

\medskip

The figure below gives an illustration of the curves and domains defined above.

\medskip

{\hfill {\centerline {\psfig {figure=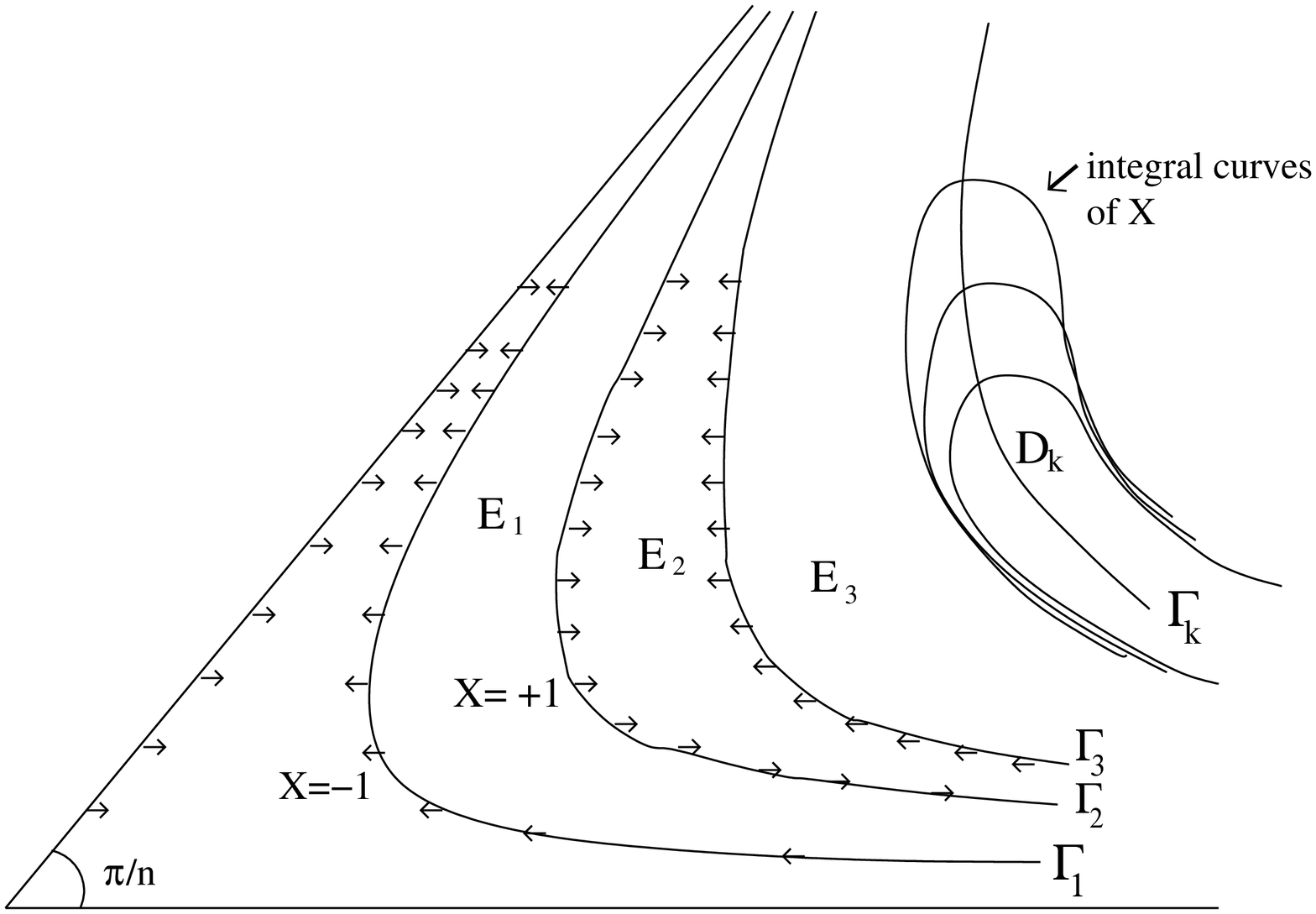,height=4cm}}}}

\bigskip

\eno{ Lemma II.5.1.2} { For $k \geq 1$, the domains $E_{2k}$ and $E_{2k+1}$ are 'trapping regions' for integral
curves of $X$ for positive and negative times respectively, ie $$ \hbox{If } z_0 \in E_{2k} \hbox{ then } Z(t;z_0)
\in E_{2k} \hbox{ for all } t > 0, \hbox{ and }$$ $$ \hbox{If } z_0 \in E_{2k+1} \hbox{ then } Z(t;z_0) \in
E_{2k+1} \hbox{ for all } t < 0. $$ }

\stit {Proof.}

The vector field $X = +1$ on $\Gamma_{2k}$ and -1 on
$\Gamma_{2k+1}$, and hence, at all points on the two curves,
points into the domain $E_{2k}$; so any integral curve starting
in the region $E_{2k}$ at time $t=0$ must stay in it for all
$t>0$ (whenever it gets near one of the boundary curves
$\Gamma_{2k}$, $\Gamma_{2k+1}$, it must flow back into the region
since $X$ points into $E_{2k}$ on and near
the boundary curves).

A similar argument shows that $E_{2k+1}$ is a trapping region for negative times. $\diamondsuit$

\medskip

\eno{ Lemma II.5.1.3} { The domains $D_k, k \geq 1$, are pairwise disjoint.}

\stit {Proof.}

First we observe that $D_{2k}$ and $D_{2k+1}$ are disjoint, ie no
integral curve starting at a point on $\Gamma_{2k}$ can intersect
an integral curve starting from a point on
$\Gamma_{2k+1}$; for, if two such curves were to meet, then it would
be possible to flow along in positive time from the starting point of
one to that of the other, ie from one
boundary point to another, which is impossible for a trapping region.

By a similar argument, domains $D_{2k-1}$ and $D_{2k}$ are disjoint;
so for $k\geq1$ the domains $D_{1,k}$
and $D_{1,k+1}$ are disjoint, and for the general case $j < k$, with
$j - k > 1$ say, the domains $D_{1,j}$
and $D_{1,k}$ are separated by the curve $\Gamma_{j+1}$ for example,
and hence disjoint. $\diamondsuit$

\medskip

\eno{ Lemma II.5.1.4} { Each curve $(Z(t;z_0))_{t \in \dd{R}}$ for $z_0 \in \Gamma_k$ gets mapped by $F$ to a full
horizontal line, ie $$F(\{ Z(t;z_0) : t \in \dd{R} \}) = \{ \hbox{ Im } w = \hbox{ Im } F(z_0) \}.$$ }

\stit {Proof.}

Let $z_0 \in \Gamma_{2k}$. Then for $t>0$,
$Z(t;z_0) \in E_{2k} = \{ z \in \Pi : 2k\pi < \hbox{Im }z^d <
(2k+1)\pi \}$, so we have
Im ${d \over dt}Z(t;z_0) = - \sin $ Im $(Z(t;z_0)^d) < 0$,
so Im $Z(t;z_0) < \hbox{Im } z_0$. From this and the fact that
$Z(t;z_0) \in E_{2k}$, it follows that arg $Z(t;z_0) < \hbox{arg
}z_0$. We can write arg $z_0 = \pi/d - \epsilon$ for some $\epsilon > 0$. Then
$$ \hbox{arg }Z(t;z_0)^d < d \, \hbox{arg }z_0 = \pi - d\epsilon, $$ so
$$\eqalign{
{d \over dt} \hbox{Re } F(Z(t;z_0)) = \hbox{Re }Z(t;z_0)^d & > {\hbox{Im }Z(t;z_0)^d \over - \tan
(d\epsilon)} \cr
                     & > {(2k+1)\pi \over - \tan (d\epsilon)} \cr
}$$
(since $2k\pi < \hbox{Im }Z(t;z_0)^d < (2k+1)\pi$ for $Z(t;z_0) \in
E_{2k}$). Since the lower bound is independent of $t>0$,
Re $F_n(Z(t;z_0)) \to +\infty$ as $t \to +\infty$, for any
$z_0 \in \Gamma_{2k}$.  Similarly one can show that
Re $F_n(Z(t;z_0)) \to -\infty$ as $t \to -\infty$ for
$z_0 \in \Gamma_{2k}$, and for $z_0 \in \Gamma_{2k+1}$ that Re $F_n(Z(t;z_0)) \to
+\infty$ and $-\infty$ as $t \to +\infty$ and $-\infty$ respectively,
 from which the result follows. $\quad \diamondsuit$

\medskip

\eno{Lemma II.5.1.5} { Let $(\gamma (t))_{t \in \dd{R}}$ be a curve such that $|\gamma(t)| \to \infty$ as $t \to
+\infty$. If for some $\epsilon > 0$ we have $(\pi/2 + \epsilon)/d \leq$ arg $\gamma(t) \leq (3\pi/2 -\epsilon)/d$
for all $t$ sufficiently large, then $$ F(\gamma(t)) \to a_1 = e^{i \pi/d} \int_0^\infty e^{-s^d} ds \, \hbox{ as
} t \to +\infty$$ }

\stit{Proof.}

Write $\gamma(t) = r e^{i \theta}$, where $r \to \infty$ as $t \to \infty$. Then
$$\eqalign{
F(r e^{i \theta}) & = F(r e^{i \pi/d}) + \bigl(F(r e^{i \theta}) - F(r e^{i \pi/d}) \bigr) \cr
 & = \int_{[0,r e^{i \pi/d}]} e^{z^d} dz + \int_{C} e^{z^d} dz \cr
}$$
where $C$ is the shorter arc of the circle joining $r e^{i \pi/n}$
and $r e^{i \theta}$. The first of the two integrals converges to $a_1$,
$$\eqalign{
\int_{[0,r e^{i \pi/n}]} e^{z^d} dz & = \int_0^r e^{{(se^{i \pi/d})}^d} e^{i \pi/d} ds \cr
                                    & = e^{i \pi/d} \int_0^r e^{-s^d} ds \cr
                                    & \to e^{i \pi/d} \int_0^\infty e^{-s^d} ds = a_1 \hbox{ as } t \to
+\infty \cr
}$$
while the second one tends to 0:
$$\eqalign{
| \int_{C} e^{z^d} dz | & \leq \bigl( \Max_{z \in C} |e^{z^d}| \bigr) \cdot (\hbox{length of } C) \cr
& \leq \bigl( \Max_{ (\pi/2 + \epsilon)/d \leq \phi \leq (3\pi/2 - \epsilon)/d } e^{-r^d \cos(n\phi)} \bigr)
\cdot 2\pi r  \cr
& \leq e^{-r^d \sin(\epsilon)} \cdot 2\pi r \to 0 \, \hbox{ as } t \to +\infty.  \quad \diamondsuit \cr
}$$

\medskip

\eno{ Lemma II.5.1.6} { There is a sequence $(\gamma_k)_{k \geq 1}$ of distinct integral curves of $X$ such that
$$
\partial D_k = \gamma_k \cup \gamma_{k+1} \, , \, k \geq 1$$ }

\stit {Proof.}

For $k \geq 1$, let ${D_k}^+ = \{ Z(t;z_0) : t>0, z_0 \in \Gamma_k \} ( \, \subset D_k)$.
Then the ${D_k}^+$'s are disjoint open sets, and moreover
${D_{2k}}^+, {D_{2k+1}}^+ \subseteq E_{2k}$ ($E_{2k}$ is a
trapping region for positive times); since $E_{2k}$ is connected,
it follows that $E_{2k} - ({D_{2k}}^+ \cup
{D_{2k+1}}^+)$ is nonempty.

So let $z^* \in E_{2k} - ({D_{2k}}^+ \cup {D_{2k+1}}^+)$,
and define $\gamma_{2k+1}$ to be the integral curve
through $z^*$.

\medskip

{\bf Claim.} $\gamma_{2k+1} = E_{2k} - ({D_{2k}}^+ \cup {D_{2k+1}}^+)$.

\medskip

{\bf Proof of Claim.} Since $z^*$ doesn't lie on any of the integral
curves starting from points on
$\Gamma_{2k}$ or $\Gamma_{2k+1}$, $\gamma_{2k+1}$ doesn't
intersect either ${D_{2k}}^+$ or ${D_{2k+1}}^+$, hence
$\gamma_{2k+1} \subseteq E_{2k} - ({D_{2k}}^+ \cup
{D_{2k+1}}^+)$ (note $\gamma_{2k+1} \subseteq E_{2k}$
since it can't intersect either boundary curve
$\Gamma_{2k}$ or $\Gamma_{2k+1}$ of $E_{2k}$).

\medskip

We make some observations on $\gamma_{2k+1}$ :

\medskip

{\bf (i) } Since $\gamma_{2k+1} \subseteq E_{2k}$, we have, as in the proof of Lemma II.5.1.4, that Im ${d \over
dt} \gamma_{2k+1} < 0$ for all $t$.

{\bf (ii) } By {\bf Obsvn 3.} made earlier,
$|\gamma_{2k+1}(t)| \to \infty$ as $|t| \to \infty$; since every
set of the form $E_{2k} \cap \{ m < \hbox{ Im } z < M \}, m,M>0$ is bounded,
$\gamma_{2k+1}(t)$ must leave
every such set as $|t| \to \infty$. In particular, it follows from this and
{\bf (i) } that Im
$\gamma_{2k+1}(t) \to +\infty$ as $t \to -\infty$.

{\bf (iii) } From {\bf (ii)} it follows that arg $\gamma_{2k+1}(t) \to \pi/d$ as $t \to -\infty$; so, by Lemma
II.5.1.5, $F(\gamma_{2k+1}(t)) \to a_1$ as $t \to -\infty$, and hence Im $F \equiv $ Im $a_1$ on $\gamma_{2k+1}$.

\medskip

We prove the inclusion $\gamma_{2k+1} \supseteq E_{2k} - ({D_{2k}}^+ \cup {D_{2k+1}}^+)$
 by contradiction; so let $z^{**} \in E_{2k} - ({D_{2k}}^+ \cup {D_{2k+1}}^+)$
 such that $z^{**} \notin \gamma_{2k+1}$. Then
$Z(\cdot; z^{**}) \subseteq E_{2k} - ({D_{2k}}^+ \cup {D_{2k+1}}^+)$ by the same
argument as for $\gamma_{2k+1}$; since $\gamma_{2k+1}, Z(\cdot;z^{**})$ are
simple disjoint curves contained in $E_{2k} - ({D_{2k}}^+ \cup {D_{2k+1}}^+)$
which both escape to infinity as $|t| \to \infty$, we can consider the
region $U \subseteq E_{2k} - ({D_{2k}}^+ \cup {D_{2k+1}}^+)$ bounded by these two
curves.

\medskip

For any $z \in U$, if $\gamma$ is the integral curve to $X$ through $z$, then
by the same arguments as in the case of $\gamma_{2k+1}$, we must
have Im $F \equiv $ Im $a_1$ on $\gamma$. But then Im $F \equiv $ Im
$a_1$ in all of $U$, a contradiction since $F$ is a nonconstant analytic function.
This proves the claim.

\medskip

Similarly we can define curves $\gamma_{2k-1}$ in the domains $E_{2k-1}$,
such that $\gamma_{2k-1} = E_{2k-1}
- ({D_{2k-1}}^- \cup {D_{2k}}^-)$. It is then straightforward to show
that $\partial D_k = \gamma_k \cup
\gamma_{k+1}, k \geq 1$, as required. $\quad \diamondsuit$

\medskip

{\bf Remark.} $\quad$ Similarly to the above we can show that $\gamma_1$
is a boundary curve of $D_0$. Thus we have the following complete
foliation of the sector $\Pi$ by integral curves of X :

$$ \Pi = \overline{D_0} \cup \bigcup_{k = 1}^\infty (D_k \cup \gamma_k) $$

(except for the 2 other boundary curves $\{ $ arg $z = 0 \}$ and $\{ $ arg $z = \pi/n \}$ of $D_0$, all other
curves above are integral curves of $X$.)

\medskip

\eno{ Lemma II.5.1.7} { For $k \geq 1$,
\medskip
(1) $F$ maps the curves $\gamma_{2k-1}$ and $\gamma_{2k}$ to the
half-lines $\{ $ Re $w > $ Re $a_1 $, Im
$w = $ Im $a_1 \}$ and $\{ $ Re $w < $ Re $a_1 $,
Im $w = $ Im $a_1 \}$ respectively.
\medskip
(2) $F$ maps the domains $D_{2k-1}$ and $D_{2k}$ univalently to
the half-planes $\{ $ Im $w > $ Im $a_1 \}$ and $\{ $ Im $w < $ Im $a_1 \}$
 respectively.
}

\stit{Proof.}

$(1)$. We observed in the proof of Lemma II.5.1.6 that $F(\gamma_{2k-1}(t)) \to a_1$
 as $t \to -\infty$; for $t > 0$, the proof of Lemma 3 applies to $\gamma_{2k-1}$
 to show that Re $F(\gamma_{2k-1}(t)) \to +\infty$ as $t
\to +\infty$. Since Re $F$ is strictly increasing and Im $F$ constant
on $\gamma_{2k-1}$, it follows that $F(\gamma_{2k-1}) =
\{ $ Re $w > $ Re $a_1 $, Im $w = $ Im $a_1 \}$.

Similarly one can show that
$F(\gamma_{2k}) =  \{ $ Re $w < $ Re $a_1 $, Im $w = $ Im $a_1 \}$.

$(2)$. It follows from Lemma II.5.1.4 that each domain $D_k$ gets mapped to a connected union of full horizontal
lines, thus either to the whole plane, or a half-plane or a horizontal strip. In either case, $F(D_k)$ is simply
connected; since $F$ is locally univalent ($F'(z) = e^{z^d} \neq 0$ everywhere), this implies that $F$ is in fact
univalent on $D_k$.

$D_k$ is not the whole plane, hence $F$ cannot map $D_k$ to the
whole plane. Considering the images of the boundary curves $\gamma_k$
and $\gamma_{k+1}$ as described in $(1)$ above, we see that
$F$ cannot map to a strip either, but must map to one of the two
half-planes $\{ $ Im $w < $ Im $a_1 \}$ or $\{ $ Im $w > $ Im
$a_1 \}$; exactly which of these two half-planes follows from
considering the orientation of the boundary curves $\gamma_k, \gamma_{k+1}$
 with respect to $D_k$ (since $F$ is orientation preserving). $\quad \diamondsuit$

\medskip

{\bf Remark.} $\quad$ Similarly one can show that $F$ maps $D_0$ univalently
to the domain bounded by the three straight lines $\{$ Re $w \geq 0$, Im $w = 0 \}$,
$\{$ Re $w \geq $ Re $a_1$, Im $w = $ Im $a_1 \}$ and
$\{$ arg $w = $ arg $\pi/n , 0 \leq |w| \leq |a_1| \}$.

\medskip

Define
$$ C_k = D_{2k-1} \cup \gamma_{2k} \cup D_{2k} \quad , k \geq 1 $$
It is immediate from the previous Lemma that $F$ maps each domain
$C_k$ univalently to the slit plane $\dd
C - \{ w : $ Im $ w = $ Im $ a_1$, Re $w \geq $ Re $a_1 \}$.
The figure below illustrates the domains $D_0, D_k, C_k, k \geq 1$
and their images under $F$. This gives a complete description of the mapping $F$ in
the sector $\Pi$.

\bigskip

{\hfill {\centerline {\psfig {figure=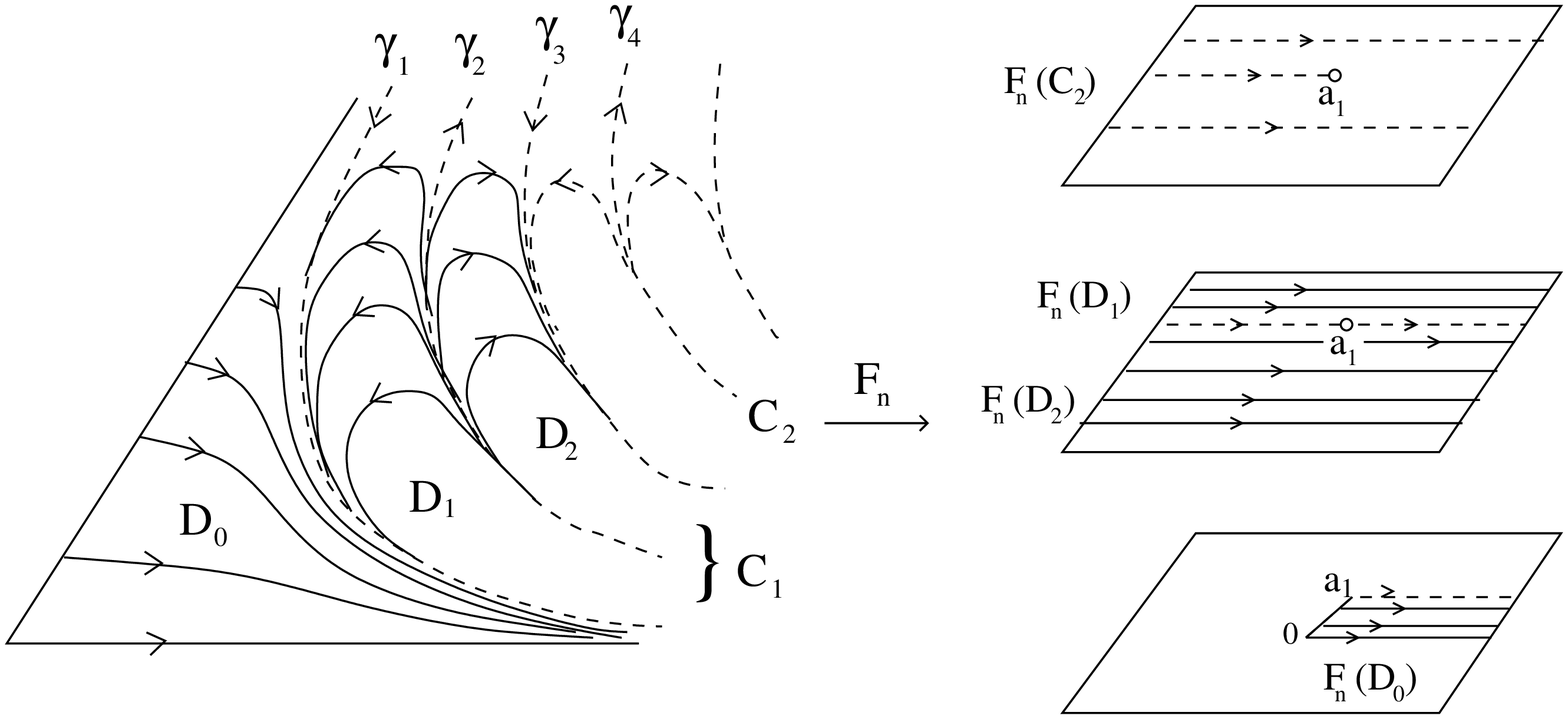,height=5cm}}}}

\bigskip

Clearly if we shift the bundaries of the domains $D_0$ and
$C_k, k \geq 1$, appropriately, we can obtain instead domains
$D^*_0$ and $C^*_k, k \geq 1$, such that

\medskip

$\bullet \quad D^*_0$ has 3 boundary curves, $\{ $ arg $z = 0 \},
\{ $ arg $z = \pi/d \}$, and a curve ${\alpha}^*_1$; $F$ maps
these to the straight lines $\{ $ arg $w = 0 \}, \{ $
arg $w = \pi/d, |w| < |a_1| \}$, and $\{ $ arg $w = \pi/d, |w| > |a_1 \}$
 respectively, and $D^*_0$ univalently to the sector $\{ 0 < $
arg $w < \pi/d \}$.

\medskip

$\bullet \quad $ Each $C^*_k, k \geq 1$, has 2 boundary curves,
${\alpha}^*_k$ and ${\alpha}^*_{k+1}$; $F$
maps both to the 'slit' $\{ $ arg $w = \pi/d , |w| > |a_1| \}$,
and $C^*_k$ univalently to the slit-plane
$\dd C - \{ $ arg $w = \pi/d, |w| \geq |a_1| \}$.

\medskip

Let $D^{**}_0, C^*_{-k}, {\alpha}^*_{-k+1}, k \geq 1$, be the
reflections through the line $\{ $ arg $z =
\pi/d \}$ of $D^*_0, C^*_k, {\alpha}^*_k, k \geq 1$, respectively.

\medskip

Let $A = \{ $ arg $z = 0 \} \cup D^*_0 \cup \{ $ arg $z = \pi/d \} \cup D^{**}_0 \{$ arg $z = 2\pi/d \}$.

\medskip

Define
$$\eqalign{
A_{0,0} & =  \bigl( \bigcup_{j=1}^d \omega^{j-1} \, A \bigr) \cup \{ 0 \}  \cr
A_{j,k} & = \omega^{j-1} \, C^*_k \quad \quad \quad \quad \quad \quad \quad 1 \leq j \leq d, \quad k \in \dd
Z - \{ 0 \}  \cr
\alpha_{j,k} & = \omega^{j-1} \, {\alpha}^*_k \quad \quad \quad \quad \quad \quad \quad 1 \leq j \leq d,
\quad k \in \dd Z  \cr
}$$

This gives the desired partition of the plane mentioned earlier:

$$ \dd C = A_{0,0} \cup \bigl( \bigcup_{j=1}^d \bigcup_{k \in \dd Z - \{0\} } A_{j,k} \bigr) \cup \bigl(
\bigcup_{j=1}^d \bigcup_{k \in \dd Z } \alpha_{j,k} \bigr) $$

\stit {Proof of Theorem II.5.1.1} $F$ is univalent on each domain in the above partition, mapping $A_{0,0}$ to the
trace of the 'base sheet', the $d$ families of domains $A_{j,k}$ to the traces of the $d$ families of clean
sheets, and the boundaries $\alpha_{j,k}$ to the traces of the cuts joining the sheets. Since $\pi$ is also
univalent
 on each sheet of the minimal atlas, $F$ has a unique lift
$\tilde{F} : \dd C \to \cl S$ and the lift is biholomorphic. $\quad \diamondsuit$

The figure below illustrates the correspondence between
the domains $A_{j,k}$ and the sheets of the minimal atlas under
$\tilde{F}$.

\bigskip

{\hfill {\centerline {\psfig {figure=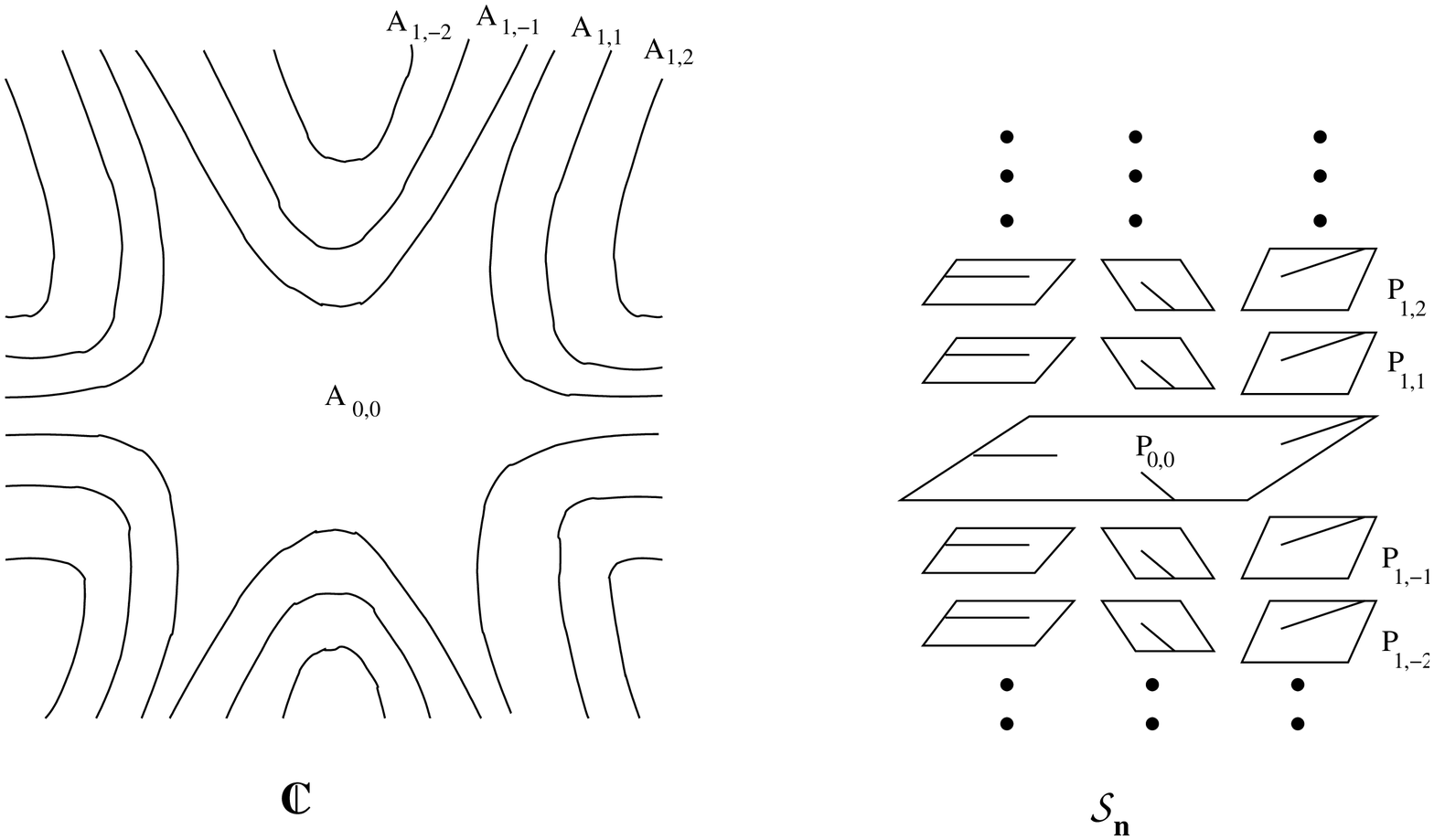,height=6cm}}}}

\bigskip

\stit{II.5.2) On the general case.}

\medskip

Let $P(z) = a_d z^d + \dots + a_0$ be a polynomial of degree $d$. We can use the same methods as above to analyze
the general case of the entire function $$ F(z) = \int_{0}^{z} e^{P(t)} \ dt $$ The Riemann surface $\cl S$ of the
inverse of $F$ can again be given a log-Riemann surface structure which admits a description very similar to that
considered above.

\medskip

\eno{Theorem II.5.2.1}{ There exists a log-Riemann surface $\cl S$ such that the function $F$ lifts to a
biholomorphic map $\tilde{F} : \dd C \to \cl S$ such that $\pi \circ \tilde{F} = F$. The log-Riemann surface $\cl
S$ is simply connected and contains exactly $d$ ramification points $w_1, \dots, w_d$, all of infinite order.
These project onto the points $$ w_j' = \pi(w_j) = \int_{0}^{\rho_j \cdot \infty} e^{P(z)} \ dz \ , \ j = 1,
\dots, d $$ where $\rho_1, \dots, \rho_d$ are the $d$ values of $(-a_d)^{-1/d}$. }

The main difference with the previous case lies in the locations of the $d$
ramification points as points in the surface $\cl S$; for
example, they need not all lie in a single base sheet, or be symmetrically
placed around $0$, but instead may be spread out over different sheets at
arbitrary positions in the sheets. In fact, and this is the content of
the following section II.5.3, any arbitrary arrangement of the $d$
ramification points may be achieved by a suitable choice of the polynomial
$P$.

\medskip

The proof of the above theorem is very similar to the previous case when
$P(z)$ was equal to $z^d$. In this case we consider the vector field
$$
X_P (z)  = e^{-i \, \hbox{\sevenrm Im} \, P(z)} \, \, \, ,z \in \dd{C}
$$
whose integral curves get mapped to horizontals by $F$.
This vector field is in general {\bf not} conformally
conjugate on any neighbourhood of $\infty$ to the vector field
$X(z) = e^{-i \, \hbox{\sevenrm Im} \, z^d}$. However on a neighbourhood
of infinity, we can make a change of variables $z = h(\xi) = c_1 \xi + c_0 +
c_{-1}/{\xi} + \dots$, such that $P(z) = \xi^d$. For convenience we
assume $P(z) = z^d + a_{d-1} z^{d-1} + \dots + a_0$ has leading coefficient
$1$ (this can always be achieved by a suitable of variables in the
integral defining $F$), so $h$ can be taken to be of the form
$h(\xi) = \xi + \dots$. We assume $h$ is defined for $|\xi| > R_0$.
We use this change of variables to study the vector field $X_P$ as follows.

\medskip

We define:

\medskip

{\bf 1.} $Z(.;z_0)$ to be the integral curve of $X_P$ starting from $z_0 \in \dd C$, $Z(0;z_0) = z_0$. As
before, all solutions are defined for all time $t \in \dd R$, and must escape to infinity as $|t| \to
\infty$.

\medskip

{\bf 2.} The sectors
$$
\Pi_j = \{ (j-1)\pi/d < \arg \xi < j\pi/d, \ |\xi| > R_0 \}, \ \Pi_j ' = h(\Pi_j), \ j=1, \dots, 2d.
$$

\medskip

{\bf 3.} The curves
$$
\Gamma_{j,k}(\alpha) = \{ \hbox{ Im } \xi^d = k\pi - \alpha, \ \xi \in \Pi_j \}, \ \Gamma_{j,k}'(\alpha) =
h(\Gamma_{j,k}(\alpha)),
$$
for $0 < \alpha < 2\pi$, $j = 1,\dots, 2d$ and $k \geq k_0$ is positive for $j$ odd, $k \leq -k_0$ is
negative for $j$ even, where $k_0$ is chosen large enough so that $\{ $ Im $\xi^d = k_0\pi - 2\pi \} \subset
\{ |\xi| > R_0 \}$. We note that $X_P(z) = X(\xi) = \pm e^{i\alpha}$ for $z \in \Gamma_{j,k}'(\alpha), \xi
\in
\Gamma_{j,k}(\alpha)$.

\medskip

{\bf 4.} The regions $E_{j,k}(\alpha)$ (resp. $E_{j,k}'(\alpha) = h(E_{j,k}(\alpha))$) to be the regions
bounded by $\Gamma_{j,k}(\alpha)$ and $\Gamma_{j,k+1}(\alpha)$ (resp. $\Gamma_{j,k}'(\alpha)$ and
$\Gamma_{j,k+1}'(\alpha)$) for $j$ odd, and by $\Gamma_{j,k}(\alpha)$ and $\Gamma_{j,k-1}(\alpha)$ (resp.
$\Gamma_{j,k}'(\alpha)$ and $\Gamma_{j,k-1}'(\alpha)$) for $j$ even.

\medskip

{\bf 5.} The domains
$$
D_{j,k}'(\alpha) = \{ Z(t;z_0) : t \in \dd R, \ z_0 \in \Gamma_{j,k}'(\alpha) \ \}
$$

\medskip

\eno{Lemma II.5.2.2}{ Fix $0 < \delta < \pi/2d$ and $\alpha_1, \alpha_2, \dots , \alpha_{2d}$ such that
$(j-1)\pi/d + \delta < \alpha_j < j\pi/d - \delta$. Then for $k_0$ large enough (depending on $\delta$), in each
sector $\Pi_j '$ the curves $\Gamma_{j,k}'(\alpha_j)$ are transverse to the vector field $X_P$.}

\medskip

\stit{Proof.} Let $(\xi(t))$ be a parametrization of a
curve $\Gamma_{j,k}(\alpha)$ and
$z(t) = h(\xi(t))$ the corresponding parametrization of
$\Gamma_{j,k}'(\alpha)$. It is easy
to verify that, independently of $\alpha$ and $k$, the
tangent vectors $\xi'(t), -\xi'(t)$ satisfy
$$
\arg(\xi'(t)),\ \arg(-\xi'(t)) \notin ((j-1)\pi/d, j\pi/d) \cup ((j-1)\pi/d + \pi, j\pi/d + \pi)
$$
Since $z'(t) = h'(\xi(t)) \cdot \xi'(t)$ and $h'(\xi) = 1 + O(1/\xi^2)$, it follows
that given $\delta > 0$ by choosing $k_0$ large enough
(and hence $|\xi|$ large enough),
the tangent vectors $z'(t), -z'(t)$ to $\Gamma_{j,k}'(\alpha)$ satisfy
$$
\arg(z'(t)),\ \arg(-z'(t)) \notin ((j-1)\pi/d + \delta, j\pi/d - \delta) \cup ((j-1)\pi/d + \pi + \delta,
j\pi/d + \pi - \delta).
$$
Since $X_P = \pm e^{i\alpha}$ on $\Gamma_{j,k}'(\alpha)$, it
follows that by choosing $\alpha_j$ such that
$(j-1)\pi/d + \delta < \alpha_j < j\pi/d - \delta$, the curves
 $\Gamma_{j,k}'(\alpha_j)$ are transverse to the vector field $X_P$.
$\diamondsuit$.

\medskip

Using this lemma we can now proceed as in the previous section. The following lemmas
follow from similar arguments as before:

\medskip

\eno{Lemma II.5.2.3}{ The domains $E_{j,2k}'(\alpha_j)$ and $E_{j,2k+1}'(\alpha_j)$ are 'trapping regions' for
integral curves of $X_P$ for positive and negative times respectively.}

\medskip

It follows that

\medskip

\eno{Lemma II.5.2.4}{ The domains $D_{j,k}'(\alpha_j)$ are pairwise disjoint.}

\medskip

\eno{Lemma II.5.2.5} { Each curve $(Z(t;z_0))_{t \in \dd{R}}$ for $z_0 \in \Gamma_{j,k}'(\alpha_j)$ gets mapped by
$F$ to a full horizontal line, ie $$F(\{ Z(t;z_0) : t \in \dd{R} \}) = \{ \hbox{ Im } w = \hbox{ Im } F(z_0) \}.$$
}

\medskip

\eno{Lemma II.5.2.6}{ Let $(\gamma (t))_{t \in \dd{R}}$ be a curve such that $|\gamma(t)| \to \infty$ as $t \to
+\infty$. If for some $\epsilon > 0$ and some odd $j$, $1 \leq j \leq 2d$, we have $j\pi/d - \pi/2d + \epsilon
\leq$ arg $\gamma(t) \leq j\pi/d + \pi/2d - \epsilon$ for all $t$ sufficiently large, then $$ F(\gamma(t)) \to
w_p' = \int_0^{e^{ij\pi/d} \cdot \infty} e^{P(s)} ds \, \hbox{ as } t \to +\infty$$ where $p = (j+1)/2$. }

\medskip

\eno{Lemma II.5.2.7}{ Let $1 \leq j \leq 2d$. Then

(1) For $j$ odd, $k \geq k_0$, there are distinct integral curves $\gamma_{j,k}'$ of $X_P$ such that
$$
\partial D_{j,k}'(\alpha_j) = \gamma_{j,k}' \cup \gamma_{j,k+1}', \ \ \partial D_{j,k}'(\alpha_j) \cap
\partial D_{j,k+1}'(\alpha_j) = \gamma_{j,k+1}'.
$$

\medskip

(2) For $j$ even, $k \leq -k_0$, there are distinct integral curves $\gamma_{j,k}'$ of $X_P$ such that
$$
\partial D_{j,k}'(\alpha_j) = \gamma_{j,k}' \cup \gamma_{j,k-1}', \ \ \partial D_{j,k}'(\alpha_j) \cap
\partial D_{j,k-1}'(\alpha_j) = \gamma_{j,k-1}'.
$$
}

\medskip

\eno{Proposition II.5.2.8}{ Let $1 \leq j \leq 2d$ and let $j=2p-1$ be odd. Then

(1) $F$ maps $\gamma_{j,k}'$ to the half-line $]w_p',
w_p'+1\cdot\infty[$ for $k$ odd and to $]w_j' - 1\cdot\infty, w_j'[$ for $k$ even. $F$ maps
$D_{j,k}'(\alpha_j)$ univalently to the half plane $\{ $ Im $w > $ Im $w_p' \}$ for $k$ odd and
to $\{ $ Im $w < $ Im $w_p' \}$ for $k$ even.

(2) $F$ maps $\gamma_{j+1,k}'$ to the half-line $]w_j' - 1\cdot\infty, w_j'[$ for $k$ odd
and to $]w_p', w_p'+1\cdot\infty[$ for $k$ even. $F$ maps $D_{j+1,k}'(\alpha_{j+1})$ univalently to the half
plane $\{ $ Im $w < $ Im $w_p' \}$ for $k$ odd
and to $\{ $ Im $w > $ Im $w_p' \}$ for $k$ even.
}

\medskip

We define $2d$ families of domains $C_{j,l}'$ for $j = 1, \dots, 2d$ as follows:

\medskip

1. For $j$ odd: Fix $l_0$ such that $2l_0 - 1 \geq k_0$. We define for $l \geq l_0$,
$$
C_{j,l}' = D_{j,2l-1}'(\alpha_j) \cup \gamma_{j,2l}' \cup D_{j,2l}'(\alpha_j)
$$

2. For $j$ even: We define for $l \leq -l_0$,
$$
C_{j,l}' = D_{j,2l+1}'(\alpha_j) \cup \gamma_{j,2l}' \cup D_{j,2l}'(\alpha_j)
$$

\medskip

The domains $C_{j,l}$ are disjoint, and $F$ maps each univalently to a slit plane,
$$\displaylines{
F(C_{j,l}') = \dd C - [w_p', w_p'+1\cdot\infty[ , \ \hbox{ for } j = 2p-1, \hbox{ odd, and} \cr
F(C_{j,l}') = \dd C - ]w_p'-1\cdot\infty, w_p'] , \ \hbox{ for } j = 2p, \hbox{ even.} \cr
}$$

Thus as before we have $2d$ families of domains $(C_{1,l}')_{l \geq l_0}, \, (C_{2,l}')_{l \leq -l_0}, \dots,
\, (C_{2d-1,l}')_{l \geq l_0}, \, (C_{2d,l}')_{l \leq -l_0}$  which correspond under $F$ to families of
planes in $\cl S$, slit and pasted around the ramification points $w_1', \dots, w_d'$, with two families
$(C_{2p-1,l}), (C_{2p,l})$ for each ramification point $w_p'$.

\medskip

It remains to understand the mapping $F$ in the region outside the domains $C_{j,l}'$.

\medskip

Let
$$
D = \dd C - \overline{\cup_{j,l} C_{j,l}' }
$$

\medskip

\eno{Proposition II.5.2.9}{ There are only finitely many integral curves $\beta_1, \dots, \beta_n$ of $X_P$ within
$D$ which get mapped to either horizontal half-lines or line segments but not to full horizontal lines.}

\medskip

We need the following lemma, which is straightforward to prove by estimating the integral defining $F$ as in the
proof of Lemma II.5.2.6 above.

\medskip

\eno{Lemma II.5.2.10}{ Fix $\varepsilon > 0$. If $z \to \infty$ in $\dd{C}$ through the union of the $d$ sectors
given by $Q_{\varepsilon} := \{ \  |\arg z - \arg((-a_d)^{-1/d}) | > \pi/2d + \epsilon \ \}$ (where the inequality
holds for all the $d$-th roots), then $|w| = |F(z)| \to \infty$. }

\medskip

\stit{Proof of Proposition II.5.2.9} Take $R>0$ large enough so that the circle $\{ |z| = R \}$ meets each of the
domains $C_{j,l}$ for $j = 1, \dots, 2d$ and $l = +l_0$ (for $j$ odd), $l = -l_0$ (for $j$ even). Then removing
$B(0,R)$ from $D$ disconnects it into $2d$ components, ie $D - B(0,R)$ consists of $2d$ connected components $T_1,
\dots, T_{2d}$. In each $T_j$, we have $$ \arg z \to j\pi/d \hbox{ as } z \to \infty \hbox{ within } T_j. $$ and
hence $$\displaylines{ F(z) \to w_p' \hbox{ as } z \to \infty \hbox{ within } T_j \hbox{ for } j = 2p-1 \hbox{
odd}, \cr F(z) \to \infty \hbox{ as } z \to \infty \hbox{ within } T_j \hbox{ for } j \hbox{ even}. \cr }$$ Hence
an integral curve $\gamma$ of $X_P$ which is contained in $D$ and whose image under $F$ is not a full horizontal
line must escape to infinity when either $t \to +\infty$ or $t \to -\infty$ through one of the components $T_j$
for an odd $j$. It suffices to show that in each such component there can only be finitely many curves that escape
to infinity through it.

\medskip

Suppose a component $T_j$, for $j = 2p-1$ odd, contained infinitely many such curves. Then these curves must have
an accumulation point $z_0 \in B(0,R) \cap D$. Each curve is mapped to a segment contained in the same horizontal
line $\{ \Im w = \Im w_p' \}$ and ending at $w_p'$. On the other hand, $F'$ is nonzero, so $F$ is univalent in a
neighborhood of $z_0$, leading to a contradiction. $\diamond$

\medskip

$D - (\beta_1 \cup \dots \cup \beta_n)$ thus has only a finite number of connected components. The image of
each of them under $F$ is, being a union of full horizontal lines, simply connected, hence $F$ is univalent
on each connected component, mapping it to either a half-plane or a horizontal strip.

\medskip

The region $D$ thus contributes only finitely many sheets of $\cl S$, so we see that $\cl S$ has precisely
$d$ ramification points lying above the points $w_1', \dots, w_d'$. We note that the difference between the
analysis for a general $P$ here and that for $P(z) = z^d$ lies in the partitioning of the region $D$; the
domains of the partition correspond to sheets in $\cl S$, each containing one or more ramification points,
but unlike the previous case we cannot now identify a distinguished sheet in which all of them lie, indeed
there need not be such a sheet.

\medskip

\stit {II.5.3) Uniformization theorem.}

\eno {Theorem II.5.3.1}{ Let $\cl S$ be a simply connected log-Riemann surface of finite type of log-degree $d\geq
0$ without finite ramification points.

The Riemann surface $\cl S$ is bi-holomorphic to $\dd C$ and
the uniformization mapping
$$
F: \dd C \to \cl S
$$
is the primitive of a polynomial $P=P_{\cl S} \in \dd C[z]$ of degree $d\geq 0$,
$$
F(z)=\int_0^z \ e^{P(z)}\ dz \ .
$$
Conversely, for each polynomial $P \in \dd C[z]$ there exits
a a log-Riemann surface of finite
type of log-degree $d\geq 0$ without finite ramification points
for which the primitive $F$ of $\exp (P)$ realizes the uniformization.

The correspondence
$$
\cl S \mapsto P_{\cl S}
$$
is bijective.

}

We assume in the proof of the direct part the converse
result for the polynomial $P(z)=z^d$.

\eno {Lemma II.5.3.2}{ Let $\cl S_1$ and $\cl S_2$ be log-Riemann surfaces without finite ramification points and
with the same log-degree. Then there exists a quasi-conformal homeomorphism $$ \varphi : \cl S_1 \to \cl S_2 \ .
$$ We can add that $\pi_2\circ \varphi=\pi_1$ at infinite in the charts (i.e. $\varphi$ is the identity on
charts), that is out of $\pi_1^{-1} (K)$ where $K\subset \dd C$ is a large compact ball.

\medskip

Conversely, two such Riemann surfaces that are quasi-conformally
homeomorphic do have the same log-degree.
}

\eno {Corollary II.5.3.3}{ All Riemann surfaces of the previous lemma are bi-holomorphic to the complex plane.}

\stit {Proof of the Corollary.}

Let $\dd E_i$ be the disk $\dd D$ or the plane $\dd C$. Consider the uniformizations $$ \eqalign { &F_1 : \dd E_1
\to \cl S_1 \cr &F_2 : \dd E_2 \to \cl S_2 \cr } $$ and let $\varphi :\cl S_1 \to \cl S_2$. Then $$ \psi: \dd E_1
\to \dd E_2 $$ constructed as $\psi=F_1\circ \varphi \circ F_2^{-1}$ is a quasi-conformal homeomorphism. Thus both
$\dd E_1$ and $\dd E_2$ are the disk or the plane. But we know one example (the Riemann surface for the primitive
of $\exp (z^d)$) for which $\dd E$ is the complex plane $\dd C$. Thus they are all bi-holomorphic to $\dd C$.
$\diamond$

\stit {Proof of the Lemma.}

We proceed by induction on $d$. For $d=0$ the result is clear
(there is only one such Riemann surface).
Each such log Riemann surface $\cl S$
of log-degree $d\geq 1$ is a log Riemann surface $\tilde \cl S$
of log-degree $d-1$ with one infinite ramification point added
as we have observed before.
Thus by induction there is a quasi-conformal homeomorphism
$$
\tilde \varphi : \tilde S_1 \to \tilde S_2
$$
Let $z_1\in \tilde \cl S_1$ and $z_2\in \tilde \cl S_2$ be the
points where the infinite ramification points are added in
order to get $\cl S_1$ and $\cl S_2$.  Let $z_1'=\tilde \varphi^{-1}
(z_2)$. By smooth deformation of $z_1$ into $z'_1$, we can
construct a diffeomorphism $\psi : \tilde \cl S_1 \to \tilde \cl S_1$
such that $\psi(z_1)=z'_1$ and $\psi$ is the identity in the
$\pi_1$-pre-image of small neighborhoods of the infinite
ramification points in the charts, and in a neighborhood of
infinite on charts (one has to bend the cuts which gives an equivalent
log-Riemann surface structure).
Thus $\psi$ is quasi-conformal (by compactness and
continuity of the differential). Now $\tilde \varphi
\circ \psi : \tilde \cl S_1 \to \tilde \cl S_2$ defines a
map on the charts
that extends to a quasi-conformal homeomorphism $\varphi: \cl S_1
\to \cl S_2$ (we only need to extend it to the plane sheets
attached to $z_1$ and $z_2$ which is straightforward).

 For the converse, just observe that a q.c. homeomorphism does
preserve the infinite ramification points.
$\diamond$

\stit {Proof of the direct part of the theorem.}

We start considering two log-Riemann surfaces $\cl S_1$ and
$\cl S_2$, of finite type and the same log-degree.
Let $F_1 : \dd C \to \cl S_1$ and   $F_2 : \dd C \to \cl S_2$
be the uniformizations. Since $(\pi_i\circ F_i)' \not= 0$ we can
write
$$
\pi_i \circ F_i=\int e^{h_i}
$$
where $h_i$ is an entire function. Let $\varphi :\cl S_1 \to
\cl S_2$ be the quasi-conformal homeomorphism given by the lemma. Let
$\psi : \dd C \to \dd C$ be the quasi-conformal homeomorphism defined
by
$$
\psi=F_2^{-1} \circ \varphi \circ F_1\ .
$$

Note that any quasi-conformal homeomorphism $\psi: \dd C
\to \dd C$ extends to a quasi-conformal homeomorphism of the Riemann
sphere, thus it is H\"older at $\infty$ for the chordal metric
(any quasi-conformal homeomorphism is H\"older).

Now using
$$
F_2\circ \psi =\varphi \circ F_1
$$
and
$$
\pi_2 \circ F_2 \circ \psi =\pi_2 \circ \varphi \circ F_1
=\pi_1 \circ F_1
$$
(the last equality holds "near infinite"), we have that
the growth at infinite of $\pi_2\circ F_2$ is H\"older
equivalent to
the one of $\pi_1 \circ F_1$. Thus the order of  $\pi_1 \circ F_1$
is the same as the one of $\pi_2\circ F_2$ (note that the notion
of order is well defined for non-holomorphic functions).

\medskip

But now the order of an entire function $f$ is equal to the
order of its derivative $f'$ (this can be proved directly using the
mean value theorem, or the computation of the order using the
coefficients of a power series expansion).
Now we consider $\cl S_1$ to be the log-Riemann surface $\cl S$
of finite type and log-degree $d\geq $ of the theorem,
and $\cl S_2$ the log Riemann of the primitive of $\exp(z^d)$.
Thus the order of $\pi_2\circ F_2$ is $d$. So the order of
$\pi_1\circ F_1$ is also $d$ as well as the order of
$$
(\pi_1 \circ F_1)'=e^{h_1} \ .
$$
Thus $\Re h_1$ has a growth at infinite that is polynomial of degree $d$,
thus the same holds for $|h_1|$ and by Liouville theorem $h_1$
is a polynomial $P_{\cl S}$ of degree $d$.
$\diamond$

\medskip

\stit {II.5.4) Schwarz-Christoffel formula.}

\medskip

In this section we describe how the Schwarz-Christoffel formula for planar polygons generalizes to the case of
log-polygons. This will be used in the following section, where we sketch another approach to the uniformization
of log-Riemann surfaces of finite type, which is very close to the method originally used by Nevanlinna ([Ne1]).
One approximates a given surface $\cl S$ by log-polygons embedded in the surface, and obtains using Caratheodory
Kernel Convergence a uniformization for $\cl S$ as the limit of the uniformizations of the approximating
log-polygons.

\medskip

 The classical Schwarz-Christoffel formula gives a formula for
 the uniformization of a planar polygon. Indeed, there are two
 versions, the first for planar polygons with sides which are Euclidean line
 segments, which asserts that the {\it nonlinearity} $F''/F'$ of
uniformization $F$ is rational (see for eg.[Ne-Pa] p.330 or [Ah1] p.236 ), while the second, for planar polygons
with sides which are either Euclidean line segments or circular arcs, asserts that the {\it Schwarzian derivative}
$\{ F, z \}$ is rational (see for eg. [Hi] p.379 ) . The same assertions do in fact generalize to the case of
log-polygons made up of either only Euclidean segments or of both Euclidean segments and circular arcs; here
vertices are allowed to be at ramification points and hence angles greater than $2\pi$ are allowed as well. It
will be useful to also have the formula for a class of log-domains slightly more extended than the class of
log-polygons, for which we need the following definitions:

\medskip

\eno{Definition II.5.4.1}{ Let $D \subset \cl S$ be a log-domain in a log-Riemann surface $\cl S$ with projection
$\pi$. An end at infinity $e$ of $D$ is given by a family of nonempty sets $e = (U_R)_{R > 0}$ such that, for each
$R > 0$, $U_R$ is a connected component of $D - \pi^{-1}(\{ |w| \leq R \})$, and such that $U_{R_1} \subseteq
U_{R_2}$ whenever $R_1 \geq R_2$.

We say that ''$w \to \infty$ through $e$'' if for every $R > 0$ eventually $w$ lies in $U_R$.

}

\medskip

\eno{Definition II.5.4.2}{ A log-polygon with ends at infinity is a log-domain $P \subset \cl S$ in a log-Riemann
surface $\cl S$ with projection $\pi$ such that:

(1) $P$ is simply connected.

(2) The boundary $\partial P \subset \cl S^*$ of $P$ in the completion $\cl S^*$ is a union of finitely many
Euclidean segments, $\partial P = \gamma_1 \cup \dots \cup \gamma_n, \, n \geq 2$, where each $\gamma_k$ is either
a finite Euclidean segment or a Euclidean half-line. The $\gamma_k$'s are called the sides of $P$ and their
end-points the finite vertices of $P$.

(3) Each side which is a finite Euclidean segment intersects exactly two other sides at its two endpoints, while
each side which is a Euclidean half-line intersects exactly one other side at its one endpoint. There are thus two
sides meeting at each finite vertex $v$, and we define the interior angle at $v$ to be the angle $\theta$ in $P$
between these two sides, where $\theta \in (0, 2\pi)$ if $v \in \cl S$ is a regular point, and $\theta \in (0,
+\infty)$ if $v \in \cl S^* - \cl S$ is a ramification point.

\medskip

If all sides of $P$ are finite Euclidean segments then $P$ is a log-polygon as defined previously; if not, we
assume further that:

\medskip

(4) For each end at infinity $e = (U_R)_{R > 0}$ of $P$, there is an $R_0 > 0$ such that for $R \geq R_0$, $U_R$
is bounded by two sides of $P$ which are half-lines and an arc of $\pi^{-1}(\{ |w| = R \})$, and there is an
isometric embedding of $U_R$ into the surface of the logarithm $\cl S_{\log}$ which is the identity on charts.
This embedding followed by the automorphism $w \mapsto 1/w$ of $\cl S_{\log}$ maps $U_R$ to an angular sector
$V_R$ bounded by two curves (each of which is either a Euclidean segment or a circle arc) meeting at $0 \in \cl
S^*_{\log}$, and an arc of the form $\{ |w| = 1/R, a \leq \arg w \leq b \}$. We define the interior angle at the
end $e$ to be the angle $\theta \in [0, +\infty)$ in $V_R$ between the two boundary curves meeting at $0 \in \cl
S^*_{\log}$. }

\medskip

We note that such a log-polygon with ends at infinity has only finitely many ends at infinity.

\medskip

\eno{Proposition II.5.4.3}{ Let $P \subset \cl S$ be a log-polygon with ends at infinity, and $K : P \to \dd D$ a
conformal representation of $P$. Then the ends at infinity of $P$ correspond to points on the boundary of the unit
disk under $K$, more precisely for any end $e = (U_R)_{R > 0}$ there is a unique point $z_e \in \partial \dd D$
such that when $w \to \infty$ through $e$, then $K(w) = z \to z_e, \, z \in \dd D$.}

\medskip

\stit{Proof:} Let $e = (U_R)_{R > 0}$ be an end at infinity of $P$, and $\theta$ the angle at $e$. Let $R_0 > 0$
be as given by condition $(4)$ above, and $h : U_{R_0} \to V_{R_0} \subset \cl S_{\log}$ be the isometric
embedding into $\cl S_{\log}$ followed by inversion, of $U_{R_0}$. The domain $V_{R_0}$ can be mapped conformally
to a Jordan domain $V \subset \dd C$ by a map $g : V_{R_0} \to V$ with expression in log-coordinates of the form
$g(w) = w^{\pi / \theta}$ (or $g(w) = e^{c / w}$ if $\theta = 0$, for some constant $c$), and $g$ extends to a
homeomorphism of the closed domains $g : \overline{V_{R_0}} \subset \cl S^*_{\log} \to \overline{V}$, mapping $0
\in \cl S^*_{\log}$ to $0 \in \dd C$.

\medskip

Let $w_0, w_1 \in \partial P$ be the points on the boundary of $U_{R_0}$ where the two half-lines bounding
$U_{R_0}$ meet the arc of $\pi^{-1}(\{ |w| = R \})$ bounding $U_{R_0}$. Since the boundary of $P$ is locally
connected, $K$ extends continuously to the points $w_0, w_1$ and one sees that $K$ maps $U_{R_0}$ to a Jordan
domain $W \subset \dd D$, bounded by one of the arcs of $\partial \dd D$ joining $K(w_0)$ to $K(w_1)$, and a curve
in $\dd D$ joining $K(w_0)$ to $K(w_1)$.

\medskip

The conformal map $\phi = K \circ h^{-1} \circ g^{-1} : V \to W$ between Jordan domains extends to a homeomorphism
of the closed domains $\phi : \overline{V} \to \overline{W}$. Let $z_e \in \partial W$ correspond under $\phi$ to
$0 \in \partial V$.

\medskip

Now as $w \to \infty$ through $e$, it is clear that $h(w) \to 0 \in \cl S^*_{\log}$, so $g(h(w)) \to 0 \in
\partial V$, hence $z = K(w) = \phi(g(h(w))) \to z_e$ as required. It is not hard to see as well that $z_e \in \partial \dd
D$. $\diamondsuit$

\medskip

\eno{Theorem II.5.4.4 (Generalized Schwarz-Christoffel formula 1).}{ Let $P \subset \cl S$ be a log-polygon with
ends at infinity, embedded in a log-Riemann surface $\cl S$ with projection mapping $\pi$. Suppose $P$ has $n$
finite vertices $w_1, \dots, w_n$ with interior angles $\pi\alpha_1, \dots, \pi\alpha_n$, and $m$ ends at infinity
$e_1, \dots, e_m$ with interior angles $\pi\beta_1, \dots, \pi\beta_m$, where $\alpha_1, \dots, \alpha_n > 0$,
$\beta_1, \dots, \beta_m \geq 0$. Then for any uniformization $\tilde{F} : \dd D \to P$ that maps the unit disk
$\dd D$ conformally onto $P$, with expression in log-coordinates $ F(z) := \pi \circ \tilde{F}(z) $, its
nonlinearity $F''/F'$ is a rational function $$ {F'' \over F'} = \sum_{k=1}^n {\alpha_k - 1 \over z - z_k} +
\sum_{j=1}^m {(-\beta_j - 1) \over z - z'_j} + C $$ where $z_1, \dots, z_n \in
\partial \dd D$ and $z'_1, \dots, z'_m \in \partial \dd D$ are the points on the boundary of the unit disk
that correspond to the finite vertices and ends at infinity respectively of $P$ and $C$ is a constant depending on
$\tilde{F}$. Since $F''/F' = {d \over dz} \log F'$, one can also solve for $F$ from the above formula, to write
$F$ in integral form as $$F(z) = A \int_{0}^{z} (t - z_1)^{\alpha_1 - 1} \dots (t - z_n)^{\alpha_n - 1} (t -
z'_1)^{-\beta_1 - 1} \dots (t - z'_m)^{-\beta_m - 1} \ dt + B \ , \ z \in \dd D $$ where $A, B$ are constants
depending on $\tilde{F}$.}

\medskip

\stit{Proof :} The proof follows the same lines as the classical case. The uniformization $\tilde{F}$ extends
continuously to $\partial \dd D - \{ z_1, \dots, z_n, z'_1, \dots, z'_m \}$, which is a disjoint union of $(n +
m)$ circular arcs, each of which is mapped one-to-one onto the corresponding side of $P$. By the Schwarz
reflection principle, the function $F = \pi \circ \tilde{F}$ can be analytically continued to any point
 $z \in \overline{\dd C} - \overline{\dd D}$ along any curve $\gamma$ that starts
 from $0 \in \dd D$ and passes through exactly one of these arcs, via
 the equation
 $$
 F(z) = S (F(1/\overline{z}))
 $$
 where $S$ denotes the reflection through the straight line in $\dd C$ containing
 the $\pi$-projection of the corresponding side of $P$ .

\medskip

The key observation here is that while the branch of $F$ obtained depends on the path $\gamma$, any two branches
$F_1$ and $F_2$ are related by a product of two reflections through straight lines, and hence by an affine linear
transformation, $F_1 = aF_2 + b$. Since the nonlinearity is invariant
 under affine linear transformations of the
dependent variable, we have $F''_1/F'_1 = F''_2/F'_2$, and it follows
 that the nonlinearity $F''/F'$ extends to a
single-valued function on $\overline{\dd C} - \{z_1, \dots, z_n, z'_1, \dots, z'_m \}$.

\medskip

A local analysis near the points $z_k, w_k$ shows that near each $z_k$, the function $F$ can be written in the
form $$ F(z) = \pi(w_k) + H_k(z) (z - z_k)^{\alpha_k} \ , \ |z - z_k| < \epsilon, z \in \dd D $$ where $H_k$ is a
function regular in a full neighbourhood $\{ |z - z_k| < \epsilon \}$ of $z_k$, and $H_k(z_k) \neq 0$. It follows
that $$ {F'' \over F'}(z) = {\alpha_k - 1 \over z - z_k} + G_k(z) \ , \ |z - z_k| < \epsilon, z \in \dd D $$ for a
function $G_k$ regular in $\{ |z - z_k| < \epsilon \}$. Since both sides of the above equation are defined on a
punctured neighbourhood $\{ 0 < |z - z_k| < \epsilon \}$ of $z_k$, they agree there as well; it follows that
$F''/F'$ has a simple pole with residue $\alpha_k - 1$ at $z_k$.

\medskip

Similarly, near each point $z'_j$ corresponding to an end $e_j$, $F$ can be written in the form $$ F(z) = (z -
z'_j)^{-\beta_j} P_j(z) \ , \ |z - z'_j| < \epsilon, z \in \dd D $$ where $P_j$ is a function regular in a full
neighbourhood $\{ |z - z'_j| < \epsilon \}$ of $z_j$, and $P_j(z'_j) \neq 0$. As above it follows that $F''/F'$
has a simple pole with residue $-\beta_j - 1$ at $z'_j$.

\medskip

Thus $F''/F'$ is regular everywhere in the extended plane $\overline{\dd C}$ except for simple poles at the points
$z_1, \dots, z_n, z'_1, \dots, z'_m$, hence is a rational function and can be written in the form given in the
theorem. $\diamondsuit$

\medskip

We also have a version of the formula for log-polygons with sides which are either finite Euclidean segments or
circular arcs.

\medskip

\eno{Theorem II.5.4.5 (Generalized Schwarz-Christoffel formula 2).}{ Let $P \subset \cl S$ be a log-polygon with
sides that are either finite Euclidean segments or circular arcs, embedded in a log-Riemann surface $\cl S$ with
projection mapping $\pi$, and suppose $P$ has $n$ vertices $w_1, \dots, w_n$ with interior angles $2\pi\alpha_1,
\dots, 2\pi\alpha_n$, where $\alpha_1, \dots, \alpha_n > 0$. Then for any uniformization $\tilde{F} : \dd D \to P$
that maps the unit disk $\dd D$ conformally onto $P$, with expression in log-coordinates $ F(z) := \pi \circ
\tilde{F}(z) $, its Schwarzian derivative $\{ F, z \}$ is a rational function $$ \{ F, z \} = \left( {w'' \over
w'} \right)' - {1 \over 2} \left( {w'' \over w'} \right)^2 =  {1 \over 2} \sum_{k=1}^n \left [{1 - {\alpha_k}^2
\over (z - z_k)^2 } + {\beta_k \over z - z_k} \right] $$ where $z_1, \dots, z_n \in
\partial \dd D$ are the $n$ points on the boundary of the unit disk that correspond to the vertices $w_1, \dots,
w_n$ respectively of $P$ and $\beta_1, \dots, \beta_n$ are constants depending on $\tilde{F}$. These constants
satisfy the relations $$\displaylines{ \sum_{k=1}^n \beta_k = 0 \ , \ \sum_{k=1}^n (2\beta_k z_k + 1 -
{\alpha_k}^2 ) = 0 \cr \sum_{k=1}^n [ \beta_k {z_k}^2 + (1 - {\alpha_k}^2)z_k ] \cr }$$ }

\medskip

\stit{Proof : } We give here only a sketch of the proof, which follows the same lines as that of the preceding
theorem. As above, by Schwarz reflection principle $F = \pi \circ \tilde{F}$ can be continued analytically along
curves which start in $\dd D$ and end in $\overline{\dd C} - \overline{\dd D}$ passing through an arc of $\partial
\dd D - \{ z_1, \dots, z_n \}$, by the formula $$ F(z) = S (F(1/\overline{z})) $$
 where now $S$ denotes either a reflection through a straight line or through
 a circle, depending on whether the corresponding side of $P$
 is a Euclidean line segment or a circular arc.

\medskip

It follows that any two branches $F_1$ and $F_2$ of $F$ are related to one another by a fractional linear
transformation, $F_1 = (aF_2 + b)/(cF_2 + d)$, and hence, since the Schwarzian derivative is invariant under
fractional linear transformations of the dependent variable, that $\{ F , z \}$ can be extended to a single-valued
function regular on all of $\overline{\dd C} - \{ z_1, \dots, z_n \}$. Local analysis near the points $z_k$ shows
that in fact $\{ F , z \}$ has double poles at these points, with principal parts of the form $${1 - {\alpha_k}^2
\over (z - z_k)^2 } + {\beta_k \over z - z_k} \ , \ k=1,\dots,n $$ for some constants $\beta_1, \dots, \beta_n$.
Moreover, any branch of $F$ is regular at infinity, from which one can show that $\{ F , z \}$ must vanish to the
fourth order at infinity (ie $z^4 \{ F , z \}$ is holomorphic at infinity), so $\{ F , z \}$ is indeed equal to
the sum of its principal parts. The conditions given on the $\beta_k$'s express the fact that when $\{ F , z \}$
is expanded in powers of $1/z$ near $z = \infty$, the terms in $1/z^m$ are missing for $m = 1,2,3$. $\diamondsuit$

\medskip

\stit {II.5.5) Uniformization via Schwarz-Christoffel formula.}

\medskip

Let $\cl S$ be a simply connected log-Riemann surface of finite log-degree $d$, and let $w^*_1, \dots, w^*_d \in
\cl S^* - \cl S$ be the $d$ infinite order ramification points. With the theorems of the previous section in hand
we may now attempt to obtain a uniformization of $\cl S$ as the limit of uniformizations of approximating
log-domains, either log-polygons with ends at infinity, or log-polygons with circular arcs, that converge to $\cl
S$ in the sense of Caratheodory.

\medskip

If one takes log-polygons with all sides finite Euclidean segments, then the number of vertices must necessarily
increase without bound; if one allows log-polygons with ends at infinity however, it is then possible, as we will
see below, to construct an approximating sequence with a uniformly bounded number of vertices plus ends at
infinity. This has the advantage that the uniformizations of these log-polygons with ends at infinity have
rational nonlinearities of bounded degree, and hence any limit of their uniformizations must have rational
nonlinearity. If one uses log-polygons with circular arcs, then it is also possible to bound uniformly the number
of vertices needed (one needs to take circular arcs which spiral around many sheets), but in this case one obtains
only that the Schwarzian of the limit uniformization is rational, and one cannot directly integrate as in the case
of rational nonlinearity to obtain a formula for the limit uniformization. The approximating sequence of
log-polygons with ends at infinity is constructed as follows:

\medskip

Consider a minimal atlas for $\cl S$, given as in section I.3.1 by taking the cells $(U(w_i))$ of the fiber $(w_i)
= \pi^{-1}(z_0)$ of a generic point $z_0 \in \dd C$. For $j=1, \dots, d$, in a neighborhood of $w^*_j$ one can
define a well-defined argument function $\arg(w - w^*_j)$; there is an angle $\theta_j$ such that in each sheet of
the minimal atlas containing $w^*_j$, we have $\theta_j + 2\pi N < \arg(w - w^*_j) < \theta_j + 2\pi(N+1)$, with
$N \in \dd Z$ being an integer depending on the sheet. Since $\cl S$ has finitely many ramification points, there
is an integer $N_0 \geq 1$ such that for $|N| \geq N_0, \ j=1,\dots, d$, any sheet containing $w^*_j$ with
$\theta_j + 2\pi N < \arg(w - w^*_j) < \theta_j + 2\pi(N+1)$ is a clean sheet (see definition I.3.1.3), containing
only $w^*_j$. We define the following sequence of log-domains $(D_N)_{N \geq N_0}$:

\medskip

For $N \geq N_0$, let $D_N$ be the interior of the closed log-domain given by the closure of the union of all
sheets $U(w_i)$ such that for all $j=1,\dots,d$, $\theta_j - 2\pi N < \arg(w - w^*j) < \theta_j + 2\pi(N+1)$ in
$U(w_i)$. The boundary $\partial D_N \subset \cl S$ of $D_N$ in $\cl S$ consists of the $2d$ Euclidean half-lines
given by $\{ \ \arg(w - w^*_j) = \theta_j + 2\pi(N+1) \ \}, \ \{ \ \arg(w - w^*_j) = \theta_j - 2\pi N \ \}, \
j=1, \dots, d$. It is straightforward to check that the $D_N$'s are log-polygons with ends at infinity. Each $D_N$
has $d$ finite vertices, namely the $d$ ramification points of $\cl S$, and has $d$ ends at infinity; the angles
at the finite vertices and at the ends at infinity are all equal to $2\pi(2N+1)$.  The uniformizations of the
$D_N$'s all have rational nonlinearities of degree $2d$. We observe that the log-domains $D_N$ converge in the
sense of Caratheodory to $\cl S$.

\medskip

Let $R(N)$ be the conformal radius of $D_N$, and $\tilde{F_N} : \dd D_{R(N)} \to D_N$ the uniformization of $D_N$
normalized so that $F_N(0) = z_0, F'_N(0) = 1$, where $F_N = \tilde{F_N}$. Let $\tilde{F} : \dd C \to \cl S$ be
the uniformization of $\cl S$ normalized so that $F(0) = z_0, F'(0) = 1$, where $F = \tilde{F}$. Then we have:

\medskip

 \eno{Theorem II.5.5.1}{

(1) $R(N) \to +\infty$ as $N \to +\infty$.

(2) $\tilde{F_N} \to \tilde{F}$ uniformly on compacts, in the sense that $d(\tilde{F_N}, \tilde{F}) \to 0$
uniformly on compacts of $\dd C$, where $d(.,.)$ is the log-euclidean metric on $\cl S$.

(3) The nonlinearity $F''/F'$ of $F$ is a polynomial $P$ of degree at most $2d$. Hence letting $Q$ be a primitive
of $P$, for some constant $A$ we have $$\eqalign{ F(z) & = A \int_{0}^{z} e^{Q(t)} \ dt + z_0 \cr & = \int_{0}^{z}
e^{P_0(t)} \ dt + z_0 \cr}$$ where $P_0 = Q + \log A$.

}

\medskip

\stit{Proof :} (1) follows from the fact that $D_N \to \cl S$ in the sense of Caratheodory Kernel Convergence and
the continuity of the conformal radius, Theorem *.*.*. For (2), consider the functions $G_N = \tilde{F}^{-1} \circ
\tilde{F_N} : \dd D_{R(N)} \to \dd C$. Since $G_N(0) = 0, G'_N(0) = 1$ and $G_N$ is univalent, the $G_N$'s form a
normal family on any disk of fixed radius $R$; any limit of this sequence must be univalent on $\dd C$, hence
affine linear, and hence by virtue of the normalizations must be the identity. It follows that $\tilde{F}^{-1}
\circ \tilde{F_N} \to id$ uniformly on compacts of $\dd C$, from which (2) follows easily.

\medskip

It follows that the functions $F''_N / F'_N$ converge normally to $F''/F'$. Since these are rational functions of
bounded degree $2d$, $F''/F'$ must be a rational function of degree at most $2d$. Each $F''_N / F'_N$ has $2d$
simple poles on the boundary of the disk $\dd D_{R(N)}$ of radius $R(N)$ and no other poles; from (1) it follows
that for any fixed compact $K \subset \dd C$, eventually none of the functions $F''_N / F'_N$ have poles on $K$.
It follows that $F'' / F'$ has no poles in the finite plane, and is hence a polynomial as stated in (3).
$\diamondsuit$

\medskip

We note that we obtain here a polynomial $P_0$ of degree at most $2d+1$; a more detailed analysis, which we forego
here, can show that in fact $P_0$ must have degree exactly $d$, a result which was already known from section
II.5.3.

\medskip

\stit {II.5.6) General uniformization theorem.}

Let $P(z) = a_d z^d + \dots + a_0$ and $Q(z) = b_m z^m + \dots + b_0$
be two polynomials of degrees $d$ and $m$ respectively. Let $F$ be the
entire function
$$
F(z) = \int_{0}^{z} Q(t) e^{P(t)} \ dt
$$
Generalizing the results of the previous sections, we have:

\medskip

\eno{Theorem II.5.6.1}{ Let $A = Q^{-1}(0)$ be the zeroes of $Q$. There exists a log-Riemann surface $\cl S$ such
that the map $F : \dd C - A \to \dd C$ lifts to a biholomorphism $\tilde{F} : \dd C - A \to \cl S$ such that $\pi
\circ \tilde{F} = F$. The surface $\cl S$ contains exactly $d$ ramification points of infinite order, and $m$
ramification points of finite order (counting multiplicities). The finitely completed Riemann surface $\cl
S^{\times}$ is simply connected, and the map $\tilde{F}$ extends to a biholomorphism of Riemann surfaces
$\tilde{F} : \dd C \to \cl S^{\times}$. The infinite ramification points $w_1, \dots, w_d$ project onto the points
$$ w_j' = \pi(w_j) = \int_{0}^{\rho_j \cdot \infty} Q(z) e^{P(z)} \ dz \ , \ j = 1, \dots, d $$ where $\rho_1,
\dots, \rho_d$ are the $d$ values of $(-a_d)^{-1/d}$. }

\medskip

We only give a sketch of the proof, which follows the same lines as
in the previous sections. We assume for convenience again that $a_d = 1$.

\medskip

We consider the vector field
$$
X_{P,Q}(z) = e^{-i \, \hbox{\sevenrm Im} \, ( P(z) + \log Q(z)) } \, \, \, ,z \in \dd{C} - A
$$
whose integral curves get mapped to horizontals by $F$ (note that $X_{P,Q}$ is well-defined
independently of the choice of $\log Q$). For large $z$ we have
$P(z) + \log Q(z) = z^d (1 + O((\log z)/z))$, and analysis of $X_{P,Q}$ can be
carried out similarly as for $X_P$. The function
$$
\xi = (P(z) + \log Q(z))^{1/d} = z (1 + O((\log z)/z))^{1/d}
$$
is well-defined and univalent in any slit domain $\{ |z| > R, \ z \notin [R, +1\cdot\infty] \}$
for $R > 0$ large enough. Hence the inverse function $z = h(\xi)$ is
a change of variables such that
$$
P(z) + \log Q(z) = \xi^d
$$
As before, we construct families of transversals
$\Gamma_{j,k}'(\alpha_j), j=1,\dots,2d, |k| > k_0$
to $X_{P,Q}$. Each $\Gamma_{j,k}'(\alpha_j)$ is a connected component
of $\{ \hbox{ Im } ( P(z) + \log Q(z) ) = k\pi - \alpha_j \}$,
and $X_{P,Q} = \pm e^{i\alpha_j}$ on $\Gamma_{j,k}'(\alpha_j)$.
Using the fact (which is easily checked) that
$$
h'(\xi) = 1 + o(1)
$$
it is possible, taking $k_0$ large enough and choosing the $\alpha_j$'s
appropriately, to ensure that the curves $\Gamma_{j,k}'(\alpha_j)$ are
transversal to $X_{P,Q}$.

\medskip

The families of transversals can then be used to construct $2d$ families
of disjoint domains $C_{j,l}', \, j = 1,\dots,2d$, which correspond under
$F$ to families of planes in $\cl S$, slit and pasted around the
ramification points $w_1', \dots, w_d'$, with two families
$(C_{2p-1,l}), (C_{2p,l})$ for each ramification point $w_p'$.

\medskip

The region complementary to these domains
$$
D = \dd C - \overline{\cup_{j,l} C_{j,l}' }
$$
is simply connected, and we have as before

\medskip

\eno{Proposition II.5.6.2}{ There are only finitely many integral curves $\beta_1, \dots, \beta_n$ of $X_{P,Q}$
within $D$ which get mapped to either horizontal half-lines or line segments but not to full horizontal lines.}

\medskip

We need the following Lemma

\eno{Lemma II.5.6.3}{ Let $z_0 \in A$ be a finite singularity of $X_{P,Q}$, ie a zero of $Q$. If the order of the
zero is $r$ then there are exactly $2(r+1)$ integral curves of $X_{P,Q}$ which accumulate at $z_0$.}

\medskip

\stit{Proof.} $z_0$ is a zero of order $r$ for $F'(z) = Q(z)e^{P(z)}$, so there
exists a local change of variables $\zeta(z) = \lambda(z-z_0)+O((z-z_0)^2)$
near $z_0$ such that
$$
F(z) = F(z_0) + \zeta^{r+1}
$$
Thus near $z_0$ there are exactly $2(r+1)$ curves terminating at $z_0$ which
get mapped by $F$ to horizontal segments. $\diamond$

\medskip

\stit{Proof of Proposition.} Any integral curve of $X_{P,Q}$
must either escape to infinity when
$|t| \to +\infty$, or otherwise accumulate one of the finite
singularities of $X_{P,Q}$. For the curves which escape to infinity, the same
compactness argument as before shows that there can only be finitely many such curves
within $D$ whose images are not full horizontal lines (note that the integral defining
$F$ converges to finite values $w_1',\dots,w_d'$ and diverges to $\infty$ in the
same angular sectors as before).

\medskip

For the other curves, which accumulate at the finite singularities,
the above Lemma shows that there can only be finitely many such curves
at each zero of $Q$, and since $Q$ has finitely many zeroes,
the result follows. $\diamond$

\medskip

\stit{Proof of Theorem.} Considering the connected components of
$D - (\beta_1 \cup \dots \cup \beta_n)$,
on each of which $F$ is univalent, we can partition the region $D$ into domains
corresponding under $F$ to sheets of $\cl S$, each one being a plane minus a
finite number of horizontal slits ending either at a point $w_j'$ or at a
critical value $F(z_0) \in F(A)$ of $F$.

\medskip

These finitely many sheets, along with those corresponding via $F$
to the domains $C_{j,l}'$, allow us to build simultaneously the
log-Riemann surface $\cl S$ of the Theorem as well as the lift
$\tilde{F} : \dd C - A \to \cl S$.

\medskip

It is straightforward to see that the surface $\cl S$ contains exactly
$d$ ramification points $w_1,\dots,w_d$ of infinite order, and finitely
many ramification points of finite orders adding up to the degree $m$ of
$Q$ $\ \diamondsuit$.

\medskip

Thus given a primitive of the form $\int Qe^{P}$, where $Q$ and $P$ are
polynomials, we can associate to it a log-Riemann surface $\cl S$ such
that the uniformization of $\cl S$ is given by this primitive, and such
that the numbers of finite and infinite ramification points correspond
exactly to the degrees of $Q$ and $P$ respectively.

\medskip

Conversely, we have the following Theorem:

\medskip

\eno{Theorem II.5.6.4}{ Let $\cl S$ be a log-Riemann surface of finite type of log-degree $d \geq 0$ and $m \geq
0$ finite ramification points (counting multiplicities), such that the finite completion $\cl S^{\times}$ is
simply connected.

Then the surface $\cl S$ is biholomorphic to $\dd C$ and the
uniformization mapping $F : \dd C \to \cl S^{\times}$ is given by a
primitive of the form
$$
F(z) = \int_{0}^{z} Q(z) e^{P(z)} \ dz,
$$
where $P, Q \in \dd C[z]$ are polynomials of degrees $d$ and $m$
respectively.
}

\medskip

The proof proceeds along lines similar to the proof of Theorem II.5.3.1 in section II.5.3.

\medskip

\eno{Lemma II.5.6.5}{ Let $\cl S$ be a log-Riemann surface of finite type of log-degree $d \geq 0$ and $m \geq 0$
finite ramification points (counting multiplicities), such that the finite completion $\cl S^{\times}$ is simply
connected. Then there exists a simply connected log-Riemann surface $\cl S_1$ of log-degree $d$ and with no finite
ramification points, such that there is a quasi-conformal homeomorphism $$ \phi : \cl S^{\times} \to \cl S_1. $$
Moreover, we can add that $\phi$ satisfies, for a constant $C$, the inequality $$ |\pi(\phi^{-1}(w))| \leq C
|\pi_1(w)| \ , \ w \in \cl S_1 $$ where $\pi, \pi_1$ denote the projection mappings of $\cl S$ and $\cl S_1$
respectively. }

\medskip

\stit{Proof of Lemma:} The proof is by induction on $m \geq 0$. For
$m = 0$ there is nothing to prove. For $m \geq 1$, by induction its
enough to construct a quasi-conformal homeomorphism
$\phi : \cl S^{\times} \to \cl S_1$ to a log-Riemann surface $\cl S_1$
of log-degree $d$ but with a strictly smaller number $m_1 < m$
of finite ramification points (counting multiplicities).

\medskip

Consider a finite ramification point $w_0 \in \cl S^{\times}$ of order
$n \geq 2$ say. The point $w_0$ appears in exactly $n$ sheets of any
minimal atlas of $\cl S$. Some of these may be 'clean' sheets, containing
no ramification points other than $w_0$, while others may contain
other ramification points as well. However, by
quasi-conformally deforming the surface $\cl S$, rotating around $w_0$
all the ramification points other than $w_0$ (along with the planes
attached to them), we may assume that $(n-1)$ of these sheets are 'clean',
and all the other ramification points are connected to $w_0$
through a single sheet. Assume that $\pi(w_0) = 0$. Let the regions
$A, B$ and $C$ be defined as follows (see the figure below):

\smallskip

$A = $ Union of the  $(n-1)$ clean sheets containing $w_0$ and the half-plane in the $n$th sheet $\{ \hbox{
Re } w < 0 \ \}$

\smallskip

$B = $ Region in the $n$th sheet bounded by the lines $\{ \hbox{ Re } w = 0 \ \}, \ \{ \ \arg(w - 1) = {\pi
\over 2n} \ \} \hbox{ and } \{ \ \arg(w - 1) = -{\pi \over 2n} \ \}$

\smallskip

$
C = \cl S^{\times} - (A \cup B)
$

\smallskip

By further deformation we may assume that
all ramification points other than $w_0$ lie in the region
$C$ as shown in the figure:

\medskip

{\hfill {\centerline {\psfig {figure=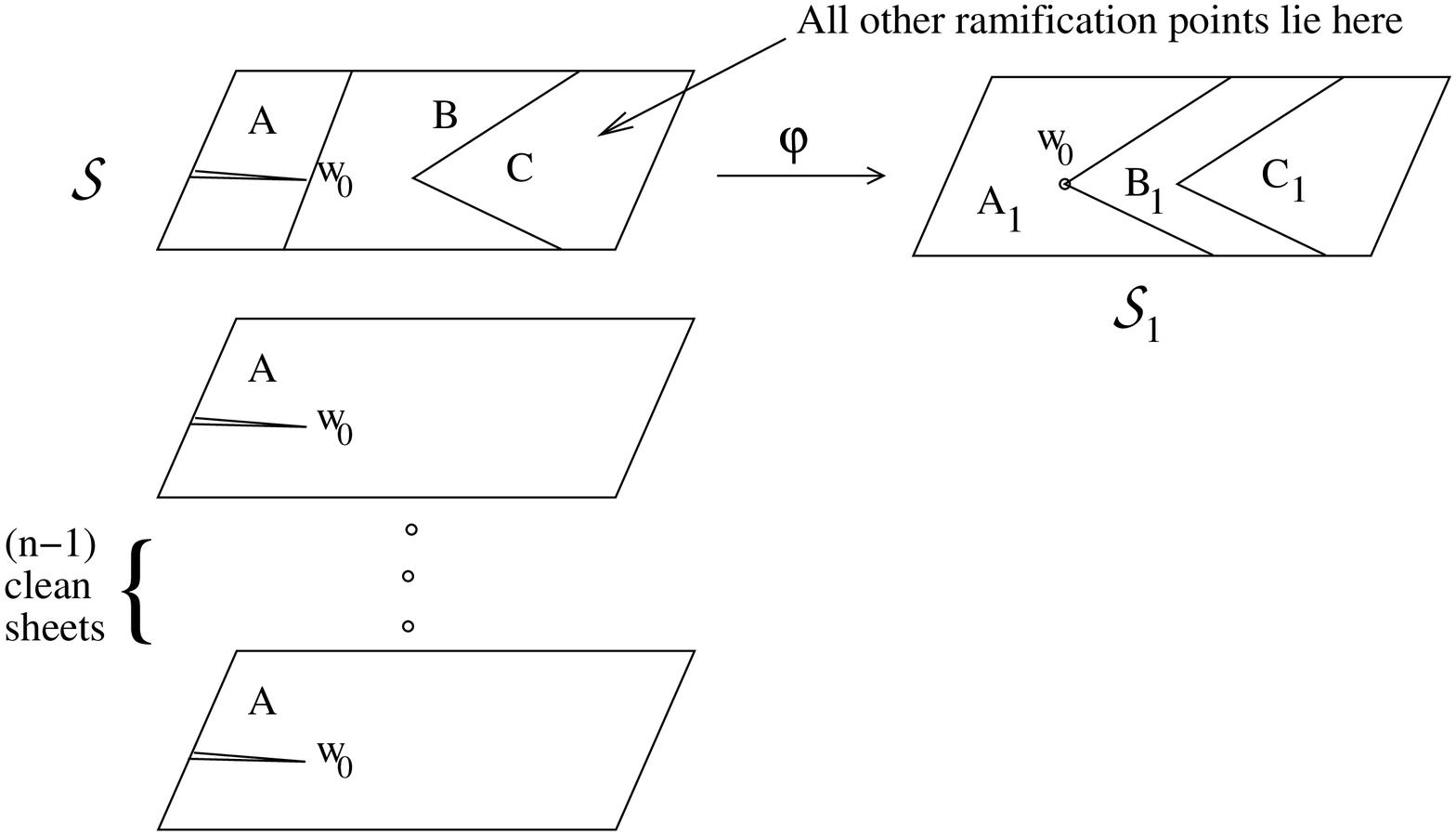,height=6cm}}}}

\medskip

Let $\cl S_1$ be the log-Riemann surface shown in the figure, given by pruning
the ramification point $w_0$ from the surface $\cl S$ (see section I.5.2). We
can define a quasi-conformal homeomorphism $\phi :\cl S^{\times} \to \cl S_1^{\times}$
as follows:

\medskip

1. Let $\theta = \arg w$ be an argument function defined in $A$ taking
values in the intervals $(\pi/2, n\pi)$ and $(-n\pi, -\pi/2)$. Define
$$
\phi(re^{i\theta}) := re^{i\theta/n} \ , \ w = re^{i\theta} \in A
$$
This maps $A \subset \cl S$ quasi-conformally onto the region
$A_1 \subset \cl S_1$.

\medskip

2. In the region $C$ define $\phi$ to be the identity in log-charts,
$$
\phi(w) := w \ , \ w \in C
$$
The region $C \subset \cl S$ corresponds isometrically to
the region $C_1 \subset \cl S_1$.

\medskip

3. Extend $\phi$ continuously to $B$ so
that it agrees on the two boundary components of $B$ with
the maps defined above, and so that $B \subset \cl S$ is
mapped quasi-conformally to the region $B_1 \subset \cl S_1$,
which is bounded by the lines $\{ \ \arg w = {\pi \over 2n} \ \},
\{ \ \arg w = -{\pi \over 2n} \ \}, \{ \ \arg (w-1) = {\pi \over 2n} \ \}$
 and $\{ \ \arg (w-1) = -{\pi \over 2n} \ \}$.

\medskip

Since $\cl S_1$ has a strictly smaller number of ramification points than $\cl S$,
the result follows by induction. We observe that the estimate in the statement
of the Lemma follows from the above construction. $\diamondsuit$

\medskip

\stit{Proof of Theorem.} It follows from the above Lemma and the main Theorem
of section II.5.3 that $\cl S^{\times}$ is parabolic.
Let $F:\dd C \to \cl S^{\times}$ be the uniformization. Since the projection
$\pi : \cl S^{\times} \to \dd C$ has critical points precisely at the finite
ramification points (and of the same orders), the entire function
$\pi \circ F : \dd C \to \dd C$ has precisely $m$ critical points (counting
multiplicities). Hence we can factor its derivative as
$$
(\pi \circ F)'(z) = Q(z) e^{h(z)}
$$
where $Q \in \dd C[z]$ is a polynomial of degree $m$ with zeroes at these $m$
critical points, and $h$ is an entire function.

\medskip

Now let $\phi : \cl S^{\times} \to \cl S_1$ be a quasi-conformal homeomorphism
as given by the Lemma to a log-Riemann surface $\cl S_1$ of log-degree $d$
and without finite ramification points. We know from section II.5.3 that
$\cl S_1$ has a uniformization $F_1 : \dd C \to \cl S_1$ given by a primitive
$\int e^{P_1}$, for some polynomial $P$ of degree $d$.

\medskip

Let
$\psi : \dd C \to \dd C$ be the quasi-conformal homeomorphism defined
by
$$
\psi = F_1^{-1} \circ \phi \circ F
$$
We can then write $\pi \circ F$ in the form
$$
\pi \circ F = \pi \circ \phi^{-1} \circ (F_1 \circ \psi)
$$
It follows from the estimate on $\phi$ given by the Lemma that $\pi \circ F$
has the same order as $\pi_1 \circ (F_1 \circ \psi)$ (note that the notion
of order is well-defined for non-holomorphic functions). Since $\psi$ is
is H\"older at $\infty \in \overline{\dd C}$ for the chordal metric, and
$\pi_1 \circ F_1$ is of finite order, it follows that $\pi_1 \circ (F_1 \circ \psi)$
and hence $\pi \circ F$ is of finite order.

\medskip

Thus $(\pi \circ F)'$ is of finite order as well, which implies that $h$ is equal to a polynomial $P \in \dd
C[z]$. Since the surface $\cl S$ has $d$ infinite ramification points it follows from Theorem II.5.6.1 at the
beginning of this section that $P$ is of degree $d$. $\diamondsuit$

%% file: ii6.tex
\input defv7.tex
\input psfig.sty

\stit {II.6) Cyclotomic log-Riemann surfaces.}

\stit {II.6.1) Definition.}

\eno {Definition II.6.1.1}{ Let $d\geq 0$ and $n\geq 1$ be integers. The cyclotomic log-Riemann surface $\cl
S_{n,d}$ of log-degree $d$ and pol-degree $n$ is the unique log-Riemann surface with uniformization given by $$
F_{n,d} : \dd C \to \cl S_{n,d} $$ with $$ \pi\circ F_{n,d} (z)=\int_0^z t^n e^{-t^d / d} \ dt \ , \ d \geq 1, $$
and $$ \pi\circ F_{n,0} (z)=\int_0^z t^n  \ dt \ . $$

}

\stit {Examples.}

\medskip

Recall the notation $a\cl S+b$ introduced in section I.1.

\medskip

{\bf 1.} For $d=0$ the log-Riemann surface $(n+1) \cl S_{n,0}$ is
$\cl S_n$, the log-Riemann surface
of $\root n \of{z}$ (see example 3 in section I.1.)

\medskip

{\bf 2.} For $d=1$ and $n=0$ the log-Riemann surface $\cl S_{0,1}$
is $\cl S_{\log}-1$.

In general for $d\geq 1$ and $n=0$ the log-Riemann surface $\cl S_{0,d}$ is
a normalized Gauss
surface of degree $d$ (see examples 7 and 8 in section I.1.)

\medskip

{\bf 3.} When $n=d-1$ the integral is computable,
$$
\int_0^z t^{d-1} e^{-t^d / d} \ dt = - e^{-z^d} -1  \ ,
$$
and the cyclotomic log-Riemann $\cl S_{d-1,d}$ is easy to describe. Just
mate together $d$ copies of $(\cl S_{\log}-1) = \cl S_{0,1}$ by making
the same straight cut at $0$ in each copy, as shown in the
figure.

\medskip

{\hfill {\centerline {\psfig {figure=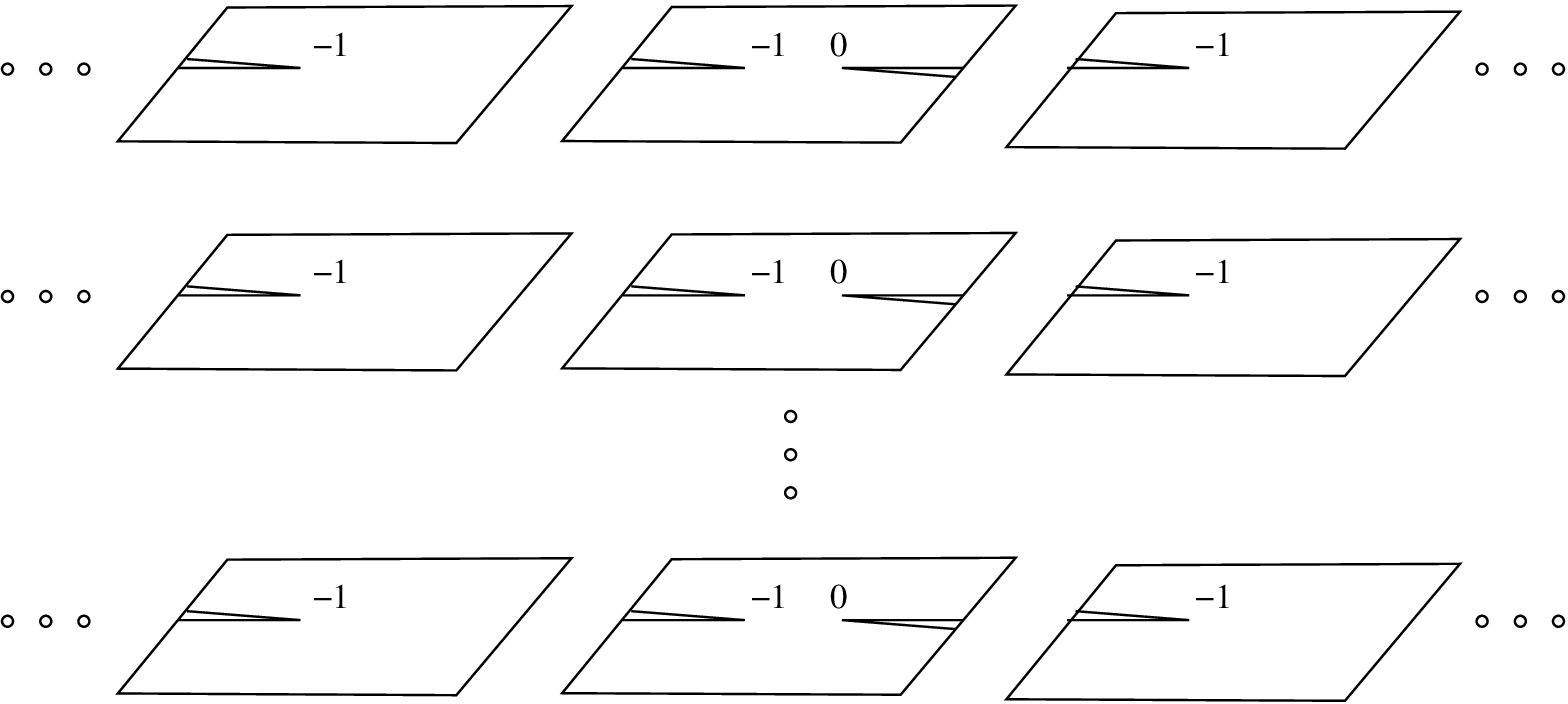,height=5cm}}}}

\medskip

\stit {II.6.2) Ramification values.}

The cyclotomic log-Riemann surface $\cl S_{n,d}$ has a unique finite ramification
point of order $n+1$ located above $0$.
The locations of the infinite ramification points are given by the evaluation
of the $\G$ function at rational values.

\eno {Theorem II.6.2.1}{ The projections of the $d$ infinite ramification points of $\cl S_{n,d}$ are given by $$
\pi (w_j^*) = \o_j^{n+1} d^{(n+1)/d - 1} \Gamma\left ( {n+1 \over d}\right ) \ , $$ where $\o_j$ is the $d$-th
root of $1$ $$ \o_j=e^{i{2j\pi \over d}} \ , $$ for $j=1,\ldots , d$. }

\stit {Proof.}

Each $\o_j$ gives a ramification point $w_j^*$ such that
$$\eqalign{
\pi(w_j^*) &= \int_0^{+\infty . \o_j} t^n e^{-t^d / d} \ dt \cr
           &= {\o_j}^{n+1} \int_0^{+\infty} s^n e^{-s^d / d} \ ds \cr
}$$
The change of variables $u=s^d / d$ gives
$$\eqalign{
{\o_j}^{n+1} \int_0^{+\infty} s^n e^{-s^d / d} \ ds &= \o_j^{n+1} d^{(n+1)/d - 1}
\int_0^{+\infty } u^{{n+1 \over d}-1} e^{-u}\ du \cr
     &=\o_j^{n+1} d^{(n+1)/d - 1}\Gamma\left ( {n+1 \over d}\right ) \ . \cr
}$$
$\diamond$

\stit {Remark.}

Very little is known about the transcendental character of the values of the $\Gamma$ function at rational
arguments. F. Lindemann's proof of the transcendence of $\pi$ ([Li]) proves the transcendence of $\G (1/2)=\root
\of {\pi}$. G.V. Chudnovsky proved the transcendence (and algebraic independence with $\pi$) of $\Gamma (1/3)$ and
$\Gamma (1/4)$ (see [Chu] and [Wa1]). Rohrlich's conjecture states that there are no multiplicative relations $$
\prod_{p/q \in \dd Q} \Gamma \left (p/q\right )^{n(p/q)} \in {\overline {\dd Q}}\ , $$ where the exponents
$n(p/q)$ are almost all zero, other than the trivial ones obtained from the basic functional equation, the
complement formula and the multiplicative formula of Euler-Legendre-Gauss (see the survey [Wa2] for more
information on this conjecture and other open problems.)

\bigskip

We can give a description of where the infinite branching takes place.

\eno {Theorem II.6.2.2}{ Let $0^*\in \cl S^*_{n,d}$ be the finite ramification point. All infinite ramification
points $w_j^*$ are at the same euclidean distance from $0*$, $$ d(0^*, w_j^*)= d^{(n+1)/d - 1} \Gamma\left ({n+1
\over d} \right ) \ , $$ and there is a unique geodesic $[0^*, w_j^*]$ realizing this distance which is an
euclidean segment. Two consecutive geodesics $[0^*, w^*_j]$ and $[0^*, w^*_{j+1}]$ form an angle of $2\pi (n+1)/d$
at their common vertex $0^*$.

The log-Riemann surface $\cl S_{n,d}$ can be obtained by grafting on
the log-Riemann surface of $\root n \of {z}$, $\cl S_n=\cl S_{n,0}$,
$d$ infinite ramification points regularly distributed and
at the distance from $0$ given above.
In general for $k\geq 1$, we can get $\cl S_{n,kd}$ from $\cl S_{n,d}$
by grafting $(k-1)d$ regularly distributed infinite ramification points
and rescaling the log-Riemann surface structure.
}

\stit {Proof.}
These statements follow from the symmetry  of the uniformization
$$
F_{n,d}\left (e^{2\pi i /d} z\right ) =e^{2\pi i (n+1)/d} F_{n,d} (z) \ ,
$$
that we obtain by a simple change of variables.$\diamond$

\stit {II.6.3) Subordination of cyclotomic log-Riemann surfaces.}

\eno {Theorem II.6.3.1}{ Let $m\geq 1$ be a divisor of $n+1$ and $d$. We have that $m^{(n+1)/d - 1} \cl S_{{n+1
\over m}-1, {d\over m}}$ is subordinate to $\cl S_{n,d}$, $$ \cl S_{n,d}\geq m^{(n+1)/d - 1} \cl S_{{n+1 \over
m}-1, {d\over m}} \ . $$ }

%% figure here

\stit {Proof.}

The change of variables $u=m^{-m/d}t^m$ in the integral
$$
F_n,d (z)=\int_0^z t^n e^{-t^d / d}\ dt
$$

gives the functional equation
$$
F_{n,d} (z)=m^{(n+1)/d - 1} F_{{n+1 \over m}-1, {d\over m}}(m^{-m/d}z^m) \ ,
$$
this gives the following commutative diagram which proves the theorem,

%%Diagram here

$\diamond$

%% figure here

\stit {II.6.4) Caratheodory limits of cyclotomic log-Riemann surfaces.}

\eno {Theorem II.6.4.1}{ We consider the normalized cyclotomic log-Riemann surfaces $$ \hat \cl S_{n,d} = {1 \over
d^{(n+1)/d - 1}\G\left ({n+1 \over d} \right )} \cl S_{n,d} \ . $$ Any Caratheodory limit of a pointed sequence
$(\hat \cl S_{n,d} , z_n )$ is either a planar log-Riemann surface $\dd C_l$, a translation of the normalized
Gauss log-Riemann surface of degree $2$, ${1\over 2} \cl S_{{\hbox {\rm Gauss}}}+{1\over 2} =\hat \cl S_{0,2}$, or
a translation of the log-Riemann surface of the logarithm, $\cl S_{\log} =\cl S_{0,1}+1$ }

\stit {Proof.}

If  $d(z_n, 0^*) \to +\infty$ we are in the first case. Otherwise, from the base point $z_n$ we can measure the
angle at $0*$ from $z_n$ to a ramification point $w_{n,j}^*$. Only one of these angular measures can stay bounded.
If there is one such ramification point, it gives one in the limit and we are in the second case (note that $0^*$
becomes the other infinite ramification point in the Gauss log-Riemann surface). Finally if all the angular
measures are unbounded and we have a Caratheodory limit, it has to be a log-Riemann surface as in the third case
(the only infinite ramification point is generated by $0^*$). $\diamond$

\stit {II.6.5) Continued fraction expansion of the uniformization.}

Consider the uniformization
$$
F_{n,d}(z) = \int_{0}^{z} t^n e^{-t^d/d} \ dt \ .
$$
Fix one of the $d$ half-lines
$$
\arg z = {2j\pi\over d}
$$
with $j=1,\dots,d$ along which this integral converges, and the
sector
$$
\left |\arg z - {2j\pi \over d }\right | < \pi/d
$$
of angle $2\pi/d$ and centered around this direction.

With the same notation as in II.6.2, we can write
$$\eqalign{
F_{n,d}(z) & = \int_{0}^{z} t^n e^{-t^d/d} \ dt \cr
           & = \int_{0}^{\omega_j\cdot\infty} t^n e^{-t^d/d} \ dt -
           \int_{z}^{\omega_j\cdot\infty} t^n e^{-t^d/d} \ ds \cr
           & = \pi(w_j^*) - \int_{z}^{\omega_j\cdot\infty} t^n e^{-t^d/d} \ dt \cr
}
$$
Repeated integrations by parts give the following
asymptotic series for the last integral appearing above:
$$
\int_{z}^{\omega_j\cdot\infty} t^n e^{-t^d/d} \ dt = z^n e^{-z^d/d} \cdot S,
$$
where $S$ has the asymptotic series
$$
S = {1 \over z^{d-1}} + {(n-d+1) \over z^{2d-1}} + {(n-d+1)(n-2d+1) \over z^{3d-1}} +
{(n-d+1)(n-2d+1)(n-3d+1) \over z^{4d-1}} + \dots
$$
Unless $n+1 \equiv 0 $ (mod $d)$, in which case
the series terminates, this series is divergent. One may try to convert
the divergent series $S$ into a convergent continued fraction with
polynomial coefficients. A direct way of doing this would go as follows: Put $S_0 = S$, and, considering the
first term of $S_0$,
$$
S_0 = {1 \over dz^{d-1} + S_1},
$$
where $S_1$ is small. Put $S_1 = a_1/z + a_2/z^2 + \dots$, substitute this series into the above equation,
and expand the fraction into a series in negative powers of $z$; equating this series to the series above
determines the $a_n$'s uniquely. If $a_m$ is the first nonzero $a_n$, then we can put
$$
S_1 = {a_m \over z^m + S_2},
$$
and $S_2 = b_1/z+b_2/z^2+\dots$, and now repeat the same procedure, applied this
time to $S_2$.

\medskip

There is no formal obstruction to continuing this procedure indefinitely, which it is clear leads to a
continued fraction representation for $S$. However, it is computationally intensive; we describe instead a
more elegant classical method for computing the continued fraction, due to Euler and Lagrange.

\medskip

The method is applicable to functions which satisy a Riccati equation, that is an equation of the form
$$
y' + A(z)y^2 + B(z)y + C(z) = 0,
$$
where $A,B,C$ are rational functions of $z$, and is based on the fact that the family of Riccati equations is
closed under Moebius transformations in the dependent variable $y$. Given a solution $y$ of such an equation,
we consider its asymptotics as $z$ approaches a point through a fixed set of directions, say for example as
in our case when $z \to \infty$ within the sector $|\arg z - {2j\pi / d }| < \pi/d$. If this is known, of the
form $y \sim a_1/z^{m_1}$ say, then we put
$$
y= {a_1 \over z^{m_1} + y_1}.
$$
Then $y_1$ satisfies a Riccati equation with new coefficients $A_1,B_1,C_1$; taking the asymptotics of $y_1$
to be of the form $y_1 \sim a_2/z^{m_2}$, if we substitute for $y_1$ an asymptotic series $y_1 = a_2/z^{m_2}
+ p_{m_2+1}/z^{m_2 + 1} + \dots$ in this equation, then the constants $a_2$ and $m_2$ are uniquely
determined by the Riccati equation. Hence we can put
$$
y_1 = {a_2 \over z^{m_2} + y_2}.
$$
As before, $y_2$ satisfies a Riccati equation, now with coefficients $A_2,B_2,C_2$, which, taking the
asymptotics of $y_2$ to be of the form $y_2 \sim a_3/z^{m_3}$, allows us as before to determine the constants
$a_3$ and $m_3$, so we proceed by putting
$$
y_2 = {a_3 \over z^{m_3} + y_3},
$$
and so on.

\medskip

In our case, we see by differentiating both sides of the expression
$$
\int_{z}^{\omega_j\cdot\infty} t^n e^{-t^d/d} \ dt = z^n e^{-z^d/d} \cdot S,
$$
that $S$ satisfies the Riccati equation
$$
S' + \left( {n \over z} - z^{d-1} \right) S + 1 = 0
$$
We apply the method described above, though for convenience we make a slight modification, that if a
coefficient $a_k$ is rational, say $a_k = p_k/q_k$, (where $S_{k-1} \sim a_k/z^{m_k}$) then we put $S_{k-1} =
p_k/(q_k z^{m_k} + S_k)$ instead of $S_{k-1} = (p_k/q_k) / (z^{m_k} + S_k)$. The general expression for the
constants $a_k, m_k$ is computable, and we have the following theorem:

\medskip

\eno{Theorem II.6.5.1}{ For $j=1,\dots,d$, in the sector $\left |\arg z - {2j\pi \over d }\right | < \pi/d$ the
function $F_{n,d}$ has a convergent continued fraction representation $$\eqalign{
 F_{n,d}(z) &= \pi(w_j^*) - z^n e^{-z^d/d} \cdot
{1\over\displaystyle z^{d-1} +
                             {\strut (d-1)-n\over\displaystyle z+
                                 {\strut d\over\displaystyle z^{d-1} +
                            {\strut (2d-1)-n\over\displaystyle z +
                                {\strut 2d\over\displaystyle z^{d-1} +
                            {\strut (3d-1)-n\over\displaystyle z +
                                {\strut 3d\over\displaystyle z^{d-1} + \dots
}}}}}}} \cr
&=\pi(w_j^*) - z^{n+1} e^{-z^d/d} \cdot
{1\over\displaystyle z^d +
                             {\strut (d-1)-n\over\displaystyle 1+
                                 {\strut d\over\displaystyle z^d +
                            {\strut (2d-1)-n\over\displaystyle 1 +
                                {\strut 2d\over\displaystyle z^d +
                            {\strut (3d-1)-n\over\displaystyle 1 +
                                {\strut 3d\over\displaystyle z^d + \dots.
}}}}}}} \cr
}$$
}

There are two parts to prove, first that the formal computations lead to the formula above, and second that
the formula converges to $F_{n,d}$ in the domains stated above. The first part follows from the following:

\medskip

\eno{Proposition II.6.5.2}{ For $k \geq 1$ we have the following relations: $$\leqalignno{ &S'_{2k} - z^{d-2} S^2
_{2k} + \left( {n \over z} - z^{d-1} \right) S_{2k} + kd = 0 & (1) \cr &S_{2k} = {kd \over z^{d-1} + S_{2k+1}} &
(2) \cr &S'_{2k+1} - S^2 _{2k+1} - \left( {n \over z} + z^{d-1} \right)S_{2k+1} - (n - ((k+1)d-1))z^{d-2} = 0 &
(3) \cr &S_{2k+1} = {((k+1)d -1) - n \over z + S_{2k+2}} & (4) \cr } $$ }

\medskip

\stit{Proof:} By induction on $k \geq 1$.

Carrying out the first few steps of the computation one can easily verify that for $k=1$ these relations hold
between the functions $S_{2k} = S_2, S_{2k+1} = S_3$ and $S_{2k+2} = S_4$.

\medskip

So assume that they hold for a $k \geq 1$. We show that they hold for $k+1$ :

\medskip

Multiplying both sides of $(4)$ by $(z+S_{2k+2}) z^{d-2}$ and rearranging terms gives
$$
z^{d-2} S_{2k+1}S_{2k+2} + z^{d-1} S_{2k+1} + (n - ((k+1)d-1))z^{d-2} = 0
$$
Adding this equation to $(3)$ and dividing the result by $S_{2k+1}$ gives
$$
{S'_{2k+1} \over S_{2k+1}} - S_{2k+1} - {n \over z} + z^{d-2} S_{2k+2} = 0
$$
Substituting for $S_{2k+1}$ using $(4)$ in the above equation gives
$$
{-(1+S'_{2k+2}) \over z + S_{2k+2} } + {(n - ((k+1)d-1)) \over z + S_{2k+2}} - {n \over z} + z^{d-2} S_{2k+2}
= 0
$$
Multiplying by $(z + S_{2k+2})$ and simplifying gives
$$
S'_{2k+2} - z^{d-2} S^2 _{2k+2} +  \left( {n \over z} - z^{d-1} \right) S_{2k+2} + (k+1)d = 0,
$$
thus relation $(1)$ holds for $k+1$.

\medskip

Substituting for $S_{2k+2}$ in the above equation a formal series in powers of $1/z$, $S_{2k+2} =
a_{2k+3}/z^{m_{2k+3}} + p_{m_{2k+3} + 1}/z^{m_{2k+3} + 1} + \dots$, where $m_{2k+3} \geq 0$, gives
$$
{ -a_{2k+3} \over z^{m_{2k+3}-(d-1)} } + (k+1)d + O \left( {1 \over z^{m_{2k+3} - d} } \right) = 0
$$
from which it follows that
$$
m_{2k+3} = d-1 \ , \ a_{2k+3} = (k+1)d,
$$
and hence
$$
S_{2k+2} = {(k+1)d \over z^{d-1} + S_{2k+3}}.
$$
Thus relation $(2)$ holds for $k+1$.

\medskip

Using relations $(1)$ and $(2)$ for $k+1$ it is straightforward to derive, in a manner similar to that above,
the relations $(3)$ and $(4)$ for $k+1$. \quad $\diamond$

\bigskip

The convergence of the continued fraction is proved in the following section.

\medskip

\stit {II.6.6) Relation to the incomplete Gamma function.}

The locations of the infinite ramification points of the
cyclotomic log-Riemann surfaces are given by values of the Gamma
function at rational arguments. More generally, the
uniformizations $F_{n,d}$ can be expressed in terms of the incomplete Gamma
function defined by (see [MOS] chapter IX)
$$
\Gamma(a,z) = \int_{z}^{+\infty} t^{a-1} e^{-t} \ dt
$$
where $a$ is a parameter. Thus $\Gamma (a,0)=\Gamma (a)$.

\medskip

Consider the uniformization
$$
F_{n,d} = \int_{0}^{z} t^n e^{-t^d/d} \ dt \ .
$$
As before, fixing one of the $d$ half-lines
$$
\arg z = {2j\pi\over d}
$$
with $j=1,\dots,d$ along which this integral converges, and the
sector
$$
\left |\arg z - {2j\pi \over d }\right | < \pi/d
$$
of angle $2\pi/d$ and centered around this direction, we can write
$$
F_{n,d}(z) = \pi(w_j^*) - \int_{z}^{\omega_j\cdot\infty} t^n e^{-t^d/d} \ dt
$$
Now making the change of variables $s = t^d/d$ in the above integral
gives
$$
F_{n,d}(z) = \pi(w_j^*) \ - \ {\omega_j}^{n+1} d^{(n+1)/d - 1}
\int_{z^d/d}^{+1\cdot\infty} s^{(n+1)/d - 1} e^{-s} \ ds
$$
thus, in terms of the incomplete Gamma function:

\medskip

\eno {Proposition II.6.6.1}{ We have $$ F_{n,d}(z) = \pi(w_j^*) \ - \ {\omega_j}^{n+1} d^{(n+1)/d - 1}
\Gamma\left({n+1 \over d},\ {z^d \over d} \right) $$ }

Note that putting $z = 0$ on both sides above gives the result in
section II.6.2 on the locations $\pi(w_j^*)$ of the infinite
ramification points.

\medskip

We can now prove the convergence of the continued fraction
for $F_{n,d}$ as follows:

\medskip

\stit{Proof of Theorem II.6.5.1 :} We recall the following continued fraction for the incomplete Gamma function
(see [Wall], pg.356):

\medskip

$$
\int_{u}^{+\infty} t^{a-1} e^{-t} \ dt = u^a e^{-u} \cdot
{1\over\displaystyle u +
                             {\strut 1-a\over\displaystyle 1+
                                 {\strut 1\over\displaystyle u +
                            {\strut 2-a\over\displaystyle 1 +
                                {\strut 2\over\displaystyle u +
                            {\strut 3-a\over\displaystyle 1 +
                                {\strut 3\over\displaystyle u + \dots.
}}}}}}}
$$
This continued fraction is convergent for $u$ in the slit plane
$\dd C - ]-1\cdot\infty, 0]$.

\medskip

For $z$ in a sector $\left |\arg z - {2j\pi \over d }\right | < \pi/d$,
the variable $u = z^d/d$ lies in this slit plane. From the
preceding proposition and the above formula with $u = z^d/d$, $a = (n+1)/d$,
we have
$$\eqalign{
F_{n,d}(z) &= \pi(w_j^*) - {\omega_j}^{n+1}
d^{{n+1 \over d} - 1} \left({z^d \over d}\right)^{n+1 \over d} e^{-z^d/d} \cdot
{1\over\displaystyle z^d/d +
                             {\strut 1-{n+1 \over d}\over\displaystyle 1+
                                 {\strut 1\over\displaystyle z^d/d +
                            {\strut 2-{n+1 \over d}\over\displaystyle 1 +
                                {\strut 2\over\displaystyle z^d/d +
                            {\strut 3-{n+1 \over d}\over\displaystyle 1 +
                                {\strut 3\over\displaystyle z^d/d + \dots.
}}}}}}} \cr
&=\pi(w_j^*) - z^{n+1} e^{-z^d/d} \cdot
{1\over\displaystyle z^d +
                             {\strut (d-1)-n\over\displaystyle 1+
                                 {\strut d\over\displaystyle z^d +
                            {\strut (2d-1)-n\over\displaystyle 1 +
                                {\strut 2d\over\displaystyle z^d +
                            {\strut (3d-1)-n\over\displaystyle 1 +
                                {\strut 3d\over\displaystyle z^d + \dots.
}}}}}}} \cr
}$$
(noting that
$$
{\omega_j}^{n+1} (z^d)^{n+1 \over d} = z^{n+1}
$$
in the sector $\left |\arg z - {2j\pi \over d }\right | < \pi/d$,
where the $d$-th root above is positive for $z^d$ real and positive).

\medskip

The above continued fraction is the same as the one appearing in Theorem II.6.5.1. $\diamondsuit$

\stit {II.6.7) Relation to Hermite polynomials.}

%% file: ii7.tex
\input defv7.tex
\input psfig.sty

\stit {II.7) Uniformization of infinite log-Riemann surfaces.}

\medskip

Let $\cl S$ be a log-Riemann surface such that its finite
completion $\cl S^{\times}$ is simply connected.

\medskip

We have seen in section II.5 that if the ramification set $\cl R$
is finite then the surface $\cl S^{\times}$ is parabolic, and the uniformization
is given, if $\cl S$ has no finite ramification points, by a primitive of the form
$$
\int \ e^{P(z)} \ dz
$$
where $P$ is a polynomial,
or, more generally, when $\cl S$ has both finite and infinite ramification points,
 by a primitive of the form
$$
\int \ Q(z) e^{P(z)} \ dz
$$
where $P$ and $Q$ are polynomials. The degrees of the polynomials $P$ and $Q$
correspond exactly to the numbers of infinite and finite ramification points
(counted with multiplicity) of $\cl S$ respectively.

\medskip

One may try and extend this correspondence between primitives and log-Riemann
surfaces to the case where the surfaces $\cl S$ have an infinite number of
ramification points, possibly by considering primitives of the form
$$
\int \ e^{h(z)} \ dz
$$
or more generally of the form
$$
\int \ g(z) e^{h(z)} \ dz,
$$
where $g$ and $h$ are no longer polynomials but instead entire functions. This
general setting poses considerably more difficulties however. For example,
considering surfaces with an infinite number of ramification points allows
for surfaces which are not necessarily parabolic but instead hyperbolic.
On the other hand, considering primitives as above with entire functions
may give rise to surfaces with more general log-Riemann surface structures
than that considered here, namely allowing charts with non-locally finite sets
of cuts. These and other questions are considered in the forthcoming [Bi-PM1].
For the moment we restrict ourselves to describing a few examples.

\medskip

{\bf 1)} The primitive $\int e^{e^z} \ dz$ :

\medskip

Let $F$ be the entire function
$$
F(z) = \int_{0}^{z} e^{e^t} \ dt
$$
The function $F$ defines a uniformization $\tilde{F} : \dd C \to \cl S$ to
a simply connected log-Riemann surface $\cl S$ such that
$F = \pi \circ \tilde{F}$, where $\pi : \cl S \to \dd C$ is the
projection mapping. The log-Riemann surface $\cl S$ is shown in
the figure below.

\medskip

%{\hfill {\centerline {\psfig {figure=figures/ii7fig1.eps,height=6cm}}}}

\medskip

The surface $\cl S$ has an infinite number of ramification points
$(w_n^*)_{n \in \dd Z}$ all of infinite order. These ramification points
are placed in a common base sheet at the points $(a_n)_{n \in \dd Z}$
given by
$$\eqalign{
a_0 &= \pi(w_0^*) = \int_{0}^{i\pi + 1\cdot\infty} e^{e^t} \ dt \cr
a_n &= \pi(w_n^*) = \int_{0}^{(2n+1)i\pi + 1\cdot\infty} e^{e^t} \ dt = a_0 + 2n\pi i \ , \ n \in \dd Z. \cr
}$$

\medskip

{\bf 2)} The primitive $\int e^{e^z + e^{-z}} \ dz$:

\medskip

Let $F$ be the entire function
$$
F(z) = \int_{0}^{z} e^{e^t + e^{-t}} \ dt
$$

The function $F$ defines a uniformization $\tilde{F} : \dd C \to \cl S$ to
a simply connected log-Riemann surface $\cl S$ such that
$F = \pi \circ \tilde{F}$, where $\pi : \cl S \to \dd C$ is the
projection mapping. The log-Riemann surface $\cl S$ is shown in
the figure below.

\medskip

%{\hfill {\centerline {\psfig {figure=figures/ii7fig2.eps,height=6cm}}}}

\medskip

The surface $\cl S$ has an infinite number of ramification points
$(v_n^*)_{n \in \dd Z}, \, (w_n^*)_{n \in \dd Z}$, all of infinite order.
These ramification points are placed in a common base sheet at the
points $(a_n)_{n \in \dd Z}, \, (b_n)_{n \in \dd Z}$ given by
$$\eqalign{
a_0 &= \pi(v_0^*) = \int_{0}^{i\pi + 1\cdot\infty} e^{e^t + e^{-t}} \ dt = F(\pi i) + \int_{0}^{\infty}
e^{-(e^s + e^{-s})} \ ds \cr
a_n &= \pi(v_n^*) = \int_{0}^{(2n+1)i\pi + 1\cdot\infty} e^{e^t + e^{-t}} \ dt = a_0 + nF(2\pi i) \ , \ n \in
\dd Z \cr
b_0 &= \pi(w_0^*) = \int_{0}^{i\pi - 1\cdot\infty} e^{e^t + e^{-t}} \ dt = F(\pi i) - \int_{0}^{\infty}
e^{-(e^s + e^{-s})} \ ds \cr
b_n &= \pi(w_n^*) = \int_{0}^{(2n+1)i\pi - 1\cdot\infty} e^{e^t + e^{-t}} \ dt = b_0 + nF(2\pi i) \ , \ n \in
\dd Z \cr
}$$
(note that $F(z + 2\pi i) = F(z) + F(2\pi i)$).

%% file: iii1.tex
\input defv7.tex
\bigskip

{\centerline {\tir {III. Algebraic theory of log-Riemann surfaces.}}}

\bigskip

\stit {III.1) A ring of special functions.}

\stit {III.1.1) Definition.}

Let $P_0(z)\in \dd C[z]$ be a polynomial of degree $d\geq 1$
$$
P_0(z)=a_d z^d+a_{d-1} z^{d-1}+\ldots +a_1 z+a_0 \ .
$$
We consider the entire functions
$$
\eqalign {
F_0(z)&=\int_0^z e^{P_0(t)} \ dt \cr
F_1(z)&=\int_0^z t\  e^{P_0(t)} \ dt \cr
\ldots & \cr
F_{d-1}(z)&=\int_0^z t^{d-1} e^{P_0(t)} \ dt \cr
}
$$

Observe that a $\dd C$-linear combination of these special functions
and the constant unit function
generate $e^{P_0}$

$$
e^{P_0}=e^{P_0(0)}.1+a_1 F_0+2a_2 F_1+\ldots (d-1)a_{d-1} F_{d-2}+
da_d F_{d-1} \ .
$$

\stit {III.1.2) Asymptotics.}

The following asymptotic estimate is a basic tool in the proofs
of the algebraic results.

\eno {Proposition III.1.2.1}{ For $j=0,1,\ldots, d-1$ we have $$ F_j(z)\sim {z^j\over P_0'(z)} e^{P_0(z)} $$ when
$z\to +\infty . a_d^{-1/d}$, that is when $z\to \infty$ in a direction given by a $d$-root of $a_d^{-1}$. }

\stit {Proof.}

The asymptotics in these directions is $+\infty$, thus we can assume that $P_0$ is non zero at $0$ by changing the
origin of integration (i.e. by a translation change of variables in the integrals).

By two integration by parts we get
$$
\eqalign {
&F_j(z) =\int_0^z t^j e^{P_0(t)} \ dt  =\int_0^z {t^j\over P'_0(t)}\ P'_0(t)
e^{P_0(t)} \ dt \cr
&=\left [{t^j\over P'_0(t)}\right ]_0^z -\int_0^z \left (
{jt^{j-1} P'_0(t) -t^j P''_0(t) \over \left (P'_0(t)\right )^2}\right )\ e^{P_0(t)}
\ dt \cr
&={z^j \over P'_0(z)}e^{P_0(z)}-\int_0^z \cl O \left({1\over
(a_d^{1/d} t)^{d-j}}\right )\ e^{P_0(t)}
\ dt \cr
&={z^j \over P'_0(z)}e^{P_0(z)}-\left [\cl O \left ({1\over (a_d^{1/d} t)^{d-j}}\right )
{1\over P'_0(t)} e^{P_0(t)}\right ]_0^z +\int_0^z \cl O \left ({1\over
(a_d^{1/d} t)^{2d-j-1}}\right )\ e^{P_0(t)} dt \cr
}
$$
Now the second and last term in the last equation are neglectable with
respect to the first one.$\diamond$

\stit {III.1.3) Linear independence.}

\eno {Proposition III.1.3.1}{ The special functions $F_0, F_1, \ldots , F_{d-1}$ and the constant function $1$ are
linearly independent over $\dd C$.}

\stit {Proof.}

Consider a non-trivial linear combination
$$
b_{-1}+b_0 F_0+b_1 F_1+\ldots +b_{d-1} F_{d-1}=0 \ ,
$$
and take one derivative. Dividing by $e^{P_0}$ we get
$$
b_0+b_1 z+\ldots +b_{d-1} z^{d-1}=0 \ .
$$
Thus $b_0=b_1=\ldots =0$ and then $b_{-1}=0$ also.$\diamond$

We give another two proofs one analytic and another more algebraic.

Take a non-trivial linear combination
$$
b_{-1}+b_0 F_0+b_1 F_1+\ldots +b_{d-1} F_{d-1}=0
$$
and let $0\leq k\leq d-1$ be the largest index such that $b_k\not=0$.
If $k=-1$ we are done. If not,
when $z\to + \infty . a_d^{-1/d}$ we have
$$
b_{-1}+b_0 F_0+b_1 F_1+\ldots +b_{d-1} F_{d-1} \sim b_k{z^k \over
P'_0(z)} e^{P_0(z)}\to \infty
$$
Contradiction.$\diamond$

\medskip

We give a more algebraic proof. First we show that $F_0, \ldots , F_{d-1}$
are $\dd C$-linearly independent.
By contradiction choose $d$ distinct points
$z_0, z_1, \ldots , z_{d-1} \in \dd C$. For $k=0,\ldots , d-1$ we have
$$
b_0 F_0 (z_k) +b_1 F_1(z_k)+\ldots +b_{d-1} F_{d-1} (z_k) =0 \ .
$$

Therefore
$$
\D(z_0, \ldots , z_{d-1})=\left |
\matrix {F_0(z_0) & F_0 (z_1) & \ldots & F_0(z_{d-1}) \cr
F_1(z_0) & F_1(z_1)& \ldots & F_1(z_{d-1})\cr
\vdots & \vdots & \ddots &\vdots \cr
F_{d-1}(z_0) & F_{d-1}(z_1) &\ldots &F_{d-1}(z_{d-1})\cr
}\right | =0
$$
But  we have
$$
\partial_{z_{d-1}} \ldots \partial_{z_1} \partial_{z_0} \ \D =
e^{P_0(z_0)} . e^{P_0(z_1)} \ldots e^{P_0(z_{d-1})}.
\left | \matrix {1 & 1 & \ldots & 1 \cr
z_0 & z_1& \ldots & z_{d-1}\cr
\vdots & \vdots & \ddots &\vdots \cr
z_0^{d-1} & z_1^{d-1} &\ldots &z_{d-1}^{d-1}\cr
} \right |
$$
and the Vandermonde determinant is not zero, thus
$\partial_{z_{d-1}} \ldots \partial_{z_1} \partial_{z_0} \ \D \not= 0$
and $\D$ is not identically $0$. Contradiction.

In order to show that $1, F_0, \ldots , F_{d-1}$ are $\dd C$-linearly
independent we proceed in a similar way evaluating the linear combination
at $d+1$ points $z_0, z_1, \ldots , z_d$. Next we apply the
differential operator $\partial_{z_1, z_2, \ldots , z_d}$ to the
Cramer determinant and develop through the first column.
$\diamond$

\medskip

We can prove more.

\eno {Proposition III.1.3.2}{ The special functions $F_0, F_1, \ldots , F_{d-1}$ and the constant function $1$ are
linearly independent over $\dd C [z]$.}

\stit {Proof.}

 By contradiction consider a non-trivial linear combination with polynomial
coefficients
$$
A_{-1}(z)+A_0(z)F_0(z)+\ldots +A_{d-1}(z) F_{d-1}(z) =0 \ .\leqno {(*)}
$$
Taking one derivative we get
$$
A'_{-1}(z)+A'_0(z)F_0(z)+\ldots +A'_{d-1}(z) F_{d-1}(z)=Q_1(z)e^{P_0(z)} \ ,
$$
where $Q_1(z)=-A_0(z)-zA_1(z)-\ldots -z^{d-1} A_{d-1}(z)$.
Iterating this procedure and taking $k$ derivatives, we get
$$
A_{-1}^{(k)}(z)+A^{(k)}_0(z)F_0(z)+\ldots
+A^{(k)}_{d-1}(z) F_{d-1}(z)=Q_k(z)e^{P_0(z)} \ ,
$$
where $Q_k(z)\in \dd C [z]$. Choose $k\geq 0$ minimal such that all $A^{(k)}_j$
are constant but not all $0$.

Let $-1\leq l_0 \leq d-1$ be the largest index such that $A_{l}^{(k)}\not= 0$.
Then if $l_0\geq 0$, when $z\to +\infty . a_d^{-1/d}$ we have the asymptotics
$$
A_{-1}^{(k)}+A_0^{(k)}F_0(z)+\ldots +A_{d-1}^{(k)} F_{d-1}(z)\sim A_{l_0}^{(k)}{z^{l_0}
\over P'_0(z)} e^{P_0(z)} \ .
$$
If $l_0=-1$, when $z\to +\infty . a_d^{-1/d}$ we have the asymptotics
$$
A_{-1}^{(k)}+A_0^{(k)}F_0(z)+\ldots +A_{d-1}^{(k)} F_{d-1}(z)\sim A_{-1}^{(k)}
$$
Therefore, in both cases, if $l_0<d-1$, we must have $Q_k\equiv 0$, thus
$$
A_{-1}^{(k)}+A_0^{(k)}F_0(z)+\ldots +A_{d-1}^{(k)} F_{d-1}(z) =0
$$
which is a non-trivial $\dd C$-linear combination of $1, F_0, \ldots , F_{d-1}$
which contradict the previous proposition. Thus $l_0=d-1$, and the degree
of $A_j$ is at most the degree of $A_{d-1}$. When $z\to +\infty . a_d^{-1/d}$
we have that $A_{d-1} F_{d-1}$ dominates $A_j F_j$ for $j< d-1$. Thus if $c$ is
the leading coefficient of
$A_{d-1}(z)$ and $m$ its degree then, when $z\to +\infty . a_d^{-1/d}$, we have
$$
A_{-1}(z)+ A_0(z)F_0(z)+\ldots +A_{d-1}(z) F_{d-1}(z)
\sim c {z^{n+d-1} \over P'_0(z)} e^{P_0(z)} \ .
$$
On the other hand $A_{-1}+A_0F_0+\ldots +A_{d-1} F_{d-1}\equiv 0$,
so $c$ must be $0$. Contradiction.$\diamond$

\stit {III.1.4) Algebraic independence.}

\eno {Theorem III.1.4.1}{ Let $\dd K_{P_0}$ be the field generated by $F_0, \ldots , F_{d-1}$ from $\dd C(z)$, $$
\dd K_{P_0}=\dd C(z)(F_0 ,\ldots , F_{d-1}) =\dd C (z, F_0, \ldots , F_{d-1}) =\dd C[z](F_0, \ldots , F_{d-1}) \ .
$$ The field $\dd K_{P_0}$ is the field of fractions of the ring $$ \dd A_{P_0}=\dd C[z] [F_0, \ldots F_{d-1}] \ .
$$ The field $\dd K_{P_0}$ has transcendence degree $d$ over $\dd C(z)$. }

It is clear that the transcendence degree is at most $d$. That it is
exactly $d$ follows from the next result:

\eno {Proposition III.1.4.2}{ For $k=0,\ldots , d-1$, $F_k$ is transcendental over $\dd C (z, F_0, \ldots ,
F_{k-1})$.}

%We have the following equivalent statement:
%
%\eno {Proposition.}{For $k=0,\ldots , d-2$, $F_k$ is transcendental
%over $\dd C (z, e^{P_0},  F_0, \ldots , F_{k-1})$.}
%
%We use the following standard lemma in differential algebra (see
%[Kaplanski] Lemma 3.9 p. 23).
%
%
%\eno {Lemma.}{Let $K\subset L$ be two differential fields of characteristic
%$0$. Let $y\in L$ such that $y'=x\in K$ but $x$ is not a derivative in $K$.
%Then $y$ is transcendental over $K$.
%}
%
%
%\stit {Proof of the Lemma.}
%
%If $y$ is algebraic over $K$, let
%$$
%y^n+by^{n-1}+\ldots =0
%$$
%be the irreducible equation. Differentiating this equation we get
%$$
%n y^{n-1} x+b' u^{n-1}+\ldots =0 \ .
%$$
%It follows that $x=-b'/n=(-b/n)'$ contradicting that $x$ is not
%a derivative in $K$.$\diamond$
%

\eno {Definition III.1.4.3}{ The exponential degree and the polynomial degree of a monomial expression $$ z^m
F_0^{n_0} F_1^{n_1} \ldots F_{d-1}^{n_{d-1}} $$ are respectively $| \dd n |=n_0+n_1+\ldots +n_{d-1}$ and
$m+n_1+2n_2+\ldots +(d-1) n_{d-1}=m+\dd {(d-1)}.\dd n$ where $\dd {(d-1)}$ denotes the vector $(0,1,2,\ldots ,
d-1)$, and $\dd n$ the vector $(n_1, \ldots , n_{d-1})$. }

\eno {Lemma III.1.4.4}{ In a vanishing linear combination of monomials in $z, F_0, \ldots F_{d-1}$ each sub-linear
combination of monomials with the same exponential and polynomial degree must vanish.}

\stit {Proof.} Note the asymptotics when $z\to +\infty .
a_d^{-1/d}$,

$$
\eqalign{
z^m F_0^{n_0} F_1^{n_1} \ldots F_{d-1}^{n_{d-1}} & \sim {z^{m+n_1+2n_2
+\ldots +(d-1)n_{d-1} }\over (P'_0(z))^{n_0+n_1+\ldots +n_{d-1}} }
e^{(n_0+n_1+\ldots +n_{d-1})P_0(z)} \cr
& \sim z^{m+\dd {(d-1)}.\dd n-|\dd n |}
e^{|\dd n |.P_0(z)} \cr
}
$$
Now consider a vanishing $\dd C$-linear combination of monomials
$$
\eqalign {
&\sum_{m, \dd n} a_{m, \dd n} \
z^m F_0^{n_0} F_1^{n_1} \ldots F_{d-1}^{n_{d-1}} \cr
&=\sum_{N\geq 0} \sum_{m, \dd n \atop |\dd n |=N} a_{m, \dd n } \
z^m F_0^{n_0} F_1^{n_1} \ldots F_{d-1}^{n_{d-1}} \cr
&=0 \cr
}
$$
The different exponential asymptotics show that for each $N\geq 0$,
$$
\eqalign{
0&=\sum_{m, \dd n \atop |\dd n |=N} a_{m, \dd n} \
z^m F_0^{n_0} F_1^{n_1} \ldots F_{d-1}^{n_{d-1}} \cr
&=\sum_{m\geq 0} \sum_{\dd n \atop |\dd n|=N} a_{m, \dd n} \
z^m F_0^{n_0} F_1^{n_1} \ldots F_{d-1}^{n_{d-1}} \ . \cr
}
$$
Again the same argument using the different asymptotics for monomials
with the same exponential degree but different polynomial degree show the
result, that is, for each $N\geq 0$ and $m\geq 0$,
$$
\sum_{\dd n \atop |\dd n |=N} a_{m,\dd n} \ z^m F_0^{n_0} F_1^{n_1} \ldots F_{d-1}^{n_{d-1}}=0 \ .
$$
$\diamond$

\eno {Lemma III.1.4.5}{ Let $N\geq 1$. The monomials $F_0^{n_0} F_1^{n_1}\ldots F_k^{n_k}$ of exponential degree
$N$ are $\dd C[z]$-linearly independent. }

\stit {Proof.}

We prove the result by induction on $N\geq 1$.

For $N=1$ the result is given by Proposition III.1.3.2. Assume the result for $N-1$ and consider, by
contradiction, a non-trivial $\dd C[z]$ linear dependence relation $$ \sum_{\dd n} A_{\dd n} (z) \  F_0^{n_0}
F_1^{n_1} \ldots F_k^{n_k} =0 \ . $$ We can assume using the previous lemma that each term in this sum has the
same polynomial degree (we could also assume with the same reasoning that each polynomial $A_n(z)$ is a monomial,
but we don't need that.) This means that there exists a constant $K$ such that for each $n$ $$ \deg A_{\dd n}
+{\dd k}.{\dd n}=K $$ where $\dd k =(0,1,2,\ldots , k)$.

Apply one more derivative to the precedent relation to get
$$
\sum_{\dd n} A'_{\dd n} (z) \ F_0^{n_0} F_1^{n_1} \ldots F_k^{n_k}
=\sum_{{\dd n} \atop j=0,1,\ldots, k}
z^j A_n(z)
F_0^{n_0}  \ldots F_j^{n_j-1}\ldots  F_k^{m_k} e^{P_0} \ .
$$
Note that the exponential degree of the terms on the right side
remain the same as the one on the left side, but the polynomial degrees
are greater by $1$, thus
$$
\sum_{\dd n} A'_{\dd n} (z) \ F_0^{n_0} F_1^{n_1} \ldots F_k^{n_k} =0 \ .
$$
We can continue taking derivatives stopping one step before all
$A^{(l+1)}$ vanish, that is when
$$
\sum_{\dd n} A^{(l)}_{\dd n}  \ F_0^{n_0} F_1^{n_1} \ldots F_k^{n_k} =0 \ ,
$$
is a non-trivial $\dd C$-linear combination of homogeneous monomials
on the $F_j$'s.
Observe now that taking one more derivative in this last relation
and dividing by $e^{P_0}$ gives
$$
\sum_{\dd n \atop j=0,1,\ldots , k} A^{(l)}_{\dd n} \  z^j F_0^{n_0}
\ldots F_{j-1}^{n_{j-1}}  F_j^{n_j-1}  F_{j+1}^{n_{j+1}} \ldots F_k^{n_k} =0 \ .
$$
Remark that each monomial in $z, F_0, \ldots , F_k$ appearing in this
sum comes from exactly one monomial in $F_0, \ldots , F_k$ of the
relation before differentiation. Thus this last relation is a non-trivial
$\dd C[z]$-linear combination between monomials of exponential degree
$N-1$. By induction assumption this is impossible.$\diamond$.

\stit {Proof of the theorem.}

If $F_k$ is not transcendental over $\dd C(z, F_0, \ldots , F_{k-1})$, then
we have a non-trivial polynomial relation between $z, F_0, \ldots F_k$.
Isolating parts of the same exponential degree we are lead to a
non-trivial $\dd C[z]$-linear relation between homogeneous monomials
in $F_0, \ldots , F_k$ which contradicts the previous Proposition.
$\diamond$.

\stit {III.1.5) Integrals.}

Using the special functions $F_0, F_1, \ldots, F_{d-1}$ we can compute a large
family of integrals.

\eno {Theorem III.1.5.1}{ We consider the $\dd C$ vector space

$$
\eqalign {
\dd V_{P_0}=\dd V_{P_0}(\dd C) &=\dd C[z] . e^{P_0(z)} \oplus \dd C . 1 \oplus
\dd C . \ F_0 \oplus \ldots
\oplus \dd C . \ F_{d-2} \cr
&=z\dd C[z] .  e^{P_0(z)} \oplus \dd C . 1 \oplus
\dd C . \ F_0 \oplus \ldots
\oplus \dd C . \ F_{d-1}\cr
}
$$
For $Q(z)\in \dd C[z]$, any primitive
$$
\int_0^z Q(t) \ e^{P_0(t)} \ dt
$$
is in the vector space $V_{P_0}$.

Conversely, any element of $\dd V_{P_0}$  vanishing at $0$  is such a primitive.
}

\stit {Proof.}

The equality of the two direct sums results from the fact that
$e^{P_0}$ is a $\dd C$-linear combination of $F_0, \ldots , F_{d-1}$.

We prove the result by induction on the degree of $Q$.
The result is clear for $\deg Q \leq d-2$ because then $\int Q e^{P_0}$ is
a linear combination of $1, F_0, \ldots , F_{d-2}$.

For $\deg Q\geq d-1$, we consider the euclidean division of $Q$ by $P'_0$,
$$
Q=AP'_0 +B
$$
where $A,B\in \dd C[z]$ and $\deg B <d-1$.
Then, by integration by parts,
$$
\eqalign {
&\int_0^z Q(t) \ e^{P_0(t)} \  dt \cr
&= \int_0^z (A(t)P'_0(t) +B(t)) e^{P_0(t)} \ dt \cr
&=\left [A(t)e^{P_0(t)}\right ]_0^z -\int_0^z A'(t)  e^{P_0(t)} \ dt
+\int_0^z B(t)\
e^{P_0(t)} \ dt \cr
&=A(z)e^{P_0(z)} -A(0)e^{P_0(0)} -\int_0^z A'(t)  e^{P_0(t)} \ dt
+\int_0^z B(t)\
e^{P_0(t)} \ dt \ .\cr
}
$$
Now $A(z)e^{P_0(z)} \in \dd C e^{P_0(z)}$, $-A(0)e^{P_0(0)} \in \dd C$,
the primitive
$$
\int_0^z B(t)\ e^{P_0(t)} \ dt
$$
is a linear combination of $1, F_0, \ldots F_{d-2}$.
Moreover we have  $\deg A' < \deg Q$. Therefore
the result follows by induction.

\medskip

Now we prove the converse.
Let $F\in \dd V_{P_0}$ vanishing at $0$. Write
$$
F(z)=zP(z)e^{P_0(z)} +c_0+c_1 F_0+\ldots c_{d} F_{d-1} \ ,
$$
where $P(z)\in \dd C[z]$ and $c_0, c_1, \ldots c_{d} \in \dd C$.
Since $F(0)=0$ we have $c_0=0$. Also
$$
c_1 F_0+\ldots c_{d} F_{d-1} =\int_0^z
(c_1 +c_2 t+\ldots +c_{d} t^{d-1} )\ e^{P_0(t)} \ dt \ ,
$$
and
$$
zP(z) e^{P_0(z)} =\int_0^z (P(t)+tP'(t)+tP(t)P'_0(t) ) \ e^{P_0(t)} \ dt  \ .
$$

$\diamond$

\stit {Addenda.}

{\bf 1.} Let $\dd K\subset \dd C$ be a field. If
$P_0(z)\in \dd K[z]$ and $P_0$ is
normalized such that $P_0(0)=0$, then any primitive
$$
\int_0^z Q(t) \ e^{P_0(t)} \ dt
$$
where $Q(z)\in \dd K[z]$ belongs to the $\dd K$-vector space
$$
\dd V_{P_0}(\dd K)=z\dd K[z] e^{P_0(z)} \oplus \dd K \oplus \dd K \ F_0
\oplus \ldots
\oplus \dd K \ F_{d-1} \ .
$$
This results from the above proof since the Euclidean division
of polynomials is well defined in the ring $\dd K [z]$, and $e^{P_0(0)}=1$.
The proof of the converse statement is also the same.

\bigskip

{\bf 2.} In general, let $\dd K$ be a field and consider the differential
field $\dd K[z]$. For $P_0\in \dd K[z]$, $\deg P_0 =d$,
we define $e^{P_0}$ as generating
the Liouville extension defined by the differential equation
$$
y'-P_0 y=0 \ .
$$
We consider the extension $\dd K_0$ generated by
$$
\eqalign {
y' &= e^{P_0} \cr
y' &= z e^{P_0} \cr
&\vdots \cr
y' &= z^{d-1} e^{P_0} \ ,\cr
}
$$
and denote by $F_0, F_1, \ldots , F_{d-1}$ these primitives.
Then the $\dd K$-vector space
$$
M_{P_0}=z\dd K[z] e^{P_0} \oplus \dd K . 1\oplus \dd K . F_0 \oplus \ldots
\oplus \dd K . F_{d-1}
$$
coincides with the set of
all primitives $\int Q e^{P_0}$ modulo constants.

\medskip

\stit {III.1.6) Differential ring structure.}

We denote by $D$ the differentiation in the ring $\dd A_{P_0}$.
Let $\dd A^{N,n}_{P_0}$ be the $\dd C$-module generated
by those monomials of exponential degree $N$ and polynomial
degree $n$. We have the graduation
$$
\dd A_{P_0}=\bigoplus_{N,n\geq 0} \dd A_{P_0}^{N,n} \ .
$$
The following proposition is immediate.

\eno {Proposition III.1.6.1}{ We have $$ D \dd A_{P_0}^{N,n} \subset \dd A_{P_0}^{N,n-1} \oplus (\dd
A_{P_0}^{N-1,n} \oplus \dd A_{P_0}^{N-1,n+1}\oplus \ldots \oplus \dd A_{P_0}^{N-1,n+d-1}) e^{P_0} \ . $$ In
particular, the principal ideal $(e^{P_0})$ generated by $e^{P_0}$ is absorbing for the derivation, i.e. any
element of $\dd A_{P_0}$ falls into $(e^{P_0})$ after a finite number of derivatives. }

\eno {Proposition III.1.6.2}{ The only elements in $\dd A_{P_0}$ with no zeros are $$ \dd C^* \cup \{ e^{nP_0} ;
n\geq 1\} \ , $$ that is, the non-zero constant functions and $e^{P_0}, e^{2P_0}, \ldots$.

The group of units in $\dd A_{P_0}$ is composed of the non-vanishing constant functions $$ \dd A_{P_0}^{\times}
=\dd C^*\ . $$ }

\stit {Proof.}
Let $F\in \dd A_{P_0}$ with no zeros. Since $\dd A_{P_0}$ is a ring of
entire functions of order at most $p$, and $F$ is zero free, we can
find a polynomial of degree less than $d$ such that
$$
F=e^Q \ .
$$
Now, when $z\to +\infty . a_d^{-1/d}$ the asymptotics of each
element in $F\in \dd A_{P_0}$ is of the form
$$
F(z)\sim c z^{a}e^{bP_0(z)}
$$
where $c\in \dd C$, and $a,b \in \dd N$, $b\geq 0$. Therefore we must have
$Q=nP_0$ for some $n\geq 1$ or $Q$ is a constant polynomial. This proves the first
statement.

\medskip

For the second statement, let $F\in \dd A_{P_0}^{\times}$ be invertible.
Then $1/F$ belongs to the ring, so it is holomorphic. Thus $F$ has no zeros.
Moreover $F$ cannot be of the form $e^{nP_0}$ for $n\geq 0$ since
$$
e^{-nP_0(z)} \to 0
$$
when $z\to +\infty . a_d^{-1/d}$ and we know that for any element
$G$ in the ring $\dd A_{P_0}$
$$
G(z)\to +\infty
$$
when $z\to +\infty . a_d^{-1/d}$. $\diamond$

\stit {III.1.7) Picard-Vessiot extensions.}

We recall that a Picard-Vessiot extension of a differential ring $\dd A$
is a differential ring extension $\dd A<u_1, \ldots , u_n>$ generated
by $u_1, \ldots , u_n$ fundamental solutions of a
homogeneous linear differential
equation of order $n$
$$
y^{(n)} +b_{n-1} y^{(n-1)}+\ldots +b_1 y'+b_0 y=0 \ ,
$$
where $b_j \in \dd A$ and the ring of constants of the extension coincides
with the ring of constants of $\dd A$.

We remind also that a {\it Liouville extension} is a
Picard-Vessiot extension generated by successive adjunctions of
integrals or exponential of integrals (see [Ka] chapter III.12
p.23, and [Rit2]). These have a solvable differential Galois group
([Ka] chapter III.13 p.24).

\eno {Theorem III.1.7.1}{ The field $\dd K_{P_0}=\dd C (z, F_0, \ldots, F_{d-1})$ and the ring $\dd A_{P_0}=\dd C
[z, F_0, \ldots, F_{d-1}]$ are Picard-Vessiot extensions of $\dd C[z]$, i.e. they are generated by the fundamental
solutions of a linear homogeneous differential equation with polynomial coefficients. Moreover these extensions
are Liouville extensions.}

The ring of constants are the constant functions. We only need to find
the homogeneous linear differential equation satisfied by $F_0, \ldots , F_{d-1}$.
We construct a homogeneous linear differential equation satisfied by $F'_0,
\ldots , F'_{d-1}$.

\medskip

We define a double sequence of functions $(y_{n,m})_{n\in \dd Z
\atop m \geq 0}$ by

\medskip

\item {$\bullet$} $y_{0,0}=e^{P_0}$,
\item {$\bullet$} For $n>m$, $y_{n,m}=0$,
\item {$\bullet$} For $n<0$, $y_{n,m}=0$,
\item {$\bullet$} For $n \in \dd N$, $m\geq 0$,
$$
y_{n,m+1}=y_{n-1,m}+y'_{n,m} \ .
$$
(Pascal's triangle rule with one derivative)
\medskip

The first proposition is straightforward.

\eno {Proposition III.1.7.2}{ We have

\item {$\bullet$} For $n\geq 0$, $y_{n,n}=e^{P_0}$.

\item {$\bullet$} For $m\geq 0$, $y_{0,m}=\left ( e^{P_0} \right )^{(m)}$.

\item {$\bullet$} For all $n\in \dd N$, $m\geq 0$, $y_{n,m}=Q_{n,m} e^{P_0}$
where $Q_{n,m}$ is a universal polynomial with
positive integer coefficients on $P'_0, P''_0, P_0^{(3)}, \ldots$
}

\eno {Proposition III.1.7.3}{ We define for $k\geq0$, $y_k(z)= z^k e^{P_0(z)}=z^k y_{k.k}$. Then we have

\item {$\bullet$} For $0\leq l\leq k$,
$$
y_k^{(l)}=z^k y_{0,l}+kz^{k-1} y_{1,l}+k(k-1) z^{k-2} y_{2,l}+
\ldots +{k!\over (k-l)!} z^{k-l} y_{l,l} \ .
$$

\item {$\bullet$} For $k\leq l$,
$$
y_k^{(l)}=z^k y_{0,l}+kz^{k-1} y_{1,l}+k(k-1) z^{k-2} y_{2,l}+
\ldots +{k! \over 1} z y_{k-1,l}+k! y_{k,l} \ .
$$
}

\stit {Proof.}

It results from a direct induction on $l$ observing
that $y'_{0,l}=y_{0,l+1}$,
$y_{0,l}+y'_{1,l}=y_{1,l+1}$, etc.$\diamond$

\medskip

\stit {Proof of the Theorem.}

We look for polynomials $b_0, b_1, \ldots , b_{d-1}$ such that $y_0, y_1, \ldots , y_{d-1}$ are solutions of $$
y^{(d)}+b_{d-1} y^{(d-1)}+\ldots b_1 y' +b_0 y =0 \ . $$ They will form a fundamental set of solutions since these
functions are $\dd C$-linearly independent. Once we find these polynomial coefficients, the special functions $1,
F_0, F_1, \ldots , F_{d-1}$ will form a fundamental set of solutions of $$ y^{(d+1)}+b_{d-1} y^{(d)}+\ldots b_1
y'' +b_0 y' =0 \ . $$ We can plug $y_k$ into the differential equation and compute $y_k^{(l)}$ using the
proposition. Then grouping together the factors of $z^j$, $j=0, \ldots , d-1$, we get a triangular system $$ b_j
y_{j,j} +b_{j+1} y_{j,j+1} +\ldots +b_{d-1} y_{j,d-1} +y_{j,d} =0 \ . $$ Thus, since $y_{j,j}=e^{P_0}$, we get $$
b_j=-b_{j+1} y_{j,j+1} e^{-P_0}-\ldots -b_{d-1} y_{j,d-1} e^{-P_0} -y_{j,d} e^{-P_0} \ , $$ and the result follows
using Proposition III.1.7.3.

\medskip

Note that the extension is a Liouville extension as
announced since each
$F_0$ is the exponential of an integral followed by an integral, and
for $j\geq 1$ the special function $F_j$ is an integral over the field
generated by $e^{P_0}$.
$\diamond$

\stit {Remark.}

The Wronskian of $F_0, F_1, \ldots, F_{d-1}$ satisfies the differential equation $$ W'-d P'_0 \ W=0 \ , $$ and is
equal to $W(z)=e^{dP_0(z)}$.

\medskip

\stit {Examples.}

{\bf 1.} For $d=1$, the equation is
$$
y'-P_0' y=0 \ .
$$

\medskip

{\bf 2.} For $d=2$, the equation is
$$
y''-2P'_0\ y'+\left [ \left (P'_0\right )^2 -P''_0\right ]\ y =0 \ .
$$
In particular, for $P_0(z)=z^2$,
$$
y''-4z\ y'+(4z^2-2) \ y =0 \ .
$$

\stit {III.1.8) Liouville classification.}

Between 1830 and 1840 J. Liouville developed a classification of transcendental functions generated by algebraic
expressions, logarithms and exponentials, and proved the non-elementary character of some natural integrals and
solutions of some differential equations. Later he noticed that his classification can be extended by allowing
integrations instead of using the logarithm function, which constitutes a particular case since any expression
$\log f$ is the primitive of $f'/f$.

\medskip

We recall Liouville's classification. Functions of order $0$ are
algebraic functions of the variable $z$, that is those functions
satisfying a polynomial equation with polynomial coefficients on
$z$. Assume by induction that order $n$ functions have been
defined. Functions of order $n+1$ are those functions that are not
of order $n$ and that can be obtained by taking an exponential or
a primitive of order $n$ functions of that satisfy an algebraic
equation with such coefficients.

\medskip

We refer to J.F. Ritt's book on elementary integration [Rit1] for more information on this subject, the precursor
of modern differential algebra.

\medskip

 Note that Liouville classification only concerns functions that
are multivalued in the complex plane, i.e. except for isolated
singularities and ramifications they can be continued
holomorphically through all the complex plane (these are called
"fluent" functions in Ritt's terminology [Rit1]).

From this classification we have:

\eno {Proposition III.1.8.1}{ Entire functions in the ring $A_{P_0}$ are functions of order at most $2$. Moreover,
if $d\geq 2$, we have that $F_0$ is of order $2$.}

For the proof of the non-elementarity of the integral giving $F_0$
see [Rit1] p.48.

%% file: iii2.tex
\input defv7.tex
\stit {III.2) Refined analytic estimates.}

\stit {III.2.1) Decomposition of the end at infinite.}

We consider the simply connected log-Riemann surface of finite type
$\cl S$.
We study the geometry of the infinite end of $\cl S$.

Note that removing $\pi^{-1} (\bar B(0,R_0))$ from $\cl S$, where $\bar B(0,R_0)$ is a closed  ball of large
radius $R_0\geq 1$, large enough so that $\pi^{-1} (B(0,R_0))$ contains all infinite ramification points, leaves
$d$ connected components $U_1(R_0), \ldots , U_d(R_0)$. Each $U_j(R_0)$ is a family of pasted planes that can be
embedded isometrically inside the log-Riemann surface of the logarithm, thus a $\log$ function, still denoted
$\log$, is well defined in each of these connected components. These function will be used in the first condition
of the Liouville theorem in section III.3.

Any ball centered at
an infinite ramification point $w^*$ of small radius, small enough not
to contain any other ramification point, can be
isometrically embedded inside the log-Riemann surface
of the logarithm. Again in such neighborhood
a $\log$ function branched at $w^*$, denoted by $\log_{w^*}$,
is well defined in the charts by $\log_{w^*} (w)=\log (\pi(w)-\pi(w^*))$.
We use these function in the
second condition of Liouville theorem.

\medskip

The next lemma describes a decomposition of a neighborhood of the infinite
end.

\eno {Lemma III.2.1.1}{ There are a finite number of sheets composing $\cl S$ intersecting more than one component
of $U_j(R_0)$. Thus there exists $M_0>0$ so that each sheet intersecting $\{w\in U_j(R_0); |\Arg w-\Im \log w| >
M_0\}$ does not intersect any other $U_i\not= U_j$. Each one of these plane sheets contains in its closure exactly
one infinite ramification point $w^*$. This ramification point $w^*_{j+}$ is the same for all these sheets with
argument larger than $M_0$ (resp.  $w^*_{j-}$ for those with argument less than $-M_0$). We denote by $\hat
U_{j+}$ (resp. $\hat U_{j-}$) the union of these sheets in $\cl S$. On $\hat U_{j\pm}$ the logarithm funcion
$\log_{w^*_{j\pm}}$ does extend holomorphically to a function denoted by $\log_{j\pm}$ (since $\hat U_{j\pm}$ can
be fully isometrically embedded in the surface of the logarithm branched at $w^*_{j\pm}$). We can also define an
argument function $\Arg_{j\pm} =\Im \log_{j\pm}$. These definitions do not depend on $R_0$, large enough so that
$\pi^{-1} (B(0,R_0))$ contains all ramification points.

Now we define for $0<r<R$,
$$
V_{j\pm} (M_0, r, R)=
\{w\in \hat U_{j\pm} ; \Arg_{j\pm} w > M_0 \ , \
\ r < |\pi(w)-\pi(w^*_{j\pm})| < R \} \ .
$$

Then the complement in $\cl S$ of the set
$$
 \bigcup_j \left (V_{j+}(M_0, r_0/2, 2R_0) \cup V_{j-}(M_0, r_0/2, 2R_0)
\right )
\cup \bigcup_{w^*}
B(w^*, r_0) \cup \left (\cl S -\pi^{-1} (B(0, R_0)) \right )
$$
is a compact set of $\cl S$.
}

\stit {III.2.2) Analytic estimates.}

We consider the simply connected log-Riemann surface of finite type
$\cl S$ whose uniformization is given by the integral $F_0$.

\eno {Definition III.2.2.1}{ Let $k_0:\cl S \to \dd C$ be the inverse of the uniformization of $\cl S$ given by
$F_0$. We define the functions on $\cl S$ $$ \eqalign{ f_0&=F_0\circ k_0 =\pi \cr f_1&= F_1\circ k_0 \cr &\vdots
\cr f_{d-1} &=F_{d-1} \circ k_0 \cr } $$ }

\eno{Definition III.2.2.2}{ We define the $\dd C$-vector space $\dd {\cl V}_{\cl S}$ of holomorphic functions
$f:\cl S\to \dd C$ of the form $$ f=F\circ k_0 \ , $$ where $F\in \dd V_{P_0}$.}

Note that $f_0, \ldots , f_{d-1}
\in \dd {\cl V}_{\cl S}$ and
$$
\dd {\cl V}_{\cl S} =k_0 \dd C[k_0 ] \ \left (e^{P_0}\circ k_0 \right )
\oplus \dd C . 1 \oplus
\dd C . f_0 \oplus \ldots \oplus \dd C . f_{d-1} \ .
$$

\eno {Proposition III.2.2.3}{ Any function $f \in \cl V_{\cl S}$ not belonging to the subspace $\dd C .1\oplus \dd
C . f_0$ has a Stolz continuous extension to $\cl S^*$ but not a continuous extension. In particular, the
functions $f_1, \ldots , f_{d-1}$ do extend Stolz continuously to $\cl S^*$ but not continuously. The function
$f_0$ also extends continuously to $\cl S^*$ for the metric topology. }

This proposition will result from the refined estimates that we
prove in what follows.

\eno {Theorem III.2.2.4}{ For $f\in \dd {\cl V}_{\cl S}$ there exist $\kappa =\kappa (f) \geq 0$ such that

{\parindent=2cm
\item {(i)} There exists $R_0 \geq 1$, such that for $w\in \cl S$,
$|\pi (w)| > R_0$,
$$
|f(w)|\leq C_0 |\pi (w) | \left | \log w\right |^{\kappa}  \  .
$$

\item {(ii)}  There exists $r_0>0$, such that
if $w^*\in \cl S^*-\cl S$ is an infinite ramification
point and $w\in B(w^*, r_0)$,
$$
|f(w)-f(w^*)| \leq C_0 |\pi (w)-\pi (w^*)|
\left | \log_{w^*} (w) \right |^{\kappa} \ .
$$

\item {(iii)} For $w\in V_{j\pm}(M_0, r_0, R_0)$ we
have
$$
|f(w)| \leq C_0 \left | \log_{j\pm} (w)\right |^{\kappa} \ .
$$
}

If $f=F\circ K_0$ and $F=\int Q e^{P_0}$ then $\kappa (f)=\deg Q /d$.
In particular
$$
\kappa (f_k)=k/d \ .
$$
}

We prove first a refinement of the asymptotic estimate in section
III.1.2.

\eno {Lemma III.2.2.5}{ Let $F\in \dd V_{P_0}$ with $$ F(z)=\int_0^z Q(t) e^{P_0(t)} \ dt \ , $$ and $Q\in \dd
C[z]$. Let $w\in \cl S$ with $\pi(w) \to \infty$. Then $z=k_0(w) \to \infty$ and we have $$ F(z)\sim {Q(z) \over
P_0'(z)} e^{P_0(z)} \ . $$

In particular putting $F=F_0$,
when $\pi (w) \to \infty$ then if $w=F_0(z)$
$$
|z|\sim |a_d|^{-1/d} |\log w|^{1/d} \ .
$$
}

\stit {Proof.}

It is clear that when $\pi (w)\to \infty $ then $z\to \infty$.
The same integration by parts as in 7.1.2 gives
$$
F(z)={Q(z) \over P'_0(z)} e^{P_0(z)} \left ( 1+\cl O (|z|^{-d})\right )
$$
and the result follows.

$\diamond$

\stit {Proof of (i).}

The lemma proves estimate (i) of the theorem. If $f=F\circ k_0$ then
when $\pi(w)\to \infty$,
$$
F(z)\sim {Q(z)\over P'_0(z)} e^{P_0(z)}\sim Q(z) F_0(z) \ .
$$
Since $w=F_0(z)$, if $\deg Q=k$, this gives that there exists $R_0$
such that for $|\pi(w)|>R_0$
$$
|F(z)|\leq C |\pi(w)||z|^k \leq C_0 |\pi(w)| |\log w|^{k/d} \ ,
$$
for some positive constants $C, C_0 >0$ where we used $|z|\sim
|a_d|^{-1/d} |\log w |^{1/d}$.
$\diamond$

We have similar asymptotics when $w\in \cl S$ approaches an infinite
ramification point.

\eno {Lemma III.2.2.6}{ Let $F\in \dd V_{P_0}$ with $$ F(z)=\int_0^z Q(t) e^{P_0(t)} \ dt \ , $$ and $Q\in \dd
C[z]$. Let $w^*\in \cl S^*-\cl S$ and $w\in \cl S$ with $w \to w^*$. Then $z=k_0(w) \to \infty$ in a sector
centered around the direction $ +\infty . a_d^{-1/d}$ for the appropriate $d$-th root, and we have $$
F(z)-\lim_{\xi\to +\infty . a_d^{-1/d}} F(\xi ) \sim {Q(z) \over P'_0(z)} e^{P_0(z)} \ . $$

In particular putting $F=F_0$, if $w=F_0(z)$
$$
F_0(z)-\lim_{\xi \to +\infty . a_d^{-1/d}} F_0(\xi)=\pi(w)-\pi(w^*) \ ,
$$
and when $w \to w^*$ then
$$
|z|\sim |a_d|^{-1/d} |\log_{w^*} w|^{1/d} \ .
$$
}

\stit {Proof.}

The argument follows the same lines as before. If $w\to w^*$ then
if $w=F_0(z)$ we have that $z\to \infty$ as described. We can
write integrating by parts
$$
\eqalign {
F(z)-\lim_{\xi\to +\infty . a_d^{-1/d}} F(\xi )
&=\int_{+\infty . a_d^{-1/d}}^z Q(t) e^{P_0(t)} \ dt \cr
&={Q(z)\over P'_0(z)} e^{P_0(z)} -\int_{+\infty . a_d^{-1/d}}^z
{Q'(t)P'_0(t) -Q(t) P''_0(t) \over (P'_0(t))^2} e^{P_0(t)}\ dt  \ .\cr
}
$$
Note that if $R=(Q'P'_0-QP''_0)/(P'_0)^2$, when $z\to \infty$
$$
R(z)=\cl O (|z|^{-(d-k)})  \ .
$$
Applying this to $F=F_0$, $Q=1$, $k=0$,
we get
$$
\pi(w)-\pi(w^*)= {1\over P'_0(z)} e^{P_0(z)} -\int_{+\infty . a_d^{-1/d}}^z
R(t) e^{P_0(t)}\ dt \ .
$$
And making the change of variables $u=F_0(t)$ in the integral, we get
$$
\pi(w)-\pi(w^*)= {1\over P'_0(z)} e^{P_0(z)} -\int_{w^*}^w R\circ k_0(u)\ du \ .
$$
There exists $r>0$ small and $C>0$ such that for $|\pi(w)-\pi(w^*)|<r$ we have
$$
\left |\int_{w^*}^w R\circ k_0(u)\ du \right | \leq  C |\pi(w)-\pi(w^*)| |z|^{-d} \ .
$$
Thus this integral is neglectable in the previous formula and we get that
$$
\pi(w)-\pi(w^*) \sim {1\over P'_0(z)} e^{P_0(z)}
$$
when $w\to w^*$. This proves that
$$
|z|\sim |a_d|^{-1/d} |\log_{w^*} w|^{1/d} \ .
$$
and the lemma for $F=F_0$.

Now in the general case we proceed by induction on the degree of $Q$.
For $k=\deg Q \leq d$ we get in the same way
$$
\left |\int_{w^*}^w R\circ k_0(u)\ du \right |
\leq  C |\pi(w)-\pi(w^*)| |z|^{-d+k} \ .
$$
Moreover when $w\to w^*$
$$
{Q(z)\over P'_0(z)} e^{P_0(z)} \sim Q(z) |\pi(w)-\pi(w^*)|
=\cl O (|\pi(w)-\pi(w^*)| |z|^k)
$$
thus  again in the integration by parts formula, the integral is neglectable
with respect to $Q/P'_0 e^{P_0}$, and
$$
F(z)-\lim_{z\to +\infty . a_d^{-1/d}} F(z)
\sim {Q(z) \over P'_0(z)} e^{P_0(z)} \ .
$$
If $k=\deg Q \geq d$ we perform the euclidean division of $Q$ by $P'_0$
$$
Q=AP'_0+B,
$$
with $\deg B \leq d-1$. Then by integration by parts, using the
induction hypothesis and observing that $\deg A' < \deg Q$,
$$
\eqalign {
\int_{+\infty . a_d^{-1/d}}^z Q(t) e^{P_0(t)} \ dt &=
\int_{+\infty . a_d^{-1/d}}^z B(t) e^{P_0(t)} \ dt + A(z)e^{P_0(z)}
-\int_{+\infty . a_d^{-1/d}}^z A'(t) e^{P_0(t)} \ dt \cr
&\sim {B(z)\over P'_0(z)} e^{P_0(z)} + {A(z) P'_0(z) \over P'_0(z)}
e^{P_0(z)}-{A'(z)  \over P'_0(z)}
e^{P_0(z)} \cr
&\sim {B(z)+A(z) P'_0(z) \over P'_0(z)} e^{P_0(z)} \cr
&\sim {Q(z)\over P'_0(z)} e^{P_0(z)} \cr
}
$$

$\diamond$

\stit {Proof of (ii).}

Writing $f=F\circ F_0$, there exists $r-0 >$ such that
for $w\in B(w^*, r_0)$,
$$
\eqalign {
|f(w)-f(w^*)| &\leq C |z|^k |\pi(w)-\pi(w^*)|\cr
&\leq C_0 |\pi(w)-\pi(w^*)| | \log (\pi(w)-\pi(w^*))|^{k/d} \cr
&\leq C_0 |\pi(w)-\pi(w^*)| | \log_{w^*} w|^{k/d} \ .\cr
}
$$
$\diamond$

The last estimate (iii) depends on the following proposition.

\eno{Proposition III.2.2.7}{ Let $w \in V_{j\pm}(M_0,r,R)$ and $w \to \infty$. Then $$ \log_{j\pm}(w) \sim P_0(z)
$$ }

\stit{Proof.}

\medskip

We recall the discussion in section II.5.2, and the existence of $2d$ families of domains $(C_{1,l}')_{l \geq
l_0}, \, (C_{1,l}')_{l \leq -l_0}, \dots, \, (C_{d,l}')_{l \geq l_0}, \, (C_{d,l}')_{l \leq -l_0}$  which
correspond under $F_0$ to families of planes in $\cl S$, slit and pasted around the ramification points, with two
families for each ramification point. Letting $$ A_ {j +} = \bigcup_{2l-1 > k_0} (C_ {j,l}' \cup \gamma_{j,2l+1}')
\ , \ A_ {j-} = \bigcup_{2l-1 < -k_0} (C_ {j,l}' \cup \gamma_{j,2l-1}') \ , $$ and choosing $M_0$ appropriately,
each $\hat{U}_{j+}$ is the image under $\tilde{F_0}$ of some $A_{j'+}$ and each $\hat{U}_{j-}$ the image of some
$A_{j'-}$ (with the $j'$ corresponding to $j$ not necessarily the same for $\hat{U}_{j+}$ and $\hat{U}_{j-}$).

\medskip

Now suppose $w = \tilde{F_0}(z)$ and $w \to \infty$ in $\cl{S}$ through some $V_{j\pm}(M_0,r,R)$, say through
$V_{j+}(M_0,r,R)$. Then $z \to \infty$ in $\dd{C}$ through $A_{j'+}$. Let $l = l(z)$ be such that $z \in C_
{j',l}' \cup \gamma_{j',2l+1}'$ (so $w$ lies in sheet $l$ of $\hat{U}_{j+}$). Since $w$ does not converge to
$w_{j+}$ or go to infinity in a sheet, $l$ must go to infinity as $w \to \infty$; we observe that $\Arg_{j+}(w)
\sim 2\pi l$, and, since Re $\log_{j+}(w) = \log|w-w_{j+}|$ is bounded on $V_{j+}(M_0,r,R)$, that therefore $$
\log_{j+}(w) \sim i \Arg_{j+}(w) \sim 2\pi i l. $$ Since $z \in C_ {j',l}' \cup \gamma_{j',2l+1}'$, $z$ lies in
the region between the curves $\Gamma_{j',2l-2}'$ and $\Gamma_{j',2l+1}'$, so $(2l-2)\pi \leq \hbox{ Im }P_0(z)
\leq (2l+1)\pi$. Also, by the propositions at the end of section II.5.2, $\arg z$ must converge to
$\arg((-a_d)^{-1/d}) + \pi/2d$ for the $d$-th root of $(-a_d)$ corresponding to $j'$ (otherwise $w$ would leave
$V_{j+}(M_0,r,R)$ eventually). Therefore $$\eqalign{ \arg P_0(z) \to (4j'-3)\pi/2 \ & \Rightarrow \hbox{ Re
}P_0(z) = o(\hbox{ Im }P_0(z)) \cr
                               & \Rightarrow P_0(z) \sim i \hbox{ Im }P_0(z) \sim 2\pi i l \sim \log_{j+}(w)
\cr
}$$
A similar argument works for the case when $w \to \infty$ through $V_{j-}(M_0,r,R)$. $\ \diamondsuit$

\medskip

\eno{Corollary III.2.2.8} { Let $w \in V_{j\pm}(M_0,r,R)$ and $w \to \infty$. Then for $F\in \dd V_{P_0}$, $F=\int
Qe^{P_0}$ with $k=\deg Q$, $$ |F(z)| = O(|\log_{j\pm}(w)|^{k/d}) $$ }

\medskip

\stit{Proof.}

An integration by parts gives
$$
F(z) = \left [Q(u) F_0(u) \right ]_0^z - \int_0^z Q'(u) F_0(u) du \ ;
$$
since $F_0(z) = w'$ is bounded for
$w \in V_{j\pm}(M_0,r,R)$, both terms on the right hand side are of the order of
$z^k$, and by the previous proposition
$$
z^d = O(\log_{j\pm}(w))
$$
so
$$
F_k(z) = O(z^k) = O(|\log_{j\pm}(w)|^{k/d}). \qquad \diamondsuit
$$

\stit {Proof of proposition III.2.2.3.}

Consider $f\in \cl V_{\cl S}$ as in the proposition. Write
$f=F\circ k_0$ with $F=\int Qe^{P_0}$ with $Q$ non-constant
polynomial. Then when $w\to w^*$
$$
f(w)-f(w^*)\sim Q(z) |\pi(w) -\pi(w^*)| \ .
$$
Consider a spiraling path $\g$ in $\cl S$ converging to
$w^*$ very slowly so that when $w\in \g$, $w\to w^*$, if
$w=F_0(z)$
$$
|Q(z)||\pi(w)-\pi(w^*)| \to +\infty \ .
$$
Note that this is always possible since the path $k_0(\g)$ in
$\dd C$ tends  to infinite thus $|Q(z)|\to +\infty$.
Then along this path $|f(w)-f(w^*)|\to +\infty$ proving that
$f$ does not have a continuous extension to $\cl S^*$. On the other
hand, in any Stolz angle when $w\to w^*$,
$$
\log_{w^*} w =\log (\pi(w)-\pi(w^*)) =\cl O (\log |\pi(w)-\pi(w^*)|) \ ,
$$
and the estimate (ii) in
the theorem shows that $f(w)\to f(w^*)$, thus
$f$ has a Stolz continuous extension to $\cl S^*$.$\diamond$

%% file: iii3.tex
\input defv7.tex
\stit {III.3) Liouville theorem on log-Riemann surfaces.}

Our goal in this section is to prove a Liouville theorem.
By Liouville theorem we mean the classical characterization
of polynomials as the only entire functions with
polynomial growth at infinite.
We would like to extend this result to simply connected
log-Riemann surfaces $\cl S$ more general than the
complex plane.

The weaker result that any bounded function is constant
does extend directly to a simply connected
log-Riemann surface $\cl S$ of finite log-degree since
such $\cl S$ is bi-holomorphic to $\dd C$.

We seek growth conditions
at infinite that identify the functions in a simple class.
It is important to note that in a log-Riemann
surface the infinite locus can be reached in different ways.
As observed before we can go to infinite in $\cl S$, that is, leave
any given compact set, by having the $\pi$-projection diverging
to $\infty$,
converging to an infinite ramification point, or "spiraling around"
without approaching ramification points or $\infty$ on the sheets.

\medskip

We prove the converse of the Theorem III.2.2.4 in the previous section. The growth conditions in Theorem III.2.2.4
do characterize the functions in the vector space $\dd {\cl V}_{\cl S}$.

\bigskip

\eno {Theorem III.3.1 (Liouville theorem).}{ Let $\cl S$ be a simply connected log-Riemann surface of log-degree
$1\leq d<+\infty$ with uniformization given by the primitive of $e^{P_0}$ with $P_0 \in \dd C[z]$ of degree $d$.
Let $f:\cl S\to \dd C$ be a holomorphic function which has a Stolz continuous extension to $\cl S^*$. We assume
that there are constants $C_0>0$ and $\kappa \geq 0$ such that

\item {(i)} There exists $R_0 \geq 1$, such that for $w\in \cl S$,
$|\pi (w)| > R_0$,
$$
|f(w)|\leq C_0 |\pi(w)| \left | \log w\right |^\kappa  \  .
$$

\item {(ii)}  There exists $r_0>0$,
if $w^*\in \cl S^*-\cl S$ is an infinite ramification
point and $w\in B(w^*, r_0)$,
$$
|f(w)-f(w^*)| \leq C_0 |\pi (w)-\pi (w^*)|
\left | \log_{w^*} (w) \right |^\kappa \ .
$$

\item {(iii)} For $w\in V_{j\pm}(M_0, r_0/4, 4R_0)$ we
have
$$
|f(w)| \leq C_0 \left | \log_{j\pm} (w)\right |^\kappa \ .
$$

\medskip

Then $f\in \cl V_{\cl S}$, that is
there exists $F\in \dd V_{P_0}$ such that $f=F\circ k_0$.
More precisely, there exists
a polynomial $Q\in \dd C[z]$ with $\deg Q \leq \kappa d$ such that
$$
F(z)=F(0)+\int_0^z Q(t)\ e^{P_0(t)} \ dt \ .
$$

In conclusion, we can identify the $\dd C$-vector space of
functions in $\cl S$ with growth conditions (i), (ii) and (iii)
with $\dd V_{P_0}(\cl S)$.
}

Observe that for a holomorphic function defined on $\cl S$, the derivative
$$
f'(w)=\lim_{w'\to w} {f(w')-f(w) \over \pi (w')-\pi (w)}
$$
is a well defined holomorphic function on $\cl S$.

\medskip

\eno {Proposition III.3.2}{ Let $f$ be as in the theorem. There exists a constant $C_2 >0$ such that we have the
following estimates for the derivative:

\item {(i)} For $w\in \cl S$, $|\pi (w)| >2R_0$,
$$
|f'(w)| \leq C_2 |\log w|^\kappa \ .
$$

\item {(ii)} If $w^*\in \cl S^*-\cl S$ and $w\in B(w^*, r_0/2)$, then
$$
|f'(w)| \leq C_2 |\log_{w^*} (w)|^\kappa \ .
$$

\item {(iii)} If $w\in V_{j\pm}(M_0, r_0/2, 2R_0)$ then
$$
|f'(w)| \leq C_2 |\log_{j\pm} (w)|^\kappa \ .
$$
}

\stit {Proof.} These estimates on the derivative $f'$ result from the combination of the estimates on $f$ given as
hypotheses in the theorem and Cauchy integral formula.

\medskip

{\bf (i)} For $w\in \cl S$, $|\pi (w)| >2R_0$, the set $\Gamma =
\{ \ \xi \in \cl S : d(w,\xi) = {1 \over 2}|\pi(w)| \}$ is a
Euclidean circle contained in a single log-chart, and in the
region $\{ |\pi (w)| > R_0| \}$ so we can estimate $f'$ using
Cauchy integral formula (we write $w, |w|$ instead of $|\pi(w)|$,
etc; and the $C$'s are positive constant):
$$\eqalign{
|f'(w)| & = \left| {1 \over 2\pi i} \int_{\Gamma} {f(\xi) \over
(\xi - w)^2} \ d\xi \right| \cr
        & \leq {1 \over 2\pi} \int_{\Gamma} {|f(\xi)| \over |\xi - w|^2} \ |d\xi| \cr
        & \leq {1 \over 2\pi} {4 \over |w|^2} \int_{\Gamma} |f(\xi)| \ |d\xi| \cr
        & \leq {1 \over 2\pi} {4 \over |w|^2} \int_{\Gamma} C_0 |\xi| \left | \log \xi\right |^\kappa \
|d\xi| \cr
        & \leq {4C_0 \over 2\pi} {1 \over |w|^2} {3|w| \over 2} \int_{\Gamma} \left | \log(w) + \log \left(1
+ {\xi - w \over w} \right)  \right |^\kappa \ |d\xi| \cr
        & \leq {6C_0 \over 2\pi} {1 \over |w|} \int_{\Gamma} \left( | \log(w) | + \left|\log \left(1 + {\xi -
w \over w} \right) \right|  \right )^\kappa \ |d\xi| \cr
        & \leq {6C_0 \over 2\pi} {1 \over |w|} \int_{\Gamma} \left( | \log(w) | + C \right )^\kappa \ |d\xi|
\ \ \cr
        & \leq {6C_0 \over 2\pi} {1 \over |w|} \, C | \log(w) |^\kappa \, 2\pi {1 \over 2}|w| \ \  \cr
        & = C | \log(w) |^\kappa \ \cr
}$$

\medskip

{\bf (ii)} For the case $w^*\in \cl S^*-\cl S$ and $w\in B(w^*,
r_0/2)$, the set $\Gamma = \{ \ \xi \in \cl S : d(w,\xi) = {1
\over 2}|\pi(w) - \pi(w^*)| \}$ is a Euclidean circle contained in
a single log-chart and in $B(w^*, r_0)$ so we can follow the same
steps as above:
$$\eqalign{
|f'(w)| & = \left| {1 \over 2\pi i} \int_{\Gamma} {f(\xi) - f(w^*)
\over (\xi - w)^2} \ d\xi \right| \cr
        & \leq {1 \over 2\pi} \int_{\Gamma} {|f(\xi) - f(w^*)| \over |\xi - w|^2} \ |d\xi| \cr
        & \leq {1 \over 2\pi} {4 \over |w - w^*|^2} \int_{\Gamma} |f(\xi) - f(w^*)| \ |d\xi| \cr
        & \leq {1 \over 2\pi} {4 \over |w - w^*|^2} \int_{\Gamma} C_0 |\xi - w^*| \left | \log(\xi - w^*)
\right |^\kappa \ |d\xi| \cr
        & \leq {4C_0 \over 2\pi} {1 \over |w - w^*|^2} {3|w - w^*| \over 2} \int_{\Gamma} \left | \log(w -
w^*) + \log \left(1 + {\xi - w \over w - w^*} \right)  \right
|^\kappa \ |d\xi| \cr
        & \leq {6C_0 \over 2\pi} {1 \over |w - w^*|} \int_{\Gamma} \left( | \log(w - w^*) | + C \right
)^\kappa \ |d\xi| \ \  \cr
        & \leq {6C_0 \over 2\pi} {1 \over |w - w^*|} \, C | \log(w - w^*) |^\kappa \, 2\pi {1 \over 2}|w -
w^*| \ \ \cr
        & = C | \log(w - w^*) |^\kappa \  \cr
}$$

\medskip

{\bf (iii)} Finally, when $w\in V_{j \pm}(M_0, r_0/2, 2R_0)$  if
we take $\Gamma = \{ \ \xi \in \cl S : d(w,\xi) = r \}$ where $r$
is a constant smaller than $r_0/4$ and $2R_0$, then $\Gamma$ is a
Euclidean circle contained in a log-chart and in the region $V_{j
\pm}(M_0, r_0/4, 4R_0)$, so we can estimate as before (noting that
$\log_{j\pm}$ has a uniformly bounded derivative in $V_{j
\pm}(M_0, r_0/4, 4R_0)$ ):

$$\eqalign{
|f'(w)| & = \left| {1 \over 2\pi i} \int_{\Gamma} {f(\xi) \over
(\xi - w)^2} \ d\xi \right| \cr
        & \leq {1 \over 2\pi} \int_{\Gamma} {|f(\xi)| \over |\xi - w|^2} \ |d\xi| \cr
        & \leq {1 \over 2\pi} {1 \over r^2} \int_{\Gamma} C_0 \left | \log_{j\pm} (\xi)\right |^\kappa \
|d\xi| \cr
        & \leq {C_0 \over 2\pi} {1 \over r^2} \int_{\Gamma} \left(\left | \log_{j\pm} (w) \right | + C
\right)^\kappa \ |d\xi| \ \ \cr
        & \leq C \left | \log_{j\pm} (w)\right |^\kappa \ \ \cr
}$$

$\diamond$

\bigskip

\eno {Proposition III.3.3}{ Let $F=f\circ F_0$. For $z\in \dd C$ and $z\to \infty$ we have $$ F'(z)e^{-P_0(z)}
=\cl O (|z|^{\kappa d}) \ . $$ }

We recall the estimates on the asymptotics of $F_0$ that we are going to use and have been established in the
previous section.

\bigskip

\eno {Lemma III.3.4}{

\item {$\bullet$} Let $w^*\in \cl S^*-\cl S$ and $w\in \cl S$, $w=F_0(z)$.
When $w\to w^*$ we have
$$
z\sim a_d^{-1/d} \left (\log_{w^*} (w)\right )^{1/d} \ ,
$$
for the same $d$-th root of $a_d$ that corresponds to $w^*$ by
$$
w^*=\int_0^{-\infty . a_d^{-1/d}} e^{P_0(t)} \ dt \ .
$$
\medskip

\item {$\bullet$} Let $w\in U_j (R)$, $w\to \infty$, then if $w=F_0(z)$,
$$
z\sim a_d^{-1/d} \left ( \log w \right )^{1/d} \ ,
$$
for the same $d$-th root of $a_d$ that corresponds to $F_0(z) \to U_j$
when $z\to +\infty . a_d^{-1/d}$.

\medskip

\item {$\bullet$} Let $w\in V_{j\pm} (M_0, r, R)$ and $w\to \infty$, then
$$
\log_{j\pm} (w)= \cl O (|z|^d) \ .
$$
}

\stit {Proof of the Proposition.}

Note that
$$
f'(w)=F'\circ k_0 (w) k'_0(w)=F'(z) {1\over F_0\circ k_0 (w)}
=F'(z) {1\over e^{P_0(z)}} \ .
$$
Thus
$$
|F'(z) e^{-P_0(z)}|=|f'(w)|
$$
and using the previous estimates for $|f'(w)|$ in each
region and the asymptotic relation between $z$ and $w$
in each region given by the lemma, the result follows.
$\diamond$

\stit {Proof of the Theorem.}

Let $Q(z)=F'(z)e^{-P_0(z)}$. The function $Q$ is an entire function
and
$$
Q(z) =\cl O (|z|^{\kappa d}) \ ,
$$
thus $Q$ is a polynomial with $\deg Q \leq \kappa d$. Moreover
$$
F(z)=F(0)+\int_0^z Q(t)e^{P_0(t)} \ dt \ ,
$$
and the result follows.$\diamond$

\bigskip
\bigskip

%% file: iii4.tex
\input defv7.tex
\input psfig.sty
\font\gothic=eufm10

\bigskip
\bigskip

\stit {III.4) The structural ring.}

\stit {III.4.1) Definition.}

\eno{Definition III.4.1.1.}{We consider the ring of entire
functions $\dd A_{P_0}$ generated by the functions of $\dd
V_{P_0}$, which is given by
$$
\dd A_{P_0} =z\dd C[z,F_0,F_1,\ldots , F_{d-1}] \oplus \dd
C[F_0,F_1,\ldots ,F_{d-1}] \ .
$$
}

\eno{Definition III.4.1.2.}{ The structural ring $\cl A_{\cl S}$
of the log-Riemann surface $\cl S$ is the ring of holomorphic
functions $f$ on $\cl S$ of the form
$$
f=F\circ k_0 \ ,
$$
where $F\in \dd A_{P_0}$.

We define the structural field $\cl K_{\cl S}$ to be the field of
fractions of $\cl A_{\cl S}$. Thus $\cl A_{\cl S} \approx \dd
A_{P_0}$.}

Observe that functions on the structural ring
do have a Stolz continuous extension
to $\cl S^*$, i.e. they have Stolz limits at
infinite ramification points.

\eno {Definition III.4.1.3.}{The coordinate ring $\dd C [\pi]$,
resp. field $\dd C(\pi )$,  is the subring of the structural ring
$\cl A_{\cl S}$, resp. subfield of the structural field $\cl
K_{\cl S}$, generated by the coordinate funcion $\pi$.}

Observe that we have
$$
\eqalign { \dd C (\pi ) & \approx \dd C (F_0) \subset \dd K_{P_0}
\ ,\cr \dd C [\pi] &\approx \dd C[F_0] \subset \dd A_{P_0} \ , \cr
}
$$
since elements $f$ of the coordinate ring are of the form
$$
f=F\circ k_0 \ ,
$$
where $F\in \dd C [F_0]$.

Observe that functions in the coordinate ring do have a continuous
extension to $\cl S^*$ for the log-euclidean topology (not just a
Stolz extension), and can be characterized by that property
according to Proposition III.2.2.3 in section III.2.2.

\bigskip

The number of infinite ramification points in the log-Riemann
surface $\cl S$ can be read algebraically as the transcendence
degree of $\cl K_{\cl S}$ over $\dd C (\pi)$.

\eno {Theorem III.4.1.4.}{The transcendence degree of $\cl K_{\cl
S}$ over $\dd C (\pi)$ is
$$
\left [\cl K_{\cl S}:\dd C (\pi)\right ]_{tr} =d \ .
$$
}

\stit {Proof.}

We have that $[\dd K_{P_0}: \dd C[F_0]]_{tr}=d$ because $z,F_0,
\ldots , F_{d-1}$ are algebraically independent.$\diamond$

\bigskip

\stit {III.4.2) Points of $\cl S^*$ as maximal ideals.}

\bigskip

Recall that to each point on $z_0\in \dd C$ we can associate a
maximal ideal ${\hbox {\gothic M}}_{z_0}$ of $\dd C[z]$, namely
the ideal of functions vanishing at $z_0$. Conversely, any maximal
ideal {\gothic M} of $\dd C[z]$ is of this form since the residual
field is $\dd C$ $$ \dd C[z]/{\hbox{\gothic M}} \approx \dd C $$
and $z$ is mapped by this quotient into some $z_0\in \dd C$, thus
${\hbox {\gothic M}}= {\hbox {\gothic M}}_{z_0}$. In that way the
points of the complex plane $\dd C$ can be reconstructed
algebraically from the ring of polynomials $\dd C[z]$, each point
corresponding to a maximal ideal. The ring is of dimension $1$ and
any prime ideal is maximal. In the same way we can reconstruct the
Riemann sphere identifying points with discrete valuation rings in
the field of fractions $\dd C (z)$.

\medskip

As in the case of the polynomial ring $\dd C[z]$ on $\dd C$, to each
point of $\cl S^*$ we can associate a maximal ring of $\cl A_{\cl
S}$

\eno {Theorem III.4.2.1.}{There is an injection of $\cl S^*$ into
the space of maximal ideals of $\cl S_{\cl S}$, $$ \eqalign { \cl
S^* &\hookrightarrow {\hbox {\rm Max}} \cl A_{\cl S} \cr w_0
&\mapsto {\hbox {\gothic M}}_{w_0} \cr } $$ where ${\hbox {\gothic
M}}_{w_0}=\{f\in \cl A_{\cl S} ; f(w_0)=0\}$.

More precisely, the ring $\cl A_{\cl S}$ separates points in $\cl
S^*$.
}

\stit {Proof.}

First observe that any ideal ${\hbox{\gothic M}}_{w_0}$ is maximal
because the kernel of
$$
\ap {\cl A_{\cl S}}{\dd C}{f}{f(w_0)}
$$
is ${\hbox{\gothic M}}_{w_0}$ and therefore
$$
\cl A_{\cl S} / {\hbox{\gothic M}}_{w_0} \approx \dd C \ ,
$$
is a field and ${\hbox{\gothic M}}_{w_0}$ is maximal.

We only need to show that two distinct points $w_1,w_2\in \cl S^*$
give two distinct ideals ${\hbox{\gothic M}}_{w_1}\not=
{\hbox{\gothic M}}_{w_2}$. This is proved in the next
lemma.$\diamond$

\eno {Lemma III.4.2.2.}{The ring $\cl A_{\cl S}$ separates the
points of $\cl S^*$.}

\stit {Proof.}

Let $w_1, w_2\in \cl S^*$ with $w_1\not=w_2$. If both points are
regular points, $w_1,w_2 \in \cl S$, let $z_1, z_2\in \dd C$ be
such that $z_i=k_0(w_i)$. Then the function $f\in \cl A_{\cl S}$,
$f=F\circ k_0$, with $$ F(z)=(z-z_1)e^{P_0(z)} \ , $$ vanishes at
$z_1$ but not at $z_2$. If one of the points, say $w_1$, is a
ramification point, then taking $f=F\circ k_0$ with $$
F(z)=e^{P_0(z)} $$ the function $f$ will vanish at $w_1$ but not
at $w_2$. The remaining case is when both points are ramification
points $w_1,w_2 \in \cl S^*-\cl S$. This case is more subtle and
is handled by the next theorem.$\diamond$

\bigskip

\eno {Theorem III.4.2.3}{Let $w_1^*$ and $w_2^*$ be two distinct
ramification points of $\cl S^*$.

We cannot have for all $j=0,1,\ldots ,d-1$,

$$
f_j(w_1^*) =f_j(w_2^*) \ .
$$
}

\stit {Proof.}

We can normalize $P_0$ such that its leading coefficient is
$-1/d$, $P_0(z)=-z^d/d+a_{d-1}z^{d-1}+\ldots +a_1 z+a_0$. Assume
by contradiction that
$$
\lim_{z\to +\infty .\omega_1} F_j(z)=f_j(w_1^*)=f_j(w_2^*)=\lim_{z\to +\infty .\omega_2} F_j(z) \ .
$$
By Theorem III.1.5.1 we have for any polynomial
$Q(z)\in \dd C[z]$,
$$
\int_0^z Q(t) e^{P_0(t)} \ dt =zA(z)e^{P_0(z)}+b_0 F_0(z)+\ldots +b_{d-1} F_{d-1}(z) \ ,
$$
where $A(z)\in \dd C[z]$ and $b_0,\ldots , b_{d-1} \in \dd C$ are
constants depending on the polynomial $Q(z)$. Therefore
$$
\int^{+\infty .\omega_1}_{0} Q(z)e^{P_0(z)} \ dz =b_0 F_0(+\infty
.\omega_1)+\ldots +b_{d-1} F_{d-1}(+\infty .\omega_1) \ ,
$$
and
$$
\int^{+\infty .\omega_2}_{0} Q(z)e^{P_0(z)} \ dz =b_0 F_0(+\infty .\omega_2)+\ldots +b_{d-1} F_{d-1}(+\infty .\omega_2)
\ .
$$
Therefore for any polynomial $Q(z)\in \dd C[z]$,
$$
\int_{+\infty .\omega_1}^{+\infty .\omega_2} Q(z)e^{P_0(z)} \ dz =0 \ .
$$
Now consider the following integral depending on the
coefficients of $P_0$ :
$$
G(u_0,u_1,\ldots , u_{d-1})=\int_{+\infty .\omega_1}^{+\infty .\omega_2} e^{-z^d/d+u_{d-1}z^{d-1}+\ldots +u_1z+u_0} \ dz \ .
$$
By uniform convergence of the integral, the function $G$ is an entire
function of $d$ complex variables defined in $\dd C^d$. We have
$$
G(a_0,a_1,\ldots ,a_{d-1})=0 \ .
$$
And also, by differentiation under the integral and using the previous property, for any
$n_0,n_1,\ldots ,n_{d-1}\geq 0$,
$$
\partial_0^{n_0}\partial_1^{n_1}\ldots \partial_{d-1}^{n_{d-1}}
G_{|(a_0,\ldots ,a_{d-1})} =\int_{+\infty .\omega_1}^{+\infty .\omega_2} z^{n_1+2n_2+\ldots +(d-1) n_{d-1}} e^{P_0(z)} \ dz =0 \ .
$$
Therefore the power series expansion of $G$ at the point
$(a_0,\ldots ,a_{d-1})$ has all coefficients equal to $0$. Thus the
entire function $G$ is identically $0$. But this contradicts the fact that
the value $G(0,\ldots ,0)$ corresponding to the cyclotomic
log-Riemann surface is non-zero, because by Theorem II.6.2.1 (with
$n=0$), we have
$$
G(0,\ldots ,0)=(\o_1-\o_2) d^{{1\over d}-1} \Gamma \left ( {1\over
d}\right ) \ .
$$
$\diamond$

\bigskip

\stit {Observation.}

The preceding argument is powerful and serves to establish much
stronger results in the next section.

\bigskip

\stit {III.4.3) The ramificant determinant.}

\bigskip

To each $f \in \cl V_\cl S$ we can associate the vector in $\dd C^d$
of its values at the infinite ramification points $(f(w^*_1),
f(w^*_2),\ldots , f(w_d^*))$.

\eno {Definition III.4.3.1.}{A normal base for the vector space
generated by $f_0,f_1,\ldots , f_{d-1}$ is a base $(g_1,\ldots ,
g_d)$ such that
$$
g_i(w_j^*) =\d_{ij} \ .
$$
}

\bigskip

\eno {Definition III.4.3.2 (Ramificant determinant)}{The determinant
$$
\Delta_{P_0}=\Delta (f_0, f_1,\ldots , f_{d-1})=\left |
\matrix{f_0(w_1^*) & f_1(w_1^*) & \ldots &f_{d-1}(w_1^*) \cr
f_0(w_2^*) &f_1(w_2^*) & \ldots &f_{d-1}(w_2^*) \cr \vdots &\vdots
&\ddots & \vdots \cr f_0(w_d^*) & f_1(w_d^*) & \ldots
&f_{d-1}(w_d^*) \cr} \right |
$$
is called the ramificant of the functions $f_0,\ldots , f_{d-1}$. }

\bigskip

\eno {Theorem III.4.3.3. (Non-vanishing of the ramificant).}{The
ramificant is never $0$, therefore there exists a normal base.}

\bigskip

Observe that this result implies Theorem III.4.2.3 from the previous
section.

\medskip

We can be more precise:

\bigskip

\eno{Theorem III.4.3.4.}{ We normalize $P_0$ to have leading
coefficient $-1/d$. For each $d\geq 0$, there exists a universal
polynomial of $d$ variables with rational coefficients
$$
\Pi_d (X_0,X_1, \ldots , X_{d-1})\in \dd Q[X_0,\ldots , X_{d-1}]
$$
with $\Pi_d (0,\ldots 0)=0$ and such that the ramificant is given by
$$
\Delta (a_0,a_1,\ldots , a_d)=\left | \matrix{f_0(w_1^*) &
f_1(w_1^*) & \ldots &f_{d-1}(w_1^*) \cr f_0(w_2^*) &f_1(w_2^*) &
\ldots &f_{d-1}(w_2^*) \cr \vdots &\vdots &\ddots & \vdots \cr
f_0(w_d^*) & f_1(w_d^*) & \ldots &f_{d-1}(w_d^*) \cr} \right | =
{(-1)^{d-1} \over \sqrt {\pi}} \left ( \pi \over d\right )^{d\over
2} V_d \ e^{\Pi_d (a_0 ,a_1 ,\ldots , a_{d-1}) } \ \ .
$$
where $V_d$ is the Vandermonde determinant of the $d$-roots of unity
$\omega_1, \ldots , \omega_d$,
$$
V_d=\left | \matrix{ 1 &\omega_1 &\omega_1^2 &\ldots &\omega_1^{d-1}
\cr 1 &\omega_2 &\omega_2^2 &\ldots &\omega_2^{d-1} \cr \vdots
&\vdots&\vdots &\ddots &\vdots \cr 1 &\omega_d &\omega_d^2 &\ldots
&\omega_d^{d-1} \cr} \right | =\prod_{i\not= j} (\omega_i-\omega_j)
\not= 0 \ .
$$
In particular this ramificant is never $0$.

Moreover the Vandermonde determinant $V_d$ can be computed
$$
V_d =(-1)^{d-1} d^d \ ,
$$
and therefore
$$
\Delta (a_0,a_1,\ldots , a_d)=\left | \matrix{f_0(w_1^*) &
f_1(w_1^*) & \ldots &f_{d-1}(w_1^*) \cr f_0(w_2^*) &f_1(w_2^*) &
\ldots &f_{d-1}(w_2^*) \cr \vdots &\vdots &\ddots & \vdots \cr
f_0(w_d^*) & f_1(w_d^*) & \ldots &f_{d-1}(w_d^*) \cr} \right | = {1
\over \sqrt {\pi}} \left ( \pi  d\right )^{d\over 2}  \ e^{\Pi_d
(a_0 ,a_1 ,\ldots , a_{d-1}) } \ \ .
$$

}

\bigskip

\stit {Proof.}

 The proof is similar to the proof of Theorem
III.4.2.3. Consider the functions of several complex variables
$(a_0,a_1,\ldots , a_{d-1})$,
$$
\Delta (a_0,a_1,\ldots ,a_{d-1})=\left | \matrix{\int_0^{+\infty .
\omega_1} e^{P_0(z)}\ dz & \int_0^{+\infty .\omega_1} z e^{P_0(z)}\
dz & \ldots &\int_0^{+\infty .\omega_1} z^{d-1} e^{P_0(z)} \ dz \cr
\int_0^{+\infty .\omega_2} e^{P_0(z)} \ dz &\int_0^{+\infty
.\omega_2} z e^{P_0(z)} \ dz & \ldots &\int_0^{+\infty .\omega_2}
z^{d-1} e^{P_0(z)} \ dz \cr \vdots &\vdots &\ddots & \vdots \cr
\int_0^{+\infty .\omega_d} e^{P_0(z)} \ dz & \int_0^{+\infty
.\omega_d} z e^{P_0(z)} & \ldots &\int_0^{+\infty .\omega_d}
z^{d-1}e^{P_0(z)} \cr} \right |
$$

Each integral is an entire function on the several complex variables
$(a_0,a_1,\ldots , a_{d-1})$, therefore the determinant is also an
entire function. Observe also that by Theorem III.1.5.1 we have that
each integral
$$
\int_0^{+\infty . \omega_i} z^n e^{P_0(z)} \ dz \ ,
$$
is a linear combination with coefficients polynomial integer
coefficients on the $(a_j)$ of the integrals for $j=0,1,\ldots
,d-1$,
$$
\int_0^{+\infty . \omega_i} z^j e^{P_0(z)} \ dz \ .
$$
Therefore, differentiating column by column, we observe that for
each $j=0,1,\ldots , d-1$,
$$
\partial_{a_j} \Delta = c_j \Delta \ ,
$$
where $c_j$ is a polynomial on the $(a_j)$ with integer
coefficients. We conclude that the logarithmic derivative of
$\Delta$ with respect to each variable is a universal polynomial
with integer coefficients on the variables $(a_j)$. This gives the
existence of the universal polynomial $\Pi_d$ such that
$$
\Delta (a_0,a_1,\ldots ,a_{d-1})=c.e^{\Pi_d(a_0,a_1,\ldots ,
a_{d-1})} \ ,
$$
with $\Pi_d (0,\ldots , 0)=0$ and  $c=\Delta (0,\ldots , 0) \in \dd
C$. It remains to prove that $c$ is not $0$. The parameter value
$(a_0,a_1,\ldots ,a_{d-1})=(0,0,\ldots ,0)$ corresponds to the case
of the cyclotomic log-Riemann surface studied in section II.6. In
this case, the computation in Theorem II.6.2.1 gives
$$
\eqalign { &\Delta (0,\ldots , 0) =\left | \matrix{d^{{1\over d}-1}
\Gamma \left ({1\over d}\right ) \omega_1& d^{{2\over d}-1} \Gamma
\left ({2\over d}\right ) \omega_1^2 & \ldots &d^{{d\over d}-1}
\Gamma \left ({d\over d}\right ) \omega_1^d \cr d^{{1\over d}-1}
\Gamma \left ({1\over d}\right ) \omega_2 &d^{{2\over d}-1} \Gamma
\left ({2\over d}\right ) \omega_2^2 & \ldots &d^{{d\over d}-1}
\Gamma \left ({d\over d}\right ) \omega_2^d \cr \vdots &\vdots
&\ddots & \vdots \cr d^{{1\over d}-1} \Gamma \left ({1\over d}\right
) \omega_d & d^{{2\over d}-1} \Gamma \left ({2\over d}\right )
\omega_d^2 & \ldots &d^{{d\over d}-1} \Gamma \left ({d\over d}\right
) \omega_d^d \cr} \right | \cr &=d^{{1\over d} (1+2+\ldots +d)
-{1\over d}} \Gamma \left ({1\over d}\right ) \Gamma \left ({2\over
d}\right )\ldots \Gamma \left ({d\over d}\right ) \left | \matrix{
\omega_1& \omega_1^2 & \ldots & \omega_1^d \cr \omega_2 &\omega_2^2
& \ldots & \omega_2^d \cr \vdots &\vdots &\ddots & \vdots \cr
\omega_d & \omega_d^2 & \ldots & \omega_d^d \cr} \right | \cr
&=d^{{1-d\over 2}} (2\pi)^{{d-1\over 2}} d^{{1\over 2} -d{1\over d}}
\Gamma (1) \left | \matrix{ \omega_1& \omega_1^2 & \ldots &
\omega_1^d \cr \omega_2 &\omega_2^2 & \ldots & \omega_2^d \cr \vdots
&\vdots &\ddots & \vdots \cr \omega_d & \omega_d^2 & \ldots &
\omega_d^d \cr} \right | \cr &={1\over \sqrt {\pi}} \left ( \pi
\over d\right )^{d\over 2} \left | \matrix{ \omega_1& \omega_1^2 &
\ldots & \omega_1^d \cr \omega_2 &\omega_2^2 & \ldots & \omega_2^d
\cr \vdots &\vdots &\ddots & \vdots \cr \omega_d & \omega_d^2 &
\ldots & \omega_d^d \cr} \right | \ ,\cr}
$$
where we have used Gauss multiplication formula
$$
\Gamma (z ).\Gamma \left ( z+{1\over d}\right )\ldots \left (
z+{d-1\over d}\right )=(2\pi)^{d-1\over 2} d^{{1\over 2}-dz}
\Gamma(dz) \ .
$$
Since $\omega_1,\omega_2,\ldots ,\omega_d$ are the $d$ roots of $1$,
we have that $\omega_j^d=1$ and the last determinant is equal to
$(-1)^{d-1} V_d$ where $V_d$ is the Vandermonde determinant
$$
V_d=\left | \matrix{ 1 &\omega_1 &\omega_1^2 &\ldots &\omega_1^{d-1}
\cr 1 &\omega_2 &\omega_2^2 &\ldots &\omega_2^{d-1} \cr \vdots
&\vdots&\vdots &\ddots &\vdots \cr 1 &\omega_d &\omega_d^2 &\ldots
&\omega_d^{d-1} \cr} \right | =\prod_{i\not= j} (\omega_i-\omega_j)
\not= 0 \ .
$$
Now the next lemma applied to the polynomial $Q(X)=X^d-1$, shows that
$$
V_d=\prod_i (d\omega_i^{d-1})=d^d \left ( \prod_i \omega_i \right
)^{d-1}=(-1)^{d-1} d^d \ .
$$

$\diamond$

\bigskip

\eno{Lemma III.4.3.5}{If $\xi_1 ,\ldots , \xi_d$ are the $d$ roots
of a monic polynomial $Q(X)$, then we can compute the following
Vandermonde determinant
$$
\left | \matrix{ 1 &\xi_1 &\xi_1^2 &\ldots &\xi_1^{d-1} \cr 1 &\xi_2
&\xi_2^2 &\ldots &\xi_2^{d-1} \cr \vdots &\vdots&\vdots &\ddots
&\vdots \cr 1 &\xi_d &\xi_d^2 &\ldots &\xi_d^{d-1} \cr} \right |
=\prod_{i\not= j} (\xi_i-\xi_j)=\prod_i Q'(\xi_i )
$$
}

\medskip

\stit {Proof.}

We have $Q'(\xi_i)=\prod_{j\not=i} (\xi_i-\xi_j)$ and the result
follows.$\diamond$

\bigskip

It is interesting to study the combinatorial properties of the
family of universal polynomials $(\Pi_d)$. We can compute the first
few polynomials $\Pi_d$ before proceeding to the general
computation.

\bigskip

\eno {Theorem III.4.3.6.}{We have
$$
\eqalign{ \Pi_1 (X_0)&=X_0 \ ,\cr \Pi_2 (X_0,X_1) &= 2X_0+{1\over 2}
X_1^2 \ ,\cr \Pi_3 (X_0,X_1,X_2) &= 3X_0 + 2 X_1 X_2 +{4\over 3}
X_2^3 \ ,\cr}
$$
for $d=4$
$$
\Pi_4(X_0,X_1,X_2,X_3)=4X_0 +3X_3 X_1+ 2
X_2^2+9X_3^2X_2+\ldots \ ,
$$
where the remaining term is a polynomial in $X_3$, and for $d\geq
5$,
$$
\Pi_d (X_0,X_1,X_2,\ldots , X_{d-1})=dX_0 +(d-1)X_{d-1} X_1+ \left (
2 (d-2) X_{d-2} + (d-1)^2 X_{d-1}^2\right ) X_2 +\ldots
$$
where the remaining terms are independent of $X_0$, $X_1$ and $X_2$.

More generally, $\Pi_d$ is of degree $1$ in $X_k$ for $k\leq d/2$.

}

\stit {Proof.}

For $d\geq 1$ the dependence of the ramificant $\Delta$ on $a_0$ is
straightforward by direct factorization of $e^{a_0}$ in the
integrals, which gives
$$
\Pi_d (X_0,\ldots ,X_{d-1}) =dX_0+\ldots
$$
with remaining terms independent of $X_0$. Also this can be seen by
differentiation column by column of $\Delta$,
$$
\partial_{a_0} \Delta =d\Delta \ ,
$$
which also gives the result.

\medskip

For the dependence on $a_1$ we use this last approach. For $d\geq
2$, we have
$$
\partial_{a_1} \Delta = (d-1) a_{d-1} \Delta \ ,
$$
because the differentiation of the first $d-1$ columns yields $0$,
and for the last column we have
$$
z^d =-zP'_0(z)+(d-1)a_{d-1}z^{d-1}+(d-2)a_{d-2}z^{d-2}+\ldots +a_1 z
\ ,
$$
and the integrals corresponding to the term $-zP'_0(z)$ contribute
$0$ because
$$
\int -z P'_0(z)e^{P_0(z)} \ dz =[-ze^{P_0}]+\int e^{P_0(z)}\ dz \ ,
$$
and by linearity of the integrals in the last column the lower order
terms $(d-2)a_{d-2}z^{d-2}+\ldots +a_1 z$ contribute $0$. Thus
the only contribution comes from the term $(d-1)a_{d-1}z^{d-1}$
which gives $(d-1)a_{d-1} \Delta $.

Now this last equation gives for $d=2$,
$$
\partial_{a_1} \Delta =  a_{1} \Delta \ ,
$$
thus $\Pi_2 (X_0,X_1)=2X_0 +{1\over 2}X_1^2 $.

For $d\geq 3$ we get
$$
\Pi_d (X_0,X_1,\ldots , X_{d-1})=dX_0+(d-1)X_{d-1}X_1+\ldots \ ,
$$
where the remaining terms are independent of $X_0$ and $X_1$.

\medskip

Now we assume $d\geq 3$ and we determine the dependence on $a_2$.

We proceed as before and differentiate column by column
$\partial_{a_2} \Delta$. Only the last two columns give a
contribution. The last but one contributes $(d-2)a_{d-2} \Delta$
because
$$
z^d =-zP'_0(z)+(d-1)a_{d-1}z^{d-1}+(d-2)a_{d-2}z^{d-2}+\ldots +a_1 z
\ ,
$$
and the last one contributes
$[(d-2)a_{d-2}\Delta+(d-1)^2a_{d-1}^2] \Delta$ because
$$
z^{d+1}=-z^2 P_0'(z)+(d-1)a_{d-1}z^{d}+(d-2)a_{d-2}z^{d-1}+\ldots
+a_1 z^2 \ ,
$$
and modulo $P'_0$ we have
$$
z^{d+1}=[(d-2)a_{d-2}\Delta+(d-1)^2a_{d-1}^2] z^{d-1}+\ldots [P_0']
$$
where the dots denote lower order terms. Thus we have
$$
\partial_{a_2} \Delta=\left (2(d-2)a_{d-2}+(d-1)^2a_{d-1}^2\right )
\Delta \ .
$$
When $d=3$ this gives
$$
\partial_{a_2} \Delta=\left (2a_1+4 a_2^2\right )
\Delta \ ,
$$
therefore
$$
\Pi_3 (X_0,X_1,X_2)=3X_0 +2X_2X_1 +{4\over 3} X_2^3 \ .
$$
When $d=4$ we get
$$
\partial_{a_2} \Delta=\left (4 a_2+9 a_3^2\right )
\Delta \ .
$$
So
$$
\Pi_4(X_0,X_1,X_2,X_3)=4X_0 +3X_3 X_1+ 2 X_2^2+9X_3^2X_2+\ldots \ ,
$$
where the remaining term is a polynomial in $X_3$.

When $d\geq 5$ we get
$$
\Pi_d (X_0,X_1,X_2,\ldots , X_{d-1})=dX_0 +(d-1)X_{d-1} X_1+ \left (
2 (d-2) X_{d-2} + (d-1)^2 X_{d-1}^2\right ) X_2 +\ldots
$$
where the remaining terms are independent of $X_0$, $X_1$ and $X_2$.

A close inspection of the procedure (for a complete analysis see
what follows next) shows that if $k\leq d/2$ then
$$
\partial_{a_k} \Delta =c \Delta \ ,
$$
where $c$ is a polynomial on $a_{d-1},a_{d-2},\ldots , a_{d-k}$ thus
the last result follows.$\diamond$

\bigskip

%\eno {Theorem III.4.3.7.}{Let $d\geq 2$. For $n\geq 0$ we define
%$(A_{n,k})_{0\leq k\leq d-1}$ as follows:
%
%For $n\leq d-1$, $A_{n,k}=0$ for $k\not=n$, and $A_{n,n}=1$.
%
%For $n=d$,
%$$
%A_{d,k}= ka_k \ .
%$$
%For $n\geq d+1$, the sequence $(A_{n,k})$ is defined
%inductively by
%$$
%A_{n+1,k}=(d-1)a_{d-1} A_{n,k}+(d-2) a_{d-2} A_{n-1,k}+\ldots + a_1
%A_{n-d+2,k} \ .
%$$
%Then modulo $P_0'$ we have, for all $n\geq 0$,
%$$
%z^n =A_{n,d-1} z^{d-1}+A_{n,d-2}z^{d-2}+\ldots +A_{n,1} z+A_{n,0} \
%[P_0'] .
%$$
%}
%
%\bigskip
%
%\stit {Proof.} For $n \leq d-1$ this is clear from the definition
%of the $A_{n,k}$'s. For $n = d$, we have
%$$\eqalign{
%z^d & = -zP_0'(z) + (d-1)a_{d-1}z^{d-1}+(d-2)a_{d-2}z^{d-2}+\ldots+a_1
%z \cr
%    & = -zP_0'(z) +
%    A_{d-1,d-1}z^{d-1}+A_{d-1,d-2}z^{d-2}+\ldots+A_{d-1,1}a_1 +
%    A_{d-1,0} \cr
%}$$
%(since $A_{d-1,k} = ka_k$), so modulo $P_0'$ the required equality holds.
%For $n \geq d+1$, the result follows by induction from the equality
%$$
%z^{n+1} =-z^{n-d+2} P'_0 +(d-1)a_{d-1} z^n+ (d-2) a_{d-2}
%z^{n-1}+\ldots +a_1 z^{n-d+2} \ .
%$$
%and the induction relation defining the $A_{n,k}$'s.
%$\diamond$

\eno {Theorem III.4.3.7.}{Let $d\geq 2$. For $n\geq 0$ we define
$(A_{n,k})_{0\leq k\leq d-1}$ to be the coefficients of the
remainder on dividing $z^n$ by $zP_0'$:
$$
z^n =A_{n,d-1} z^{d-1}+A_{n,d-2}z^{d-2}+\ldots +A_{n,1} z+A_{n,0} \
[zP_0'] .
$$

For $n\leq d-1$, $A_{n,k}=0$ for $k\not=n$, and $A_{n,n}=1$.

For $n=d$,
$$
A_{d,k}= ka_k \ .
$$
And for $n\geq d+1$, we can compute the sequence $(A_{n,k})$ by
induction by
$$
A_{n+1,k}=(d-1)a_{d-1} A_{n,k}+(d-2) a_{d-2} A_{n-1,k}+\ldots + a_1
A_{n-d+2,k} \ .
$$
}

\bigskip

\stit {Proof.}

Everything is clear except for the induction relation where we use
$$
z^{n+1} =-z^{n-d+2} P'_0 +(d-1)a_{d-1} z^n+ (d-2) a_{d-2}
z^{n-1}+\ldots +a_1 z^{n-d+2} \ .
$$
$\diamond$

\bigskip

\eno {Corollary III.4.3.8.}{For $d\geq 2$, $0\leq k\leq d-1$, and
$n\geq d$, $A_{n,k}$ is a polynomial with integer coefficients on
$a_0,a_1,\ldots , a_{d-1}$ of total degree $n-d+1$ }

\bigskip

\stit {Proof.}

This is straightforward from the induction relations.$\diamond$

\bigskip

Now we can compute the polynomial $\Pi_d$ using the polynomials
$(A_{n,k})$

\eno {Corollary III.4.3.9.}{For $d\geq 2$, the polynomial $\Pi_d$ is
uniquely determined by the equations, for $0\leq k\leq d-1$,
$$
\partial_{a_k} \Pi_d (a_0,\ldots , a_{d-1}) =A_{d-1+k,d-1}+A_{d-2+k,d-2}+\ldots
A_{d,d-k} \ .
$$
}

\stit {Proof.}

By differentiation column by column we get (as is clear from the
first computations above)
$$
\partial_{a_k} \Delta =\left ( A_{d-1+k,d-1}+A_{d-2+k,d-2}+\ldots
A_{d,d-k}\right ) \ \Delta \ ,
$$
and the result follows.$\diamond$

\bigskip

The non-vanishing of the ramificant has several corollaries.

\eno {Corollary III.4.3.10}{We consider the locus ramification
mapping $\Upsilon: \dd C^d \to \dd C^d$,
$$
\Upsilon (a_0,a_1 ,\ldots , a_{d-1})=(f_0(w_1^*),f_0(w_2^*),\ldots ,
f_0(w_d^*)), \ .
$$
Then $\Upsilon$ is a local diffeomorphism everywhere.}

\stit {Remark}

The ramification locus is not a global diffeomorphism as is easily
seen constructing two distinct log-Riemann surfaces with $d$
ramification points with the same images by the projection mapping
$\pi$.

\bigskip

\stit {Proof.} The computation of the differential at a point gives
the value of the ramificant at this points,
$$
D_{a_0,\ldots ,a_{d-1}} \Upsilon =\Delta (a_0,\ldots , a_{d-1}) \ ,
$$
and the result follows from the non-vanishing of the
ramificant.$\diamond$

\bigskip

The right philosophy is to think of $(f_i(w_j^*))$ as transalgebraic
numbers when $P_0(z)\in {\overline  {\dd Q }} [z]$. It is then
natural to ask if we have some relation between the $(f_i(w_j^*))$
and the coefficients of $P_0$ similar to the fundamental symmetric
formulas. We have the following:

\eno {Theorem III.4.3.11}{For $j=1,\ldots , d-1$ (note that $j=0$ is
excluded), we have that $e^{-a_0} a_j$ is a universal rational
function on $(f_k(w_l^*))_{k=0,\ldots , d\atop l=1,\ldots ,d}$.

More precisely, $\Delta e^{-a_0} a_j$, where $\Delta$ is the
ramificant, is a universal polynomial function of degree $d-1$ on
$(f_k(w_l^*))_{k=0,\ldots , d\atop l=1,\ldots ,d}$. }

\stit {Proof.}

Observe that for $l=1,\ldots ,d$ we have
$$
\eqalign { &d \ F_{d-1}(+\infty .\omega_l) + (d-1) a_{d-1}
F_{d-2}(+\infty . \omega_l)+\ldots +a_1 F_0(+\infty .\omega_l) \cr
&=\int_0^{+\infty . \omega_l} P'_0(z) e^{P_0(z)} \ dz \cr &=\left [
e^{P_0(z)}\right ]_0^{+\infty .\omega_l} \cr &=0-e^{a_0}=-e^{a_0} \
.\cr }
$$
Therefore if we consider the matrix
$$
M =\pmatrix{f_0(w_1^*) & f_1(w_1^*) & \ldots &f_{d-1}(w_1^*) \cr
f_0(w_2^*) &f_1(w_2^*) & \ldots &f_{d-1}(w_2^*) \cr \vdots &\vdots
&\ddots & \vdots \cr f_0(w_d^*) & f_1(w_d^*) & \ldots
&f_{d-1}(w_d^*) \cr} \ ,
$$
we have
$$
M.\pmatrix { a_1 \cr 2 a_2 \cr \vdots \cr (d-1)a_{d-1} \cr d \cr} =
-e^{a_0} \pmatrix { 1 \cr 1 \cr \vdots \cr 1 \cr 1 \cr} \ .
$$
Thus
$$
\pmatrix { a_1 \cr 2 a_2 \cr \vdots \cr (d-1)a_{d-1} \cr d
\cr}=-e^{a_0} M^{-1}.\pmatrix { 1 \cr 1 \cr \vdots \cr 1 \cr 1 \cr}
\ ,
$$
and the coefficients of $M^{-1}$ are polynomials on the entries of
$M$ divided by the ramificant $\Delta =\det M$.$\diamond$

\bigskip

As we have observed, just the location of the ramification points,
i.e. the values $(f_0(w_k^*))$ are not enough to characterize the
polynomial $P_0$ (or the log-Riemann surface). This changes if we
consider all values $(f_j(w_k^*))$ as the next corollary shows.

\bigskip

\eno {Corollary III.3.4.12.}{ Let $P_0$ and $Q_0$ be two normalized
polynomials,
$$
\eqalign {P_0(z)&={1\over d} z^d +a_{d-1} z^{d-1} +\ldots +a_1 z
+a_0\ , \cr Q_0(z)&={1\over d} z^d +b_{d-1} z^{d-1} +\ldots +b_1 z
+b_0\ . \cr}
$$
Consider the associated functions,
$$
\eqalign {F_j(z) &=\int_0^z t^j e^{P_0(t)} \ dt \ ,\cr G_j(z)
&=\int_0^z t^j e^{Q_0(t)} \ dt \ .\cr}
$$
If for $j=0,\ldots , d-1$ and $k=1,\ldots , d$,
$$
F_j(+\infty . \omega_k)=G_j (+\infty . \omega_k) \ ,
$$
then for $j=1,\ldots d-1$, we have $e^{a_0}a_j=e^{b_0} b_j$, i.e.
$$
e^{P_0(0)} (P_0(z)-P_0(0))=e^{Q_0(0)} (Q_0(z)-Q_0(0)) \ .
$$

In particular, if the polynomials have no constant term, then
$$
P_0=Q_0 \ .
$$
}

\stit {III.4.4) Infinite ramification points.}

The next step consists in distinguishing algebraically regular
points from infinite ramification points.

\eno {Theorem III.4.4.1.}{Consider a point $w_0 \in \cl S^*$ and let
${\hbox {\gothic M}} = {\hbox {\gothic M}}_{w_0}$ be the associated
maximal ideal in the ring $\cl A = \cl A_{\cl S}$. Let $\cl A_{\hbox
{\gothic M}}$ be the localization of $\cl A$ at the maximal ideal
${\hbox {\gothic M}}$, and let $\widehat{{\hbox {\gothic M}}}
\subset \cl A_{\hbox {\gothic M}}$ be the image of ${\hbox {\gothic
M}}$ in $\cl A_{\hbox {\gothic M}}$.

\medskip

{\parindent=2cm

 \item {$\bullet$} If $w_0 \in \cl S$ is a regular point, we have
$$
\widehat{{\hbox {\gothic M}}} / \widehat{{\hbox {\gothic M}}}^2
\approx \dd C^{d+1} \ .
$$

\item {$\bullet$} If $w_0 \in \cl S^*-\cl S$ is an infinite ramification point, we have
$$
\widehat{{\hbox {\gothic M}}} / \widehat{{\hbox {\gothic M}}}^2
\approx \dd C[z] \oplus \dd C^d \ .
$$

} }

\stit {Proof.}

Consider $f\in {\hbox {\gothic M}}$. For the corresponding $F\in
A_{P_0}$, we can write
$$
F(z)=z A(z,F_0, \ldots , F_{d-1}) e^{P_0(z)} + B(F_0,\ldots ,
F_{d-1}) \ ,
$$
where $A$ and $B$ are polynomials.

\medskip

If $w_0$ is a regular point of $\cl S^*$, then from the Taylor
expansions of $A$ and $B$ around the points $(z_0, F_0(z_0), \ldots
, F_{d-1}(z_0))$ and $(F_0(z_0), \ldots , F_{d-1}(z_0))$
respectively (where $z_0 = k_0(w_0) \in \dd C$), modulo ${\hbox
{\gothic M}}^2$ we have
$$
f = \left( \sum_{i = 1}^N b_i (z - z_0)^i e^{P_0(z)}  + c_0 (F_0(z)
- F_0(z_0)) + \dots + c_{d-1} (F_{d-1}(z) - F_{d-1}(z_0)) \right)
\circ k_0 \quad (mod \ {\hbox {\gothic M}}^2)
$$
for some constants $b_1, \dots, b_N, c_0, \dots, c_{d-1}$. Since the
function $e^{P_0(z)}$ doesn't vanish at $z_0$, the corresponding
element of $\cl A$ doesn't belong to ${\hbox {\gothic M}}$, and is
hence invertible in the localization $\cl A_{\hbox {\gothic M}}$.
Since $((z - z_0) e^{P_0(z)} ) \circ k_0 \in \widehat{ {\hbox
{\gothic M}} }$, it follows that $(z - z_0) \circ k_0 \in \widehat{
{\hbox {\gothic M}} }$ and so $((z - z_0)^i e^{P_0(z)}) \in
\widehat{ {\hbox {\gothic M}} }^2$ for $i \geq 2$. So modulo
$\widehat{ {\hbox {\gothic M}} }^2$,
$$
f = \left( b_1 (z - z_0) e^{P_0(z)}  + c_0 (F_0(z) - F_0(z_0)) +
\dots + c_{d-1} (F_{d-1}(z) - F_{d-1}(z_0)) \right) \circ k_0 \quad
(mod \ \widehat{{\hbox {\gothic M}}}^2)
$$
and so
$$
\widehat{{\hbox {\gothic M}}} / \widehat{{\hbox {\gothic M}}}^2
\approx \dd C^{d+1} \ .
$$

\medskip
If $w_0$ is an infinite ramification point, then expanding $A$ and
$B$ around the points $(0, f_0(w_0), \dots, f_{d-1}(w_0)),
(f_0(w_0), \dots, f_{d-1}(w_0))$ respectively, modulo ${\hbox
{\gothic M}}^2$ we have
$$
f = \left( \sum_{i = 1}^N b_i z^i e^{P_0(z)}  + c_0 (F_0(z) -
f_0(w_0)) + \dots + c_{d-1} (F_{d-1}(z) - f_{d-1}(w_0)) \right)
\circ k_0 \quad (mod \ {\hbox {\gothic M}}^2)
$$
for some constants $b_1, \dots, b_N, c_0, \dots, c_{d-1}$. Since the
function $e^{P_0} \circ k_0$ vanishes at all the infinite
ramification points, it remains noninvertible in the localization
$\cl A_{\hbox {\gothic M}}$, all the terms $(z^i e^{P_0(z)}) \circ
k_0$ are in $\widehat{{\hbox {\gothic M}}} - \widehat{{\hbox
{\gothic M}}}^2$, and are linearly independent in $\widehat{{\hbox
{\gothic M}}} / \widehat{{\hbox {\gothic M}}}^2$, so
$$
\widehat{{\hbox {\gothic M}}} / \widehat{{\hbox {\gothic M}}}^2
\approx \dd C[z] \oplus \dd C^d \ .
$$

$\diamond$

\null
\vfill
\eject

%% file: bib.tex
\input defv7.tex
\stit {Bibliography.}

\bigskip

\ref {[Ab]}{N.H. ABEL}{M\'emoire sur une propri\'et\'e g\'en\'erale d'une
classe tr\`es \'etendue de fonctions transcendantes}
{M\'emoires pr\'esent\'ees par divers savants, t. VII, Paris, 1841= Oeuvres compl\`etes,
edited by Sylow and Lie, Tome I, Christiana, 1881, p. 145-211.}

\ref {[Ah1]}{L.V. AHLFORS}{Complex analysis}
{Third Edition, McGraw-Hill, 1979.}

\ref {[Ah2]}{L.V. AHLFORS}{Lectures on quasiconformal mappings}
{Van Nostrand, 1966.}

\ref {[Ap-Go]}{P. APPELL, E. GOURSAT}{Th\'eorie des fonctions alg\'ebriques d'une
variable et des transcendances qui s'y rattachent}
{Tomes I et II; Gauthier Villars, 1930.}

\ref {[BBIF]}{Y. BORISOVICH, N. BLIZNYAKOV, Y. IZRAILEVICH, T. FOMENKO}
{Introduction to topology}
{MIR Publishers, Moscow, 1985.}

\ref {[Be]}{V. BELYI}{On the Galois extension of the maximal cyclotomic
field}
{Translated by A. Koblitz, Math. USSR-Izv., {\bf 14}, 2, 1980, p.247-253.}

\ref {[Bi]}{G.D. BIRKHOFF}{Dynamical Systems}
{American Mathematical Society, AMS Colloquium Publications, 1927.}

\ref {[Bi-PM1]}{K. BISWAS, R. P\'EREZ-MARCO}{Log-Riemann surfaces and entire functions}
{In preparation.}

\ref {[Bi-PM2]}{K. BISWAS, R. P\'EREZ-MARCO}{Tube-log Riemann surfaces}
{In preparation.}

\ref {[Bi-PM3]}{K. BISWAS, R. P\'EREZ-MARCO}{Transalgebraic
Kronecker-Weber theory} {In preparation.}

\ref {[Bo]}{J.B. BOST}{Introduction to compact Riemann surfaces,
Jacobians, and Abelian varieties}{ {\it From Number Theory to
Physics}, Les Houches 1989, Springer0Verlagm Berlinm 1992,
p.64-211.}

\ref {[BMS]}{I. BENJAMINI, S. MERENKOV, O. SCHRAMM}{A negative
answer to Nevanlinna's type question and a parabolic surface with
a lot of negative curvature} {arXiv: math.CV/0209334v2, 2002.}

\ref {[Che]}{C. CHEVALLEY}{Introduction to the theory of algebraic
functions of one variable} {Mathematical Surveys, {\bf 6},
American Mathematical Society, 1951.}

\ref {[Chu]}{G.V. CHUDNOVSKY}{Algebraic independence of constants connected with the
functions of analysis}
{Dokl. Ukrainian SSR Academy of Sciences, {\bf 8}, 1976.
Noticees of the AMS, {\bf 22}, 1975.}

\ref {[De-We]}{R. DEDEKIND, H. WEBER}{Theorie der algebraischen Funktionen einer
Ver\"anderlichen}
{ J. De Crelle, {\bf 92}, 1882, p.181-290.}

\ref {[Du]}{P.L. DUREN}{Univalent functions}
{Springer-Verlag, 1983.}

\ref {[Er1]}{A. EREMENKO}{Ahlfors contribution to the theory of meromorphic functions}
{ Lectures in the Memory of Lars Ahlfors, R. Brooks and M. Sodin Editors, Israel Math.
Conf. Proc., vol. 14, AMS, Providence RI, 2000.}

\ref {[Er2]}{A. EREMENKO}{Geometric theory of meromorphic functions}
{ "In the tradition of Ahlfors-Bers", Proc. Ahlfors-Bers colloquium,
AMS, Providence RI, 2003.}

\ref {[Eu1]}{L. EULER}{Opera Matematica}
{Teubner and O. F\"ussli, Leipzig-Berlin-Z\"urich, 1911-1985.}

\ref {[Eu2]}{L. EULER}{Analysin Infinitorum}
{.}

\ref {[FK]}{H.M. FARKAS, I. KRA}{Riemann surfaces}
{Graduate Texts in Mathematics, {\bf 72}, 2nd Edition, Springer Verlag, 1992.}

\ref {[Ga1]}{\'E. GALOIS}{Oeuvres math\'ematiques} {Publi\'ees par
\'Emile Picard, Gauthier-Villars, Paris, 1897. \'Editions J.
Gabay, 2001.}

\ref {[Ga2]}{\'E. GALOIS}{\'ecrits et m\'emoires math\'ematiques}
{\'Editeurs R. Bourgne et J.-P. Azra, Gauthier-Villars, 1976.}

\ref {[Gam]}{T. GAMELIN}{Complex analysis}
{Undergraduate Texts in Mathematics, Springer Verlag, 2001.}

\ref {[Gr]}{M. GROMOV}{Metric structures for Riemannian and non-Riemannian
spaces}
{Progress in Mathematics, {\bf 152}, Birkh\"auser, Boston MA, 1999.}

\ref {[He]}{M. HEINS}{Selected topics in the classical theory of functions of a
complex variable}
{Holt, Rinehart and Winston, New York, 1962.}

\ref {[Hi]}{E. HILLE}{Ordinary Differential Equations in the Complex Domain}{John Wiley \& Sons, Inc., New York,
1976}

\ref {[Ho]}{C. HOUZEL}{La g\'eom\'etrie alg\'ebrique. Recherches historiques}
{Blanchard, Paris, 2002.}

\ref {[Ho-Ri]}{H. HOPF, W. RINOW}{Ueber den Begriff der
vollst\"andigen differentialgeometrischen Fl\"ache}
{Commentarii, {\bf 3}, 1931, p.209-225.}

\ref {[Iwa]}{K. IWASAWA}{Collected Papers}
{Vol I and II, Springer, 2001.}

\ref {[Jo]}{C. JORDAN}{Cours d'analyse}
{Tomes I, II et III, Gauthier-Villars, Paris, 1909.}

\ref {[Ka]}{I. KAPLANSKI}{An introduction to differential algebra}
{Actualit\'es scientifiques et Industrielle, Hermann, 1957.}

\ref {[Ko]}{Z. KOBAYASHI}{Theorems on the conformal representation of Riemann
surfaces}
{Sci. Rep. Tokyo Bunrika Daigaku, sect. A, {\bf 39}, 1935.}

\ref {[La]}{J.L. LAGRANGE}{Oeuvres}
{14 volumes, Gauthier-Villars, Paris, 1867-1892.}

\ref {[Le-Vi]}{O. LEHTO, K. VIRTANEN}{Quasiconformal mappings in the plane}
{Springer-Verlag, 1973.}

\ref {[Li]}{F. LINDEMANN}{\"Uber die Zahl $\pi$}
{M.A. {\bf 20}, 1882, p.213-225.}

\ref {[Ma]}{A. MARKUSHEVICH}{Teor\'\i a de las funciones anal\'\i ticas}
{Volumen I and II, Translated from the russian by E. Aparicio Bernardo,
ediciones MIR, 1978.}

\ref {[Mi1]}{J. MILNOR}{Morse theory}
{Annals of Mathematics Studies, {\bf 51}, Princeton University Press, 1963.}

\ref {[Mi2]}{J. MILNOR}{Dynamics in one complex variable}
{Friedr. Vieweg \& Sohn, Braunschweig, 1999; arXiv:math.DS/9201272v1, 1990.}

\ref {[MOS]}{W. MAGNUS, F. OBERHETTINGER, R.P. SONI}{Formulas and
theorems for the special functions of Mathematical Physics} {Die
Grundlehren der Mathematischen Wissenschaften in
Einzeldarstellungen, {\bf 52}, Springer-Verlag, 1966.}

\ref {[Mu]}{J. MU\~NOZ D\'IAZ}{Curso de teor\'\i a de funciones}
{Editorial Tecnos, Madrid, 1978.}

\ref {[Na]}{R. NARASIMHAN}{Several complex variables}{Chicago Lectures in Mathematics,
University of Chicago Press, 1971.}

\ref {[Ne1]}{R. NEVANLINNA}{\"Uber Riemannsche Fl\"ache mit endlich vielen Windungspunkten}
{Acta Mathematica, {\bf 58}, 1932.}

\ref {[Ne2]}{R. NEVANLINNA}{Analytic functions}
{Grundlehren der Mathematischen Wissenschaften in Einzeldarstellungen,
{\bf 162}, 2nd Edition, Springer-Verlag, 1953.}

\ref {[Ne3]}{R. NEVANLINNA}{Le Th\'eor\`eme de Picard-Borel et la th\'eorie
des fonctions m\'eromorphes}
{Collection Monographies sur la th\'eorie des fonctions, Gauthier-Villars,
Paris, 1929.}

\ref {[Ne-Pa]}{R. NEVANLINNA, V. PAATERO}{Introduction to Complex Analysis}{Addison-Wesley Publishing Company,
1969}

\ref {[Oe]}{J. OESTERL\'E}{Polylogarithmes}
{S\'eminaire Bourbaki {\bf 762}, Ast\'erisque, {\bf 216}, 1993, p.49-67.}

\ref {[Pe]}{M. PERDIGAO DO CARMO}{Geometria Riemanniana}
{2nd edition, IMPA, CNPq, Rio de Janeiro, 1988.}

\ref {[Pi]}{\'E. PICARD}{Trait\'e d'analyse}
{Tomes I, II et III; Gauthier-Villars, 3rd Edition, 1926.}

\ref {[PM1]}{R. P\'EREZ-MARCO} {Sur les dynamiques holomorphes non lin\'earisables et une conjecture de V.I.
Arnold} {Ann. Scient. Ec. Norm. Sup. 4 serie, {\bf 26}, 1993, p.565-644; (C.R. Acad. Sci. Paris, {\bf 312}, 1991,
p.105-121).}

\ref {[PM2]}{R. P\'EREZ-MARCO} {Uncountable number of symmetries for non-linearisable holomorphic dynamics}
{Inventiones Mathematicae, {\bf 119}, {\bf 1}, p.67-127,1995; (C.R. Acad. Sci. Paris, {\bf 313}, 1991, p.
461-464).}

\ref {[PM3]}{R. P\'EREZ-MARCO} {Functions with exponential singularities on compact Riemann surfaces} {Manuscript,
2003.}

\ref {[PM4]}{R. P\'EREZ-MARCO} {Log-euclidean geometry and "Grundlagen der Geometrie"} {Manuscript, 2003.}

\ref {[PMBBJM]}{R. P\'EREZ-MARCO, K. BISWAS, D. BURCH, N. JONES,
E. MU\~NOZ-GARCIA} {Transalgebraic Number Theory} {Manuscript,
2002.}

\ref {[Ri1]}{B. RIEMANN}{Theorie der Abel'schen Functionen}
{Journal de Crelle, {\bf 54}, 1857; Oeuvres Math\'ematiques de B. Riemann, Blanchard,
1968, p.89.}

\ref {[Ri1]}{B. RIEMANN}{Theorie der Abel'schen Functionen}
{Journal de Crelle, {\bf 54}, 1857; Oeuvres Math\'ematiques de B. Riemann, Blanchard,
1968, p.89.}

\ref {[Ri2]}{B. RIEMANN}{Sur les hypoth\`eses qui servent de fondement \`a la g\'eom\'etrie}
{M\'emoire de la Soci\'et\'e Royale des Sciences de G\"ottingue, 1867;
Oeuvres Math\'ematiques de B. Riemann, Blanchard,
1968, p.280.}

\ref {[Rit1]}{J.F. RITT}{Integration in finite terms. Liouville's
theory of elementary methods} {Columbia University Press, New
York, 1948.}

\ref {[Rit2]}{J.F. RITT}{Differential Algebra} {AMS, 1950.}

\ref {[Ta1]}{M. TANIGUCHI}{Explicit representation of structurally finite
entire functions}
{Proc. Japan Acad., {\bf 77}, 2001, p.69-71.}

\ref {[Ta2]}{M. TANIGUCHI}{The size of the Julia set of a structurally finite
transcendental entire function}
{Math. Proc. Camb. Phil. SOc., {\bf 135}, 2003, p.181-192.}

\ref {[Ta3]}{M. TANIGUCHI}{Synthetic deformation space of an entire function}
{"Value distribution theory and complex dynamics, Hong Kong, 2000, Contemp. Math.,
{\bf 303}, Amer. math. Soc., Providence, RI, 2002, p.107-136.}

\ref {[Tsu]}{M. TSUJI}{Potential theory in modern function theory}
{Maruzen, Tokio, 1959.}

\ref {[Val]}{G. VALIRON}{Cour d'analyse. Equations fonctionnelles}
{Masson, 1945.}

\ref {[Var]}{V.S. VARADARAJAN}{Algebra in ancient and modern times}
{Mathematical World, {\bf 12}, AMS, 1991.}

\ref {[Wa1]}{M. WALDSCHMIDT}{Open Diophantine Problems} {Preprint,
www.math.jussieu.fr /~miw/articles/pdf/odp.pdf, 2004.}

\ref {[Wa2]}{M. WALDSCHMIDT}{Les travaux de Chudnovsky sur les nombres transcendants}
{S\'eminaire Bourbaki, $28$\`eme ann\'ee, 1975/76, {\bf 488}, June 1976, Lecture Notes
in Mathematics, {\bf 567}, p.274-292.}

\ref {[Wall]}{H.S. WALL}{Analytic theory of continued fractions}
{Chelsea Publishing Company, Bronx, N.Y., 1967.}

\ref {[Was]}{L.C. WASHINGTON}{Introduction to cyclotomic fields}
{Graduate Texts in Mathematics, {\bf 83}, Second Edition, Springer-Verlag, 1997.}

\ref {[We]}{H. WEYL}{The concept of a Riemann surface}
{Addisson-Wesley, 1964.}